\documentclass[11pt,a4paper]{amsart}
\usepackage[T1]{fontenc}
\usepackage{lmodern}
\usepackage{amsmath,amssymb,amsthm,mathtools}
\usepackage[a4paper,margin=27mm]{geometry}
\usepackage[expansion=false]{microtype}
\usepackage{booktabs,longtable,array,enumitem,needspace}
\usepackage[colorlinks=true,linkcolor=blue,citecolor=blue,urlcolor=blue]{hyperref}
\usepackage[capitalise,nameinlink,noabbrev]{cleveref}

\AddToHook{env/center/begin}{%
  \setlength{\topsep}{8pt plus 2pt minus 1pt}%
  \setlength{\partopsep}{0pt}%
}
\hypersetup{pdftitle={Classification of quasi-lisse affine vertex algebras},bookmarksdepth=3}
\allowdisplaybreaks[2]
\numberwithin{equation}{section}
\newtheorem{theorem}{Theorem}[section]
\newtheorem{proposition}[theorem]{Proposition}
\newtheorem{lemma}[theorem]{Lemma}
\newtheorem{corollary}[theorem]{Corollary}
\theoremstyle{definition}

\theoremstyle{remark}
\newtheorem{remark}[theorem]{Remark}
\newtheorem{example}[theorem]{Example}
\newcommand{\C}{\mathbb C}
\newcommand{\Z}{\mathbb Z}
\newcommand{\Q}{\mathbb Q}
\newcommand{\g}{\mathfrak g}
\newcommand{\h}{\mathfrak h}
\newcommand{\uu}{\mathfrak u}
\newcommand{\lev}{\mathfrak l}
\newcommand{\N}{\mathcal N}
\newcommand{\OO}{\mathbb O}
\newcommand{\clO}[1]{\overline{\mathbb O}_{(#1)}}
\newcommand{\spg}{\mathfrak{sp}}
\newcommand{\sog}{\mathfrak{so}}
\newcommand{\slg}{\mathfrak{sl}}
\newcommand{\glg}{\mathfrak{gl}}
\newcommand{\SF}{\operatorname{SF}}
\newcommand{\Sp}{\operatorname{Sp}}
\newcommand{\Vir}{\operatorname{Vir}}
\newcommand{\ch}{\operatorname{ch}}
\newcommand{\Hom}{\operatorname{Hom}}
\newcommand{\Ext}{\operatorname{Ext}}
\newcommand{\im}{\operatorname{im}}
\newcommand{\id}{\operatorname{id}}
\newcommand{\tr}{\operatorname{tr}}
\newcommand{\ad}{\operatorname{ad}}
\newcommand{\rank}{\operatorname{rank}}
\newcommand{\Sym}{\operatorname{Sym}}
\newcommand{\Find}{\operatorname{F-ind}}

\newcommand{\arxiv}[1]{\href{https://arxiv.org/abs/#1}{arXiv:#1}}

\newcommand{\wt}{\operatorname{wt}}
\newcommand{\Com}{\operatorname{Com}}

\newcommand{\WW}{\mathcal W}

\newcommand{\file}[1]{\nolinkurl{#1}}

\usepackage[width=.95\textwidth]{caption}

\AddToHook{env/proposition/begin}{\crefalias{theorem}{proposition}}
\AddToHook{env/lemma/begin}{\crefalias{theorem}{lemma}}
\AddToHook{env/corollary/begin}{\crefalias{theorem}{corollary}}
\AddToHook{env/definition/begin}{\crefalias{theorem}{definition}}
\AddToHook{env/remark/begin}{\crefalias{theorem}{remark}}
\AddToHook{env/example/begin}{\crefalias{theorem}{example}}

\title[Classification of quasi-lisse simple affine vertex algebras]{Classification of quasi-lisse simple affine vertex algebras}

\author{Dra\v{z}en Adamovi\'{c}}
\email{adamovic@math.hr}
\address{Department of Mathematics, Faculty of Science, University of Zagreb, Bijeni\v{c}ka 30, 10 000 Zagreb, Croatia}

\author{Xuanzhong Dai} 
\email{xzdai@nju.edu.cn}
\address{School of Mathematics, Nanjing University,
Nanjing 210093, China}

\author{Kaiwen Sun}
\email{kwsun@ustc.edu.cn}

\address{School of Mathematical Sciences, University of Science and Technology of China, Hefei 230026, Anhui, China}

\date{}

\begin{document}
\begin{abstract}
We classify all noncritical levels $k$ for which the simple affine vertex algebra $L_k(\g)$ is quasi-lisse, for every finite-dimensional complex simple Lie algebra $\g$. We prove that Fiebig's equivalence between the corresponding noncritical blocks of category $\mathcal O$ identifies the relevant parabolic Wakimoto modules and the homomorphisms from the corresponding universal affine vertex algebras. Combining this compatibility with a characterization of semisimple points in associated varieties, we establish that for coprime positive integers $p,q$, $L_{-h^\vee+p/q}(\g)$ is quasi-lisse if and only if $L_{-h^\vee+p/\gcd(q,r)}(\g)$ is quasi-lisse, where $r$ is the lacing number of $\g$. The classification at the resulting levels is obtained using Drinfeld--Sokolov reduction, singular vectors and relations in the reductions of the maximal ideals of $V^k(\g)$. We also determine the associated varieties at the quasi-lisse levels for all $\mathfrak{g}$ of rank at most six and for $E_7$ and $E_8$.
\end{abstract}
\maketitle
\enlargethispage{3pt}
\setcounter{tocdepth}{1}\tableofcontents\medskip
\section{Introduction}\label{sec:classification}

Let $\g$ be a finite-dimensional complex simple Lie algebra, and let $L_k(\g)$ denote the simple affine vertex algebra at level $k$. In this paper we classify all noncritical levels $k$ for which $L_k(\g)$ is quasi-lisse. Our main result, stated in
\cref{thm:classification}, gives the classification in every simple type, including both admissible and nonadmissible levels.

Recall that a vertex algebra $V$ is called quasi-lisse if its associated variety $X_V$ has finitely many symplectic leaves \cite{AK}.
For a simple affine vertex algebra $L_k(\g)$, this condition is equivalent to $X_{L_k(\g)}\subseteq\mathcal N(\g)$, 
where $\mathcal N(\g)$ denotes the nilpotent cone of $\g$.
In this case, $L_k(\g)$ has only finitely many irreducible ordinary modules, and their normalized characters satisfy modular linear differential equations \cite{AK}.
At admissible levels, the associated varieties are closures of nilpotent orbits, determined in \cite{AAV,AVEM}. The representation theory of admissible affine vertex algebras was studied in \cite{AdamovicMilas1995,ArakawaRational2016}.
Examples at nonadmissible levels satisfying
$X_{L_k(\g)}\subseteq\mathcal N(\g)$ were obtained in
\cite{ADS2026, AVD42026,  ArakawaEtAl2025, AK,  ArakawaMoreau2018,ArakawaMoreau2017,ArakawaFutornyKrizka2024}.

Two earlier developments are particularly relevant to our proof.
The construction of ideals by singular vectors in
\cite{Adamovic2003} provides explicit relations in universal affine
vertex algebras. Conformal embeddings and collapsing levels studied in
\cite{AdamovicEtAl2018,AdamovicEtAl2017,AKMPP2020} identify simple quotients
of $W$-algebras and relate their representation theory to affine vertex
algebras. We use both kinds of results, while keeping track of which
relations hold in the reduction of the maximal ideal.

The classification relies on the fact that $L_k(\g)$ is quasi-lisse if and only if its associated variety contains no nonzero semisimple point. A differential equation satisfied by the polynomial images of affine singular vectors reduces this question to finitely many conditions on semisimple elements. We first determine which representative levels are quasi-lisse, using singular vectors,
Drinfeld--Sokolov reduction, known descriptions of associated
varieties and estimates of character growth. We then compare
parabolic Wakimoto homomorphisms at different levels to extend
the classification to all denominators of $k+h^\vee$.

Recently, Shan, Yan and Zhao conjectured a description of associated
varieties at positive rational values of $k+h^\vee$ in simply-laced types
\cite{ShanYanZhao2025,ShanYanZhao2026}. Our classification confirms their
criterion for quasi-lisseness, see \cref{cor:SYZ}. Their description of
the entire associated variety as a generalized sheet closure is stronger.

Four-dimensional $\mathcal N=2$ superconformal field theories (SCFT) provide further motivation. The Schur-sector construction \cite{BeemEtAl2015}
associate a vertex algebra to every 4D $\mathcal N=2$ SCFTs. The associated variety of the vertex algebra is conjecturally the Higgs branch of the 4D SCFT \cite{BeemRastelli2018}. Several simple affine vertex algebras in our classification occur
as the vertex algebras of Argyres--Douglas theories
\cite{XieYanYau2021,SongXieYan2017,WangXie2019,
XieYanLisse2021,LiXieYan2023}.
Free-field realizations of affine examples are constructed in
\cite{BeemMeneghelliRastelli2019}. 
Affine VOAs with negative integer levels also appear in non-critical string theories \cite{DZL,JaramilloKashaniPoor2023}.

Related connections also arise from three-dimensional $\mathcal N=4$ gauge theories and their boundary vertex algebras
\cite{Gai2019,CG2019,CCG2019,BF2025,BCDN2023,AFM2026}.
For families of $A$-twisted boundary vertex algebras associated
with classical gauge groups, Arakawa, Dai and Song identify the associated varieties with the Hamiltonian reductions describing the Higgs branches
and prove quasi-lisseness \cite[Section~6]{ArDS2026}. Their construction in particular, includes the quasi-lisse simple affine vertex algebras
$L_{-2m}(\mathfrak{so}_{4m+4})$, $m\ge1$, of type~$D$.

The classification in this paper also provides examples for the study of modularity of vertex-algebra characters. Arakawa and Kawasetsu \cite{AK} proved that the
normalized characters of ordinary modules over a quasi-lisse VOA satisfy
a common modular linear differential equation (MLDE). 
Beem and Rastelli \cite{BeemRastelli2018} relate these equations to null vectors and Schur indices of four-dimensional $\mathcal N=2$ SCFTs. Flavored characters also keep track of Cartan weights and depend on flavor parameters as well as the modular parameter.
Null vectors can give flavored modular differential equations in these variables, see  \cite{ZhengPanWang2022,PanWang2023,PanYang2024}. The modular invariance of characters of quasi-lisse vertex algebras has also recently been studied in \cite{ArakawaVanEkerenLi2026}. 

\subsection{Notation and scope}\label{prop:critical}
For the affine vertex-algebra construction and Zhu algebras we refer to
\cite{FrenkelZhu1992,Zhu1996}. The Poisson algebra $R_V$ below is the
degree-zero component associated with Li's canonical filtration
\cite{Li2005}.

For a vertex algebra $V$, 
its \emph{Poisson algebra $R_V$} and \emph{associated variety $X_V$} are defined by
$$
 \begin{gathered}
 R_V=V/C_2(V),\qquad
 X_V=(\text{Spec} \, R_V)_{\mathrm{red}},
 \end{gathered}
$$
where $C_2(V)=\operatorname{span}\{a_{(-2)}b:a,b\in V\}.$ The $R_V$ is also called the $C_2$ algebra of $V$. 
A vertex algebra is \emph{lisse} if its Poisson algebra is finite-dimensional,
and \emph{quasi-lisse} if its associated variety has finitely many symplectic leaves.

Let $\g$ be a finite-dimensional complex simple Lie algebra and normalize
its invariant form so that long roots have squared length two. Let
$\widehat\g=\g\otimes\C[z,z^{-1}]\oplus\C K$, with
$$
 [x(m),y(n)]=[x,y](m+n)+m(x,y)\delta_{m+n,0}K.
$$
The universal affine vertex algebra is
$$
 V^k(\g)=U(\widehat\g)
 \otimes_{U(\g\otimes\C[z]\oplus\C K)}\C_k,
$$
where $\g\otimes\C[z]$ acts trivially on $\C_k$ and $K$ acts by $k$.
Let $I_k$ be its maximal proper graded ideal and define the simple affine
vertex algebra by $L_k(\g)=V^k(\g)/I_k$.
Write $t=k+h^\vee$. For a noncritical level that is $t\ne0$,
the Sugawara vector associated with dual bases $u_a,u^a$ of $\g$ is
$$
 \omega=\frac1{2(k+h^\vee)}\sum_a u_a(-1)u^a(-1)\mathbf1.
$$
Its conformal grading is nonnegative and integral, with finite-dimensional
homogeneous subspaces. We call a conformal module \emph{ordinary} if $L_0$
acts diagonally, its eigenspaces are finite-dimensional, and its spectrum
is contained in finitely many sets $a+\Z_{\ge0}$.

We say that homogeneous vectors $u^1,\ldots,u^r$ strongly generate a
vertex algebra $V$ if $V$ is spanned by $\mathbf1$ and the vectors
$u^{i_1}_{(-n_1-1)}\cdots u^{i_s}_{(-n_s-1)}\mathbf1$ with $n_j\ge0$. 
These are iterated normally ordered products of derivatives of the
$u^i$. Strong generation is a spanning assertion. It does not assert
linear independence. The generators are freely strongly generating if
the products with the factors in a fixed order form a basis, with the
usual restriction on repeated odd factors in the superalgebra case.
For positively graded vertex algebras, strong generation can also be
checked on the images in $R_V$, see \cite[Section~2.3]{AAV}.
Whenever we use a basis theorem below, we specify the universal
enveloping algebra or universal $W$-algebra to which it applies.

The invariant form identifies $\g$ with $\g^*$. All orbit closures and
associated varieties below are reduced. Unless explicitly stated otherwise, $t\ne0$ and the conformal vector
is the Sugawara vector.

The lacing number $r$ is the ratio of the squared lengths of long and short
roots. In a simply-laced root system $r=1$. The constants are
\begin{center}
\begin{tabular}{@{}lrrrrrrrrr@{}}\toprule
Type&$A_n$&$B_n$&$C_n$&$D_n$&$E_6$&$E_7$&$E_8$&$F_4$&$G_2$\\\midrule
$h^\vee$&$n+1$&$2n-1$&$n+1$&$2n-2$&12&18&30&9&4\\
$h$&$n+1$&$2n$&$2n$&$2n-2$&12&18&30&12&6\\
$r$&1&2&2&1&1&1&1&2&3\\\bottomrule
\end{tabular}
\end{center}
Here $A_n$ means $\slg_{n+1}$, $B_n$ means $\sog_{2n+1}$,
$C_n$ means $\spg_{2n}$, and $D_n$ means $\sog_{2n}$.

We use the standard criterion \cite[Lemma~6.1]{AK}
\begin{equation}\label{eq:nilpotent-criterion}
 L_k(\g)\text{ is quasi-lisse}
 \quad\Longleftrightarrow\quad X_{L_k(\g)}\subseteq\N(\g).
\end{equation}
It is already known that $L_k(\g)$ is lisse precisely when
$k\in\Z_{\ge0}$ \cite{AAV}. Our interest is therefore in the additional
quasi-lisse levels. At the critical level, the simple
graded quotient obtained from the augmentation character of the center
has $C_2$ algebra $\C[\N(\g)]$ and associated variety $\N(\g)$
\cite[Section~2 and Theorem~1.2]{Critical}. It has no Sugawara conformal
vector and is outside the classification theorem below.

\subsection{Classification of quasi-lisse levels}
\begin{theorem}\label{thm:classification}\label{thm:classical-classification}
Let $\g$ be a finite-dimensional complex simple Lie algebra and let
$k\ne-h^\vee$. The vertex algebra
$L_k(\g)$ is quasi-lisse if and only if $t=p/q\in\Q_{>0}$,
with $p,q\in\Z_{>0}$ and $\gcd(p,q)=1$, and the corresponding
condition in the following table holds.
\begin{center}\small
\renewcommand{\arraystretch}{1.24}
\begin{tabular}{@{}p{.13\textwidth}p{.77\textwidth}@{}}\toprule
Type&Condition on $k+h^\vee=p/q$\\\midrule
$A_n$, $n\ge1$&$p\ge n+1$.\\
$B_n$, $n\ge2$&
 $2n\le3(p-1)$ if $p,q$ are both odd, and otherwise $p>n$ and $p\ne2n-1$.\\
$C_n$, $n\ge2$&
 $p\ge\lceil n/2\rceil+1$ and $p\ne n$ if $q$ is odd, and
 $p\ge n+2$ if $q$ is even.\\
$D_n$, $n\ge4$&
 $p\ge2\lceil n/2\rceil$, or $p$ is even and $2n/3<p\le n-2$.\\
$E_6$&$p\in\{6,7\}$ or $p\ge9$.\\
$E_7$&$p\in\{6,8\}$ or $p\ge10$.\\
$E_8$&$p=6$ or $p\ge9$.\\
$F_4$&$p\ge5$, or $p=3$ and $q$ is odd.\\
$G_2$&$p\ge2$ if $3\nmid q$, and $p\ge4$ if $3\mid q$.\\\bottomrule
\end{tabular}
\end{center}
Every condition is imposed on the \emph{reduced} fraction. In particular,
$q$ even forces $p$ odd, and $3\mid q$ forces $3\nmid p$.
\end{theorem}

The admissible ranges are
\begin{equation}\label{eq:admissible-general}
 p\ge h^\vee\quad\text{if }\gcd(q,r)=1,
 \qquad
 p\ge h\quad\text{if }\gcd(q,r)=r>1.
\end{equation}
Their quasi-lisseness and associated varieties are known
\cite{AAV,AVEM}. The additional entries in the table are established
in the type-specific sections.

For a positive integer $m$, let $\OO_{\g}(m)$ be the dense orbit in
$\{x\in\N(\g):(\ad x)^{2m}=0\}$, as in
\cite[Theorem~2.1.6]{ShanYanZhao2025}. Recall that a nilpotent orbit
is \emph{distinguished} if it meets no proper Levi subalgebra. Comparing our
table with the orbit descriptions in
\cite[Sections~2.1.9, 2.1.12 and 2.1.13]{ShanYanZhao2025}
gives the following consequence.
\begin{corollary}\label{cor:SYZ}
Let $\g$ be simple and simply-laced, and write
$k+h^\vee=p/q>0$ with $p,q$ coprime positive integers. Then
$$
 L_k(\g)\text{ is quasi-lisse}
 \quad\Longleftrightarrow\quad
 \OO_{\g}(p)\text{ is distinguished}.
$$
In particular, quasi-lisseness is independent of $q$ for fixed reduced
numerator $p$.
\end{corollary}
For the conjectural description of associated varieties underlying
this criterion, see \cite{ShanYanZhao2026}.

\subsection{How we prove quasi-lisseness}\label{sec:organization}
At the representative levels we ask whether $X=X_{L_k(\g)}$ contains
a nonzero semisimple element. Section~\ref{intro:level-comparison}
extends the answer to all denominators. The basic differential equation
is as follows. Let
$\sigma:V^k(\g)\to R_{V^k(\g)}=\Sym\g$ be the quotient map, and
normalize the root vectors by $(e_\alpha,f_\alpha)=1$.
For $u\in\h$, regular $h\in\h$, and $v\in V^k(\g)$, the identity
proved in \cref{com:radial} is
$$
 (k+h^\vee)\partial_u\sigma(v)(h)=
 \sum_{\alpha>0}\frac{\alpha(u)}{\alpha(h)}
       \sigma\bigl(f_\alpha(0)e_\alpha(0)v\bigr)(h)
       +\sigma\bigl(u(1)v\bigr)(h).
$$
This is equation~\eqref{com:radial-equation} in the body of the paper.
For the finite-dimensional $\g$-submodule $E$ generated by an affine
singular vector, define $F(h)(v)=\sigma(v)(h)$ for $v\in E$, and let
$A_\alpha$ be the transpose of $f_\alpha(0)e_\alpha(0)$ on $E$.
The positive-mode term vanishes, and the equation becomes
$$
 (k+h^\vee)\,dF=\sum_{\alpha>0}A_\alpha F\,\frac{d\alpha}{\alpha}.
$$
This is the rational Casimir equation on the zero-weight subspace of
$E^*$. The Casimir connection was discovered by De Concini around 1995
in the study of braid-group monodromy and quantum Weyl groups
\cite[Sections~1 and~5]{ToledanoLaredo2002}, and was mentioned in print
by Procesi in 1996 \cite{Procesi1996}. It was also obtained
independently by Millson--Toledano Laredo
\cite{MillsonToledanoLaredo2005} and by
Felder--Markov--Tarasov--Varchenko \cite{FMTV2000}. The normally
ordered form is given in \cite[Theorem~2.1]{FMTV2000}.
Here the variables are Cartan coordinates, and the unknowns are
polynomial images of vectors in an affine vertex algebra. The modular
linear differential equations for characters mentioned earlier have a
different role.

For the whole ideal $I_k$, the term $\sigma(u(1)v)$ must be retained.
It relates consecutive conformal degrees and allows the same argument
to apply to all defining relations of $L_k(\g)$. The resulting
geometric statement is simple: if a semisimple element $s$ belongs
to $X$, then the entire center $\mathfrak z(\g^s)$ belongs to $X$.
Thus it suffices to exclude one nonzero element $h_i$ from the center
of each standard maximal Levi subalgebra
(\cref{com:finite-centers}). The differential equation also controls
orders of vanishing of individual singular vectors. We use this more
explicitly in type $B$ (Section~\ref{B:covariant-restriction}) and
at $k=-8/3,-7/3$ in type $G_2$ (see \cref{vtwo:one-return} and the discussion following it).

We use several methods to show that these semisimple elements lie outside the associated variety. In some
families the singular vectors already give sufficiently strong
polynomial equations on $X$. Their images impose matrix rank bounds,
determinantal or Pfaffian relations, and restrictions on the Cartan
subalgebra. For example, in Section~\ref{sec:D-upper}, a quadratic
matrix identity together with a rank bound forces every point to be
nilpotent. The corresponding direct arguments in type $B$ are given
in Sections~\ref{B:relations}--\ref{B:half-integral}.

When these equations do not settle the question, we use
Drinfeld--Sokolov reduction. For a nilpotent element $f$, write
$H_f=H_f^0(L_k(\g))$. The basic geometric statements are
$$
 H_f\ne0\quad\Longleftrightarrow\quad f\in X,
 \qquad X_{H_f}=X\cap\mathcal S_f
$$
when $H_f\ne0$, where $\mathcal S_f$ is the Slodowy slice.
Consequently, a vanishing reduction excludes a nilpotent orbit.
Parabolic Wakimoto homomorphisms relate this information to semisimple
points. We assume, for a contradiction, that $h_i\in X$. The
factorization criterion of \cref{com:point-map} then gives an embedding
of $L_k(\g)$ and forces a family of adjoint orbits into $X$, together
with its nilpotent boundary (\cref{com:levi-support}). If that boundary
contains an orbit independently shown to lie outside the associated variety, the assumption $h_i\in X$
is impossible. This argument is used in type $C$
(Section~\ref{seed:C-positive}) and in several exceptional cases
(\cref{vtwo:exceptional-caps}).

A further argument is needed when the boundary orbit itself is allowed.
We then determine relations in the reduction $H_f^0(I_k)$ of the
maximal ideal. These can force $X\cap\mathcal S_f$ to be a single
point, whereas a family containing semisimple elements would give a
positive-dimensional intersection. The Jantzen filtration helps locate these relations. Invariant bilinear
forms and vertex products identify the relations themselves. Affine vertex subalgebras
at nonnegative integral levels often account for part of the reduction,
while relations among the remaining invariant generators complete the argument. This is
used in Sections~\ref{B:lower-integral}, \ref{seed:D-positive} and
\ref{exc:seeds}. In other cases the class of the conformal vector
in $R_{H_f}$ is a polynomial in the classes of the weight-one
generators. Evaluating this identity gives $X\cap\C h=\{0\}$ for
the semisimple element $h$ of the chosen $\slg_2$-triple, see
\cref{vtwo:sugawara-axis}.
The essential point throughout is that these are relations in
$H_f^0(L_k(\g))$ itself. A relation in some other simple quotient of
$\mathcal W^k(\g,f)$ would not establish the required bound on $X$.
We also use previously known associated varieties and intersections
of sheets where they already give the answer, see
\cref{exc:positive-geometric}.

\subsection{How we prove that a vertex algebra is not quasi-lisse}
At the representative levels, we either find a nonzero semisimple
point in $X_{L_k(\g)}$ or contradict quasi-lisseness directly.
A semisimple point already proves non-quasi-lisseness. A Wakimoto
embedding deduced from it gives an additional realization. The use
of these embeddings to compare levels is explained below.

Drinfeld--Sokolov reduction is often more accessible than $L_k(\g)$
itself. Conformal-weight inequalities, Jantzen identities or associated varieties establish $H_f\ne0$, sometimes conditionally on quasi-lisseness (Section~\ref{com:reduction}). We identify the simple graded quotient $B$ of $\mathcal W^k(\g,f)$. Uniqueness gives $H_f\twoheadrightarrow B$.
It suffices to show that $B$ is not quasi-lisse, since nonzero reductions and quotients preserve quasi-lisseness in this setting.

Several concrete arguments are used for this last step. If $B$
contains a Heisenberg vertex algebra as a tensor factor, or has a
realization in a Heisenberg vertex algebra in which its conformal
vector has a nonzero quadratic image in the $C_2$ algebra, then its
conformal class is not nilpotent.
This contradicts quasi-lisseness by
\cref{com:conformal-obstruction}. Miura homomorphisms and tensor
product decompositions make this argument applicable in many of the
classical and exceptional cases, see Sections~\ref{sec:C-exclusions},
\ref{sec:D-exclusions}, \ref{seed:FG} and \ref{sec:exceptional-e}.
Another argument constructs infinitely many simple ordinary modules
for the quotient itself, as in the commutant of the vertex algebra generated by $J$ considered in
Section~\ref{seed:B-coset}. Finally, the remaining odd numerators in
type $D$ require a character-growth argument
(Section~\ref{seed:D-growth}). There, quasi-lisseness would force the even part of a certain commutant
to be lisse. The upper bound on character growth obtained from its strong
generators then contradicts the modular lower bound. Some exceptional
cases follow directly from published descriptions of associated varieties
as sheet closures containing nonzero semisimple elements, see
Section~\ref{sec:exceptional-e}.

\subsection{Parabolic Wakimoto realizations of the simple affine vertex algebra}\label{intro:wakimoto-embedding}
Let $t=k+h^\vee\in\Q_{>0}$ and fix a maximal parabolic with Levi
subalgebra $\lev=\C h_i\oplus\bigoplus_a\lev_a$.
Section~\ref{com:parabolic} constructs a homomorphism
$\Phi_{h_i}(t):V^k(\g)\to T_{h_i}(t)$, where
$$
 T_{h_i}(t)=\mathcal D^{\mathrm{ch}}(U_i)\otimes
 \mathcal H_{t(h_i,h_i)}\otimes
 \bigotimes_a L_{d_at-h_a^\vee}(\lev_a).
$$
Here $\mathcal D^{\mathrm{ch}}(U_i)$ is the $\beta\gamma$ vertex
algebra associated with the nilradical, $\mathcal H_{t(h_i,h_i)}$
is the rank-one Heisenberg vertex algebra generated by $H$ with
$[H_\lambda H]=t(h_i,h_i)\lambda\mathbf1$, and $d_a$ is the index
of the invariant form on $\lev_a$. As a $\widehat\g$-module, this is
the parabolic Wakimoto module induced from
$\bigotimes_a L_{d_at-h_a^\vee}(\lev_a)$, regarded as a module over
itself, and the Heisenberg Fock module on which $H(0)$ acts as zero
\cite{dBF,ABGG}.

The homomorphism is initially defined on $V^k(\g)$. Write
$B_i(t)=\operatorname{im}\Phi_{h_i}(t)$. To realize the simple quotient,
we must show that $\Phi_{h_i}(t)(I_k)=0$.
Theorem~\ref{com:point-map} gives the criterion
$$
 h_i\in X_{L_k(\g)}
 \quad\Longleftrightarrow\quad
 \Phi_{h_i}(t)(I_k)=0
 \quad\Longleftrightarrow\quad B_i(t)\text{ is simple}.
$$
Thus the occurrence of a prescribed semisimple point is expressed as
the simplicity of a particular Wakimoto image. Once these conditions
hold, the induced nonzero homomorphism from $L_k(\g)$ is injective.
It gives an embedding of vertex algebras into $T_{h_i}(t)$ and
identifies $L_k(\g)$ with a $\widehat\g$-submodule of the
corresponding parabolic Wakimoto module.

For methods that determine maximal ideals by proving simplicity of
quotients at nonadmissible levels, see \cite{APVMaximal2026}.
Explicit generators of $I_k$ are not needed for the criterion used here.
The property of singular vectors established in \cref{com:leading-fiber}
controls every homogeneous component of the ideal $I_k$. Any nonzero
image of $I_k$ would contain,
in its lowest nonzero conformal degree, a $\g$-highest-weight vector
annihilated by every $a(n)$ with $a\in\g$ and $n>0$. The lemma
detects such a vector by a nonzero polynomial obtained from its
$C_2$ image. On a dense open subset, the relevant Wakimoto coordinates
describe elements conjugate to nonzero scalar multiples of $h_i$.
If $h_i$ lies in the associated variety, then every polynomial
relation from $I_k$ vanishes there. This forces the image
of the ideal to be zero without a separate calculation for each generator. See the proof of \cref{com:point-map}.

If $h_i\in X_{L_k(\g)}$, the criterion shows that the parabolic
Wakimoto homomorphism associated with $h_i$ factors through
$L_k(\g)$, giving an embedding into $T_{h_i}(t)$.
Conversely, a direct proof that $B_i(t)$ is simple establishes
$h_i\in X_{L_k(\g)}$. We first prove one of these conditions
at a representative level. The equivalence in
Section~\ref{com:transport} then identifies the corresponding
Wakimoto images at different levels and preserves their simplicity.
For example, if Drinfeld--Sokolov reduction shows that
$L_k(\g)$ is not quasi-lisse, \cref{com:finite-centers} gives
$h_i\in X_{L_k(\g)}$ for at least one $i$. The criterion therefore
gives an embedding for the maximal parabolic associated with
that $h_i$. It does not assert that such an embedding exists
for every maximal parabolic.

The embeddings also give modules. Use the simple affine vertex algebras
attached to the Levi factors as modules over themselves and vary the $H(0)$-eigenvalue of the Fock
module. If $\varpi_i$ corresponds to the omitted root, charges giving
highest weights $m\varpi_i$, $m\in\Z_{\ge0}$, yield simple affine
highest-weight subquotients \cite[Section~7.1]{ABGG}. Their lowest spaces
are finite-dimensional $\g$-modules and their conformal-weight spaces
are finite-dimensional, so these are pairwise nonisomorphic simple
ordinary $L_k(\g)$-modules. This agrees with the finiteness theorem
of \cite{AK}. The full Wakimoto modules generally are not ordinary.
When the embedding was deduced from a semisimple point, this family
is a further consequence, not an independent exclusion argument.

\subsection{Comparison of levels and completion of the classification}\label{intro:level-comparison}
The simplicity criterion of \cite[Theorem~0.2.1]{GK} first reduces
the possible levels: $V^k(\g)$ is nonsimple precisely when
$$
 r(k+h^\vee)\in\Q_{\ge0}\setminus\{1/m:m\in\Z_{>0}\}.
$$
At noncritical levels outside $k+h^\vee\in\Q_{>0}$, it follows that
$L_k(\g)=V^k(\g)$ and $X_{L_k(\g)}=\g^*$. Within the positive
rational range, the same simplicity criterion also settles $p=1$
when $r\mid q$.

Next, if $k+h^\vee=p/q$ is in lowest terms and $e=\gcd(q,r)$,
Section~\ref{com:transport} proves
\begin{equation}\label{eq:binary-transport}
 L_{-h^\vee+p/q}(\g)\text{ is quasi-lisse}
 \quad\Longleftrightarrow\quad
 L_{-h^\vee+p/e}(\g)\text{ is quasi-lisse}.
\end{equation}
Here the Wakimoto criterion is used to pass between two different
levels. Set $t_0=p/e$ and $t=p/q$. The equivalence between the
corresponding affine highest-weight categories identifies the
parabolic modules $T_{h_i}(t_0)$ and $T_{h_i}(t)$ and the
homomorphisms defining $B_i(t_0)$ and $B_i(t)$.
It therefore preserves simplicity of their images. The comparison is
$$
\begin{aligned}
 h_i\in X_{L_{t_0-h^\vee}(\g)}
 &\quad\Longleftrightarrow\quad B_i(t_0)\text{ is simple}\\
 &\quad\Longleftrightarrow\quad B_i(t)\text{ is simple}\\
 &\quad\Longleftrightarrow\quad h_i\in X_{L_{t-h^\vee}(\g)}.
\end{aligned}
$$
Section~\ref{com:transport} identifies these particular homomorphisms,
proving the middle equivalence independently of the classification.
Applying the criterion to every maximal Levi gives
\eqref{eq:binary-transport}. Thereafter no new Wakimoto calculation
is needed at $t$. It remains to treat $k+h^\vee=p$ in simply-laced
types, and $p$ or $p/r$ otherwise. The associated varieties may differ.

Example~\ref{ex:sl4-transport} illustrates this argument for
$L_{-4+3/q}(\slg_4)$, for every positive integer $q$ coprime to $3$.
It identifies the parabolic Wakimoto homomorphism carried by the
equivalence and explains how its simple image produces the semisimple
point at the new level.

Admissibility settles large numerators. Uniform arguments treat the
remaining classical ranges (Sections~\ref{sec:type-A}--\ref{sec:type-D}),
and finitely many representatives suffice for exceptional types
(Sections~\ref{seed:FG}--\ref{sec:exceptional-e}). The common results
are in Sections~\ref{com:centers}--\ref{com:reduction}.
Appendix~\ref{app:calculations} lists the calculations and gives the
bounds used to evaluate them in \cite{Certificates}.
Section~\ref{cert:orbits} determines the orbits listed in
Appendix~\ref{app:associated-varieties}.

After the initial draft of this paper was completed, the preprint
\cite{JiangWu2026} appeared, which determines the highest weights and
multiplicities of the singular vectors of lowest conformal weight
in types~$B$ and~$C$. Several of the singular subspaces described
there were obtained independently in earlier versions of the
present manuscript.
In the present version, we use the results of \cite{JiangWu2026}
to shorten the corresponding existence proofs in
Section~\ref{B:relations}.

\section{Semisimple points and parabolic Wakimoto homomorphisms}\label{com:centers}
We first reduce the study of semisimple points of the associated variety to centers of maximal Levi subalgebras, then characterize when the corresponding homomorphisms from $V^k(\g)$ into parabolic Wakimoto realizations factor through $L_k(\g)$.

\subsection{Semisimple points of the associated variety}
Fix a complex simple Lie algebra $\g$, with long roots of squared length
two. Recall $t=k+h^\vee$, and denote for short
\begin{equation}
    V_t=V^{t-h^\vee}(\g),\qquad L_t=L_{t-h^\vee}(\g),\qquad I_t=\ker(V_t\to L_t).
\end{equation}
Unless indicated otherwise, $t\ne0$. The invariant form identifies
$\g$ with $\g^*$, and the quotient map
$$
 \sigma:V_t\longrightarrow R_{V_t}=\Sym\g
$$
identifies $x(-1)\mathbf1$ with the linear function $(x,-)$.
By the Poincar\'e--Birkhoff--Witt (PBW) theorem, $V_t$ has a basis of
vectors $x_1(-r_1)\cdots x_m(-r_m)\mathbf1$, where the $x_j$
belong to a fixed Cartan and root-vector basis, $r_j\ge1$, and the
pairs $(x_j,r_j)$ are arranged in a fixed order.
Such a vector has zero image under $\sigma$ if some $r_j\ge2$. For a proper graded vertex ideal $I$, set
$X_I=X_{V_t/I}$. All varieties in this section are reduced.

On the zero-weight restrictions of affine singular spaces, the
homogeneous equation below is the rational Casimir equation
\cite[Theorem~2.1]{FMTV2000}, see also \cite{MillsonToledanoLaredo2005}.
We give the vertex-algebra calculation because the positive-mode term
is needed for the entire ideal, including vectors that are not singular.

\begin{lemma}\label{com:radial}
Choose $e_\alpha,f_\alpha$ with
$(e_\alpha,f_\alpha)=1$ and
$[e_\alpha,f_\alpha]=\alpha^\sharp$.
For $u\in\h$, $v\in V_t$, and regular $h\in\h$,
\begin{equation}\label{com:radial-equation}
 t\partial_u\sigma(v)(h)=
 \sum_{\alpha>0}\frac{\alpha(u)}{\alpha(h)}
 \sigma\bigl(f_\alpha(0)e_\alpha(0)v\bigr)(h)
 +\sigma\bigl(u(1)v\bigr)(h).
\end{equation}
Consequently, for every semisimple $s\in X_I$,
\begin{equation}\label{com:center-propagation}
                 \mathfrak z(\g^s)\subset X_I.
\end{equation}
\end{lemma}
\begin{proof}
First observe that if $v$ belongs to the first positive term of Li's filtration, then
$\sigma(u(1)v)(h)=0$ for $u,h\in\h$. Commute $u(1)$ through a
basis vector $x_1(-r_1)\cdots x_m(-r_m)\mathbf1$ as above.
A factor $x(-r)$ with $r\ge2$ remains,
except possibly when $[u,a](-1)$ replaces $a(-2)$. This replacement is
zero for Cartan $a$, and is a root vector otherwise, so its Cartan
evaluation is zero. Reordering negative modes introduces only further
factors $x(-r)$ with $r\ge2$. Thus the claimed evaluation depends only
on $P=\sigma(v)$.

On a symmetrized product of the operators $x(-1)$, $x\in\g$, the central contraction contributes
$k\partial_uP$. The two-factor contraction contributes
$$
 \sigma(u(1)v)(h)=k\partial_uP(h)+
 \sum_{\alpha>0}\alpha(u)\alpha(h)
                       P_{e_\alpha f_\alpha}(h).
$$
Indeed the relevant bracket is half the sum
$[[u,a],b]+[[u,b],a]$. Only opposite roots survive evaluation.
The zero-mode derivations satisfy
$$
 (f_\alpha e_\alpha P)(h)=
 \alpha(h)\partial_{\alpha^\sharp}P(h)
 -\alpha(h)^2P_{e_\alpha f_\alpha}(h).
$$
Combining these identities with
$\sum_{\alpha>0}\alpha(u)\alpha^\sharp=h^\vee u$
proves \eqref{com:radial-equation}. In another root-vector
normalization its root term is divided by $(e_\alpha,f_\alpha)$.

Conjugate $s$ into $\h$ and set $Z=\mathfrak z(\g^s)$.
For each Sugawara degree $D$, let $E_D=I\cap(V_t)_D$ and let
$F_D(h)\in E_D^*$ be evaluation of its symbols at $h$.
These spaces are finite dimensional, $\g$-stable, and
$u(1)E_D\subset E_{D-1}$. Equation
\eqref{com:radial-equation} is therefore a finite linear system for
$F_D$, with an inhomogeneous term depending on $F_{D-1}$.
For $u\in Z$, omit all roots satisfying $\alpha(s)=0$ before
restricting the equation to $Z$: their numerators $\alpha(u)$ vanish
identically. Every remaining denominator is nonzero near $s$ in $Z$.

Induct on $D$. Properness gives $E_0=0$. If $F_{D-1}|_Z=0$, then
on each short segment $s+zu$ in $Z$ the equation for $F_D$ is a
homogeneous linear ODE with analytic coefficients and initial value
$F_D(s)=0$. Its unique solution is zero. Hence $F_D$ vanishes near
$s$ in $Z$, and polynomiality makes it vanish on all of $Z$.
This proves \eqref{com:center-propagation} for the whole ideal.
\end{proof}

For each standard maximal parabolic $\mathfrak p_i=\lev_i+\uu_i$,
choose a real $h_i\in\mathfrak z(\lev_i)$ positive on $\uu_i$.
Its norm is nonzero. Its scale will not matter.
The characterization of quasi-lisseness by $X_I\subset\N(\g)$ is
\cite[Lemma~6.1]{AK}. The preceding lemma reduces that criterion to
the following finite collection of semisimple elements.
\begin{corollary}\label{com:finite-centers}
For a proper graded quotient $V_t/I$ with $t\ne0$,
$$
 X_I\subset\N(\g)
 \quad\Longleftrightarrow\quad h_i\notin X_I\text{ for every }i.
$$
These conditions also characterize quasi-lisseness of $V_t/I$.
\end{corollary}
\begin{proof}
If $X_I$ contains a nonnilpotent point, it contains its nonzero
semisimple part, which lies in its adjoint-orbit closure. The centralizer
of that part is a proper Levi. After conjugation it is contained in a
standard maximal Levi $\lev_i$, and
$\mathfrak z(\lev_i)\subset\mathfrak z(\g^s)$.
The lemma implies that $h_i$ lies in $X_I$. The reverse implication is immediate.
Affine associated varieties are closed conical unions of adjoint orbits,
which are their symplectic leaves. There are finitely many nilpotent
orbits. Conversely, conical scaling of a nonzero semisimple invariant
produces infinitely many orbits. This proves the last assertion.
\end{proof}

\subsection{The parabolic Wakimoto homomorphism}\label{com:parabolic}
In this section $t>0$ is rational. Fix a maximal parabolic with Levi
$\lev=\C h\oplus\bigoplus_a\lev_a$, where $h$ is as above and
$\lev_a$ is simple. Let $d_a$ be its form index:
\begin{equation}\label{com:levi-normalization}
 (\ ,\ )_{\g}|_{\lev_a}=d_a(\ ,\ )_a,
 \qquad k_a=d_at-h_a^\vee.
\end{equation}
Here the form on $\lev_a$ is normalized by its own longest roots.
Thus $d_a=2/(\alpha_a,\alpha_a)_\g$ for a longest root $\alpha_a$
of that component. It is one for a component containing long roots
of $\g$. An entirely short component has index equal to the lacing number.
The quantum correction to the Levi level is
$k(\ ,\ )_\g+(\kappa_\g-\kappa_\lev)/2$, which gives
\eqref{com:levi-normalization}, rather than $k_a=d_ak$.

Composing the parabolic Wakimoto homomorphism with
$V^{k_a}(\lev_a)\twoheadrightarrow L_{k_a}(\lev_a)$ gives
\begin{equation}\label{com:target}
 \Phi_h(t):V_t\longrightarrow T_h(t):=
 \mathcal D^{\mathrm{ch}}(U)\otimes
 \mathcal H_{t(h,h)}\otimes
 \bigotimes_a L_{d_at-h_a^\vee}(\lev_a).
\end{equation}
The factor $\mathcal H_{t(h,h)}$ is regarded as its own Fock
module, and its generating operator $H(0)$ acts as zero. We use the usual
polynomial big-cell realization. Its construction and conformal vector
are the parabolic Wakimoto formulas of de Boer--Feh\'er,
\cite[Theorem~3]{dBF}. Their module version is Feigin--Frenkel
induction, as in \cite[Section~7.1, equation~(50)]{ABGG}.
Write $\widetilde T_h(t)$ for the tensor product obtained by replacing
$L_{k_a}(\lev_a)$ in \eqref{com:target} by $V^{k_a}(\lev_a)$ for every $a$. Denote the corresponding homomorphism by
$\widetilde\Phi_h(t)$ and the quotient map by
$\pi_h(t):\widetilde T_h(t)\twoheadrightarrow T_h(t)$.

\begin{lemma}\label{com:universal-parabolic-injective}
The homomorphism $\widetilde\Phi_h(t)$ is injective.
\end{lemma}
\begin{proof}
Its associated graded map is the pullback of the parabolic moment map.
At a point with nonzero central coordinate and zero nilradical
coordinates, its differential is surjective: variation in the Levi
factor gives $\lev$, while the two nilradical directions are identified
with $\uu$ and $\uu^-$ using the invertibility of $\ad h$ there.
Thus the moment map is dominant and smooth on a dense open subset.
The induced maps on finite jet spaces are dominant as well.
Use Li's filtrations on the source and target. Their associated
graded algebras are the coordinate rings of the corresponding arc
spaces, and the induced homomorphism is the pullback just described.
These filtrations are separated, so injectivity of the associated
graded homomorphism implies injectivity of $\widetilde\Phi_h(t)$.
\end{proof}

The vertex algebra $\mathcal D^{\mathrm{ch}}(U)$ admits commuting
left and right actions of $\uu[z,z^{-1}]$. Let $R_x(z)$,
$x\in\uu$, be the fields defining the right action. The subspace
annihilated by the left action of $\uu[z]$ is the module generated
by $\mathbf1$ under $R_x(-m)$, $x\in\uu$, $m\ge1$.
Each $R_x$ is linear in the variables $\beta$, has conformal weight
one and positive $h$-degree. The variables $\gamma$ have negative
$h$-degree. Write $H$ for the generator of $\mathcal H_{t(h,h)}$.
The fields associated with $\lev_a$ act on $R_x$ through the adjoint
action and on the $V^{k_a}(\lev_a)$-module used in the construction. Their formulas contain no independent term proportional to $H$.
For $y\in\uu^-$, the only term involving $H$ in the image of
$y(-1)\mathbf1$ is
$$
                         :H\phi_y(\gamma):,
$$
where $\phi_y(0)=0$. The coefficient of the coordinate corresponding
to the root vector dual to $y$ is a nonzero multiple of
$\alpha(h)/(h,h)$. This follows by projecting the first commutator
in $\operatorname{Ad}(\exp\gamma)y$ to $\C h$.
The fields $R_x$ are primary of weight one for the image of the
Sugawara vector of $V_t$. The $h$-grading also gives these
properties: a correction to $R_x$ independent of $\beta$, or a
scalar coefficient in the relevant products of the left and right
fields or of the conformal vector with $R_x$, would have positive
$h$-degree while being a polynomial in coordinates of negative
$h$-degree.

For comparison, \cite[Theorem~1.2]{ArakawaJiangMoreau2021} proves that
a nonzero affine singular vector in $V^k(\g)$ has nonzero image in
$R_{V^k(\g)}$ at a noncritical level. We need a statement for vectors
in $T_h(t)$ after a further restriction of their $C_2$ images.
The following proof establishes that additional property.

\begin{lemma}\label{com:leading-fiber}
Let $0\ne v\in T_h(t)$ be a homogeneous $\g$-highest-weight
vector of dominant highest weight, annihilated by $a(n)$ for every
$a\in\g$ and $n>0$. Its $C_2$ symbol is nonzero on
the fiber obtained by setting the coordinates $\gamma$ and the
derived-Levi augmentation ideals equal to zero.
\end{lemma}
\begin{proof}
We show that a leading component of $v$ contains only the modes that
survive on the stated $C_2$ fiber.

First, we choose the leading component.
Left $\uu[z]$-invariance first places $v$ in
$$
 U\bigl(\uu^{\mathrm{op}}[z^{-1}]z^{-1}\bigr)
 \otimes\C[H(-1),H(-2),\ldots]\otimes A,
 \qquad A=\bigotimes_a L_{k_a}(\lev_a).
$$
Maximize the number $d$ of $H$ creation factors, and then the number $r$ of operators $R_x(-l)$. The resulting nonzero associated graded
component $P$ is a polynomial in $R_x(-l)$ and $H(-s)$ with
coefficients in $A$. This order of the two filtrations is essential.

Second, we remove the higher right-action modes. 
Choose root vectors $x_\alpha\in\g_\alpha$ for
$\alpha\in\Delta(\uu)$ and dual root vectors
$y_\alpha\in\g_{-\alpha}$. At the coordinate origin,
$$
 (R_{x_\beta}\phi_{y_\alpha})(0)
       =c_\alpha\delta_{\alpha\beta},\qquad c_\alpha\ne0.
$$
Indeed, projection of $[x_\beta,y_\alpha]$ to $\C h$ vanishes
unless $\alpha=\beta$. In that case its coefficient is a nonzero
multiple of $\alpha(h)/(h,h)$, as in the construction above.
Thus each opposite root vector tests one right-action variable.

For the vector $y_\alpha(n)v$ use the momentum-degree filtration of the
whole beta-gamma space, with coordinates of degree zero.
At $\gamma=0$, including all coordinate creation derivatives, the image of a product of $r$ operators $R_x(-l)$ becomes
a monomial of degree $r$ in the variables $\beta$. Examine
$y_\alpha(n)v$ in bidegree $(d+1,r-1)$, $n\ge1$.
Only $:H\phi_{y_\alpha}:$ can increase the first degree. It must contract
at least one field $R_x$ to create a negative $H$ mode. The
identity
$$
 [\,:Hf(\gamma):(n),R_x(m)\,]
                  =\pm:H(R_xf)(\gamma):(n+m)
$$
shows that each further contraction consumes another right factor.
Thus exactly one contraction contributes to the chosen momentum
degree. After coordinate evaluation, singularity gives
\begin{equation}\label{com:fiber-equation}
 \sum_{l>n}c_\alpha H(n-l)
                  \frac{\partial P}{\partial R_{x_\alpha}(-l)}=0,
 \qquad \alpha\in\Delta(\uu).
\end{equation}
Lower lexicographic components cannot contribute there. In particular,
the coordinate evaluation is a linear projection on the displayed
polynomial space. Compatibility with the module action is not assumed. If $L>1$ is the largest index of an operator $R_x(-L)$ occurring in $P$,
take $n=L-1$. Multiplication by $H(-1)$ is injective and every
$c_\alpha$ is nonzero. Applying \eqref{com:fiber-equation} for each
$\alpha$ therefore removes every $R_x(-L)$. Descending induction
leaves only the $R_x(-1)$.

Next, we reduce the Levi coefficients to scalars.
For $x\in[\lev,\lev]$ and $n>0$, the operator $x(n)$ preserves
the number of $H(-s)$ factors and does not increase the number of
$R_y(-l)$ factors. On an $R_y(-1)$ factor it produces an operator of index $n-1$,
which is nonnegative. To survive it must commute through another operator $R_x(-l)$, decreasing length. Hence at bidegree $(d,r)$ the
positive Levi modes act only on the coefficients of $P$.
In $L_{k_a}(\lev_a)$, the common kernel of $x(n)$,
$x\in\lev_a$, $n>0$, is $\C\mathbf1$: a nonzero
highest-weight vector of positive degree in this kernel would
generate a proper $\widehat{\lev}_a$-submodule. Thus the
coefficients of $P$ in $A$ are scalars.

In the end, we remove the higher Heisenberg modes and restrict to the fiber. 
Affine singularity and Sugawara imply $L_nv=0$ for $n>0$.
At the same top bidegree the contribution from the fields $R_x$ again lowers
length, and the Levi contribution annihilates $\mathbf1$. The scalar part is
$$
 \sum_{s>n}sH(n-s)\frac{\partial P}{\partial H(-s)}=0.
$$
The background-charge terms and the terms containing two positive Heisenberg
modes lower $H$ count.
Taking $n=S-1$ for the largest $S>1$ removes $H(-S)$, successively.
We obtain a nonzero polynomial
$$
                 P\in\Sym(\uu[-1])\otimes\C[H(-1)].
$$
On the stated $C_2$ fiber its variables become the independent
momentum coordinates and the Heisenberg coordinate $[H]$. Its leading bidegree cannot be
cancelled by a lower component. This proves the assertion.
\end{proof}

Parabolic realizations and their relation to Dixmier sheets are also
used in \cite[Theorem~3.24]{ArakawaFutornyKrizka2024}, under the hypotheses imposed
there on the nilradical and the induced Levi levels. The criterion below
uses the simple affine Levi factors and applies to the entire ideal
of $V_t$.

\begin{theorem}\label{com:point-map}
For \eqref{com:target}, write $B_h(t)=\im\Phi_h(t)$. Then
\begin{equation}\label{com:point-equivalence}
 h\in X_{L_t}\quad\Longleftrightarrow\quad
 \Phi_h(t)(I_t)=0\quad\Longleftrightarrow\quad B_h(t)\text{ is simple}.
\end{equation}
If these conditions hold, then also
\begin{equation}\label{com:levi-support}
 \overline{G\cdot\bigg(\C h+
       \prod_a X_{L_{k_a}(\lev_a)}\bigg)}\subset X_{L_t}.
\end{equation}
More generally, for any proper graded ideal $J\subset V_t$,
$h\in X_{V_t/J}$ if and only if $\Phi_h(t)(J)=0$.
In this case \eqref{com:levi-support} holds with $X_{V_t/J}$
in place of $X_{L_t}$.
\end{theorem}
\begin{proof}
On the fiber in \cref{com:leading-fiber}, the moment map is
$(\beta,H)\mapsto Hh/(h,h)+\uu^-$. For $H\ne0$, each such point
is conjugate to $Hh/(h,h)$: eliminate its nilradical component in
increasing center degree using the invertibility of $\ad h$.
Thus if $h\in X_{L_t}$, conicity and invariance make every ideal
symbol vanish on a dense open subset of this fiber, and hence on all
of it. If $\Phi_h(t)(I_t)\ne0$, choose its least nonzero conformal
degree and a highest-weight vector for $\g$ there. Every operator $a(n)$ with $a\in\g$ and $n>0$
annihilates that vector. Local finiteness of the $\g$-action follows because it is
an image of $V_t$. By the nonvanishing statement in \cref{com:leading-fiber},
its symbol is nonzero on the specified $C_2$-fiber, a contradiction. Therefore $\Phi_h(t)$ factors through
$L_t$, and its image is isomorphic to $L_t$.

Conversely, if $B_h(t)$ is simple, then $\ker\Phi_h(t)=I_t$. Evaluate the
induced $C_2$ map at zero beta-gamma coordinates, the Levi origins,
and $H=(h,h)$ to obtain $h\in X_{L_t}$. Allowing arbitrary points
of the associated varieties of the simple affine vertex algebras attached
to the Levi factors in this evaluation proves
\eqref{com:levi-support}, followed by $G$-invariance and closedness.
For the last assertion, the same nonvanishing statement on the specified
$C_2$-fiber, \cref{com:leading-fiber}, applies to $\Phi_h(t)(J)$: its least nonzero degree is finite-dimensional as a $\g$-module,
and a $\g$-highest-weight vector in that degree is annihilated by
$a(n)$ for every $a\in\g$ and $n>0$.
Conversely, $\Phi_h(t)(J)=0$ gives a map from $V_t/J$ to $T_h(t)$.
The same $C_2$ evaluations give the point $h$ and the support inclusion.
\end{proof}

\section{Reduction of the denominator of \texorpdfstring{$k+h^\vee$}{k+h dual}}
\label{com:transport}
Let $r$ be the lacing number: $r=1$ in types $A,D,E$,
$r=2$ in types $B,C,F_4$, and $r=3$ in type $G_2$.
To compare quasi-lisseness at two levels, we use the semisimple-point
criterion of \cref{com:point-map}. For each maximal Levi subalgebra,
we must compare the corresponding homomorphisms $\Phi_h$.
There are two steps: identify their Wakimoto targets under the block
equivalence, then identify the maps from the universal affine vertex
algebras. Exactness then identifies their images and preserves simplicity.

Root dilation supplies the data for Fiebig's block equivalence
\cite[Theorem~4.1]{Fiebig}. Lemmas~\ref{com:twisting-compatibility}
and~\ref{com:relative-limit} compare the semi-infinite inductions of the
modules $M_\lev(\mu)$ defined below.
Lemma~\ref{com:induction-fully-faithful} allows us to pass to the simple
Levi factors that define the targets. Finally, the relevant space of
vacuum homomorphisms is one dimensional, which identifies $\Phi_h$
up to a nonzero scalar. The resulting comparison gives the theorem
below. It preserves quasi-lisseness, although the associated varieties
at the two levels may differ.

\begin{theorem}\label{thm:denominator}
Let $t=p/q>0$, with $p,q$ coprime positive integers, and set
$e=\gcd(q,r)$. Then
\begin{equation}\label{com:denominator-statement}
 L_{p/q-h^\vee}(\g)\text{ is quasi-lisse}
 \quad\Longleftrightarrow\quad
 L_{p/e-h^\vee}(\g)\text{ is quasi-lisse}.
\end{equation}
\end{theorem}

\subsection{Integral root systems and Levi subalgebras}
The pairing of $k\Lambda_0+\widehat\rho$ with the coroot of
$\alpha+m\delta$ is
\begin{equation}\label{com:integral-pairing}
 \langle\rho,\alpha^\vee\rangle+
                    \frac{2mp}{q(\alpha,\alpha)}.
\end{equation}
Thus the integral periods are $q$ for long roots and $q/e$ for
short roots. At the representative $t_0=p/e$ they are $e$ and one,
respectively. Set $D=q/e$. The map
\begin{equation}\label{com:root-dilation}
 \alpha+m\delta\longmapsto\alpha+Dm\delta
\end{equation}
is an isomorphism of based integral real-root systems, preserving the
finite subgroup, positivity, the pairings with $k\Lambda_0+\widehat\rho$, and the
full stabilizer of its closed dominant representative. This includes
singular representatives. These stabilizers are finite: $k+h^\vee>0$ places the representative on a proper face of the affine alcove,
whose stabilizer is a finite parabolic subgroup. This argument does not apply at $k=-h^\vee$. Fiebig's block theorem
(\cite{Fiebig}, Theorem~4.1) gives a labelled exact equivalence
$\mathcal E$ of the corresponding positive blocks. We use their
large versions allowing algebraic coproducts, with the common dominant
upper bound. Integrability for the finite-dimensional Lie algebra $\g$ is specified
by the finite simple-root labels, so the equivalence also preserves
$V^k(\g)$ and its simple quotient at the corresponding levels.

The Levi labels agree at the corresponding levels. Since
$(\alpha,\alpha)_a=d_a(\alpha,\alpha)_\g$,
$$
 \frac{2m(d_at)}{(\alpha,\alpha)_a}
       =\frac{2mt}{(\alpha,\alpha)_\g}.
$$
Moreover $\rho_\g-\rho_\lev$ is central on $\lev$.
The derived integral subsystems and every finite Levi label therefore
correspond under \eqref{com:root-dilation}. In particular, for a short-root $B_1$ component, identified with
$A_1$ with its normalized invariant form, the induced level satisfies
$\ell+2=2t$. For a $C_1$ component it satisfies $\ell+2=t$.
For a component of type $A$ generated by short roots in types $C,F_4$,
one has $\ell+h^\vee_{A}=2t$, and in type $G_2$ one has
$\ell+2=3t$. In every case the induced level plus the dual Coxeter number is positive.

\subsection{Categories of modules for the affine Lie algebra \texorpdfstring{$\widehat\lev$}{l-hat}}
\label{com:Levi-categories}
We include the degree derivation $d$ in the affine Cartan and use
$[d,x(n)]=nx(n)$. In a positive noncritical block $\Lambda$, let
$\lambda^+$ be its closed dominant representative, with a fixed
degree coordinate. We use the large category $\mathcal O_\Lambda$
of Cartan weight modules which are locally finite over the positive
Borel and have their simple subquotients in $\Lambda$.
Every object has weights in $\lambda^+-\widehat Q_+$.
Arbitrary weight multiplicities and algebraic direct sums are allowed. Each finitely generated submodule has finite-dimensional weight spaces.
The projective cover $P(\lambda)$ is the projective cover in
$\mathcal O^{\leq\lambda^+}_\Lambda$.
By \cite[Theorem~2.2]{Fiebig}, it has a Verma filtration with
$$
 (P(\lambda):M(\mu))=[M(\mu):L(\lambda)].
$$
Only $\lambda\leq\mu\leq\lambda^+$ can occur. This interval is
finite and each multiplicity is finite. The projective covers are
finitely generated, hence small with respect to algebraic direct sums.
They generate the large block, so every object has a presentation by
direct sums of them, see \cite[Section~4]{Fiebig}.

Write $[\lev,\lev]=\bigoplus_i\lev_i$. Fix levels $\ell_i$ with
$\ell_i+h_i^\vee>0$, one positive block for each
$\widehat{\lev}_i$, and a scalar $a$ such that $h(0)=a\,\id$
on the inducing Levi modules, where $Z(\lev)=\C h$. The Heisenberg Lie algebra associated with $Z(\lev)$ has
nonzero level. On a bounded module, the oscillator projection onto
the common kernel of its positive modes is finite on every vector.
Multiplication by negative modes identifies the module with the
charge-$a$ Fock representation tensored with this kernel.
Let $L^{\rm H}_0$ be the Sugawara operator of the vertex algebra
associated with $Z(\lev)$. This operator is semisimple.

For each $\widehat{\lev}_i$ at level $\ell_i$, let $L^{(i)}_0$ be its Sugawara
operator and $D_i$ its semisimple part. The operators $L^{(i)}_0$
preserve the finite-dimensional Cartan weight spaces of every
finitely generated submodule, so they are locally finite. Their
commuting generalized eigenspace decompositions give
$$
 [D_i,x_j(n)]=-n\delta_{ij}x_j(n),\qquad x_j\in\lev_j.
$$
Every homomorphism preserves these decompositions, and taking a
generalized eigenspace is exact on locally finite modules.
Consequently
$$
 d+\sum_iD_i+L^{\rm H}_0
$$
is a semisimple central operator. Its eigenvalue is fixed by the
chosen block and degree normalization. Distinct eigenvalues give
subcategories with no extensions between them. Choose the constants
$c_i,c_{\rm H}$ according to the highest-weight normalizations of
the separate affine blocks and the Fock representation. The operators
$$
 d_i=-D_i+c_i,\qquad d_{\rm H}=-L^{\rm H}_0+c_{\rm H}
$$
then satisfy
$$
 [d_i,x_j(n)]=n\delta_{ij}x_j(n),\qquad
 d=\sum_i d_i+d_{\rm H}.
$$
These constructions also apply to extensions, since the operators
remain locally finite and the fixed generalized eigenspaces are exact.
On a module with separate derivations, each $d_i+D_i$ is a
semisimple central operator with the scalar fixed by its block.
The construction therefore recovers the given $d_i$.
Adding and forgetting the normalized separate derivations are mutually
inverse exact functors. Thus the category of Levi modules with the
common derivation is identified with the specified block for the
direct sum of the Lie algebras $\widehat{\lev}_i$ with their
separate degree derivations, tensored with
the charge-$a$ Fock representation. We denote this
category by $\mathcal C_\lev(t,a)$.

For a tuple $\mu=(\mu_i)_i$ of highest weights in the chosen affine
blocks, let $M_i(\mu_i)$ be the ordinary Verma module for
$\widehat{\lev}_i$ of highest weight $\mu_i$, and let $\mathcal F_a$
be the charge-$a$ Fock module for the Heisenberg Lie algebra associated
with $Z(\lev)$. We write
$$
 M_\lev(\mu)=\biggl(\bigotimes_i M_i(\mu_i)\biggr)\otimes\mathcal F_a.
$$
The levels $\ell_i$, the charge $a$, and the degree normalizations are
fixed as above and suppressed from this notation.
Tensor products of the projective covers are small projective
generators of $\mathcal C_\lev(t,a)$. Indeed, successive tensor--Hom
adjunctions reduce projectivity to exactness of
$\operatorname{Hom}(P_i,-)$ in each affine block. Since $P_i$ is
finitely generated, these Hom spaces remain in the bounded categories
for the other Lie algebras $\widehat{\lev}_j$. The same adjunctions give generation
and smallness. Tensoring the finite Verma filtrations gives finite
filtrations by the modules $M_\lev(\mu)$. Every Levi object therefore
has a presentation by direct sums of these projectives.

We use the principal grading of $\widehat\g$, in which every simple
root of $\widehat\g$ has degree one. It is positive on the positive
Levi roots, and
the specified categories have a common principal upper bound.
For a positive finite root $\beta$ in $\uu$, the root
$\beta+j\delta$ has positive principal degree for $j\geq0$ and
negative principal degree for $j<0$. In particular, the enveloping
algebra $U(\uu[z])$ has finite-dimensional homogeneous components.
No uniform bound on the dimensions of module weight spaces is imposed.

\subsection{Compatibility of twisting functors with the block equivalence}
\label{com:twisting-section}
To compare relative Wakimoto modules at $t_0$ and $t$ as direct limits
of twisted Verma modules, we need the block equivalence to commute with
twisting functors on modules and homomorphisms.
The lemma uses Arakawa's semiregular twisting and adjunction
\cite[Lemma~4.1]{ABGG}, his treatment of roots simple only in the
integral root system \cite[Corollaries~5.6, 5.9 and Lemma~5.8]{ABGG},
and ordinary and parabolic BGG reciprocity
\cite[Theorem~2.2]{Fiebig}, \cite[Theorem~3.4]{Liu2022}.

\begin{lemma}\label{com:twisting-compatibility}
Let $t=p/q>0$, where $p,q$ are coprime positive integers, let
$e=\gcd(q,r)$, and set $t_0=p/e$. Let
$\mathcal E:\mathcal O_{\Lambda_0}\longrightarrow\mathcal O_\Lambda$
be the labelled exact equivalence between the corresponding large
positive blocks, with the degree normalizations specified above.
Let $\iota:W(\Lambda_0)\longrightarrow W(\Lambda)$ be the isomorphism
of integral Weyl groups induced by \eqref{com:root-dilation}.
For $w\in W(\Lambda_0)$, choose a reduced expression in integral
simple reflections. Let $T_w^{(t_0)}$ and
$T_{\iota(w)}^{(t)}$ be the corresponding products of twisting
functors. For a root simple only in the integral root system,
use the conjugation construction for twisting functors in
\cite[Corollary~5.6 and Lemma~5.8]{ABGG}.
There is a natural isomorphism
$$
 \mathcal E\circ T_w^{(t_0)}
       \cong T_{\iota(w)}^{(t)}\circ\mathcal E.
$$
The assertion includes singular highest weights and is natural
with respect to all module homomorphisms in these categories.
\end{lemma}
The proof is largely based on known results, notably Fiebig's equivalence
of noncritical blocks \cite[Theorem~4.1]{Fiebig}. The description of twisting
functors on projective modules by maximal integrable quotients is known
in the finite-dimensional semisimple setting
\cite[Lemma~2.1 and Proposition~2.4(1)]{MazorchukStroppel2007}.
See also Abe's construction in \cite{Abe2012}.
Since this compatibility is essential for comparing parabolic Wakimoto
modules at different levels, we give a full proof in
Appendix~\ref{app:twisting-compatibility}, with the details needed for the
affine blocks and singular highest weights considered here.

\subsection{Semi-infinite induction of Verma modules for the affine Levi factors}
We now describe the semi-infinite induction of $M_\lev(\mu)$
as an inductive limit of twisted Verma modules for $\widehat\g$.
For ordinary Wakimoto modules, such a description is given in
\cite[Proposition~4.3]{ABGG}. The following lemma establishes the
corresponding statement for induction from $\lev$. We include
the proof to specify the transition maps and to explain why
the construction also applies to singular highest weights.

Let $F_t=\Find$ denote Feigin--Frenkel induction with the
bilinear forms on $\lev$ specified in \eqref{com:target}.
On underlying vector spaces, this functor tensors the inducing
module with the semiregular module. It is therefore exact and
commutes with algebraic direct sums. For the loop algebra of
the nilradical $\uu$, write
$$
 \mathfrak n=\uu[z,z^{-1}],\qquad
 \mathfrak n_+=\uu[z],\qquad
 \mathfrak n_-=z^{-1}\uu[z^{-1}],\qquad
 \mathcal D=U(\mathfrak n_+)^*,
$$
where the star denotes the restricted dual with respect to the
principal grading. For the highest-weight label $\mu$ of
$M_\lev(\mu)$, let $\widetilde\mu$ denote the corresponding
$\widehat\g$-highest weight determined by these induction conventions.

To construct the inductive system, choose an integral coweight
$\eta$ that vanishes on the roots of $[\lev,\lev]$ and is
positive on the roots of $\uu$. After multiplying $\eta$ by a
sufficiently divisible positive integer, we may assume that it
belongs both to the coroot lattice of $\g$ and to the integral
translation lattice at $t_0$. At level $t_0$, set
$w_n=t_{-n\eta}$. At level $t$, replace $\eta$ by $D\eta$
throughout the construction, where $D=q/e$, so that
$w_n=t_{-nD\eta}$. Below, $w_n$ denotes the translation sequence
at the level under consideration.

\begin{lemma}\label{com:relative-limit}
For every label in the specified Levi block there is an isomorphism
of affine modules
\begin{equation}\label{com:relative-limit-isom}
 F_tM_{\lev}(\mu)\cong
       \varinjlim_n M^{w_n}(\widetilde\mu).
\end{equation}
The transition maps are the semiregular maps normalized on the dual
unit. Their Hom spaces are one dimensional, and all these maps are
nonzero. The assertion includes singular and nonadmissible labels.
\end{lemma}
\begin{proof}
The semiregular realization gives
$$
 F_tM_{\lev}(\mu)=
 U(\mathfrak n_-)\otimes\mathcal D\otimes M_{\lev}(\mu)
$$
as a vector space. The operators $x(n)$ with $x(n)$ in the negative Levi subalgebra
act triangularly in the degree of
$M_\lev(\mu)$, with diagonal action equal to its given Levi action.
In particular its character equals that of the ordinary
$\widehat\g$-Verma module $M(\widetilde\mu)$,
and all its weight spaces are finite dimensional.

Define the subalgebras
$$
 \mathfrak b_n=\bigoplus_{\beta\in\Delta(\uu)}
       \bigoplus_{0\le j<n\beta(\eta)}\g_\beta z^j,
 \qquad
 \mathfrak c_n=\bigoplus_{\beta\in\Delta(\uu)}
       \bigoplus_{j\ge n\beta(\eta)}\g_\beta z^j.
$$
Their closure under brackets follows by adding the two thresholds.
Choose a Poincar\'e--Birkhoff--Witt basis with the elements of
$\mathfrak c_n$ preceding those of $\mathfrak b_n$. The quotient by the
right ideal $\mathfrak c_nU(\mathfrak n_+)$ is then $U(\mathfrak b_n)$
as a right $\mathfrak b_n$-module. Dualizing yields a
$\mathfrak b_n$-linear embedding
$$
 j_n:U(\mathfrak b_n)^*\longrightarrow\mathcal D.
$$
For $m\ge n$ the embeddings are compatible. In each fixed principal
degree the tail $\mathfrak c_n$ eventually has no vectors, so
$\bigcup_n\im j_n=\mathcal D$. The quotient is a right-module quotient, so it requires no
commutativity assumption on $\uu$.

The inversion subalgebra of $w_n$, after translation, is precisely
$\mathfrak b_n$. The module $F_tM_\lev(\mu)$ is cofree over this subalgebra, and
the twisting adjunction of \cite[Lemma~4.1]{ABGG} gives
$$
 \Hom\bigl(M^{w_n}(\widetilde\mu),F_tM_{\lev}(\mu)\bigr)
 \cong
 \Hom\bigl(M(w_n^{-1}\!\cdot\widetilde\mu),
                    G_{w_n}F_tM_{\lev}(\mu)\bigr).
$$
Here $G_{w_n}$ is the right adjoint of $T_{w_n}$. In its
semiregular Hom realization, coinduction removes the restricted-dual
factors indexed by the inversion set. The dual unit contributes
$w_n^{-1}\widehat\rho-\widehat\rho$. The remaining factors give
$$
 \ch G_{w_n}F_tM_{\lev}(\mu)
       =\ch M(w_n^{-1}\!\cdot\widetilde\mu).
$$
Since $\eta$ vanishes on the derived-Levi roots, translation
preserves every root factor from the Levi module. Thus only the stated
nilradical factors are removed. The identity follows factor by factor from
$w_n^{-1}(\widehat\Delta_+\setminus\mathcal I_n)
=\widehat\Delta_+\setminus(-w_n^{-1}\mathcal I_n)$,
where $\mathcal I_n=\widehat\Delta_+\cap w_n\widehat\Delta_-$.
Writing $v_\mu$ for the Levi highest vector, that highest-weight line
in the Hom realization is spanned by
$$
 \xi_n:U(\mathfrak b_n)^*\longrightarrow
 U(\mathfrak n_-)\otimes\mathcal D\otimes M_\lev(\mu),
 \qquad \phi\longmapsto1\otimes j_n(\phi)\otimes v_\mu,
$$
with the Weyl twist just specified. Its weight is
$w_n^{-1}\!\cdot\widetilde\mu$. The character identity shows that
all positive root operators annihilate it. The Verma universal
property and adjunction give the required homomorphism.

For $m\ge n$, set $u=w_n^{-1}w_m$. Translation lengths add.
In the semiregular bimodules the map is $s\mapsto s\otimes1_u^*$,
where $1_u^*$ evaluates the constant coefficient. Consequently
\cite[\S4.1, equation~(23), and \S4.4]{ABGG} gives
$$
 \Hom\bigl(M^{w_n}(\widetilde\mu),M^{w_m}(\widetilde\mu)\bigr)
 =\C\nu_{m,n},\qquad
 \nu_{m,n}(s\otimes v_{w_n^{-1}\cdot\widetilde\mu})
 =(s\otimes1_u^*)\otimes v_{w_m^{-1}\cdot\widetilde\mu}.
$$
On restricted duals these maps are the embeddings above, so
$\nu_{l,m}\nu_{m,n}=\nu_{l,n}$. The dual unit is nonzero for every
label, including singular labels.

The resulting map from the colimit to $F_tM_{\lev}(\mu)$ is
surjective: the embedded duals exhaust $\mathcal D$, the action of $U(\mathfrak n_-)$ generates the free factor, and
the negative Levi subalgebra generates $M_\lev(\mu)$ by induction
in its Verma degree.
Every finite stage has the same finite weight multiplicities as
the module $F_tM_\lev(\mu)$. Its colimit has no larger weight spaces, so the
surjection is an isomorphism. No scalar depending on $\mu$ was
inverted. This proves the assertion also on singular walls.
\end{proof}

Root dilation identifies the entire specified set of Levi labels.
Lemmas~\ref{com:twisting-compatibility} and~\ref{com:relative-limit}
therefore give
$$
 \mathcal E\bigl(F_{t_0}M_{\lev}(\mu)\bigr)
       \cong F_tM_{\lev}(\mu').
$$
Indeed the images of transition arrows differ from the corresponding arrows
only by nonzero scalars. Rescale the stage identifications recursively
to make the sequential diagrams commute. Equivalences preserve these
algebraic colimits. Thus $\mathcal E$ identifies these semi-infinite
inductions and their defining direct systems.

\subsection{Semi-infinite induction for Levi modules with fixed \texorpdfstring{$h(0)$}{h(0)}-eigenvalue}
Fix $t$ and write $F=F_t$. Normalize $h\in Z(\lev)$ so that
its eigenvalues on $\uu$ are positive integers.
The next lemma gives the categorical property needed to pass from
the modules $F_tM_\lev(\mu)$ to the inductions of the simple quotients
of $M_\lev(\mu)$.
Throughout this subsection the inducing modules have the same
eigenvalue $a$ of $h(0)$, where $Z(\lev)=\C h$. The operator $h(0)$
need not act by one scalar on their semi-infinite inductions.

\begin{lemma}\label{com:induction-fully-faithful}
On the category of Levi modules specified in
Section~\ref{com:Levi-categories}, with $h(0)=a\,\id$ for a fixed $a$,
$F$ induces isomorphisms on $\operatorname{Hom}$ and
$\operatorname{Ext}^1$. Its essential image is the smallest full
subcategory containing the modules $FM_\lev(\mu)$ and closed under
extensions, algebraic direct sums and cokernels of morphisms between
its objects.
\end{lemma}
\begin{proof}
Let $S$ be the semiregular module of $\mathfrak n$ and use the
same embedded dual $\mathcal D$. The multiplication maps
in \cite[\S2.4, Eqs.~(5)--(6), and Lemma~2.4]{ABGG}, with
$\mathfrak n_-$ and $\mathfrak n_+$ as the two subalgebras, give
compatible left-module identifications
\begin{equation}\label{com:compatible-semiregular}
 S=\operatorname{Ind}_{\mathfrak n_+}^{\mathfrak n}\mathcal D,
 \qquad S\cong U(\mathfrak n_-)\otimes\mathcal D.
\end{equation}
The semi-infinite character of $\mathfrak n$ is zero. In the coinduced
realization
$\operatorname{Hom}_{U(\mathfrak n_-)}(U(\mathfrak n),U(\mathfrak n_-))$,
the element $\varphi\in\mathcal D$ is the map
$$
 \widetilde\varphi(v_-v_+)=v_-\varphi(v_+),
 \qquad v_\pm\in U(\mathfrak n_\pm).
$$
This embedded dual is left $\mathfrak n_+$-stable. The multiplication
map $u\otimes\varphi\mapsto u\cdot\varphi$, evaluated on
$v\in U(\mathfrak n_+)$, is $\widetilde\varphi(vu)$.
Choose root-vector Poincar\'e--Birkhoff--Witt bases of
$U(\mathfrak n_+)$ and $U(\mathfrak n_-)$.
Commuting the factors of $v$ in $\mathfrak n_+$ past the factors
of $u$ in $\mathfrak n_-$ gives $u\varphi(v)$ and terms containing brackets.
A bracket lying in $\mathfrak n_-$ replaces one negative-mode factor
by one of strictly larger $h(0)$-eigenvalue. A bracket lying in
$\mathfrak n_+$ decreases the number of negative-mode factors.
Further Poincar\'e--Birkhoff--Witt reordering among negative-mode factors
can only decrease their number. Hence every correction term has either
smaller negative-mode length or the same length and strictly larger
total $h(0)$-eigenvalue.

Let $c_{\max}$ be the largest eigenvalue of $\operatorname{ad}h$ on
$\uu$. For negative-mode length at most $L$, the total
$h(0)$-eigenvalue lies between zero and $Lc_{\max}$. Ordering these
terms by the pair consisting of their length and the negative of
their total $h(0)$-eigenvalue gives a finite filtration. In the
vector-space identification with $U(\mathfrak n_-)\otimes\mathcal D$,
the multiplication map is $1+R$, where $R$ strictly lowers this
filtration. It is therefore invertible, with inverse
$\sum_{j\ge0}(-R)^j$, a finite sum on every vector. This proves the
left-module identification using the specified embedded dual.
The maps preserve principal degree. The fixed-degree subspaces of
$U(\mathfrak n_-)\otimes\mathcal D$ are finite dimensional, so
these identifications take place in the algebraic restricted-dual
realization and require no completion.
The left coinduced realization in the cited Eq.~(6), restricted to
$\mathfrak n_+$, is cofree with multiplicity space $U(\mathfrak n_-)$.
These descriptions use the same embedded dual, so induction of a
map from that dual is compatible with the cofree structure.

For a Levi module $M$, the restriction of $F(M)$ to
$\mathfrak n$ is $S\otimes M$, with $\mathfrak n$ acting on $S$.
See \cite[Section~7.1, equation~(50)]{ABGG}. In particular, $F$ is
exact and preserves algebraic direct sums.

The $\mathfrak n_+$-invariants of $S\otimes N$ have $h(0)$-eigenvalues at least $a$, and their eigenvalue-$a$ subspace is the
copy of $N$. All generators of $U(\mathfrak n_-)$ have positive
$h(0)$-eigenvalue, while nonconstant elements of $\mathcal D$
have negative charge. A $\mathfrak n$-intertwiner of negative
charge kills this charge-$a$ invariant subspace. Induction on the magnitude of
the negative charge in $\mathcal D$ then kills its entire embedded
dual: applying $\mathfrak n_+$ reduces that magnitude, and a
remaining image would be an invariant of charge below $a$.
Freeness over $U(\mathfrak n_-)$ kills the whole map. The same
argument, after subtracting its restriction to that subspace, shows that a
charge-zero intertwiner is exactly $\id_S\otimes f$.
Compatibility with the Levi action makes $f$ Levi-linear. Such a
map is a homomorphism of $\widehat\g$-modules by functoriality of induction. This
proves full faithfulness on Hom spaces.

Consider an extension of $\widehat\g$-modules $FM$ by $FN$. We first split
it over $\mathfrak n$. Set $B=U(\mathfrak n_-)\otimes N$ and
let $b_0$ be its principal upper bound. Bounded graded coinduction
from $B$ has degree-$d$ space
$$
 \bigoplus_{0\le j\le b_0-d}
       \Hom_\C\bigl(U(\mathfrak n_+)_j,B_{d+j}\bigr).
$$
This finite sum equals the restricted-dual tensor realization,
even when the weight multiplicities of $B$ are infinite.
Coinduction is right adjoint to the exact forgetful functor to
graded Cartan weight spaces, a semisimple category. The cofree
restriction of $FN$ is consequently injective in this category.
Pull back the extension to $\mathcal D\otimes M$ and choose a
weight-preserving $\mathfrak n_+$-linear section. Inducing this
section by \eqref{com:compatible-semiregular} gives a
$\mathfrak n$-linear section of the original extension. Its
composite with the quotient map is the identity because it is
the identity on the embedded dual.

For $y\in\lev[z,z^{-1}]$, the Feigin--Frenkel action on
$S\otimes M$ is $A_y\otimes1+1\otimes y$, where the operator
$A_y$ on $S$ is independent of $M$ and includes the prescribed
change in the invariant form. In the chosen splitting, the
off-diagonal action of $y$ is a charge-zero $\mathfrak n$-homomorphism. It is therefore $\id_S\otimes f_y$.
The Lie relations for these operators are exactly the extension
cocycle relations for the Levi algebra. They define an extension
$0\to N\to R\to M\to0$ of Levi modules with the common derivation.
Since the splitting preserves the full Cartan, it also preserves
the levels and scalar $h(0)$-action. Its weights have the common upper bound,
and its simple subquotients belong to the specified Levi block.
The exact equivalence constructed in \cref{com:Levi-categories}
equips $R$ with the normalized separate affine derivations and
preserves the inclusion and quotient maps. Hence this is an extension
in $\mathcal C_\lev(t,a)$.

The actions on the given extension and on $FR$ now agree on
the subalgebra $\mathfrak n+\lev[z,z^{-1}]$. For $y(n)$ with $y$ in the opposite nilradical and $[h,y]=-cy$,
compare its two actions. Commuting
the difference with an element of $\mathfrak n$ of positive $h(0)$-eigenvalue gives
either a nilradical or Levi difference, already zero, or a difference for an element of the opposite nilradical of smaller $c$. Induction on $c$
therefore makes the difference a negative-charge nilradical
intertwiner. Such an intertwiner is zero by the preceding argument.
Only finitely many positive $h$-weights occur in the finite-dimensional
nilradical, so this induction covers every loop mode. The $\widehat\g$-actions consequently agree. Every such extension is induced, and Hom full faithfulness
identifies equivalences and split extensions. This proves the Ext
assertion.

Full faithfulness and exactness imply that the image is closed under
cokernels, and the extension argument and direct-sum realization give
the other two closures. Conversely every Levi projective is a finite
extension of modules $M_\lev(\mu)$. Every Levi object has a presentation
by direct sums of these projectives. Applying the exact functor $F$
shows that every object of its image lies in the stated closure.
\end{proof}

\begin{proof}[Proof of Theorem~\ref{thm:denominator}]
The description of the essential image uses only the modules
$F_tM_\lev(\mu)$ and their $\widehat\g$-module homomorphisms. The labelled equivalence
$\mathcal E$ preserves these modules by
Lemmas~\ref{com:twisting-compatibility} and~\ref{com:relative-limit}.
It preserves all their morphisms,
extensions, sums and cokernels because it is an exact equivalence.
Lemma~\ref{com:induction-fully-faithful} therefore identifies the
two induced categories. The unique simple quotient, in this category, of the induction of
$M_\lev(\mu)$ with finite highest weight zero and charge $a=0$ is the induction
of its simple quotient. The latter is the tensor product of the simple
affine vertex algebras $L_{\ell_i}(\lev_i)$ and the Fock module
for the Heisenberg Lie algebra of $Z(\lev)$ on which its zero modes act as zero. Its induction is $T_h(t)$.
Consequently
$$
 \mathcal E(T_h(t_0))\cong T_h(t).
$$
Here simplicity in the induced category identifies the target
$T_h(t)$. The target need not be simple as a $\widehat\g$-module.
We now identify the map into this target.

The space $\operatorname{Hom}(V_t,T_h(t))$ is one dimensional. A vector
of weight $k\Lambda_0$ has degree zero, and in that degree the
only possible additional vectors are coordinate polynomials. Every
nonconstant coordinate monomial has strictly negative $h(0)$-eigenvalue. The kernel of $h(0)$ in degree zero is therefore
$\C\mathbf1$. The equivalence sends $V_{t_0}$ to $V_t$, so
$$
 \mathcal E(\Phi_h(t_0))=c_h\Phi_h(t),\qquad c_h\ne0.
$$
Since $\mathcal E$ is exact, it identifies the images of the
corresponding homomorphisms. Since it is an equivalence, it also
preserves simplicity as affine modules. For these vacuum quotients, affine submodules are vertex-algebra
ideals: the current fields generate the vertex algebra, and the
Sugawara formula makes every affine submodule stable under translation.
Thus simplicity of the image is also simplicity as a vertex algebra.
The criterion for factorization through $L_k(\g)$ and
Corollary~\ref{com:finite-centers} now prove Theorem~\ref{thm:denominator}.
\end{proof}

\begin{remark}
The argument uses maximal Levi subalgebras, whose centers are one dimensional.
The theorem reduces types $B,C,F_4$ to denominators one or two,
and type $G_2$ to denominators one or three. It preserves the condition
$X_{L_k(\g)}\subset\N(\g)$. Nilpotent orbit labels, Drinfeld--Sokolov reductions, and conformal gradings
can depend on the denominator. The theorem assumes $k+h^\vee\in\Q_{>0}$.
\end{remark}

\begin{example}\label{ex:sl4-transport}
Consider $\g=\slg_4$, $h=\operatorname{diag}(1,1,1,-3)$, and
$k_q=-4+3/q$, where $q$ is any positive integer coprime to $3$.
We start with the known equality
$X_{L_{-1}(\slg_4)}=\overline{G\cdot(\C^\times h)}$
\cite[Theorem~1.1]{ArakawaMoreau2017}. This settles $q=1$ independently.
For the Levi subalgebra $\slg_3\oplus\C h$, the parabolic Wakimoto
homomorphism is
$$
 \Phi_q:V^{k_q}(\slg_4)\longrightarrow T_q,
 \qquad
 T_q=L_{-3+3/q}(\slg_3)\otimes\mathcal H_{36/q}
       \otimes\mathcal S(3).
$$
Here $\mathcal S(3)$ is the vertex algebra of three $\beta\gamma$
pairs. The rank-one Heisenberg vertex algebra $\mathcal H_{36/q}$
is generated by $H$ with
$[H_\lambda H]=(36/q)\lambda\mathbf1$. This follows from
$(h,h)=\tr(h^2)=12$ and $k_q+4=3/q$.
Write $B_q=\operatorname{im}\Phi_q$.

For $q=1$, the Levi factor is $L_0(\slg_3)=\C$.
The known point $h\in X_{L_{-1}(\slg_4)}$ and
Theorem~\ref{com:point-map} give $B_1\cong L_{-1}(\slg_4)$.
Thus $B_1=\operatorname{im}\Phi_1$ is a simple vertex algebra.
Section~\ref{com:transport} compares $\Phi_1$ with $\Phi_q$.
Its exact equivalence $\mathcal E$ sends $V^{-1}(\slg_4)$ to
$V^{k_q}(\slg_4)$ and $T_1$ to $T_q$, as
$\widehat{\slg_4}$-modules.
The latter identification is established in
Lemmas~\ref{com:relative-limit} and~\ref{com:induction-fully-faithful},
including the passage from $F_tM_\lev(\mu)$ to the inductions of the
simple quotients of $M_\lev(\mu)$.
The relevant space of homomorphisms is one dimensional, so
$\mathcal E(\Phi_1)=c_q\Phi_q$ for some $c_q\ne0$, after these
identifications. Exactness gives $B_q\cong\mathcal E(B_1)$.
An equivalence preserves simple modules, so $B_q$ is simple.
Consequently $\Phi_q$ annihilates the whole maximal ideal of
$V^{k_q}(\slg_4)$ and induces an embedding
$$
 L_{-4+3/q}(\slg_4)\hookrightarrow
 L_{-3+3/q}(\slg_3)\otimes\mathcal H_{36/q}\otimes\mathcal S(3).
$$

The embedding now gives a semisimple point at the new level.
It induces a homomorphism
$R_{L_{k_q}(\slg_4)}\to R_{T_q}$.
Set the six $\beta\gamma$ coordinates to zero, evaluate the
$\slg_3$ factor at the origin of its associated variety, and set
the Heisenberg coordinate $[H]$ equal to $12$.
The Wakimoto formulas then send $x\in\slg_4$ to
$(x,h)=\tr(xh)$. This is evaluation at $h$ under the identification
$\slg_4^*\cong\slg_4$.
Every polynomial relation in $R_{L_{k_q}(\slg_4)}$ must therefore
vanish at $h$, because the homomorphism is defined on the simple
quotient. Hence $h\in X_{L_{k_q}(\slg_4)}$.
Since $h$ is nonzero and semisimple, $L_{k_q}(\slg_4)$ is not
quasi-lisse. This is the reverse implication in
Theorem~\ref{com:point-map}, made explicit in this example.

At $q=1$, the semisimple point was known in advance. At the other
levels, it is obtained only after the comparison proves that
$\Phi_q$ factors through the simple algebra.
The Wakimoto realization is thus part of the proof of
non-quasi-lisseness at those levels. Preserving simple modules alone
would not give this conclusion: the particular homomorphism
$\Phi_q$ produces the semisimple point.

The individual cases $q=1$ and $q=2$ already have infinite families of
irreducible ordinary modules, by \cite[Theorem~6.1]{APFusion2014}
and \cite[Theorem~1.2(1)]{APVsl42023}, respectively.
The comparison above treats all positive $q$ coprime to $3$ uniformly.
\end{example}

\section{Drinfeld--Sokolov reduction and the Jantzen filtration}
\label{com:reduction}
We recall the exactness and associated-variety theorems for
Drinfeld--Sokolov reduction, then justify their application to completed
Jantzen character identities so that character coefficients can detect
relations in the reduction.

Let $f$ be nilpotent, and choose its dominant Dynkin element $x$, so
that $[x,f]=-f$. Write $H_f^i$ for quantized Drinfeld--Sokolov
reduction and $\mathcal W^k(\g,f)=H_f^0(V^k(\g))$. 
We use the BRST construction of quantum Hamiltonian reduction in \cite{KRW2003,KacWakimoto2004}.
The conformal degree after reduction is the eigenvalue of the operator
induced by $L_0-x(0)$ and the ghost conformal operator.
Here $\mathrm{KL}_k$ is the graded category of smooth
$\widehat\g$-modules at level $k$ on which the operators
$a(0)$, $a\in\g$, act locally finitely, with semisimple
degree operator, finite-dimensional degree spaces, and degree
support in a finite union of lower-bounded integral translates.
See \cite{AAV}, \S2.6. BRST complexes use its enlargement by algebraic filtered colimits
and direct sums, which need not have finite-dimensional pieces.

The vanishing and exactness assertions below are
\cite[Theorem~4.5.7]{AAV}, and the associated-variety statement and the
nonvanishing criterion are \cite[Theorem~4.6.1 and
Proposition~4.6.2]{AAV}. We also record the formal extension to
algebraic colimits that will be used below.

\begin{theorem}\label{com:ds}
For $M\in\mathrm{KL}_k$, one has $H_f^i(M)=0$ for $i\ne0$,
and $H_f^0$ is exact. In the algebraic BRST enlargement it
commutes with filtered unions and algebraic direct sums. In
particular higher-cohomology vanishing extends to such colimits
of objects of $\mathrm{KL}_k$. For a finitely generated vertex-algebra quotient $Q$ of $V^k(\g)$ in this category,
$X_{H_f^0(Q)}=X_Q\cap\mathcal S_f$. In particular, for the simple
affine vertex algebra $L_k(\g)$,
\begin{equation}\label{com:ds-support}
 H_f^0(L_k(\g))\ne0\quad\Longleftrightarrow\quad
 f\in X_{L_k(\g)}.
\end{equation}
When nonzero, its associated variety is the reduced intersection
$X_{L_k(\g)}\cap\mathcal S_f$, where $\mathcal S_f$ is the
Slodowy slice.
\end{theorem}
For the colimit assertion, each term of the algebraic BRST complex is a tensor product
with a fixed ghost space, and filtered colimits and direct sums of
complex vector spaces are exact. 

\subsection{Characters of Drinfeld--Sokolov reductions}
The completion used below is determined by the Sugawara grading. Assume $t=k+h^\vee>0$. Every module in the
following argument has its Sugawara grading,
nonnegative integral eigenvalues, finite-dimensional graded pieces,
and locally finite $\g$-action through the operators $a(0)$. For a dominant weight $\lambda$
set
\begin{equation}\label{com:energy-functions}
 C_\lambda=(\lambda,\lambda+2\rho),\qquad
 \Delta_t(\lambda)=\frac{C_\lambda}{2t},\qquad
 \delta_x(\lambda)=\Delta_t(\lambda)-\lambda(x).
\end{equation}
Let $\mathbb V_k(\lambda)$ denote the generalized Verma
$\widehat\g$-module induced from
the finite-dimensional $\g$-module $E_\lambda$, with
$\g\otimes z\C[z]$ acting trivially. Its top has Sugawara degree
$\Delta_t(\lambda)$. In the character formulas, $q$ denotes the formal
grading variable, which is independent of the denominator used to write
the fixed level $t=k+h^\vee$. The standard BRST character calculation gives
\begin{equation}\label{com:standard-ds-character}
 \ch H_f^0(\mathbb V_k(\lambda))
 =q^{\Delta_t(\lambda)}\ch E_\lambda(q^{-x})
       \ch\mathcal W^k(\g,f),
\end{equation}
with least reduced degree $\delta_x(\lambda)$, see
\cite{AAV}, \S4.3. 

We use the following construction throughout this subsection. For a
nonzero module $M$, let $d$ be its least Sugawara degree. The finite-dimensional $\g$-module $M_d$ is semisimple and is annihilated by $\g\otimes z\C[z]$.
Each constituent $E_\lambda$ therefore satisfies
$d=\Delta_t(\lambda)$. Let $P$ be the finite direct sum of the
corresponding generalized Verma modules, with their multiplicities in $M_d$, and let
$F$ be the image of the resulting map $P\to M$. Set
$K=\ker(P\to F)$ and $Q=M/F$. The map is an isomorphism in degree
$d$, so both $K$ and $Q$ begin strictly above $d$. The two exact
sequences give
\begin{equation}\label{com:degree-decomposition}
 [M]=[P]-[K]+[Q].
\end{equation}
Iterating the construction on the successive quotients $Q$ gives an
exhaustive filtration by quotients of finite sums of generalized Verma modules.

\begin{lemma}\label{com:tail}
For every real $B$ there is an integer $N_B$ such that, for every
module $M$ as above,
$$
 M_d=0\ (d\le N_B)\quad\Longrightarrow\quad
 H_f^0(M)_a=0\ (a\le B).
$$
The bound depends only on $\g,t,x,B$, not on $M$.
\end{lemma}
\begin{proof}
The identity
\begin{equation}\label{com:coercivity}
 2t\delta_x(\lambda)
 =\|\lambda+\rho-tx\|^2-\|\rho-tx\|^2
\end{equation}
shows that only finitely many dominant integral weights have
$\delta_x(\lambda)\le B$. Choose $N_B\ge0$ bounding all values
$\Delta_t(\lambda)$ that are nonnegative integers for these weights.

If $M$ starts above $N_B$, every generalized Verma module used in the preceding
filtration has $\delta_x(\lambda)>B$. Equation
\eqref{com:standard-ds-character}, exactness of $H_f^0$, and its
commutation with the exhaustive union prove the assertion.
\end{proof}

\begin{proposition}
\label{com:completion}
Let a character identity in the preceding Sugawara-graded category
be expressed as a locally finite, possibly signed sum of generalized Verma
characters. Then it may be reduced by
\eqref{com:standard-ds-character}, provided the expansion in generalized Verma characters
is obtained in the degree completion from exact sequences.
In particular this holds for a locally finite affine Jantzen
character identity for $V^k(\g)$.
\end{proposition}
\begin{proof}
Induction on the lowest conformal degree in
\eqref{com:degree-decomposition} expresses $[M]$, modulo characters
vanishing in degrees at most $N$, as a finite integral linear combination
of generalized Verma characters. Finiteness follows because every step strictly
increases the lowest integral degree and involves only finitely many
irreducible $\g$-summands. Keeping the full $\g$-characters makes the coefficients unique: at the first nonzero degree, the irreducible
$\g$-characters are linearly independent.

Apply exact reduction to these finitely many exact sequences.
For a bound $B$ on the conformal degree after reduction, choose $N\ge N_B$. All remainders then
vanish through degree $B$ by \cref{com:tail}. Thus the finite reduced
calculation stabilizes at every reduced degree, giving the asserted
completed identity. For a Jantzen filtration $J^1\supseteq J^2
\supseteq\cdots$, the module $\bigoplus_{a\ge1}J^a$ has
finite-dimensional homogeneous subspaces for the Sugawara grading by degreewise finiteness
of the filtration. The same argument applies to its character.
\end{proof}

\begin{lemma}
\label{vtwo:finite-filtration}
Let $M$ have the nonnegative integral Sugawara grading and finite
$\g$-modules in each degree used in \cref{com:tail}. For an integer $N$
there is a finite filtration with factors of two kinds:
\begin{enumerate}
\item simple $\widehat\g$-modules $L_\mu$ with
$\Delta_t(\mu)\le N$,
\item modules with all original pieces through $N$ zero.
\end{enumerate}
The multiplicities of the first kind are intrinsic and agree with the
simple-character multiplicities in that finite degree range. If $N$ is
a uniform bound as in \cref{com:tail} for reduced degree $B$, a positive multiplicity of
$L_\mu$ with $H_f^0(L_\mu)_B\ne0$ implies $H_f^0(M)_B\ne0$.
\end{lemma}
\begin{proof}
Induct on $b_N(M)=\sum_{d\le N}\dim M_d$. If it is zero, $M_d=0$ for $d\le N$. Otherwise apply the construction above to one irreducible
$\g$-submodule $E_\mu$ in the least degree $d\le N$.
Its generated submodule $F\subset M$ is a quotient of
$\mathbb V_k(\mu)$ in which the entire top survives, hence has simple head
$L_\mu$. The kernel of $F\twoheadrightarrow L_\mu$ starts above
$d$. Both this kernel and $M/F$ have strictly smaller $b_N$.
Their inductively constructed filtrations, together with $L_\mu$,
give the required filtration of $M$.

Characters of irreducible $\g$-modules at the least remaining degree
determine the simple multiplicities recursively. Exact
$H_f^0(-)_B$ vanishes on every factor of the second kind by \cref{com:tail}, and its
finite-dimensional values on the remaining factors add with
nonnegative multiplicities. This proves the final assertion.
No finite length of $M$ as a whole, filtration by generalized Verma modules, or
identification of Jantzen layers with radical layers is assumed.
The filtration consists of affine subquotients, not truncations
obtained by cutting off the higher graded pieces.
\end{proof}

All these character identities hold equivariantly for a reductive or
finite subgroup of the triple centralizer. The BRST differential, finite
exact sequences and generalized Verma character calculation are equivariant. The
bound for the omitted terms is uniform over the finite-dimensional
representations occurring in these character decompositions.
Taking invariants is exact in each finite reduced degree. In particular
an equivariant completed Jantzen sum computes the sum of the characters
of all the nested ideals $H_f^0(J^a)$, rather than only that of $H_f^0(J^1)$.

Let $A(\nu)=\sum_{w\in W(\g)}\det(w)e^{w\nu}$ and
$\chi(\mu)=A(\mu+\rho)/A(\rho)$. Denote
$$
 P_\g(Q)=\prod_{m\ge1}(1-Q^m)^{-\rank\g}
                 \prod_{\alpha\in\Delta}(1-e^\alpha Q^m)^{-1}.
$$
The following formula is the Jantzen sum formula obtained from the
Gorelik--Kac determinant \cite{GK}, written using finite Weyl
characters. We include the conversion to fix the root normalization
and the conformal degrees used in subsequent calculations.

\begin{proposition}\label{com:vacuum-jantzen-root-sum}
For a simple Lie algebra $\g$ and $t=k+h^\vee\in\Q_{>0}$,
the Jantzen filtration of $V^k(\g)$ obtained by deforming $k$ satisfies
\begin{equation}\label{com:vacuum-jantzen-formula}
 \sum_{r\ge1}\ch J^r
 =P_\g(Q)\sum_{\substack{a\ge1,\ \alpha\in\Delta\\
 b=2((\rho,\alpha)+at)/(\alpha,\alpha)\in\Z_{>0}}}
                   Q^{ab}\chi(-b\alpha).
\end{equation}
For a regular alternant its term is
$\varepsilon\ch\mathbb V_k(\mu)$, where
$\mu+\rho$ is the dominant translate of $\rho-b\alpha$ and
\begin{equation}\label{com:vacuum-jantzen-weight}
 \Delta_t(\mu)=ab.
\end{equation}
Singular alternants are zero.
\end{proposition}
\begin{proof}
To obtain this formula from
\cite[Theorem~1.2.8, Example~1.2.12, Sections~2.2.3 and~6.1]{GK}, the
finite denominator is replaced by its Weyl alternant. For a real root
$\alpha+a\delta$, the possible vanishing factor is proportional to
$r(at+(w\rho,\alpha)-r(\alpha,\alpha)/2)$. It vanishes when
$r=2(at+(w\rho,\alpha))/(\alpha,\alpha)$.
Its derivative is $ra\ne0$, so it contributes once to the determinant
valuation. Imaginary factors are nonzero because $t>0$.
Changing $\alpha$ to $w\alpha$ gives
\eqref{com:vacuum-jantzen-formula}. Weyl invariance of the norm gives
$C_\mu=b^2(\alpha,\alpha)-2b(\rho,\alpha)=2tab$,
proving \eqref{com:vacuum-jantzen-weight}.
\end{proof}

\begin{proposition}\label{com:vacuum-test}
Suppose $J^1=I_k$ is the maximal proper graded ideal of $V^k(\g)$ and
$J^1\supseteq J^2\supseteq\cdots$ is its locally finite Jantzen
filtration. Set
$$
 K_a=H_f^0(J^a)\subseteq\mathcal W^k(\g,f),\qquad
 R_f=\sum_{a\ge1}\dim(K_a)_0.
$$
Then $R_f$ is a finite nonnegative integer and
$$
 H_f^0(L_k(\g))\ne0\quad\Longleftrightarrow\quad R_f=0.
$$
\end{proposition}
\begin{proof}
The $K_a$ are nested ideals by exactness. The Dynkin
grading of the universal $\mathcal W$-algebra is nonnegative and
has degree-zero space $\C\mathbf1$. Hence $(K_a)_0\ne0$ exactly
when $K_a=\mathcal W^k$. If any $K_a$ has this property, so does
$K_1$. Conversely $K_1=\mathcal W^k$ contributes one to the sum.
The reduction of $0\to I_k\to V^k(\g)\to L_k(\g)\to0$
identifies its cokernel with $H_f^0(L_k(\g))$. For every $N$,
one has $(J^a)_d=0$ for $d\le N$ and all sufficiently large $a$. \Cref{com:tail}, with $B=0$,
makes all sufficiently deep $(K_a)_0$ zero.
\end{proof}

If a complete parabolic Jantzen expansion is
$\sum_\lambda c_\lambda\ch\mathbb V_k(\lambda)$ and all its
terms satisfy $\delta_x(\lambda)\ge0$, the preceding results give
the finite formula
\begin{equation}\label{com:equality-vacuum-sum}
 R_f=\sum_{\delta_x(\lambda)=0}
               c_\lambda\dim E_\lambda[\lambda(x)].
\end{equation}
If negative $\delta_x$ occur, descendants in the universal
$\mathcal W$ factor of \eqref{com:standard-ds-character} must also
be included. Neither a partial sum of root contributions nor a formal
identification of the total Jantzen sum with $J^1$ is sufficient.

\subsection{A criterion for a nilpotent element to belong to the associated variety}
The conformal-weight inequality gives a direct criterion for a nilpotent
element to belong to the associated variety. When equality occurs, the
criterion involves the entire $\g$-isotypic component of the ideal.

\begin{lemma}\label{com:energy}
Let $I$ be a proper graded ideal of $V^k(\g)$, with $t>0$.
If $f\notin X_{V^k/I}$, there are a nonzero dominant root-lattice
weight $\lambda$, a positive integer $D=\Delta_t(\lambda)$, and
$v\in I_D[\lambda]$ such that
$$
 \sigma(v)(f)\ne0,\qquad v\in I_D[\lambda][x=D],\qquad
                       \delta_x(\lambda)\le0.
$$
\end{lemma}
\begin{proof}
Choose the least $D$ at which an ideal symbol evaluates nontrivially
at $f$, then an irreducible $\g$-submodule $E_\lambda$ with
this property. For dual bases $u_a,u^a$ of $\g$, Sugawara gives
\begin{equation}\label{com:ideal-sugawara}
 (2tD-C_\lambda)\sigma(v)
       =2\sum_a u_a\,\sigma(u^a(1)v).
\end{equation}
Indeed all terms $u_a(-r)u^a(r)v$ with $r\ge2$ vanish modulo
$C_2$. The remaining positive-mode vectors belong to $I_{D-1}$,
so the right side evaluates to zero at $f$. Consequently
$D=\Delta_t(\lambda)$. The trivial $\g$-module would give
$D=0$ and is impossible in a proper ideal. Equivariance and
homogeneity, using $[x,f]=-f$, show that evaluation at $f$ is zero
on every $x$-eigenspace except the eigenvalue-$D$ space. Dominance
of $x$ bounds its eigenvalues by $\lambda(x)$, proving the claim.
\end{proof}

\begin{corollary}\label{com:equality}
To test whether $f\in X_{V^k/I}$, it suffices to evaluate the
full spaces
\begin{equation}\label{com:finite-energy-spaces}
 I_{\Delta_t(\lambda)}[\lambda][x=\Delta_t(\lambda)],
 \quad
 \lambda\in Q\cap P_+\setminus\{0\},\quad
 \Delta_t(\lambda)\in\Z_{>0},\quad\delta_x(\lambda)\le0.
\end{equation}
There are finitely many such $\lambda$. If every nonzero dominant
root-lattice weight with $\Delta_t(\lambda)\in\Z_{>0}$ has
strictly positive $\delta_x$, then
$f$ belongs to the associated variety. If all have nonnegative
$\delta_x$, only the equality spaces in
\eqref{com:finite-energy-spaces} remain.
\end{corollary}
\begin{proof}
Finiteness follows from \eqref{com:coercivity}. For instance
$\delta_x(\lambda)\le0$ implies
$\|\lambda\|\le2\|\rho-tx\|$. A nonzero evaluation in a listed
space implies $f\notin X_{V^k/I}$. Conversely, if $f\notin X_{V^k/I}$,
\cref{com:energy} guarantees the existence of such a space. For equality
$\Delta_t(\lambda)=\lambda(x)$ one must test every copy of
$E_\lambda$ and its entire maximal-$x$ eigenspace. This eigenspace
need not be one dimensional when $x$ is singular.
\end{proof}

A useful sufficient condition is that every equivariant map
$\Sym^D\g\to E_\lambda^*$ vanishes at $f$ in the relevant
eigenspace: it then applies to every copy in the actual ideal.
Nonstrict conformal-weight inequalities alone do not suffice. For
$\g=\slg_2$, $k=0$, $t=2$, and principal $f$, the types
$\lambda=m\alpha$ have $\delta_x(\lambda)=m(m-1)/2\ge0$.
For $m=1$, the symbols of $\g(-1)\mathbf1\subset I_0$ do not all
vanish at $f$, where $I_0$ is the maximal proper graded ideal of $V^0(\slg_2)$. Indeed $L_0(\slg_2)=\C$ and its principal
reduction vanishes.

\subsection{Conformal-weight inequalities in classical types}
\begin{lemma}\label{v3:master-energy}
For $B_n,D_n$ set $t=T$ and $\sigma=1,0$, respectively. For $C_n$
set $t=T/2$ and $\sigma=2$. Let $\eta=2x$ denote the nonnegative
grading eigenvalues in decreasing order, with $\eta_n=0$ in the type-$D$ applications.
Set $A_0=0$ and
$$
 A_l=l(2n+\sigma-l)-T\sum_{i\le l}\eta_i.
$$
For a dominant root-lattice weight with coordinates
$a_1\ge\cdots\ge a_n\ge0$, taking the absolute value of the last
type-$D$ coordinate, and $a_{n+1}=0$, one has
$$
 2T\delta_x(\lambda)=\sum_i a_i(a_i-1)
              +\sum_l(a_l-a_{l+1})A_l.
$$
On a consecutive block of $L$ equal grading eigenvalues starting after $b$,
$$
 A_{b+u}=\frac{L-u}{L}A_b+\frac uL A_{b+L}+u(L-u).
$$
Thus nonnegative endpoint values imply $\delta_x\ge0$. Its nonzero
equality weights are precisely columns $(1^l)$ at endpoints $A_l=0$,
subject to the root-lattice parity: even $l$ in $C,D$, unrestricted
$l$ in $B$. Both final signs are allowed in type $D$.
\end{lemma}
\begin{proof}
Insert the standard epsilon-coordinate Weyl vectors into
$C_\lambda/(2t)-\lambda(x)$ and telescope the linear terms using
$a_i=\sum_{l\ge i}(a_l-a_{l+1})$. The root lattices make all $a_i$
integers. The second identity follows by subtracting the affine
interpolation from the quadratic polynomial in $u$. Interior values
are strictly positive. Equality forces every $a_i\in\{0,1\}$,
so dominance gives exactly the stated columns. The condition
$\eta_n=0$ makes the last sign irrelevant. Nonvanishing in an equality
case still requires the ideal calculation in \cref{com:equality}.
\end{proof}

\begin{lemma}
\label{vtwo:current-descendant}
Let $\g$ be simple, $k\ne0$, and
$\Hom_\g(\g,\Sym^2\g)=0$. No proper graded ideal of $V^k(\g)$
contains an adjoint summand in Sugawara degree two. Consequently, if
\cref{com:equality} leaves only the adjoint representation in conformal weight two, that equality
case cannot make the reduction of the quotient of $V^k(\g)$ zero.
\end{lemma}
\begin{proof}
The central term $a(1)b(-1)\mathbf1=k(a,b)\mathbf1$ gives $I_1=0$.
The Poincar\'e--Birkhoff--Witt theorem for $U(\widehat\g)$ gives
the following isomorphism of $\g$-modules:
$$
 V^k(\g)_2=\g(-2)\mathbf1\oplus\Sym^2(\g(-1)).
$$
The only adjoint copy is therefore $\g(-2)\mathbf1$. If it lay in
$I_2$, then $a(1)b(-2)\mathbf1=[a,b](-1)\mathbf1$ would place all
nonzero vectors of $\g(-1)\mathbf1$ in $I_1$, a contradiction. For $\g=\sog_N$, $N\ge7$,
the hypothesis follows from the curvature decomposition
$$
 \Sym^2(\Lambda^2W)=\Lambda^4W\oplus\mathcal R_0(W)
                  \oplus\Sym^2_0W\oplus\C:
$$
none of these is $\Lambda^2W$. At $N=8$ the exterior summand splits
into two middle-degree summands, neither adjoint. At $N=7$ its Hodge
dual is $\Lambda^3W$, also distinct. At $N=6$ the exterior summand
\emph{is} adjoint, so the assertion is not applied to $D_3=A_3$.
At $k=0$, the ideal of $V^0(\g)$ generated by $\g(-1)\mathbf1$
provides a counterexample if one omits the
level hypothesis. The nonnegative conformal-weight differences and the complete list of equality cases remain separate premises of the corollary.
\end{proof}

\subsection{Comparison of nilpotent orbits at different levels}
We compare the nilpotent orbits in the associated varieties at
levels $k$ and $k+\ell$, where $\ell\in\Z_{\ge0}$.
The proof uses the diagonal $\widehat\g$-action on
$L_k(\g)\otimes L_\ell(\g)$ and the tensor-product reduction theorem
of Arakawa--Creutzig--Feigin
\cite[Theorem~1]{ArakawaCreutzigFeigin2022}.
\begin{proposition}
\label{com:translation}
Let $\ell\in\Z_{\ge0}$ and suppose $k$ and $k+\ell$ are
noncritical. Then
$$
 X_{L_{k+\ell}(\g)}\cap\N(\g)
       \subseteq X_{L_k(\g)}\cap\N(\g).
$$
\end{proposition}
\begin{proof}
The vectors $x(-1)\mathbf1\otimes\mathbf1+\mathbf1\otimes x(-1)\mathbf1$,
$x\in\g$, generate a vertex subalgebra $D$ of
$L_k(\g)\otimes L_\ell(\g)$, and $D$ surjects onto
$L_{k+\ell}(\g)$. These objects and the whole tensor product lie
in $\mathrm{KL}_{k+\ell}$: the homogeneous subspaces for the tensor-product grading are
finite-dimensional $\g$-modules and the operators $x(n)$, $x\in\g$, $n>0$, lower
degree. If $f\notin X_{L_k}$, then $H_f^0(L_k)=0$.
The tensor-product reduction theorem
(\cite[Theorem~1]{ArakawaCreutzigFeigin2022}) gives
$$
 H_f^0(L_k\otimes L_\ell)
       \cong H_f^0(L_k)\otimes L_\ell=0.
$$
Exactness first on the embedding of $D$, then on its simple
quotient, gives $H_f^0(L_{k+\ell})=0$.
Now use \eqref{com:ds-support}. 
\end{proof}

\section{Simple quotients and tensor product decompositions}
\label{com:current-factors}
\label{com:images}
We first give the elementary graded form of
\cite[Proposition~4.5.12]{AAV}, relating simple quotients of
$\mathcal W^k(\g,f)$ to $H_f^0(L_k(\g))$, and its consequence for
quasi-lisseness, with a proof for completeness. We then give criteria
for affine and Heisenberg vertex subalgebras to be tensor factors
and for equality with the affine vertex subalgebra.

\begin{lemma}\label{com:simple-image}
Let $U$ be a nonnegatively graded conformal vertex algebra with
$U_0=\C\mathbf1$. Suppose $H\ne0$ is a graded quotient of $U$
and there is a graded unital surjection $U\twoheadrightarrow B$ onto
a simple vertex algebra $B$. Then
$H$ surjects onto $B$. In particular, if
$U=\mathcal W^k(\g,f)$ and $H=H_f^0(L_k(\g))\ne0$, a
non-quasi-lisse simple quotient $B$ of $U$ proves that $L_k(\g)$ is not
quasi-lisse.
\end{lemma}
\begin{proof}
The sum of all proper graded ideals has zero degree-zero part, so
$U$ has a unique maximal graded ideal $M$. Both quotients are graded,
and the kernel defining $H$ lies in $M=\ker(U\to B)$.
Thus $H\twoheadrightarrow B$.

If $L_k$ were quasi-lisse, \cref{com:ds} would identify $X_H$ with
$X_{L_k}\cap\mathcal S_f$, whose finitely many symplectic leaves
are the connected components of its nilpotent-orbit intersections.
The closed Poisson subvariety $X_B$ would be a union of those leaves,
contradicting the hypothesis.
\end{proof}

When $B$ is realized in another vertex algebra, its grading is induced
from $U$. The containing vertex algebra may have a different conformal
vector. It is enough that $B$ is not lisse when
$X_{L_k(\g)}\subset\N(\g)$ would force
$X_{L_k(\g)}\cap\mathcal S_f$ to be zero-dimensional, as for principal
$f$, or when $G\cdot f$ is maximal among the nilpotent orbits allowed
in $X_{L_k(\g)}$.

For the two even gradings below, Arakawa's almost irreducibility
theorem \cite[Theorem~5.5.4 and Corollary~5.5.5]{Arakawa2011}
gives simplicity. We explain its application here, since the
identification of the two BRST constructions and the one-dimensional
finite $W$-algebra module are needed here.

\begin{lemma}\label{com:cominuscule-simple}
Let $\g=B_n$, $f=(3,1^{2n-2})$, or let $\g=C_n$, $f=(2^n)$.
For $k\ne-h^\vee$, every nonzero reduction
$H_f^0(L_k(\g))$ is simple.
\end{lemma}
\begin{proof}
The Dynkin grading has $\g=\g_{-1}\oplus\g_0\oplus\g_1$,
with $\g_{\pm1}$ abelian. In the two cases the deleted simple
root has index $i=1,n$, respectively, and $\alpha_i(x)=1$.
The isomorphism $y_0=w_0t_{-x}$ of
\cite[Section~4.3]{Arakawa2011} sends a root operator $e_\alpha(m)$
to a nonzero scalar multiple of $e_{-\alpha}(m+\alpha(x))$.
Its Cartan correction is $h(0)\mapsto-h(0)+k(x,h)$.
Thus the vector $\mathbf1\in L_k(\g)$ is a highest-weight vector for the
twisted $\widehat\g$-module: a positive Levi root with index zero
has image of index zero, and a positive root in $\g_1$ has image
of index one. For every root and $m\ge1$, the image has index
at least zero and annihilates $\mathbf1$. The twisted simple
$\widehat\g$-module has level $k$ and finite highest weight
$k\omega_i$, zero on the simple coroots of
$[\g_0,\g_0]$. It consequently belongs to the category
$\mathcal O_{0,k}$ of that reference.

For the charged ghosts of roots in $\g_1$, the same isomorphism
has the form $\psi_\alpha(r)\mapsto\psi_{-\alpha}(r+1)$.
The inverse shift is by $-1$. These shifts interchange the two
annihilator conditions $r\ge0$ and $r\ge1$ and preserve the
ghost vacuum. Since $\g_1$ is abelian, the cohomology differential
has no cubic ghost term. Its terms involving the action of
$\widehat\g$ and the character defined by $f$ are preserved by the
isomorphism. This gives the identification of complexes in
Section~4.3 for the quotient $L_k(\g)$.

Write $E$ for the simple highest-weight $\g$-module of weight
$k\omega_i$. It is a quotient of the scalar generalized Verma
$\g$-module $M=U(\g)\otimes_{U(\g_{\ge0})}\C_{k\omega_i}$.
As a $U(\g_{-1})$-module, $M$ is its regular module, so for the
Whittaker character $\chi$ its coinvariant space is
$$
 M/\langle(a-\chi(a))M:a\in\g_{-1}\rangle\simeq\C.
$$
By \cite[Theorem~5.5.4(ii)]{Arakawa2011}, nonvanishing of the
reduction $H_f^0(L_k(\g))$ implies that $E$ has Gelfand--Kirillov dimension
$\dim\g_{-1}$, the largest possible value in $\mathcal O_0(\g)$.
The exact finite Whittaker functor of
\cite[Theorem~5.1.1]{Arakawa2011} therefore gives a nonzero
quotient of this one-dimensional space. It is simple over the finite
$W$-algebra. Proposition~5.4.2 and Corollary~5.5.5 of the same
reference identify it with the initial space of the reduction
and prove that the reduction is simple.
\end{proof}

Nilpotence of the conformal class is
\cite[Proposition~3.3]{AK}. We recall its short geometric argument
and record the Heisenberg realizations to which it will be applied.

\begin{lemma}
\label{com:conformal-obstruction}
Let $V$ have a nonnegative conformal grading, $V_0=\C\mathbf1$,
and a finitely generated $C_2$ algebra. If $V$ is quasi-lisse,
every positive-degree Poisson-central element of $R_V$ is
nilpotent. In particular its conformal class $[\omega]$ is
nilpotent. Let $\mathcal H$ be the rank-$r$ Heisenberg vertex
algebra generated by $h_1,\ldots,h_r$, with
$[h_i{}_{\lambda}h_j]=K_{ij}\lambda\mathbf1$ and $K$ invertible.
Suppose a graded vertex-algebra homomorphism $V\to\mathcal H$ sends $\omega$ to
$$
 \frac12\sum_{i,j}(K^{-1})_{ij}:h_ih_j:
                     +\sum_i b_i\partial h_i,
 \qquad r>0.
$$
Then $[\omega]$ is not nilpotent and $V$ is not quasi-lisse.
The same conclusion holds for a homomorphism into
$\mathcal H\otimes C$, where $R_C\ne0$, whose conformal-vector image is
the displayed vector plus $1\otimes\omega_C$.
\end{lemma}
\begin{proof}
A central regular function is constant on every symplectic leaf.
Finitely many leaves give finitely many values. Positive
homogeneity and the conical action force its value to be zero
everywhere. The Nullstellensatz then makes its class nilpotent.
The conformal class is central since $\omega_{(0)}$ is the
translation operator, whose image is zero in $R_V$.

The map on $C_2$ quotients sends $[\omega]$ to the nonzero
quadratic polynomial $\frac12h^TK^{-1}h$ in
$R_{\mathcal H}=\C[h_1,\ldots,h_r]$. Every positive power is
nonzero. In the tensor-factor version, take the highest degree
in $h_1,\ldots,h_r$. Its leading term is the same
nonzero power, tensored with $\mathbf1$ in the other factor.
\end{proof}

\subsection{Invariant bilinear forms on conformal weight-one subspaces}

We use invariant bilinear forms as in \cite{LiForms1994}. Gradings
are discrete and bounded below, with finite-dimensional pieces.
Forms are homogeneous and normalized by $(\mathbf1,\mathbf1)=1$.
On a simple conformal vertex algebra they are nondegenerate because
their radicals are ideals.

The following consequence of the invariant-form theory of
\cite{LiForms1994} specifies which pairings must vanish.

\begin{lemma}\label{com:zero-current}
Let $S$ be a simple nonnegatively graded conformal vertex algebra
with $S_0=\C\mathbf1$ and a normalized invariant form. If
$j\in S_1$ satisfies $j_{(1)}v=0$ for every $v\in S_1$, then
$j=0$. Consequently a semisimple Lie ideal in $S_1$ whose
level is zero vanishes. For an abelian Lie subalgebra of $S_1$, the same conclusion requires
its pairing with all of $S_1$ to vanish.
\end{lemma}
\begin{proof}
Invariance gives $(L_1v,\mathbf1)=(v,L_{-1}\mathbf1)=0$.
Thus $L_1S_1=0$. The invariant pairing on $S_1$ is, up to the
fixed sign convention, the scalar $j_{(1)}v$. Its restriction is
nondegenerate because distinct conformal weights are orthogonal.
For a semisimple Lie ideal $\mathfrak a\subseteq S_1$, invariance
of the weight-one pairing and $[\mathfrak a,\mathfrak a]
=\mathfrak a$ show that its mixed pairing with a commuting
complement is zero. The asserted vanishing follows.
\end{proof}

The full pairing and form hypotheses are essential. Two isotropic generators of a Heisenberg vertex algebra
may have nonzero mutual pairing, while the conformal
vector $T_{\mathrm{std}}+Q\partial J$ gives
$L_1J=-2Q\kappa\mathbf1$ and need not admit a normalized form.
For the half-integral Dynkin quotients, \cref{seed:integral-form}
gives the ordinary form on the integer-weight fixed subalgebra.
The half-integral convention is as follows.

\begin{lemma}\label{v5:half-pairing}
Let $V$ be a simple conformal vertex algebra with finite-dimensional graded pieces,
$V=\bigoplus_{d\in\frac12\Z_{\ge0}}V_d$, $V_0=\C\mathbf1$ and
$L_1V_1=0$. Its normalized contragredient bilinear form is nondegenerate
on every weight space. For a quasiprimary $u\in V_d$, its pairing
with $v\in V_d$ is $e^{\mathrm i\pi d}$ times the coefficient
of $\mathbf1$ in $u_{(2d-1)}v$. Consequently, if $V_d$ is quasiprimary and its
top pairing is identically zero, then $V_d=0$.
\end{lemma}
\begin{proof}
Use the opposed field
$Y(e^{zL_1}(e^{\mathrm i\pi}z^{-2})^{L_0}u,z^{-1})$.
The half-integral contragredient construction identifies invariant
forms with $(V_0/L_1V_1)^*$, see \cite[\S2.1, Remarks 2.5, 2.8
and the proof of Proposition 2.9]{BeemGarner2026}. The degree-zero functional taking value one on $\mathbf1$ gives such
a form. Its radical is the kernel of
a $V$-module map $V\to V'$, hence zero by simplicity. The self-adjoint
$L_0$ makes unequal weights orthogonal. The displayed top coefficient
follows by $L_1u=0$. The form need not be symmetric: for even states
of weight $d$ it has symmetry $(-1)^{2d}$, allowing the alternating
pairing in weight $3/2$. No positivity or graded-unitarity hypothesis
is used.
\end{proof}

If the image of a semisimple Lie algebra $\mathfrak a$ in $V_1$ is zero,
then $\mathfrak a$ acts trivially on $V$. Every nontrivial irreducible
$\mathfrak a$-module of generators therefore has zero image.
The generators invariant under $\mathfrak a$ require separate relations.

\subsection{Tensor product decompositions with affine vertex algebras and their commutants}
For the rank-one Heisenberg vertex algebra generated by $J$, the
Fock-module decomposition when $J(0)=0$ respects the vertex-algebra
structure. This is the zero-charge case of the oscillator decomposition
in \cite[Section~3.1, (3.3) and Proposition~3.2]{CKLRHeisenberg2019}.
We give the argument without assuming that $V$ is simple or $J$ is
primary for its conformal vector.

\begin{lemma}\label{com:heisenberg}
Let $V_0=\C\mathbf1$, let $J\in V_1$ with
$J_{(1)}J=\kappa\mathbf1$, $\kappa\ne0$, and suppose $J_0=0$
on all of $V$. Then multiplication is an isomorphism
$$
       \mathcal H_\kappa\otimes\Com(\mathcal H_\kappa,V)
                  \xrightarrow{\ \sim\ } V.
$$
If $V$ is simple, its commutant is simple. No primarity of $J$
or invariant form on $V$ is needed. The same assertion holds for
several mutually commuting elements of $V_1$ with a nondegenerate Gram matrix and all zero
modes acting trivially.
\end{lemma}
\begin{proof}
On every finite graded piece the oscillator projector
\begin{equation}\label{com:oscillator-projector}
 \Pi=\prod_{r>0}\left(
   \sum_{a\ge0}\frac{(-1)^aJ_{-r}^{a}J_r^{a}}
                         {a!(r\kappa)^a}\right)
\end{equation}
is a finite expression. The oscillator commutator gives
$J_r\Pi=0$ for $r>0$, and $\Pi$ is the identity on their common
kernel. Taylor expansion in the creation operators $J_{-r}$
gives the direct vector-space decomposition
$$
 V=\C[J_{-1},J_{-2},\ldots]\otimes\bigcap_{r>0}\ker J_r.
$$
Since $J_0=0$, the second factor is the commutant. Its fields
commute with $Y(J,z)$, so the decomposition respects
all vertex products. A proper ideal of that factor would give a
proper ideal of $V$. For a nondegenerate Gram matrix, change to
dual oscillators and repeat the same calculation.
\end{proof}

For a vertex algebra $U$ with unique maximal graded ideal $I$, a nonzero
graded image $B$ on which $J_0=0$ gives
$J_0U\subseteq\ker(U\to B)\subseteq I$.
If $J_{(1)}J=\kappa\mathbf1\ne0$, \cref{com:heisenberg} therefore
applies to $U/I$. Since $J_0$ is a derivation, the condition $J_0=0$ need only
be checked on strong generators.

\begin{proposition}\label{com:neutral-split}
Suppose $V$ is simple and nonnegatively graded. Let $\g_i$ be
simple Lie algebras and let
$\bigotimes_i V^{\ell_i}(\g_i)\to V$ be a graded homomorphism.
Assume that $V$ is strongly generated by the images of
$\g_i(-1)\mathbf1$ and positive-weight vectors annihilated by
$x(0)$ for all $x\in\bigoplus_i\g_i$, and that
$\ell_i+h_i^\vee\ne0$. The latter vectors can be modified by adding expressions in
lower-weight generators so that the resulting vectors belong to the commutant
of the image of $\bigotimes_i V^{\ell_i}(\g_i)$. If $A$ denotes
the vertex algebra generated by these corrected vectors, multiplication
gives an isomorphism
\begin{equation}\label{com:neutral-tensor}
           V\cong\bigotimes_i L_{\ell_i}(\g_i)\otimes A,
\end{equation}
and $A$ is simple. The assertion also allows Heisenberg vertex
algebras with nondegenerate bilinear forms among the tensor factors,
provided their zero modes annihilate every additional strong generator.
If the images of $\g_i(-1)\mathbf1$ and the Heisenberg generators
are primary, the conformal vector of $A$ is the difference between
that of $V$ and the sum of their Sugawara conformal vectors.
\end{proposition}
\begin{proof}
Let $E$ be the sum of the Sugawara operators $L_0^{\g_i}$. Induct on the
conformal weights of the generators annihilated by all zero modes. After correcting smaller
weights, the vertex algebra they generate together with the images of
$\g_i(-1)\mathbf1$ is a homomorphic image of
$$
  \bigotimes_i V^{\ell_i}(\g_i)\otimes A_{\mathrm{old}}.
$$
Call it $U_{\mathrm{old}}$. Its $E$-spectrum is nonnegative integral,
and $E$ kills $A_{\mathrm{old}}$. For the next such generator $w$,
all $J_nw$ with $n>0$ lie in $U_{\mathrm{old}}$. The finite Casimir
term vanishes, giving
$$
 Ew=\sum_i\frac1{\ell_i+h_i^\vee}
       \sum_{n>0,a}J^{u_a}_{-n}J^{u^a}_n w.
$$
The relation $[E,J_{-n}]=nJ_{-n}$ shows that $Ew$ lies in the positive
$E$-spectrum of $U_{\mathrm{old}}$. Hence its finite spectral inverse
$u=E^{-1}Ew$ has the same conformal weight and is annihilated by all zero modes.
Then $v=w-u$ satisfies $Ev=0$, while $J_nv\in U_{\mathrm{old}}$
has $E$-eigenvalue $-n$ and is zero for $n>0$. Zero modes also kill
$v$. Induction gives a surjection
$\bigotimes_i V^{\ell_i}(\g_i)\otimes A\to V$.

Let $I_{\ell_i}(\g_i)$ be the maximal proper graded ideal of
$V^{\ell_i}(\g_i)$. The image of
$\bigotimes_{j<i}V^{\ell_j}(\g_j)\otimes I_{\ell_i}(\g_i)
\otimes\bigotimes_{j>i}V^{\ell_j}(\g_j)\otimes A$
is an ideal of $V$ concentrated in positive $E$-degrees, and is
therefore zero by simplicity. The map
therefore factors through \eqref{com:neutral-tensor}. To prove
injectivity, set $B=\bigotimes_iL_{\ell_i}(\g_i)$. Its adjoint mode
module is simple with scalar endomorphisms. For independent $b_j$,
density gives mode operators taking one $b_j$ to $\mathbf1$ and
the others to zero. Applied to a kernel element
$\sum_jb_j\otimes a_j$, these force every $a_j=0$, since $A\subset V$.
An ideal of $A$ would tensor to an ideal of $V$, proving simplicity.
For the Heisenberg vertex algebras, use \cref{com:heisenberg}.
Subtracting their Sugawara vectors and those of $V^{\ell_i}(\g_i)$
from the conformal vector of $V$ proves the conformal assertion.
\end{proof}

Since $R_{B\otimes A}=R_B\otimes R_A$, a tensor factor isomorphic to a nonzero-rank Heisenberg vertex algebra
provides a free central polynomial coordinate and obstructs
quasi-lisseness. An inclusion of a Heisenberg vertex algebra alone does not suffice:
integrable simple affine vertex algebras contain such vertex
subalgebras and are lisse.

\begin{lemma}
\label{vtwo:whole-factor}
Let $U,A$ be nonnegatively integer-graded conformal vertex algebras with
finite-dimensional graded pieces and degree-zero space $\C$, and let $A$ be simple. Suppose
$B\subseteq U\otimes A$ is a graded subalgebra containing
$\mathbf1\otimes A$ and the total conformal vector. Then
$$
 B=C\otimes A,\qquad C=\{u\in U:u\otimes\mathbf1\in B\},
$$
and its graded simple quotient is $C_{\mathrm s}\otimes A$. If
$[\omega_A]$ is not nilpotent in $R_A$, the total conformal class in
this quotient is not nilpotent.
\end{lemma}
\begin{proof}
Write $b=\sum_i u_i\otimes a_i$ with independent $a_i$. The mode
algebra of the simple countable adjoint $A$-module has scalar commutant.
By density, there is a finite mode combination taking one $a_i$ to
$\mathbf1$ and the others to zero. Its action preserves $B$, so every
$u_i\otimes\mathbf1$ lies in $B$. Conversely it generates
$u_i\otimes A$, proving the tensor statement. Subtracting
$\mathbf1\otimes\omega_A$ shows that $\omega_U$ lies in $C$.
The unique graded simple quotient of $C$ tensored with $A$ is simple,
again by density, and is consequently the unique graded simple quotient
of $B$. Finally
$$
 R_{C_{\mathrm s}\otimes A}=R_{C_{\mathrm s}}\otimes R_A
       \longrightarrow R_A
$$
using the positive-grading augmentation of the first factor sends the
total conformal class to $[\omega_A]$. Every power remains nonzero.

\end{proof}

\subsection{Quotients of Heisenberg commutants and comparison of graded dimensions}
We next treat quotients and specializations of realizations. In
\cref{com:finite-image}, finitely many graded dimensions suffice to
prove an isomorphism under the stated ideal-generation hypothesis.

\begin{lemma}
\label{com:heisenberg-quotient}
Let $V$ be a nonnegatively graded vertex algebra and let $J\in V_1$
satisfy $[J_\lambda J]=\kappa\lambda\mathbf1$, with $\kappa\ne0$.
Write $\mathcal H_\kappa$ for the vertex algebra generated by $J$,
and assume that $J(0)$ acts semisimply on $V$ with integer eigenvalues. For every charge-preserving graded
quotient $V\twoheadrightarrow\overline V$, the induced map
$$
 \Com(\mathcal H_\kappa,V)\longrightarrow
       \Com(\mathcal H_\kappa,\overline V)
$$
is surjective. It remains surjective after taking invariants of
an involution sending $J$ to $-J$. The same assertion holds for a
family over a parameter ring in which every homogeneous
$J(0)$-eigenspace is a finite free module, with fixed integer eigenvalues.
The commutants and their invariant subalgebras commute with
specialization wherever $\kappa$ is invertible.
\end{lemma}
\begin{proof}
Lift a commutant vector in charge zero and apply
\eqref{com:oscillator-projector}. This produces a commutant lift.
For the involution, average the lift. The projector commutes with
the involution. Over the parameter ring localized at $\kappa$,
each homogeneous charge-zero piece has the same finite idempotent
$\Pi$. Its image is a direct summand and base change preserves
its image. The reflection average is another commuting
idempotent. 
\end{proof}

If moreover $V$ is simple and conformal of CFT type with normalized
invariant form and primary $J$, the commutant $C$ and its conformal reflection
fixed points are simple. Indeed $J_0$ is skew, so only opposite
charges pair. The charge-zero space and its oscillator decomposition
$\mathcal H_\kappa\otimes C$ have nondegenerate forms, as do $C$
and $C^{\Z/2}$ by reflection averaging. Invariance makes every proper
graded ideal orthogonal to the whole algebra, proving simplicity.
These form hypotheses concern only this assertion.

\begin{lemma}
\label{com:finite-image}
Let $U$ be a nonnegatively graded vertex algebra with unique maximal
graded ideal $I$, and suppose that $I$ is generated as an ideal by
homogeneous vectors of degrees at most $N$. Write $A=U/I$. Let $H\ne0$ be strongly generated in degrees at most
$N$, and let $\phi:U\to H$ be a graded unital map. Assume the
compared pieces $A_d,H_d$ are finite dimensional. If
$$
                  \dim A_d=\dim H_d\qquad(0\le d\le N)
$$
at every allowed graded degree, then $\phi$ is surjective and
$\ker\phi=I$. Hence $H\cong A$.
\end{lemma}
\begin{proof}
Set $K=\ker\phi\subseteq I$. In every degree $d\le N$,
$(U/K)_d$ surjects onto $A_d$ and embeds in $H_d$. Equality of
the outer dimensions forces $(U/K)_d=H_d$ and $K_d=I_d$.
The strong generators of $H$ are therefore in the image, so the
map is onto. The stated generators of $I$ lie in $K$, giving
$K=I$. 
\end{proof}

\begin{lemma}
\label{vtwo:small-subalgebra}
Let $U_0=\C$ with nonnegative grading, $I\subset U$ a proper graded
ideal, and $H=U/I$. Let $C\subseteq U$ contain the total conformal
vector and have a graded conformal surjection onto a simple $A$.
Suppose $H$ is generated by vectors in a finite set of degrees $S$ and
$\dim H_d=\dim A_d<\infty$ for $d\in S$. Then $H$ has a conformal
surjection onto $A$.
\end{lemma}
\begin{proof}
The unique maximal graded ideal of $C$ contains $I\cap C$. Hence the
image $C_H$ surjects onto $A$. In each generator degree,
$$
 \dim A_d\le\dim(C_H)_d\le\dim H_d=\dim A_d.
$$
Thus $C_H$ contains all generators of $H$ and is $H$. Its quotient map is conformal. An isomorphism requires additional
hypotheses, such as ideal generation in \cref{com:finite-image}.
\end{proof}

\subsection{Conformal embeddings, conformal levels and collapsing levels}

An embedding $\varphi:V\hookrightarrow V'$ of vertex operator algebras
is called \emph{conformal} if it preserves the conformal vectors,
that is, $\varphi(\omega_V)=\omega_{V'}$. We use the same convention
for conformal vertex algebras and vertex superalgebras.

Let $k$ be noncritical, and let $\mathcal W_k(\g,f)$ be the simple
quotient of the universal $W$-algebra $\mathcal W^k(\g,f)$ defined
by BRST cohomology. Let $A$ be its affine vertex subalgebra generated
by the canonical generators of conformal weight one. When the Sugawara
conformal vector $\omega_A$ is defined, the level $k$ is called
\emph{conformal} if the inclusion $A\hookrightarrow\mathcal W_k(\g,f)$
is conformal, that is, $\omega_A=\omega_{\mathcal W_k(\g,f)}$.
The level $k$ is called \emph{collapsing} if
$A=\mathcal W_k(\g,f)$. At a conformal level, $A$ may still be a proper
vertex subalgebra. These conventions agree with
\cite[Sections~3--4]{AdamovicEtAl2018}. They concern the simple
$W$-algebra $\mathcal W_k(\g,f)$. An identification involving
$H_f^0(L_k(\g))$ is asserted separately, since a nonzero reduction
need not be simple.

More generally, let $A\subseteq V$ be generated by a reductive Lie
algebra $\mathfrak a\subseteq V_1$, and let $S$ be its Sugawara
conformal vector. The next lemma gives a sufficient condition for
$S=\omega_V$ and then a condition under which $V=A$. Equality of
central charges alone is not sufficient here: the lemma also requires
a condition on the $\mathfrak a$-invariant subspace of $V_2/A_2$.

The comparison of Sugawara conformal weights used below is the same
principle as the criterion for conformal embeddings in
\cite[Theorem~2.2]{APConformal2013} and its use for minimal
$W$-algebras in \cite[Section~4]{AdamovicEtAl2018}. We state the form
needed for general strong generators and give its proof. Applications
to collapsing levels and nilpotent Slodowy slices are developed in
\cite{AKMPP2020,AVEMSlodowy2024}.

\begin{lemma}\label{v3:resonance}
Let a nonzero nonnegatively graded conformal vertex algebra $V$, with finite-dimensional
homogeneous subspaces and $V_0=\C$, be strongly generated by the primary images
of a reductive Lie algebra $\mathfrak a\subseteq V_1$ and finitely
many finite-dimensional $\mathfrak a$-modules. Assume the levels
$\ell_i$ on the simple ideals of $[\mathfrak a,\mathfrak a]$ are
noncritical, the bilinear form $K$ on $Z(\mathfrak a)$ is
nondegenerate, and the $\mathfrak a$-action is semisimple.
Let $A$ be the vertex subalgebra generated by $\mathfrak a$ and
$S$ its Sugawara vector. If $\omega_V=S$ and every generating
$\mathfrak a$-module outside $A$, of highest weight $\lambda_i$
on each simple ideal, central character $q$ and conformal weight $d$, satisfies
$$
 d\ne E(\lambda,q):=
 \sum_i\frac{(\lambda_i,\lambda_i+2\rho_i)}{2(\ell_i+h_i^\vee)}
                      +\frac12q^{\mathsf T}K^{-1}q,
$$
then $V=A$. For a simple $V$ with nondegenerate normalized invariant
form, the premise $\omega_V=S$ follows from $c(V)=c(S)$ if
$(V_2/A_2)^{\mathfrak a}$ is spanned by the image of $\omega_V$.
\end{lemma}
\begin{proof}
In the least generating weight not contained in $A$, the operators
$x(n)$, $x\in\mathfrak a$, $n>0$, take values in $A$. Sugawara applied modulo $A$ therefore acts by the displayed
Sugawara conformal weight $E(\lambda,q)$, while the conformal grading acts by $d$.
The inequality $d\ne E(\lambda,q)$ makes this quotient zero.
Induction over the discrete grading proves $V=A$.
For the second assertion, $D=\omega_V-S$ belongs to $\Com(A,V)$ and
is orthogonal to $A_2$. It pairs only with $\mathfrak a$-invariant vectors, hence
only with $\omega_V$ modulo $A_2$, and
$(D,\omega_V)=(c(V)-c(S))/2=0$. Nondegeneracy gives $D=0$.
For a simple $W$-algebra with half-integral conformal grading, this last
argument uses its integer-weight subalgebra and \cref{seed:integral-form}.
The first assertion
requires no simplicity. An identification of $A$ with a tensor product of simple affine vertex
algebras requires, in addition, integrability or simplicity.
\end{proof}

\section{\texorpdfstring{Miura homomorphisms and generators of $W$-algebras}{Miura homomorphisms and generators of W-algebras}}
\label{vtwo:image-tools}
The Miura homomorphism gives concrete realizations of
$\mathcal W^k(\g,f)$. The structure theorem of Kac--Wakimoto
\cite[Theorem~4.1]{KacWakimoto2004} gives a freely strongly generating
set $w(a_i)$ indexed by a homogeneous basis $a_i$ of $\g^f$.
If $[x,a_i]=j_i a_i$, then $w(a_i)$ has conformal weight $1-j_i$.
Thus derivatives and normally ordered products of these generators,
with factors in a fixed order, form a basis of the universal
$W$-algebra. Their images strongly generate a quotient, but need not
remain linearly independent. At the levels considered here,
we also compose the Miura homomorphism with the maps from the
universal affine vertex algebras associated with $[\g_0,\g_0]$
to their simple quotients. We determine generators of the resulting image and prove that they
generate it, including the symplectic-fermion invariant
subalgebras occurring in the applications.

\subsection{Images of weight-two generators under the Miura homomorphism}\label{com:quadratic-miura-section}
Suppose $\g=\bigoplus_{j\in\Z}\g_j$ is an integral good grading
for $f\in\g_{-1}$. Set $\mathfrak c=\g_0^f$ and
$\mathfrak n=\bigoplus_{j>0}\g_j$. Choose bases $a_i$ of $\g_0$
and $b_r$ of $\g_1$. For $u\in\g_{-1}^f$, define
$$
 A_{ri}=([f,a_i],b_r),\qquad [u,b_r]=\sum_i B_{ri}a_i,
 \qquad [a_i,b_s]=\sum_r c_{is}^{r}b_r,
$$
$$
 \kappa_s=k(u,b_s)+\tr_{\mathfrak n}(\ad u\,\ad b_s),
 \qquad C_s=-\frac12\sum_{i,r}B_{ri}c_{is}^{r}.
$$
The trace is well defined since $\ad u\,\ad b_s$ preserves each
graded summand. Goodness makes $A$ surjective. Choose a right inverse
$R$, with $AR=I$, and set
\begin{equation}\label{com:quadratic-miura-formula}
 M=RB+B^{\mathsf t}R^{\mathsf t}-RBA^{\mathsf t}R^{\mathsf t},
 \qquad d=-R(\kappa+C).
\end{equation}

\begin{proposition}\label{com:quadratic-miura}
The matrix $M$ is symmetric and
$$
 W_u=J^u+\frac12\sum_{i,j}M_{ij}:J^{a_i}J^{a_j}:
          +\sum_i d_i\partial J^{a_i}
$$
is a cocycle in the small BRST complex. Its complete Miura image is
$$
 \mu(W_u)=\frac12\sum_{i,j}M_{ij}:J^{a_i}J^{a_j}:
          +\sum_i d_i\partial J^{a_i}.
$$
These formulas hold at every level in the natural grading. Every
degree-two field is a sum of these fields, normally ordered quadratic
products of the weight-one fields $J^a$, $a\in\mathfrak c$, and their
first derivatives. For a fixed $u$, the difference of two choices belongs to $\Sym^2\mathfrak c\oplus\partial\mathfrak c$.
It vanishes under a vertex algebra homomorphism annihilating
$J^a$ for every $a\in\mathfrak c$.
\end{proposition}
\begin{proof}
Invariance and $[f,u]=0$ imply
$$
 (AB^{\mathsf t})_{rs}=([f,[u,b_s]],b_r)
 =([f,[u,b_r]],b_s)=(BA^{\mathsf t})_{rs}.
$$
Thus \eqref{com:quadratic-miura-formula} gives $M=M^{\mathsf t}$
and $AM=B$.

Write $\phi^r$ for the odd charged ghosts dual to $b_r$, with ghosts
ordered to the left. The small-complex differential
\cite[equation~(2.6)]{KacWakimoto2004} reads
$$
 d_{\rm BRST}J^{a_i}=\sum_r A_{ri}\phi^r,
 \qquad
 d_{\rm BRST}J^u=-\sum_{r,i}B_{ri}:\phi^rJ^{a_i}:
                         +\sum_r\kappa_r\partial\phi^r.
$$
The possible scalar ghost from $b\in\g_2$ has coefficient
$([f,u],b)=0$, and higher grades give no term. There are no neutral fields because the grading is integral. The bracket
$[\phi^r{}_{\lambda}J^{a_i}]=\sum_s c_{is}^{r}\phi^s$
gives
$$
 :J^{a_i}\phi^r:=:\phi^rJ^{a_i}:
                   -\sum_s c_{is}^{r}\partial\phi^s.
$$
Consequently the coefficients of $:\phi^rJ^{a_i}:$ in the
differential of the quadratic term form the matrix $AM$ and derivative coefficient
$-\frac12\sum_{i,j,r}M_{ij}A_{rj}c_{is}^{r}=C_s$.
The full differential has coefficients $AM-B$ and
$Ad+\kappa+C$, both zero. The Miura projection kills $J^u$.

A ghost-number-zero vector of degree two in the small complex consists
of a $\g_{-1}$ term, symmetric quadratic $\g_0$ terms and derivatives
of $J^a$, $a\in\g_0$. The $\g_2$ ghost coefficient forces its leading
vector into $\g_{-1}^f$. The remaining equations are exactly the two
just solved, proving completeness. For two solutions with the same
$u$, their symmetric-matrix difference $N$ satisfies $AN=0$, so
$N\in\Sym^2(\ker A)=\Sym^2\mathfrak c$. Its derivative correction
is zero and the derivative difference lies in $\ker A$.
Conversely, products and derivatives of $J^a$, $a\in\mathfrak c$,
are cocycles, proving the assertion.
\end{proof}

The coefficient of the leading small-complex field is $+J^u$.
For $\slg_2$ this normalization gives
$W_f=J^f-:hh:/4-(k+1)\partial h/2$.
The six integral gradings below therefore include every degree-two
ghost coefficient. Nonvanishing of $H_f^0(L_k(\g))$, generation of the Miura image and simplicity
are proved in the applications.
\begin{center}\small
\begin{tabular}{@{}lllr@{}}\toprule
Type, $p$&Orbit&Size of $A$&$\dim\g_0^f$\\\midrule
$F_4,2$&$F_4(a_3)$&$12\times12$&0\\
$E_6,2$&$E_6(a_3)$&$12\times12$&0\\
$E_6,3$&$D_4(a_1)$&$18\times20$&2\\
$E_7,3$&$E_7(a_5)$&$21\times21$&0\\
$E_8,3$&$E_8(a_6)$&$24\times24$&0\\
$E_8,5$&$E_8(a_7)$&$40\times40$&0\\\bottomrule
\end{tabular}
\end{center}
The maps on $\g_0^f$ and the projections of these fields are given
in \hyperref[calc:QUAD]{Calculation~\ref*{calc:QUAD}},
\hyperref[calc:E6-LIFT]{Calculation~\ref*{calc:E6-LIFT}} and
\hyperref[calc:CENTRALIZER]{Calculation~\ref*{calc:CENTRALIZER}}.
They include all quadratic and derivative terms before projection.

\subsection{Strong generators of \texorpdfstring{$\mathcal W^k(\g,f)$}{Wk(g,f)} and the Lie algebra \texorpdfstring{$\g^f$}{g-f}}
For a homogeneous element $a\in\g^f$ with $[x,a]=ja$, let
$w(a)$ be the corresponding strong generator of conformal weight $1-j$ in
$\mathcal W^k(\g,f)$. Its symbol is the corresponding linear coordinate
on the Slodowy slice. The linear terms of their zeroth products are
controlled by the Lie bracket of $\g^f$. We record the consequence
of the standard Slodowy Poisson structure that will be used below
\cite[Sections~3.1--3.3 and Corollary~4.5.10(i)]{AAV}.

\begin{lemma}\label{com:linear-bracket}
If $a,b\in\g^f$ are homogeneous and
$d=\wt w(a)+\wt w(b)-1$, then
\begin{equation}\label{com:linear-bracket-formula}
 w(a)_{(0)}w(b)-w([a,b])
 \in\langle w(c):\wt w(c)<d\rangle.
\end{equation}
\end{lemma}
\begin{proof}
The Poisson algebra $R_W=\C[\mathcal S_f]$ is generated by
$\ell_a(f+y)=(a,y)$, $y\in\g^e$. Its cotangent Lie bracket at $f$
is the bracket of $\g^f$. For example, in the Dirac formula defining
the transverse bracket, each correction contains two factors
$\{\ell_a,\varphi_i\}$ which vanish at $f$, because $[a,f]=0$.
Consequently $\{\ell_a,\ell_b\}-\ell_{[a,b]}$ has no linear term.
Every nonlinear term of weight $d$ involves generators of weights
less than $d$. The same is true of a vector $u_{(-2)}v$ of weight
$d$, since $\wt u+\wt v+1=d$. Induction on conformal weight gives
strong generation and lifts this Poisson identity to the assertion.
\end{proof}

\begin{lemma}\label{vtwo:cotangent}
Let $B$ be the graded image of a universal affine
$\mathcal W^k(\g,f)$ for a discrete Dynkin good grading with all
strong generators of weights at least one. Assume $B_0=\C$ and that its homogeneous subspaces are finite dimensional.
For $\mathfrak m_B=(R_B)_{>0}$ there is a surjective graded Lie map
$$
 \g^f\longrightarrow\mathfrak m_B/\mathfrak m_B^2.
$$
Here the Lie grading is conformal weight minus one: the Poisson
bracket lowers conformal weight by one. Its kernel is a Lie ideal. If
$\dim(\mathfrak m_B/\mathfrak m_B^2)_d\le n_d$ for each weight $d$,
then $B$ has a strong generating set with at most $n_d$ vectors of
weight $d$.
\end{lemma}
\begin{proof}
By \cite[\S3.3, Corollary 4.5.10(i)]{AAV},
$R_{\mathcal W^k}=\C[\mathcal S_f]$ and its linearized bracket at
$f$ is the centralizer bracket, also for half-integral good gradings.
The quotient preserves the augmentations defined by the positive gradings
and hence induces a Lie map on cotangent spaces: the square of the
augmentation ideal is a Lie ideal in the augmentation ideal. Positive-weight induction shows that
cotangent lifts generate $R_B$. Reconstruction gives strong generation,
since $a_{(-2)}b$ uses lower-weight fields or their derivatives.
For half-integral weights, induct on twice the weight.
If $B$ is realized in a Heisenberg vertex algebra $\mathcal H$,
vanishing in $R_{\mathcal H}$ alone does not imply vanishing in $R_B$.
\end{proof}

In particular, suppose the generators $w(a)$, $a\in\g_0^f$,
and two oppositely charged
weight-two spaces vanish, and their centralizer bracket spans the unique
weight-three direction. Its class vanishes in $R_B$: every other
weight-three polynomial correction contains a vanished weight-one
coordinate. If $B_1=0$ and $B_2=\C T$, then
$C_2(B)_3=\partial B_2$, so that last field is a derivative of $T$.
Thus $B$ is a Virasoro quotient.

\begin{lemma}\label{v3:abelianization}
Write $K=\g^f=\bigoplus_{j\ge0}K_j$, where
$K_j=\g^f\cap\g_{-j}$ and the latter index is the Dynkin
$h$-weight. Thus Dynkin degree $j$ corresponds to conformal weight
$1+j/2$.
Suppose all generators $w(a)$ with $a\in K_0$ and all fields
with linear symbols in a homogeneous
subspace $R$ vanish in a graded quotient $B$ of the universal
$W$ algebra. Set $Q=K/\langle K_0,R\rangle$. Homogeneous lifts of
$$
 Q/[Q,Q]=K/([K,K]+K_0+R)
$$
Lie-generate $Q$ and vertex-generate $B$. If this quotient has one
generator of Dynkin degree two, $B$ is a Virasoro quotient whenever its conformal vector is nonzero. If it has Dynkin degrees two and six, corresponding to conformal weights two and four, and the first lift is central
in $Q$, then $B$ has strong generators only in those two conformal weights.
\end{lemma}
\begin{proof}
$Q$ is positively graded. Subtracting lifts of its abelianization leaves
brackets of strictly smaller positive degrees, so induction proves Lie
generation. Lift each bracket by a zeroth vertex product. By
\cref{vtwo:cotangent}, its discrepancy is decomposable modulo $C_2$.
Positive-weight induction expresses it using lower fields and derivatives.
This proves vertex generation also in half-integral Dynkin gradings.
A Lie algebra generated by one element is one-dimensional. Two generators
with one central commute, so their Lie algebra is abelian and
two-dimensional. In these two cases the cotangent dimension bound gives
the stated strong generation. The field in weight two is the nonzero conformal vector. The vanishing of the generators indexed by $K_0$ is necessary for
this argument, which uses a positively graded Lie algebra.
The finite abelianization rows used below are \cite[Tables 7--11]{deGraaf}.
\end{proof}

\begin{lemma}
\label{v5:equivariant-generation}
Let $K=\bigoplus_{j\ge0}K_j$ be a finite graded Lie algebra with
reductive $K_0$ acting semisimply. Let $J$ be the ideal generated by
a specified homogeneous subspace $R\subseteq K$. If brackets of already known elements
of $J$ with $K$ span a $K_0$-submodule, that module belongs to $J$.
For an isotypic component $V_\lambda\otimes M$, the bracket images must span
$M$. For an irreducible summand of multiplicity one a nonzero projection suffices.
Thus these bracket maps give quotient bounds for
\cref{vtwo:cotangent}. If their sum is itself an ideal containing $R$,
they identify $J$ exactly.
\end{lemma}
\begin{proof}
The ideal is $K_0$-stable, so each nonzero equivariant map onto a
simple $K_0$-module is surjective. Applying this to the multiplicity spaces
and inducting over the listed bracket maps proves the containment. The
reverse inclusion in the last assertion follows from the definition of the
ideal generated by $R$.
For normalized $\slg_2$ generators $e,f,h$, highest vectors
$u\in V_m$, $v\in V_n$ and $r=(m+n-s)/2\in\Z$ with
$0\le r\le\min(m,n)$, a highest-weight vector in the component of type $V_s$ is
$$
 \sum_{i=0}^r c_i[(\ad f)^iu,(\ad f)^{r-i}v],\quad
 c_0=1,\quad
 c_{i+1}=-c_i\frac{(r-i)(n-r+i+1)}{(i+1)(m-i)}.
$$
The identity $ef^iu=i(m-i+1)f^{i-1}u$ proves it is highest.
For the graded Heisenberg Lie algebra $K_1=V_1$, $K_2=\C z$,
$[u,v]=\det(u,v)z$, this expression is $2z$. Adding an independent
central line $\C w$ does not eliminate $w$: the corresponding multiplicity
space then has dimension two.
\end{proof}

\begin{lemma}
\label{v5:triangular-generation}
Let $V$ be a nonnegatively graded vertex algebra with $V_0=\C\mathbf1$,
strongly generated by homogeneous vectors $u_1,\ldots,u_r$ of positive
weights, and let $B\subseteq V$ be a vertex subalgebra. For each
weight $d$ occurring among these generators, suppose that $B$ contains
vectors
$$
 v_\alpha=\sum_{\wt u_i=d}c_{\alpha i}u_i+P_\alpha,
$$
where each $P_\alpha$ is a linear combination of derivatives and
normally ordered products of generators of weights strictly less than
$d$, and the matrix $(c_{\alpha i})$ has full column rank. Then $B=V$.
No linear independence of normally ordered products is assumed.
\end{lemma}
\begin{proof}
Induct on the finitely many weights of the $u_i$. By the induction
hypothesis, every $P_\alpha$ belongs to $B$. The full-rank condition
then expresses each $u_i$ of weight $d$ as a linear combination of the
$v_\alpha-P_\alpha$. Hence every strong generator belongs to $B$,
and $B=V$.
In applications, zeroth products can be computed using
\cref{com:linear-bracket}. Coefficients of higher products require a
calculation in the vertex algebra itself. The lemma proves generation
as a vertex algebra, not necessarily strong generation by the vectors
from which the $v_\alpha$ were constructed.
\end{proof}

\subsection{Symplectic-fermion invariant subalgebras}
\label{com:scalar}
We use two arguments to establish simplicity: restriction of a positive-definite unitary form, and restriction of the standard nondegenerate invariant bilinear form on the symplectic-fermion vertex superalgebra.

Let $V$ be a unitary vertex operator algebra with positive-definite
Hermitian form $\langle\cdot,\cdot\rangle$, linear in the second variable.
Its \emph{PCT involution} (parity, charge conjugation and time reversal)
is the conjugate-linear vertex operator algebra automorphism $\theta$
with $\theta^2=\id$ for which
$$
 B(u,v)=\langle\theta u,v\rangle
$$
is an invariant bilinear form. In particular, $\theta$ fixes the vacuum
and conformal vectors. A vertex subalgebra $A\subseteq V$ is
\emph{PCT-stable} if $\theta(A)=A$. If $A$ has the same conformal
vector as $V$, this condition ensures that the restricted Hermitian
form defines a unitary structure on $A$, see
\cite[Section~5.1 and Proposition~5.23]{CarpiEtAl2018}.

The following standard consequence of unitarity is contained in
\cite[Propositions~5.3 and~5.23]{CarpiEtAl2018}. We include the short
argument for completeness.
\begin{lemma}\label{com:scalar-simple}
A conformal subalgebra of a positive-definite unitary vertex
algebra, with the same conformal vector and stable under its PCT
involution, is simple if its degree-zero space is
$\C\mathbf1$.
\end{lemma}
\begin{proof}
The restricted Hermitian form is positive definite, and the
adjoints of the subalgebra modes are again its modes, by PCT
stability and the common conformal structure. The orthogonal
complement of a graded ideal is therefore a submodule. A proper
ideal has no degree-zero part, so this complement contains the
vacuum and hence the entire subalgebra. The ideal is zero.
\end{proof}

For the second argument, use the common odd lattice
$\Z\phi_1\oplus\cdots\oplus\Z\phi_s$, with orthonormal
$\phi_i$ and compatible cocycles, and write $h_i=\partial\phi_i$.
Set
\begin{equation}\label{com:sf-lattice}
 e_i=e^{\phi_i},\qquad f_i=\partial e^{-\phi_i},\qquad
 T_s=\frac12\sum_{i=1}^s(:h_i^2:-\partial h_i).
\end{equation}
Both $e_i,f_i$ have weight one for $T_s$. Their OPEs are
$$
 e_i(z)f_j(w)\sim\frac{\delta_{ij}}{(z-w)^2},
 \qquad f_i(z)e_j(w)\sim-\frac{\delta_{ij}}{(z-w)^2},
$$
and all other pairings have no singular part. These fields
realize the rank-$s$ symplectic-fermion algebra $\SF(s)$.
The lattice expansion
$$
 e^{\phi_i}(z)e^{-\phi_i}(w)
 =(z-w)^{-1}\left(1+(z-w)h_i
       +\frac{(z-w)^2}{2}(:h_i^2:+\partial h_i)+\cdots\right)
$$
gives $:e_if_i:=(\partial h_i-:h_i^2:)/2$, so the conformal
vector is exactly $T_s$ in \eqref{com:sf-lattice}. In particular
all ranks and coordinate subsets are realized compatibly in the
vertex superalgebra of this odd lattice.

Denote
\begin{equation}\label{com:sf-invariants}
                    \mathcal A_s=\SF(s)^{\Sp_{2s}}.
\end{equation}
The diagonal torus in $\Sp_{2s}$ makes its invariants have zero
lattice charge in every coordinate. Thus $\mathcal A_s\subseteq
\mathcal H^{\otimes s}$. It contains $T_s$ and is fixed
pointwise by coordinate permutations. The established
symplectic-fermion invariant theory says that $\mathcal A_s$ is freely
strongly generated in weights $2,4,\ldots,2s$, and for $s\ge2$
is generated as a vertex algebra by its weight-two and
weight-four fields. These results are due to Creutzig--Linshaw
\cite[Theorems~3.2 and~3.9, Corollary~5.7]{CLSF}.

The symplectic-fermion vertex superalgebra $\SF(s)$ carries a
nondegenerate invariant bilinear form induced by the symplectic
pairing on its generating space. It is nondegenerate on each
finite graded piece and preserved by $\Sp_{2s}$. Complete reducibility
makes its invariant restriction nondegenerate. The central element
$-1$ removes odd states, so $\mathcal A_s$ is an ordinary conformal vertex
algebra of CFT type. Invariance makes every proper ideal orthogonal
to the algebra, proving simplicity. The invariant form used here is the restriction of the form on
$\SF(s)$. No invariant-form assertion for $\mathcal H^{\otimes s}$
with conformal vector $T_s$ is needed. By
\cref{com:conformal-obstruction}, $[T_s]$ is not nilpotent and $\mathcal A_s$
is not quasi-lisse for $s>0$.

\begin{proposition}
\label{com:sf-recognition}\label{vtwo:scalar-recognition}
Let $r\ge1$, and use one common lattice realization of
$\mathcal A_r=\SF(r)^{\Sp_{2r}}$ and $\mathcal A_{r+1}$. Let $a$ denote an additional
generator of a rank-one Heisenberg vertex algebra, normalized by
$[a_\lambda a]=\lambda\mathbf1$, with conformal vector
$T_a=\tfrac12(:a^2:-\partial a)$. Suppose the conformal image $B$
satisfies the containment
$$
 B\subseteq\mathcal H_a\otimes \mathcal A_r,
 \qquad T=T_a+T_r\in B,
$$
and is strongly generated by at most one field in each weight
$2,4,\ldots,2r+2$, with no generators of other weights. Suppose also that it has a
weight-four field whose image in
$R_{\mathcal H_a\otimes\mathcal H^{\otimes r}}$ is
$$
 \alpha p_4+\beta p_2^2,\qquad \alpha\ne0,
 \quad p_j=a^j+\sum_{i=1}^r h_i^j.
$$
Then $B=\mathcal A_{r+1}$ as subalgebras of $\mathcal H_a\otimes\mathcal H^{\otimes r}$. In particular,
$B$ is simple and its conformal class is not nilpotent.
\end{proposition}
\begin{proof}
The weight-four primary calculation in
Appendix~\ref{app:primary-four} has two consequences. First, the map
to $R_{\mathcal H_a\otimes\mathcal H^{\otimes r}}$ is injective on
the weight-four primary subspace of $\mathcal H_a\otimes \mathcal A_r$.
Second, a primary vector with symbol $\alpha p_4+\beta p_2^2$,
$\alpha\ne0$, must have $\beta=6\alpha/(5r-6)$ and is a nonzero
multiple of the weight-four primary generator of the concrete
subalgebra $\mathcal A_{r+1}$. Thus the calculation identifies the complete
vector, including its derivative terms.

The generator bound on $B$ gives $B_1=0$, $B_2=\C T$ and
$B_3=\C\partial T$. Subtracting Virasoro descendants from its weight-four generator gives a primary vector. The weight-four descendant Gram
determinant is nonzero at $c=-2(r+1)$, since $c\ne0,-22/5$. The same correction follows directly from the $L_1,L_2$ equations. This correction changes the
quartic symbol only by a multiple of $p_2^2$. The uniqueness established in Appendix~\ref{app:primary-four}
shows that $W_{r+1}$ lies in $B$, whence $\mathcal A_{r+1}\subseteq B$ by weight-two/four
vertex generation. Since $B$ has at most one strong generator in
each of the weights $2,4,\ldots,2r+2$, spanning by normally ordered
products gives, with formal variable $z$,
$$
 \dim B_d\le [z^d]\prod_{j=1}^{r+1}\prod_{m\ge0}
                         (1-z^{2j+m})^{-1}=\dim(\mathcal A_{r+1})_d.
$$
The equality on the right uses free strong generation of $\mathcal A_{r+1}$.
Together with $\mathcal A_{r+1}\subseteq B$, this gives equality of the
homogeneous subspaces. All strong generators of $B$ belong to $\mathcal A_{r+1}$, proving
equality. Simplicity and the obstruction were proved above.
\end{proof}

For rank one, $\mathcal A_1=\Vir_{-2}$. The universal Virasoro vertex algebra at
$c=-2$ is simple, since the Virasoro minimal-model central-charge formula would require
$p/q=2$ or $1/2$, incompatible with coprime $p,q\ge2$
(\cite{Salazar2026}, Theorem~4.1). Proposition~\ref{com:sf-recognition} requires the degree-four component
to belong to $\C p_4+\C p_2^2$ with nonzero coefficient of $p_4$. A containment in the specified lattice realization,
the bound on the strong generators and the nonvanishing of
$H_f^0(L_k(\g))$ are established separately in each application. These conditions concern the image of the specified homomorphism.

\begin{lemma}
\label{vtwo:torus-detector}
Let an algebraic torus $T$ act nontrivially on the generating symplectic
space of $\SF(s)$, $s>0$. Every parity-even conformal subalgebra
$B\subseteq\SF(s)^T$ containing the same conformal vector has a
nonnilpotent conformal class. If $R_B$ is finitely generated, $B$ is
not quasi-lisse.
\end{lemma}
\begin{proof}
Choose a nontrivial one-dimensional subtorus with integral weights
$(a_i,-a_i)$ on symplectic pairs, where $a\ne0$. Use
\eqref{com:sf-lattice}. Every charge of $B$ belongs to the positive
even lattice
$$
 K=\{n\in\Z^s:a\cdot n=0,\ \textstyle\sum_i n_i\equiv0\pmod2\}.
$$
It spans $a^\perp$, so the lattice decomposition embeds
$B\hookrightarrow V_K\otimes\mathcal H_a$. Orthogonal decomposition
of the conformal vector, followed by the augmentation of $R_{V_K}$ defined
by its standard conformal grading, gives a homomorphism of commutative $C_2$ algebras
$$
 R_B\longrightarrow R_{V_K}\otimes\C[z]\longrightarrow\C[z],
 \qquad [T]\longmapsto\tfrac12z^2.
$$
Derivative background terms vanish in $C_2$. The lattice augmentation
uses its standard positive grading, not the background-charge grading.
It exists without needing lattice $C_2$-cofiniteness. Nonzero powers of
$z^2/2$ prove the assertion by \cref{com:conformal-obstruction}.
No $C_2$ injection for a subalgebra is assumed.

For $\SF(\C^d\otimes\C^{2r})^{\Sp_{2r}}$, the sublattice of zero torus charge is
$\mathcal A_{d-1}^{\oplus r}$ and the same map sends the conformal vector to
$\frac12\sum_{i=1}^r z_i^2$. At $s=1$, $K=0$ and this is the familiar
rank-one Heisenberg realization. The parity condition is essential:
for $a=(1,2)$ the even kernel lattice is generated by $(4,-2)$, not
$(2,-1)$. A trivial torus, a subalgebra missing the conformal vector, or an
arbitrary simple quotient need not preserve this argument.
\end{proof}

\subsection{Quadratic generators of symplectic-fermion invariants}
\begin{proposition}\label{v3:quadratic-generation}
For $d\ge2,r\ge1$, set
$A_{d,r}=\SF(\C^d\otimes\C^{2r})^{\Sp_{2r}}$.
It is vertex-generated by its $d(d+1)/2$ symmetric weight-two
quadratics. Identify these with self-adjoint endomorphisms for a
nondegenerate symmetric form $G$, writing $Q(I)=T$.
If $D$ has distinct eigenvalues and the off-diagonal graph of a second
self-adjoint matrix $M$ in a $D$-eigenbasis is connected, then
$T,Q(D),Q(M)$ generate $A_{d,r}$.
Consequently, if a conformal map from a weight-two-generated algebra
$U$ has these three complete images and $\dim U_2\le d(d+1)/2$,
its entire image is $A_{d,r}$, which is simple and not quasi-lisse.
\end{proposition}
\begin{proof}
Use odd pairs $e_{i,l},f_{i,l}$ with opposite double-pole pairings,
and set
$$
 q_{ij}^{a,b}=-\frac12\sum_l
 (:\partial^ae_{i,l}\partial^bf_{j,l}:
       +:\partial^be_{j,l}\partial^af_{i,l}:).
$$
The symplectic first fundamental theorem and exact invariants in the
fermion filtration give strong generation by all these contractions
\cite[\S3, (3.10)]{CLSF}. For $i\ne j$, the separate conformal vectors for individual rows
$T_i=q_{ii}^{0,0}$ translate the two indices independently, generating
every off-diagonal contraction. A single Wick contraction gives
$$
 2(q_{ij}^{m,0})_{(1)}q_{ij}^{0,0}
     =q_{ii}^{m,0}+(-1)^m(m+1)q_{jj}^{m,0}.
$$
Together with the interchanged row, this has determinant $-m(m+2)$
for $m>0$. The diagonal $m=0$ fields are already present. Total
translation then generates all diagonal derivatives. Full contractions
have pole order $m+4$ and do not enter this first product.

Wick contraction also gives $Q(A)_{(1)}Q(B)=Q(AB+BA)$.
Polynomial interpolation in $D$ produces its spectral idempotents.
Two idempotents isolate each nonzero edge of $M$. Products along paths
give every symmetric matrix unit. This proves generation. The inclusion
$A_{d,r}\subseteq\im U$ and the weight-two bound force equality in
weight two and hence equality of the vertex algebras. The invariant
oscillator form proves simplicity, as in \cref{com:scalar}. Finite strong
generation is \cite[Theorem 3.10]{CLSF}, and \cref{vtwo:torus-detector} proves nonnilpotence. Complete field images,
not only their leading quadratic symbols, are required here.
\end{proof}

\section{Affine vertex subalgebras of \texorpdfstring{$H_f^0(Q)$}{Hf0(Q)} and nilpotent orbits}
\label{com:boundary-current-section}

We use the nilpotent-orbit conventions and standard geometric results in \cite{CollingwoodMcGovern1993}.
Related questions on the local geometry of special pieces and their transverse slices are studied in \cite{FuJuteauLevySommers2023}.

Let $Q$ be a graded quotient of $V^k(\g)$. All affine vertex subalgebras
below lie in $H_f^0(Q)$, a quotient of $\mathcal W^k(\g,f)$ which we
do not assume to be simple. In the applications, $Q=L_k(\g)$ or a
specified intermediate quotient of $V^k(\g)$.

Write $\g^\natural=\g_0^f$ for the centralizer of the
$\mathfrak{sl}_2$-triple, where $\g_0=\ker\operatorname{ad}x$.
At a maximal nilpotent orbit, its action on the reduced Slodowy-slice
intersection is trivial. This forces the homomorphisms from universal
affine vertex algebras to factor through simple integrable affine vertex
algebras, imposing relations on the associated variety of the reduction.

The following adapts the integrability argument of
Dong and Mason \cite[Theorem~3.1]{DongMason2006}.
Here $C_2$-cofiniteness is replaced by triviality of the action on
the reduced $C_2$ algebra, and $V$ need not be simple. The level-zero
case is also included.

\begin{lemma}\label{com:bounded-current}
Let $V$ be nonzero and finitely strongly generated, with $V_0=\C$,
a discrete nonnegative conformal grading, finite-dimensional pieces, and
no weights strictly between zero and one. Suppose a simple Lie algebra
$\mathfrak a$ defines a graded homomorphism $V^\ell(\mathfrak a)\to V$
whose weight-one generators are primary, with long roots of squared
length two. If the operators $x(0)$, $x\in\mathfrak a$, act trivially
on $(R_V)_{\mathrm{red}}$, then $\ell\in\Z_{\ge0}$ and this
homomorphism factors through an embedding
$L_\ell(\mathfrak a)\hookrightarrow V$.
\end{lemma}
\begin{proof}
Let $\theta$ be the highest long root. Choose finitely many homogeneous
Cartan eigenvectors generating $R_V$. Each generator of nonzero
$h_\theta$-charge vanishes in the reduced algebra and is nilpotent.
Thus all $h_\theta$-charges occurring in $R_V$ are bounded in absolute
value by some $B$. Powers of neutral generators do not affect this bound.

Choose homogeneous representatives of a basis of the entire $R_V$, omitting
the vacuum vector. This basis may be infinite. The difference-one spanning
theorem of Gaberdiel--Neitzke \cite[Proposition~8]{GaberdielNeitzke2003}
(see also \cite[Proposition~5.7]{Buhl2010}) spans $V$ by
$$
 x^1_{(-n_1)}\cdots x^s_{(-n_s)}\mathbf1,
 \qquad n_1>\cdots>n_s>0.
$$
It does not require $R_V$ to be finite dimensional. The usual induction
by generator weight also applies to the positive discrete half-grading
used here. If $d_i\ge1$ and $b_i$ are the weight and charge of $x^i$,
the eigenvalue of $A=L_0-h_\theta(0)/2$ on such a word satisfies
$$
 A=\sum_i(d_i+n_i-1-b_i/2)
 \ge\frac{s(s+1)-Bs}{2}.
$$
A fixed value of $A$ bounds $s$, and the original conformal weight is at
most $A+Bs/2$. Consequently every $A$-eigenspace is finite dimensional.
The operator $e_\theta(-1)$ preserves $A$ and raises the original weight.
It is locally nilpotent. The operator $f_\theta(1)$ lowers original
weight and is locally nilpotent too. The operators $e_\alpha(0)$ and $f_\alpha(0)$ act locally
nilpotently on each finite-dimensional homogeneous $\mathfrak a$-module. These are the affine
Chevalley generators, proving integrability.

If $N$ is least with $e_\theta(-1)^N\mathbf1=0$, the real-root
$\slg_2$ identity gives
$$
 0=f_\theta(1)e_\theta(-1)^N\mathbf1
   =N(\ell-N+1)e_\theta(-1)^{N-1}\mathbf1.
$$
Hence $\ell=N-1$, and
$e_\theta(-1)^{\ell+1}\mathbf1=0$ in $V$. The homomorphism
$V^\ell(\mathfrak a)\to V$ therefore factors through
$L_\ell(\mathfrak a)$. Since $V\ne0$, the resulting unital
homomorphism from the simple vertex algebra is injective.
\end{proof}

The next lemma records the consequence of Lusztig--Spaltenstein
induction needed for conical associated varieties. We give the
geometric argument to make the one-dimensional parameter explicit
\cite{LusztigSpaltenstein}.
\begin{lemma}\label{vtwo:conic-boundary}
Let $X\subseteq\g$ be closed, conical and $G$-stable, and let
$z=s+n\in X$ be its Jordan decomposition. Set $\lev=\g^s$ and let
$L$ have Lie algebra $\lev$. Then $X$ contains the closure of
$\operatorname{Ind}_{\lev}^{\g}(L\cdot n)$, of dimension
$$
 \dim\g-\dim\lev+\dim(L\cdot n)=\dim(G\cdot z).
$$
More precisely, if $s\ne0$, the cone
$Y=\overline{G(\C^*z)}$ has nilpotent intersection equal to this
induced orbit closure and dimension $\dim(G\cdot z)+1$.
\end{lemma}
\begin{proof}
Choose a parabolic $P$ with Levi $L$ and nilradical $\uu$. The family
$$
 Z=G\times_P\bigl(\C s+\overline{L\cdot n}+\uu\bigr)
$$
is closed in $\g\times G/P$, so projection to $\g$ is proper. For
$c\ne0$, the operator $\ad(cs+y)$, with $y\in\overline{L\cdot n}$,
is invertible on each graded piece of the nilradical filtration:
$\ad(s)$ has nonzero eigenvalues and the commuting $\ad(y)$ is
nilpotent. Successive unipotent conjugations therefore remove every
$\uu$ component of $cs+y+u$. The nonzero-parameter locus is dense,
and scaling $n$ can be absorbed by its Jacobson--Morozov cocharacter
in $L$. The proper image of $Z$ is consequently
$Y=\overline{G(\C^*s+L\cdot n)}\subseteq X$.
A point with $c\ne0$ has nonzero semisimple part. Thus $Y\cap\N$
is precisely the image of the zero-parameter fiber, which is the
induced orbit closure. The usual induction dimension formula gives
the stated orbit dimension. The semisimple parameter varies freely
up to finite Weyl ambiguity, adding one dimension to $Y$.
For $s=0$, take $Y=\overline{G\cdot n}$ and $\lev=\g$.
The induction and dimension theorem is the standard
Lusztig--Spaltenstein theorem \cite{LusztigSpaltenstein}, see also
\cite[Appendix B.0.5, Theorem B.1]{AMLectures}.
\end{proof}

\begin{lemma}\label{vtwo:point-slice}
Let $X$ be as above and suppose $X\cap\mathcal S_f=\{f\}$.
There is no nonnilpotent $z\in X$ for which the induced boundary
in \cref{vtwo:conic-boundary} is $G\cdot f$. In particular, a
one-parameter semisimple family with induced dimension $\dim Gf$
and that boundary cannot lie in $X$.
\end{lemma}
\begin{proof}
Its cone $Y\subseteq X$ would have dimension $\dim Gf+1$ and contain
$f$. Slodowy transversality, or the smooth local action map
$G\times\mathcal S_f\to\g$ at $f$, gives
$\dim_f(Y\cap\mathcal S_f)\ge1$, a contradiction. It suffices that $Gf$ be the unique possible boundary of the equality
dimension. To establish the identity $X\cap\mathcal S_f=\{f\}$ by reduction,
the actual reduction must be nonzero and $C_2$-cofinite. Its simple quotient does not suffice.
\end{proof}

\begin{lemma}\label{com:maximal-boundary}
Let $Q$ be a finitely generated graded quotient of $V^k(\g)$, and let $f$
be maximal, in orbit closure order, among nilpotent orbits in $X_Q$.
Then the reduced Poisson bracket on $X_Q\cap\mathcal S_f$ is zero.
If $H_f^0(Q)$ satisfies the hypotheses of \cref{com:bounded-current},
the homomorphisms from universal affine vertex algebras associated with
the simple ideals of $\g^\natural$ factor through integrable simple
affine vertex algebras at nonnegative integral levels. At level zero,
the corresponding Lie algebra has zero image in $H_f^0(Q)$. No simplicity of $Q$ or of its reduction is required.

If, in addition,
$$
 X_Q\cap\N\subseteq\overline{Gf},\qquad
 X_Q\cap\mathcal S_f=\{f\},
$$
every nonnilpotent orbit in $X_Q$ has dimension at most $\dim Gf-2$.
If $Q$ is a quotient of $V^k(\g)$ with $k+h^\vee\in\Q_{>0}$, then for a maximal Levi $\lev_j$, a nonzero point of its center in $X_Q$ and a nilpotent
orbit of dimension $d$ in the product of associated varieties of the simple affine vertex
algebras associated with $[\lev_j,\lev_j]$ imply
\begin{equation}\label{com:center-dimension-test}
 \dim\g-\dim\lev_j+d\le\dim Gf-2.
\end{equation}
The final implication uses the support inclusion \cref{com:levi-support},
obtained from the parabolic Wakimoto homomorphism. The assertion that
the reduced Poisson bracket vanishes does not use this inclusion.
\end{lemma}
\begin{proof}
For $z\in X_Q\cap\mathcal S_f$, take the cone $Y$ of
\cref{vtwo:conic-boundary}. The contracting Kazhdan action preserves
$Y$ and places $f$ in its nilpotent boundary. Hence $f$ is in the
closure of the induced orbit. Maximality makes that orbit $Gf$.
The transverse Poisson rank at $z$ is therefore
$\dim Gz-\dim Gf=0$. This proves that the reduced slice bracket is
zero. The action of $\g^\natural$ is Hamiltonian, so its elements
act trivially on the reduced $C_2$ algebra. Apply
\cref{com:bounded-current}. At level zero, the simple affine vertex algebra is
$L_0(\mathfrak a)=\C$.

Under the extra hypotheses, \cref{vtwo:point-slice} rules out $Gf$
as the induced boundary of a nonnilpotent point. Every remaining
boundary has strictly smaller even orbit dimension, giving the drop
by two. If $X_Q$ contains a nonzero point of a maximal Levi center, \cref{com:levi-support}
gives the mixed orbit with dimension
$\dim\g-\dim\lev_j+d$, proving the last inequality.
\end{proof}

\begin{corollary}\label{com:scalar-point}
Let $H\ne0$ be a finitely strongly generated vertex algebra, with a
nonnegative discrete conformal grading, finite-dimensional graded pieces,
and $H_0=\C$. Let a semisimple Lie algebra $\mathfrak a$ act by the zero modes of
vectors $J^a\in H_1$, $a\in\mathfrak a$, satisfying
$J^a_{(0)}J^b=J^{[a,b]}$. Choose a finite-dimensional
homogeneous $\mathfrak a$-stable space of strong generators, and assume
that $[J^a]$ is nilpotent in $R_H$ for every $a\in\mathfrak a$.
Then $(R_H)_{\mathrm{red}}$ is generated by the images of the trivial
$\mathfrak a$-summands in that space. If these images are zero, $H$ is
$C_2$-cofinite.

This nilpotence hypothesis holds if the homomorphisms from the
universal affine vertex algebras associated with the simple ideals of $\mathfrak a$
factor through integrable simple affine vertex algebras. It also holds
if the reduced Poisson bracket on $R_H$ is zero. In the latter case
no affine identification is needed.

If $Q$ is a quotient of $V^k(\g)$, $H=H_f^0(Q)$, and
the images of all remaining $\mathfrak a$-invariant generators vanish, then
$X_Q\cap\mathcal S_f=\{f\}$. The same conclusion holds
if its conformal vector vanishes, in which case $H=\C$.
\end{corollary}
\begin{proof}
The nilradical is a Poisson ideal. Thus the nilpotent classes $[J^a]$
act trivially on $(R_H)_{\mathrm{red}}$, and equivariance kills every
nontrivial irreducible summand of the strong generators. The remaining
images of the $\mathfrak a$-invariant generators generate this algebra. If they vanish, finite generation
and $(R_H)_{\mathrm{red}}=\C$ imply that $R_H$ is finite-dimensional.

Integrable simple affine vertex algebras are $C_2$-cofinite, so the classes of their weight-one generators, and their images in $R_H$, are nilpotent. A zero reduced
bracket instead gives
$\overline{J^{[a,b]}}=\{\overline{J^a},\overline{J^b}\}=0$,
and perfection of $\mathfrak a$ gives $\overline{J^a}=0$ for every $a\in\mathfrak a$.
For an actual reduction, \cref{com:ds} identifies the variety with
$X_Q\cap\mathcal S_f$, whose origin is $f$.
Finally $\omega_H=0$ gives $L_0=0$, hence $H=H_0=\C$.
\end{proof}

Thus it remains to find relations in $H$ that make the classes of
the remaining $\mathfrak a$-invariant generators nilpotent in $R_H$. For a sole conformal generator, a relation expressing
it as a polynomial in $[J^a]$, $a\in\mathfrak a$, with zero constant term
suffices. Integral levels alone do not prove the nilpotence of $[J^a]$, and
no identification with the affine subalgebra is asserted.

\subsection{A criterion for \texorpdfstring{$X_Q\cap\C h=\{0\}$}{XQ intersect C h = 0}}
\label{vtwo:affirmative-tools}
\begin{lemma}
\label{vtwo:sugawara-axis}
Let $Q$ be a finitely generated quotient of $V^k(\g)$, with $k\ne-h^\vee$, and let
$(e,h,f)$ have integral Dynkin grading. Suppose $H_f^0(Q)\ne0$ and its conformal vector has $C_2$ class equal to a polynomial
in the classes $[J^a]$, $a\in\g^\natural$, with zero constant term.
Then $X_Q$ contains no nonzero point of $\C h$.
\end{lemma}
\begin{proof}
At $y=f+ce\in\mathcal S_f$, $c\ne0$, the linear functions corresponding to $J^a$,
$a\in\g^\natural$, vanish by grading. The proposed polynomial is zero there,
whereas the actual conformal invariant has value
$$
 \frac{(y,y)}{2t}=\frac{c(f,e)}t\ne0.
$$
Thus $y\notin X_Q\cap\mathcal S_f$. Within the triple's
$\mathrm{SL}_2$, it is conjugate to $\sqrt c\,h$, excluding every
nonzero point of $\C h$. For $\slg_2$ this is the direct calculation with
$f+ce=\left(\begin{smallmatrix}0&c\\1&0\end{smallmatrix}\right)$.
A relation only in the simple quotient of $H_f^0(Q)$ would not rule out
points of $X_Q\cap\mathcal S_f$.
\end{proof}

\section{The Jantzen filtration and ideals of \texorpdfstring{$\mathcal W^k(\g,f)$}{Wk(g,f)}}\label{vtwo:actual-tools}
The reduced Jantzen sum counts all ideals $H_f^0(J^a)$, whereas a
relation in the reduction of $L_k(\g)$ must belong to $H_f^0(J^1)$.
We distinguish these using nestedness, representation multiplicities
and vertex products, including proofs of the elementary subspace
lemmas to clarify the use of character coefficients.
Throughout, $t=k+h^\vee>0$, maps preserve Sugawara degrees,
and every specified weight lattice for the $\g$-action applies to all
subquotients.

\begin{lemma}
\label{vtwo:scalar-vacuum}
Let $V=V^k(\g)$ be universal, and let $s$ be a nonzero affine singular
vector of dominant $\g$-highest weight $\lambda$ and degree
$d=\Delta_t(\lambda)>0$. Let $I$ be its generated ideal. Assume, for
the chosen dominant Dynkin grading $x$, that
\begin{enumerate}
\item every allowable dominant $\mu$ with $\Delta_t(\mu)>d$ satisfies
$\Delta_t(\mu)-\mu(x)>0$,
\item $H_f^0(\mathbb V_k(\lambda))$ has no negative reduced degrees
and its degree-zero space has dimension one.
\end{enumerate}
Then $H_f^0(V/I)=0$.
\end{lemma}
\begin{proof}
The generalized Verma module $\mathbb V_k(\lambda)$ maps onto $I$. Its entire
irreducible $\g$-module in its lowest degree maps injectively, so its kernel $K$ vanishes
through degree $d$. The proof of \cref{com:tail}, using hypothesis (1),
gives $H_f^0(K)_0=0$. Exactness carries the nonzero scalar of (2) to
$H_f^0(I)_0$ and embeds that space in $H_f^0(V)_0=\C\mathbf1$.
The reduced ideal therefore contains $\mathbf1$. Its quotient
is zero. The allowable weight lattice in hypothesis (1) is essential.
\end{proof}

\begin{lemma}
\label{vtwo:first-top}
Let $U$ be a graded module for an algebra acting by
degree-preserving operators, and let $K_1\supseteq K_2\supseteq\cdots$
be graded submodules. In each degree $a$, suppose that the modules
$(K_i)_a$ have finite length and vanish for sufficiently large $i$. If their total Grothendieck class vanishes below $d$ and
equals $[R]$ in degree $d$, where $R$ is simple, then
$$
 (K_1)_d\cong R,\qquad (K_i)_d=0\ (i\ge2),
 \qquad (K_i)_a=0\ (a<d).
$$
\end{lemma}
\begin{proof}
The classes of simple objects form a free basis of the Grothendieck group. All
composition multiplicities on the left are nonnegative. Thus exactly
one nonzero term can occur in degree $d$, and it is $R$. Nestedness
forces that term to occur at $i=1$. The same positivity proves vanishing below
$d$. Simplicity of $R$ is essential: $K_1=K_2=\C$ gives total
dimension two but a one-dimensional first ideal.
\end{proof}

\begin{lemma}\label{vtwo:equivariant}
Let $U$ be a conformal vertex algebra with conformal vector $T$
and finite-dimensional homogeneous subspaces. Let a reductive
algebraic group $A$ and a finite group $\Gamma$ act on $U$ by commuting
conformal automorphisms, and let $K_1\supseteq K_2\supseteq\cdots$
be stable graded ideals, with only finitely many nonzero terms in
each degree. Suppose a $\Gamma$-equivariant linear map
$\Phi:U_2^A\to P$
identifies a subspace $P_0\subseteq(K_1)_2^A$ with $P$. Set
$K=\ker\Phi$, assume $K$ is multiplicity-free as a $\Gamma$-module,
and suppose $(K_i)_2^A\subseteq K$ for $i\ge2$. If
$$
 \sum_i[(K_i)_2^A]=[P_0]+[K],
$$
then $(K_1)_2^A=U_2^A$. The special case $U_2^\Gamma=\C T$ and
$\sum_i\dim(K_i)_2^\Gamma=1$ directly gives $T\in K_1$.
If $U/K_1$ has its inherited nonnegative conformal grading and degree-zero space $\C\mathbf1$, either conclusion containing $T$ implies $U/K_1=\C$.
\end{lemma}
\begin{proof}
The splitting gives $U_2^A=P_0\oplus K$ and
$(K_1)_2^A=P_0\oplus L$ with $L\subseteq K$. Every later term lies
in $L$. Each irreducible constituent of $K$ occurs exactly once in
the displayed sum. If absent from $L$, nestedness would make it absent
from the sum. Thus $L=K$, and later terms are zero. For an invariant
line the same argument is one-dimensional. If $T$ vanishes in the
quotient, $L_0=T_{(1)}$ vanishes there, killing every positive graded
piece. Multiplicity-freeness prevents repetition of a line across layers.
\end{proof}

For the equivariant DS characters just used, $\Gamma$ is a subgroup of
the triple centralizer and acts on the generalized Verma $\widehat\g$-modules and
Jantzen ideals. In the classical BRST calculation the finite $G$-module
bundle over $N_\infty\times\mathcal S_\infty$ is trivialized by
$$
 (n,s,v)\longmapsto(n,s,n^{-1}v),
$$
where $n$ acts on the finite module through evaluation of its loop at
zero. Conjugation of $n$ and the natural actions on $s,v$ make this
trivialization $\Gamma$-equivariant and Kazhdan compatible. The common
ghost/universal factor has constant character one. Thus finite-dimensional highest-$x$ spaces have their Levi representations, including central
determinant twists. Taking invariants of the completed exact identity
is legitimate. 

\subsection{Invariant relations in \texorpdfstring{$H_f^0(I_k)$}{Hf0(Ik)}}
We apply the next lemma to homogeneous subspaces of the reduced
Jantzen ideals $H_f^0(J^1)\supseteq H_f^0(J^2)\supseteq\cdots$,
with $J^1=I_k$. Here invariance refers to the action of the specified
reductive Lie algebra by zero modes. The lemma determines when the
character sum forces an invariant relation to belong to $H_f^0(I_k)$.
Representation multiplicities and vertex products rule out additional
occurrences in later ideals. Once we have shown that the relation belongs to $H_f^0(I_k)$, invariant bilinear forms determine the coefficients of the relation modulo $C_2$.

\begin{lemma}\label{v3:scalar-detection}
Let $\mathfrak a$ be a reductive Lie algebra, let $U_d$ be a
finite-dimensional completely reducible $\mathfrak a$-module, and let
$K_1\supseteq K_2\supseteq\cdots$ be $\mathfrak a$-submodules of
$U_d$, with $K_i=0$ for sufficiently large $i$. Let
$B\subseteq K_1$ be a $\mathfrak a$-submodule and let $\C t$ be
a specified trivial $\mathfrak a$-submodule of $U_d/B$.
Suppose absence of that line from $K_1/B$ would place the classes of
$K_1/B$ and every later $K_i$ in an additive monoid $\mathcal M$ of
the representation ring. If
$\sum_i[K_i]-[B]\notin\mathcal M$, then $t\in K_1/B$.
An additive dimension congruence annihilating $\mathcal M$ suffices.
\end{lemma}
\begin{proof}
Under absence, the displayed difference is the sum of the classes
$[K_1/B]$ and $[K_i]$ for $i\ge2$, hence belongs to $\mathcal M$.
Corrections from $B$ already belong to the first ideal, so a
specified lift of the invariant line belongs to it. The residual
condition must include any possible repetitions of $B$ in later
layers. Subtracting a known subspace once does not by itself rule out such repetitions.
For a unique invariant line and $B=0$, a positive invariant multiplicity
in the actual sum already suffices by nestedness.
\end{proof}

For Virasoro quotients, the following consequence of the
Gorelik--Kac singular-vector theorem \cite[Theorem~0.3.2]{GK}
also gives the nilpotence bound needed below. We include the short
argument to keep track of that bound.

\begin{lemma}\label{v3:vir-quotient}
Every nonzero quotient of the universal Virasoro vertex algebra $\Vir^c$ by a nonzero
proper graded ideal is $C_2$-cofinite. If the ideal meets degree $B$,
then its conformal class satisfies $[T]^{\lfloor B/2\rfloor}=0$.
The same applies to a chain $\Vir^c\twoheadrightarrow H
\twoheadrightarrow H'\ne0$ when the second map has a nonzero kernel
in degree $B$: the conclusion and nilpotence hold in $H'$.
\end{lemma}
\begin{proof}
Pull the ideal back to $\Vir^c$ and choose a vector of least positive
degree. Positive Virasoro modes kill it. By \cite[Theorem 0.3.2]{GK}
its degree is $2N$ and its $L_{-2}^N\mathbf1$ coefficient is nonzero.
The standard basis of $\Vir^c$ consists of
$L_{-n_1}\cdots L_{-n_s}\mathbf1$ with $n_1\ge\cdots\ge n_s\ge2$.
Every such vector of degree $2N$, other than $L_{-2}^N\mathbf1$,
contains a factor $L_{-j}$ with $j\ge3$ and vanishes in
$R_{\Vir^c}=\C[T]$. Thus $[T]^N=0$, with $2N\le B$.

\end{proof}

\section{Drinfeld--Sokolov reductions of \texorpdfstring{$L_k(\g)$}{Lk(g)} in types \texorpdfstring{$A,B,C,D$}{A, B, C, D}}\label{seed:section}
\label{seed:classification}
Write $H_f=H_f^0(L_k(\g))$ and $\mathcal W_k(\g,f)$ for the simple
quotient of the universal $W$-algebra. If $H_f\ne0$,
\cref{com:simple-image} makes $\mathcal W_k(\g,f)$ a quotient of $H_f$,
which we use to prove non-quasi-lisseness. Conversely, relations in
$H_f$ and Slodowy-slice geometry exclude semisimple points from
$X_{L_k(\g)}$. We begin with calculations common to types $A,B,C,D$.

\subsection{Equivariant polynomial maps and nonvanishing of Drinfeld--Sokolov reductions}\label{seed:energy}
We apply \cref{com:energy} to the maximal proper graded ideal of $V^k(\g)$. Set
$\Delta_t(\lambda)=(\lambda,\lambda+2\rho)/(2t)$. If
$\Delta_t(\lambda)>\lambda(x)$ for every nonzero dominant root-lattice
weight, then $H_f^0(L_k(\g))\ne0$. If
$\Delta_t(\lambda)\ge\lambda(x)$ for every nonzero dominant
root-lattice weight, it remains to show that, for each $\lambda$
satisfying equality, every $G$-equivariant homogeneous polynomial map
$\g\to E_\lambda^*$ of degree $\Delta_t(\lambda)$ vanishes at $f$.
Here $E_\lambda$
is the finite-dimensional irreducible $\g$-module of highest weight
$\lambda$. Such maps are also called polynomial covariants. The
following bounds apply to the entire space of these maps, as required.

\begin{lemma}\label{seed:covariants}
Let $V$ be the defining representation.
\begin{enumerate}
\item For $\sog_{2n}$, if $f^{2m+1}=0$ and
$j(2m-1)<D<n-j$, every degree-$D$ equivariant polynomial map to
$\Lambda^{2j}V$ vanishes at $f$.
\item For $\sog_N$ with $N,r$ odd,
$\Hom_{\mathrm{SO}_N}(\Sym^D\sog_N,\Lambda^rV)=0$ if $2D<N-r$.
\item For $\spg(V)$, if $f^{2s}=0$, $r$ is even and
$D>r(s-1)$, every degree-$D$ map to the primitive exterior space
$\Lambda^r_0V$ vanishes at $f$.
\end{enumerate}
\end{lemma}
\begin{proof}
Use the orthogonal and symplectic first fundamental theorems
(\cite[arXiv v1, Theorem~2.1]{LehrerZhang2017} and \cite[arXiv v1, Theorem~4.8]{LehrerZhang2015}).
In the first case $2D+2j<2n$, so no epsilon tensor occurs. Contracting
endomorphism indices gives cycles and paths joining the $2j$ exterior indices.
A cycle is a positive trace power of a nilpotent and vanishes. Exterior
projection makes every path odd, and nilpotence bounds its length by
$2m-1$. The $j$ paths cannot have total length $D>j(2m-1)$.
In the second case the total number of slots $2D+r$ is odd. Pairings
cannot fill them and an epsilon tensor needs at least $N$ slots.
Two epsilon tensors reduce to pairings, so the stated bound is absolute.
In the symplectic case all invariant tensors are pairings, in every
rank. Alternating the $r$ exterior indices leaves only paths of even length.
A length-zero path is a
symplectic factor and disappears under primitive projection. The $r/2$
paths have length at most $2s-2$, proving the bound. Relations among
contraction diagrams in small rank cannot create additional tensors.
\end{proof}

Use the common endpoint identity of \cref{v3:master-energy}, with
$\sigma=0,1,2$ for $D,B,C$, respectively, and
$k+h^\vee=T,T,T/2$, respectively.
Write $A_l=l(2n+\sigma-l)-T\sum_{i\le l}\eta_i$, $\eta=2x$.
It suffices to establish the inequalities at the endpoints of each
block on which the spectrum is constant. Interpolation
adds $u(L-u)>0$ inside a block. The zero tail has
$A_{L+u}=A_L+u(2(n-L)+\sigma-u)>0$.
Below $F(C_l)=A_l$ for an orthogonal column and $F(\omega_l)=A_l$
for a symplectic fundamental column. Equality can only be a column $(1^l)$ with $A_l=0$, subject to even
column length in the $C,D$ root lattices. The preceding lemma applies
to all equivariant polynomial maps of the required degree for each
weight satisfying equality.
\label{seed:D-columns}\label{seed:B-columns}\label{seed:C-columns}

\subsection{Symplectic-fermion invariants as Miura images and tensor factors}\label{seed:scalar}
Recall the invariant subalgebra
$\mathcal A_s=\SF(s)^{\Sp_{2s}}$ from Section~\ref{com:scalar}.
We identify certain classical Miura images with $\mathcal A_s$
and establish conditions under which it occurs as a tensor factor. 
\begin{proposition}\label{v6:scalar-partitions}
For the following pairs $(\g,f)$, the simple vertex algebras obtained
from the Miura homomorphism by replacing the universal affine
vertex algebras associated with $[\g_0,\g_0]$ by their simple
quotients are $\mathcal A_s$:
$$
\begin{array}{c|c|c}
 \g&t&f\\\hline
 C_n, n=sP&P/2&((2s)^P)\\
 \sog_N, N=(2s+1)P+2&P&((2s+1)^P,1^2)
\end{array}
$$
In the orthogonal row $P$ is odd in type $B$
and even in type $D$.
\end{proposition}
\begin{proof}
Let $R_d$ denote the $d$-dimensional irreducible
$\mathfrak{sl}_2$-module and let $U=\C^P$ carry its symmetric form.
The standard representation is $R_{2s}\otimes U$ in the symplectic
case and $(R_{2s+1}\otimes U)\oplus E$, with $E=\C^2$ orthogonal,
in the orthogonal case. We use the Dynkin grading. It is integral,
and its degree-zero subalgebra is
$$
 \g_0=\glg_P^s\quad\hbox{or}\quad
 \g_0=\glg_P^s\oplus\sog_{P+2},
$$
respectively. The Miura bilinear form is
$\tau_k=t(\cdot,\cdot)_\g-\frac12\kappa_{\g_0}$, where
$\kappa_{\g_0}$ is the Killing form. Its restriction to each
$\slg_P$ is $(2t-P)\operatorname{tr}$ in the symplectic case and
$(t-P)\operatorname{tr}$ in the orthogonal case. In the latter case
its restriction to $\sog_{P+2}$ is also zero. Thus all derived
factors have level zero, including both simple factors when $P=2$.
Their simple quotients are trivial. Denote the resulting homomorphism
by $\overline\mu$ and its image by $B_s$.

Let $I_i$ be the identity matrix in the $i$th $\glg_P$ summand.
In both cases $\tau_k(I_i,I_j)=P^2\delta_{ij}$, so
$h_i=I_i/P$ satisfies $[h_i{}_{\lambda}h_j]=\delta_{ij}\lambda$.
The image of the conformal vector is
$$
 T_s=\frac12\sum_{i=1}^s(:h_i^2:-\partial h_i).
$$
For completeness, this follows directly by projecting the BRST
conformal vector. Its quadratic part on the centers is
$\frac12\sum_i:h_i^2:$, and its derivative term is the image of
$\partial J^{x-(\rho-\rho_0)^\sharp/t}$, where $\rho_0$ is the
half-sum of the positive roots of $\g_0$ and $\sharp$ uses the
normalized form on $\g$. If $x_i$ is the grading value on the $i$th
positive block and $q_i$ is the mean of the coordinates of
$\rho-\rho_0$ on that block, then
$$
\begin{array}{c|c|c|c}
 &x_i&q_i&(D_j,D_j)_\g\\\hline
 C_{sP}&s-i+\frac12&P(s-i+\frac12)+\frac12&2\\
 \sog_{(2s+1)P+2}&s-i+1&P(s-i+1)+\frac12&1.
\end{array}
$$
Here $D_j$ has eigenvalues $1,-1$ on a paired pair of standard
coordinates. In either row, the coefficient of $I_i$ in
$x-(\rho-\rho_0)^\sharp/t$ is $-1/(2P)$, as required.

The bound on the strong generators follows from equivariance and the strong
generators of \cite[Theorem~4.1]{KacWakimoto2004}. In the
symplectic case use the reductive group $K=O(U)$ centralizing the
triple. The strong generating space indexed by $\g^f$ has
$K$-modules
$$
 \operatorname{Sym}^2U\quad\hbox{in weights }2,4,\ldots,2s,
 \qquad
 \Lambda^2U\quad\hbox{in weights }1,3,\ldots,2s-1.
$$
These follow by decomposing
$\operatorname{Sym}^2(R_{2s}\otimes U)$ under the triple.
In the orthogonal case use $K=S(O(U)\times O(E))$. Decomposing
$\Lambda^2((R_{2s+1}\otimes U)\oplus E)$ gives
$$
\begin{array}{c|c}
 K\text{-module}&\text{conformal weights}\\\hline
 \operatorname{Sym}^2U&2,4,\ldots,2s\\
 \Lambda^2U&1,3,\ldots,2s+1\\
 U\otimes E&s+1\\
 \Lambda^2E&1.
\end{array}
$$
The BRST construction and the Miura homomorphism are $K$-equivariant.
One may choose the strong generators equivariantly: their prescribed
symbols have homogeneous lifts, and reductivity splits the
corresponding surjections in each finite-dimensional conformal-weight
space. The quotient $\overline\mu$ has trivial $K$-action because
$K$ fixes every $I_i$ and all other factors have been killed.
It therefore annihilates the nontrivial $K$-summands of this generating
space. The invariant part of $\operatorname{Sym}^2U$ is its trace
line, and $(\Lambda^2U)^{O(U)}=0$.
In the orthogonal case $SO(E)$ has no invariant vectors in $U\otimes E$, and simultaneous
reflections in $U$ and $E$ rule out contributions from $\Lambda^2E$ and also
$\Lambda^2U$ when $P=2$. For $P\ge3$ the latter has no
$SO(U)$-invariant vectors, and for $P=1$ it is zero.
Consequently $B_s$ is strongly generated by at most one field in
each weight $2,4,\ldots,2s$, and by no other fields.

For $s\ge2$, consider the invariant polynomial
$p_4(X)=\operatorname{tr}_{\mathrm{std}}(X^4)$.
Its restriction to the Slodowy slice $\mathcal S_f$ has conformal
degree four. The Li-filtration isomorphism
$R_{\mathcal W^k(\g,f)}\cong\C[\mathcal S_f]$ therefore gives
a homogeneous lift $W\in\mathcal W^k(\g,f)_4$ of this polynomial.
The induced Miura map on these $C_2$-algebras is the restriction
$$
 \C[\mathcal S_f]\simeq
 \C[f+\g_{\ge0}]^{G_{>0}}
 \longrightarrow\C[f+\g_0]
$$
by \cite[Proposition~2.8 and Section~2.5, equation~(2.8)]{NakatsukaMiura}.
For $y\in\g_0$, every term containing $f$ in the expansion of
$\operatorname{tr}(f+y)^4$ has negative grading, so
$$
 p_4(f+y)=p_4(y).
$$
After restriction to the centers, the eigenvalues are
$z_1,\ldots,z_s,-z_s,\ldots,-z_1$, each of multiplicity $P$,
and additionally zero eigenvalues in the orthogonal case. Thus the
restriction is $2P\sum_i z_i^4$. The common nonzero normalization
of the center coordinates gives
$$
 [\overline\mu(W)]=\gamma\sum_{i=1}^s h_i^4
 \quad\hbox{in }R_{\mathcal H^{\otimes s}},\qquad \gamma\ne0.
$$
Thus the component of $\overline\mu(W)$ containing four factors from $h_1,\ldots,h_s$ and no derivatives is exactly
$\gamma\sum_i:h_i^4:$. Its image in $R_{\mathcal H^{\otimes s}}$
is independent of $(\sum_i h_i^2)^2$ for $s\ge2$.

It remains to prove $B_s\subseteq\mathcal H_1\otimes B_{s-1}$,
as required for Proposition~\ref{com:sf-recognition}. Delete the first positive
grading block and its paired negative block. For $s\ge2$, the
resulting standard Levi subalgebra is
$$
 \glg_P\oplus C_{(s-1)P}
 \quad\hbox{or}\quad
 \glg_P\oplus\sog_{(2s-1)P+2}.
$$
The deleted simple root joins consecutive grading values, so it
has degree one. The restricted grading is the Dynkin grading of
$0\oplus((2s-2)^P)$ or of
$0\oplus((2s-1)^P,1^2)$, respectively.
Theorem~5.10 of \cite{Genra} applies at the stated levels and its
induced level $\ell$ on the smaller simple Lie algebra $\g'$
satisfies $\ell+h^\vee_{\g'}=P/2$ or $P$, respectively. Its compatibility with the Miura maps,
followed by the same derived-factor quotients, gives
$$
 B_s\subseteq\mathcal H_1\otimes B_{s-1}.
$$
For $s=1$, the bound on the strong generators and $T_1\in B_1$ show that
$B_1$ is generated by $T_1$ alone. It is the simple universal
Virasoro algebra at $c=-2$, hence $B_1=\mathcal A_1$ in the vertex algebra generated by $h_1$.
Inductively $B_{s-1}=\mathcal A_{s-1}$. The containment, the
bound on the strong generators and the degree-four polynomial computed above
are precisely the hypotheses of Proposition~\ref{com:sf-recognition}.
It follows that $B_s=\mathcal A_s$ for every $s\ge1$.
\end{proof}

\begin{proposition}\label{seed:scalar-split}
Suppose a conformal subalgebra $B\subset U\otimes\mathcal A_s$
has conformal vector $T_U+T_s$, with $U_0=\C$. Assume that it contains
$Q\in B_2$ and $P\in B_4$ such that
$$
 Q-1\otimes T_s\in\partial B_1,\qquad
 \overline{\operatorname{pr}_{U_0\otimes\mathcal A_s}P}
 \notin\C(\textstyle\sum_i h_i^2)^2\quad(s\ge2).
$$
The bar denotes the image in $R_{\mathcal H^{\otimes s}}$.
If $(T_U)_{(1)}$ is semisimple with nonnegative spectrum on every
finite-dimensional eigenspace of $(T_U+T_s)_{(1)}$, then
$$
 B=(B\cap(U\otimes1))\otimes\mathcal A_s.
$$
\end{proposition}
\begin{proof}
Subtracting the derivative in that hypothesis gives $1\otimes T_s\in B$.
The operator $(T_U)_{(1)}=(T_U+T_s)_{(1)}-(T_s)_{(1)}$ preserves $B$.
On each finite-dimensional eigenspace of $(T_U+T_s)_{(1)}$, the
projection to $\ker(T_U)_{(1)}$ is a polynomial in $(T_U)_{(1)}$.
This projection sends $P$ into $1\otimes(\mathcal A_s)_4$ with
symbol independent of $(\sum_i h_i^2)^2$ for $s\ge2$.
The generation of $\mathcal A_s$ gives $1\otimes\mathcal A_s\subset B$.
The density theorem for the adjoint representation of the simple
vertex algebra $\mathcal A_s$ gives
$B=(B\cap(U\otimes1))\otimes\mathcal A_s$.
The restriction of any nonzero graded quotient of $B$ to
$1\otimes\mathcal A_s$ is injective because $\mathcal A_s$ is simple.
\end{proof}

The two field hypotheses can be checked simultaneously for
$$
 \begin{array}{c|c|c|c}
 \g&f&t&U\\\hline
 \sog_N&((2s+3)^r,(2s+1)^\epsilon,1^2)&P&L_\epsilon(\slg_r)\otimes\mathcal H\\
 \spg_{2n}&((2s+2)^r,(2s)^\epsilon)&P/2&L_\epsilon(\slg_r)\otimes\mathcal H
 \end{array}
 \qquad P=r+\epsilon,\quad \epsilon=1,2.
$$
Here $s,r\ge1$, and the affine factor is omitted for $r=1$.
Parabolic induction for the outer $\glg_r$ block gives the containment
in $U\otimes\mathcal A_s$. Here $\mathcal H$ corresponds to its center,
with the induced conformal vector. After taking the simple quotients
of the inner derived factors, the shorter-chain diagonal is
$(d_1,\ldots,d_L)=(h_1,\ldots,h_s,0,-h_s,\ldots,-h_1)$ for
$L=2s+1$, and $(h_1,\ldots,h_s,-h_s,\ldots,-h_1)$ for $L=2s$.
Set $H_a=\sum_{b\le a}d_b$ and $d_{L+1}=0$. Since $H_L=0$,
\begin{equation}\label{seed:telescope}
 -\frac12\sum_{a=1}^{L}H_a(d_{a+1}-d_a)
       =\frac12\sum_{a=1}^{L}d_a^2
       =\sum_{i=1}^s h_i^2.
\end{equation}
To identify this expression with the Miura symbol of the shorter-block
trace, write $f_{\rm sh}$ for the restriction of $f$ to those blocks,
extended by zero. On a shorter chain choose $f$ to have unit
subdiagonal entries and choose the degree-one matrix $X$ with
superdiagonal entries $-H_a$. In the paired chain bases these matrices
preserve the orthogonal or symplectic form, respectively, since
$d_{L+1-a}=-d_a$ and $H_{L-a}=H_a$. Then $[X,f]=-y$ on that chain.
The degree-one part of $\exp(\operatorname{ad}X)(f+y)$ is
$\frac12[X,y]$. On this chain, the trace of its product with $f$
is the left side of \eqref{seed:telescope}. The remaining transformations
to the Slodowy slice may be chosen block diagonally. Transformations
on longer blocks have zero pairing with $f_{\rm sh}$, and subsequent
degree-two corrections contribute $[Z,f]$, whose pairing is zero
because $[f_{\rm sh},f]=0$. Higher-degree terms do not pair with
$f_{\rm sh}$. This computes the complete quadratic restriction of
the linear slice coordinate indexed by $f_{\rm sh}$. In particular,
there are no quadratic or mixed terms involving $U$.

Take a homogeneous lift of this coordinate in the universal
$W$-algebra and rescale its image $Q$ so that its quadratic symbol
is $\frac12\sum_i h_i^2$. The invariant form and the $\epsilon$
identical shorter blocks only contribute a nonzero common scalar.
In the vertex algebra $U\otimes\mathcal H^{\otimes s}$, the kernel
of the $C_2$ projection in weight two is the derivative of its
weight-one space. Hence $Q-1\otimes T_s=\partial j$ for some
vector $j\in(U\otimes\mathcal H^{\otimes s})_1$. Since $L_1T_s=0$ and
$L_1\partial j=2j$, one has
$$
 Q-1\otimes T_s=\tfrac12\partial(L_1Q)\in\partial B_1.
$$
This also shows that a quasiprimary choice of $Q$ has image
$1\otimes T_s$ exactly.

For $s\ge2$, lift $\operatorname{tr}_{\rm std}(X^4)|_{\mathcal S_f}$
as in Proposition~\ref{v6:scalar-partitions}. On the degree-zero
algebra its Miura symbol is $\operatorname{tr}_{\rm std}(y^4)$.
Setting the complementary $U$ coordinates to zero leaves each
$h_i,-h_i$ with multiplicity $P$, so its restriction is a nonzero
multiple of $\sum_i h_i^4$. Projection to $U_0\otimes\mathcal A_s$
induces on $R_U\otimes R_{\mathcal H^{\otimes s}}$ the positive-grading augmentation
of $R_U$ tensored with the identity. Thus the projected weight-four
field has this same quartic symbol, independent of
$(\sum_i h_i^2)^2$. These computations verify both field hypotheses
of Proposition~\ref{seed:scalar-split} for the stated applications.
For families III and IV of \cref{B:odd-exclusions} with $j\ge2$, take
$s=j-1$ and $r=p-\epsilon$ with $\epsilon=1,2$, respectively.
For the odd-rank numerator-two row of \cref{sec:D-exclusions},
take $r=\epsilon=1$. For its row $n=(2m+1)p/2$, $m\ge2$,
take $s=m-1$, $r=p-1$ and $\epsilon=1$. In the second row of
\eqref{seed:C-negative-list}, take $s=j$, $r=P-c$, $\epsilon=c$.

\subsection{Invariant bilinear forms on simple quotients of \texorpdfstring{$W$}{W}-algebras}\label{seed:negative-families}
\begin{lemma}\label{seed:integral-form}
Let $S=\mathcal W_k(\g,f)$, where $\g$ is a simple Lie algebra,
$k\ne-h^\vee$ and the grading is the nonnegative half-integral
conformal grading associated with the Dynkin grading of $f$.
Then its integer-weight subalgebra $S^{\mathrm{int}}$ is a simple
CFT-type VOA with a normalized nondegenerate invariant form.
Consequently the nondegeneracy argument in \cref{com:zero-current} and the
integral-weight null-vector arguments below apply to the full spaces
$S_d$, $d\in\Z_{\ge0}$.
\end{lemma}
\begin{proof}
The strong generators of $\mathcal W^k(\g,f)$ have weights
$1-j$ for $\g_j^f$, where $j\le0$. Their images strongly generate
$S$. Hence $S_0=\C$, no weights occur in $(0,1)$, every homogeneous
subspace is finite dimensional, and $S_1$ is the image of $\g^\natural$.
For $a\in\g^\natural$, the coefficient of the third-order pole in the
operator product of the conformal vector and the corresponding generator $J^a$ is
$\frac12\tr_{\g_{>0}}(\ad a)-(k+h^\vee)(a,x)$
\cite[Theorem~2.1(b)]{KacWakimoto2004}. Both terms vanish:
$a$ centralizes the triple, and in
$\g=\bigoplus_d V_d\otimes M_d$ every multiplicity space $M_d$
is self-dual, so $\tr_{M_d}a=0$. This includes the $\sog_2$ centers.
Thus $L_1S_1=0$. The automorphism $e^{2\pi iL_0}$ has order at most
two and its fixed algebra is $S^{\rm int}$. Countable quantum Galois
theory \cite[arXiv v1, Theorem~5.5]{DongRenYang2024} makes this fixed algebra simple.
Li's invariant-form theorem \cite{LiForms1994}, recalled in
\cite[Proposition~4.6]{CarpiEtAl2018}, then gives its normalized nondegenerate form. Its integral weight
spaces are the entire $S_d$. Thus a vector orthogonal to all of its
integral weight space is zero in $S$ and acts trivially also on
half-integral fields.

\end{proof}
Each vanishing assertion below for a Lie subalgebra of $S_1$ uses
its orthogonality to all of $S_1$. The lemma therefore
justifies \cref{com:zero-current}, including a central $\sog_2$
generator. Its self-pairing alone would not suffice.
We shall also use the following composition of maps between Jordan blocks.
For the partitions $((E+1)^a,E^b,1^c)$ occurring below, write the defining
module as
$$
 V=(R_{E+1}\otimes U)\oplus(R_E\otimes W)\oplus B,
 \qquad R_d=\C[z]/z^d,\quad
 \dim U=a,\quad\dim W=b,\quad\dim B=c.
$$
Let $\mathfrak a_W$ and $\mathfrak a_B$ be the orthogonal or symplectic
Lie algebras preserving the induced forms on $W$ and $B$, respectively.
For this argument, assume $\mathfrak a_W\ne0$ and that $W$ has no
$\mathfrak a_W$-invariant vectors. If the generators corresponding to
$\mathfrak a_W$ vanish, their zero modes therefore force the generators
in the modules $U\otimes W$ and $W\otimes B$ to vanish. The maps
$$
 R_1\xrightarrow{z^{E-1}}R_E\xrightarrow{z}R_{E+1},
 \qquad1\longmapsto z^E\ne0
$$
compose nontrivially. Contracting the nondegenerate form on $W$ therefore
spans the remaining module $U\otimes B$ in the centralizer bracket.
The level-independent linear part of the Slodowy Poisson bracket,
\cref{com:linear-bracket}, gives the corresponding nonzero coefficient
in the zeroth vertex product. To control its lower-weight corrections,
use the involution induced by $-\id_B$ on $B$ and the identity on the other
two summands of $V$. Its negative eigenspaces among the strong generators
are precisely the modules involving $W\otimes B$ or $U\otimes B$.
Every correction with eigenvalue $-1$ contains a lower-weight generator
with that eigenvalue, which has already vanished. Hence the generators
in $U\otimes B$ vanish as well. Apart from the generators corresponding
to $\mathfrak a_B$, all remaining strong generators are annihilated by
$J^a_{(0)}$ for $a\in\mathfrak a_B$. The Sugawara correction and the density
theorem \cref{com:neutral-split} give a tensor product decomposition
containing the simple affine vertex algebra associated with
$\mathfrak a_B$, without assuming rationality of that factor.

\begin{proposition}\label{com:orthogonal-point}
Let $f\in\sog(V)$ have partition $(3^a,2^m,1^4)$, where $a\ge3$
and $m\ge0$ is even. Write
$$
 V=(S_2\otimes U)\oplus(S_1\otimes W)\oplus B,
 \qquad \dim U=a,\quad\dim W=m,\quad\dim B=4,
$$
where $S_i$ is the $(i+1)$-dimensional $\slg_2$-module,
$U,B$ are orthogonal and $W$ is symplectic.
Let $k\ne-h^\vee$ and let $H\ne0$ be a graded quotient of
$\mathcal W^k(\sog(V),f)$ satisfying the grading hypotheses of
\cref{com:scalar-point}. Suppose the classes of the generators $J^a$, $a\in\g^\natural$,
are nilpotent in $R_H$.

If $m>0$, then $H$ is $C_2$-cofinite. If $m=0$, its reduced
$C_2$ algebra is generated by one scalar of weight two. In that case,
an actual positive-weight homogeneous relation that becomes a nonzero
polynomial in that scalar after setting the classes of $J^a$, $a\in\g^\natural$, and all
nontrivial generating modules to zero also makes
$H$ $C_2$-cofinite.
\end{proposition}
\begin{proof}
Set $\mathfrak a=\sog(U)\oplus\spg(W)\oplus\sog(B)$, omitting
the middle factor if $m=0$. Decomposing $\Lambda^2V$ under the
triple gives the following weights and $\mathfrak a$-modules of
strong generators of $\mathcal W^k(\sog(V),f)$. Their images strongly
generate $H$:
$$
\begin{array}{c|l}
\text{weight}&\mathfrak a\text{-module}\\\hline
1&\Lambda^2U\oplus\Sym^2W\oplus\Lambda^2B\\
3/2&(U\otimes W)\oplus(W\otimes B)\\
2&\Sym^2U\oplus\Lambda^2W\oplus(U\otimes B)\\
5/2&U\otimes W\\
3&\Lambda^2U.
\end{array}
$$
The only scalar summands are the trace in $\Sym^2U$ and, if
$m>0$, the trace in $\Lambda^2W$. All other generating classes
vanish in $(R_H)_{\mathrm{red}}$ by \cref{com:scalar-point}.

For completeness, the two traces when $m>0$ follow from elementary
centralizer brackets. Choose chain bases $u_+,u_0,u_-$ in $S_2$
and $w_+,w_-$ in $S_1$, with
$$
 fu_+=u_0,\quad fu_0=u_-,\quad fw_+=w_-,\qquad
 (u_+,u_-)=1,\quad(u_0,u_0)=-1,\quad(w_+,w_-)=1.
$$
Write $\beta$ for the alternating form on $W$ and identify
$\sog(V)=\Lambda^2V$. The weight-$3/2$ centralizer elements are
$$
 x(u,w)=(u_0\otimes u)\wedge(w_-\otimes w)
       -(u_-\otimes u)\wedge(w_+\otimes w),\qquad
 y(w,b)=(w_-\otimes w)\wedge b.
$$
Set
$$
 \tau_W(w,w')=(w_-\otimes w)\wedge(w_-\otimes w'),
$$
$$
 \tau_U(u,u')=(u_-\otimes u')\wedge(u_0\otimes u)
             +(u_-\otimes u)\wedge(u_0\otimes u').
$$
The orthogonal bracket gives
$$
 [y(w,b),y(w',b')]=-(b,b')\tau_W(w,w'),
$$
$$
 [x(u,w),x(u',w')]=(u,u')\tau_W(w,w')
                  +\beta(w,w')\tau_U(u,u').
$$
Thus these brackets span both trace directions. In the universal
$C_2$ algebra, the cotangent bracket is the centralizer bracket,
by \cref{com:linear-bracket,vtwo:cotangent}. Its correction to
each displayed linear term has weight two and no linear part,
so it is a quadratic polynomial in the classes of $J^a$, $a\in\mathfrak a$.
Those terms vanish in the reduced $C_2$ algebra, as do the
brackets of the vanished cross fields. Both
scalar classes are therefore zero. When $m=0$ only the first
trace remains, and a nonzero homogeneous relation is a nonzero
multiple of a positive power of that class.
\end{proof}

\subsection{Pfaffian singular vectors in \texorpdfstring{$V^k(\mathfrak{so}_N)$}{Vk(soN)}}
\label{seed:orth-positive}
For a Pfaffian of full size, the critical-level assertion below is
Molev's Segal--Sugawara construction
\cite[Theorem~2.1]{MolevPfaffian2011}. We give the commutator calculation
in $V^k(\sog_N)$, including the powers of Pfaffians needed here.
\begin{lemma}\label{lem:orthogonal-pfaffian}\label{lem:D-pfaffian}
For $\g=\sog_N$, the $\g$-submodule spanned by the Pfaffians of
size $2s$ in $V^k(\g)$ is annihilated by
$\g\otimes z\C[z]$ at $k=2-2s$. Its image in $\Sym\g$ defines
the locus $\rank X\le2s-2$.
If $2s\le\lfloor N/2\rfloor$, the square of the highest isotropic
Pfaffian is singular at $k=3-2s$. Its $G$-translates impose the same
rank bound. At a middle exterior degree both orthogonal halves are used.
\end{lemma}
\begin{proof}
Let $F_{ab}$ be the orthogonal matrix generators and let $P_I$ be the
Pfaffian in the commuting disjoint-pair operators $F_{ab}(-1)$ on
the vacuum, for $|I|=2s$. The commutator formula gives
$$
 F_{ab}(1)P_I=-(k+2s-2)\frac{\partial P_I}{\partial F_{ab}},
 \qquad F_{ab}(r)P_I=0\quad(r\ge2).
$$
The direct contraction of a matching contributes $-k$. The two
alternate matchings through each other pair contribute $-2$.
For an index outside $I$, the remaining operator with index zero
acts on disjoint factors and is zero. This proves the first assertion.
The polynomial Pfaffians cut out the stated rank locus.
For the highest isotropic Pfaffian $P_s$, the first two commutators
with $f_\theta(1)$ give $(k+2s-2)P_{s-1}$ and $-2P_sP_{s-1}$,
and the third is zero. Therefore
\begin{equation}\label{seed:pf-power}
 f_\theta(1)P_s^a
   =a(k+2s-1-a)P_s^{a-1}P_{s-1}.
\end{equation}
The finite simple-root raising operators annihilate $P_s^2$.
At $a=2$, $k=3-2s$, equation~\eqref{seed:pf-power} also gives
$f_\theta(1)P_s^2=0$. These are the simple-root raising operators
of $\widehat\g$, so $P_s^2$ is singular. Commutation with
$\g\otimes1$ then gives the assertion for the finite-dimensional
$\g$-submodule it generates.
Outside the middle exterior degree, the $G$-span of the highest
exterior vector is the space spanned by these Pfaffians. In the
middle exterior degree, this space has two irreducible summands,
interchanged by the outer reflection, and both highest-weight vectors
are used. The translated squares vanish precisely when all these
Pfaffians vanish, giving the rank bound over $\C$.
\end{proof}

\section{\texorpdfstring{Quasi-lisse simple affine vertex algebras of type $A$}{Quasi-lisse simple affine vertex algebras of type A}}\label{sec:type-A}

Explicit singular-vector and Zhu-algebra calculations at integer levels,
including $k=-1$, were carried out by Adamovi\'c and Per\v se in
\cite{APNonrational2008}. Drinfeld--Sokolov reduction shows uniformly that the remaining
nonadmissible levels are not quasi-lisse: we find a nonzero reduction
with a simple quotient in whose $C_2$ algebra the class of the
conformal vector is not nilpotent.

For the boundary-admissible family $D_q(A_n)$, the partition is $(q^a,b)$, where $n+1=aq+b$ and
$0\le b<q$, with a zero part omitted.
The examples $L_{-4/3}(A_1)$ and $L_{-3/2}(A_2)$ correspond to
$H_1=(A_1,A_3)$ and $H_2=(A_1,D_4)$, with minimal nilpotent
Higgs branches. The family $L_{-2+2/(2m+1)}(A_1)$, $m\ge1$,
corresponds to $(A_1,D_{2m+1})$, with Higgs branch
$\C^2/\mathbb Z_2$, see \cite[Sections~4 and~5.2]{BeemRastelli2018}.

For $L_k(\slg_N)$ with $k=-N+p/q$, where $p,q\in\Z_{>0}$ are
coprime, admissibility covers precisely $p\ge N$.
By \cref{thm:denominator}, it remains to prove that
$L_{p-N}(\slg_N)$ is not quasi-lisse for $1\le p<N$.

\begin{proposition}\label{seed:A}
For $1\le p<N$, write $N=ap+r$, $0\le r<p$, and set $s=p-r$ and
$f=((a+1)^r,a^s)$. Then $H_f^0(L_{p-N}(\slg_N))\ne0$, and the unique simple graded
quotient of this reduction is not quasi-lisse. For $p\ge N$, the simple affine vertex algebra $L_{p-N}(\slg_N)$ is
integrable.
\end{proposition}
\begin{proof}
We first prove that the reduction is nonzero. We then distinguish the
nonrectangular case, where a Heisenberg tensor factor survives, from the
rectangular case, where the Miura image gives the obstruction directly.
The eigenvalues of the grading element $x$, in decreasing order, have
alternating constant blocks of lengths $r$ and $s$.
Set $D_j=j(N-j)-2p\sum_{i\le j}x_i$. At the block endpoints
$bp,bp+r$ it is zero. Inside a block of length $L$ it is $u(L-u)$.
Thus $(\rho-px,\omega_j)=D_j/2\ge0$, and for every nonzero dominant
$\lambda$ and $0<t\le p$,
$$
 \Delta_t(\lambda)-\lambda(x)=
 \frac{\|\lambda\|^2}{2t}+\frac{(\rho-px,\lambda)}t
 +(p/t-1)\lambda(x)>0.
$$
The nonvanishing assertion follows from \cref{com:energy}.

If $r>0$, use the even left-aligned pyramid with columns $p^a,r$.
A commuting map from a chain of length $L$ to one of length $M$ starts
at $q\ge\max(0,M-L)$ and has degree $-q\le0$, verifying goodness.
The invariant bilinear form in the Miura homomorphism is
$p\tr-\frac12\mathrm{Kill}_{\lev}$. Composing with
$V^0(\slg_p)\twoheadrightarrow L_0(\slg_p)=\C$ in each corresponding
tensor factor gives a homomorphism with values in $V^{p-r}(\slg_r)\otimes\mathcal H$. The factor
$V^{p-r}(\slg_r)$ is omitted when $r=1$. The centralizer element acting by
$as$ on long chains and $-(a+1)r$ on short chains has scalar averages
$-rs/p$ on a full column and $as$ on the last one. Its image is a vector $J$ generating a rank-one Heisenberg
vertex subalgebra and satisfying
$$
 J_{(1)}J=p\bigg(ap\frac{r^2s^2}{p^2}+ra^2s^2\bigg)\mathbf1
 =ars^2N\mathbf1.
$$
Here $ars^2N>0$.
The operator $J_{(0)}$ annihilates $V^{p-r}(\slg_r)\otimes\mathcal H$
and hence the unique simple graded quotient of the image. Its nonzero level preserves $J$ there, so the
decomposition of \cref{com:heisenberg} gives a tensor factor
isomorphic to the rank-one Heisenberg vertex algebra generated by $J$.

If $r=0$, the rectangular Miura polynomial generators specialize,
after coefficient extraction, to the elementary symmetric polynomials
$e_d(h_1,\ldots,h_a)$, with $\sum h_i=0$ and Gram matrix
$\delta_{ij}-1/a$. Derivative terms carry an extra specialization
parameter. The generators and their formulas at arbitrary level are given in
\cite[Theorem~3.1 and Corollary~3.2]{ArakawaMolev2017}. The bilinear form in their equation~(2.1) is
$\kappa(X,Y)=k(\tr(XY)-\tr X\,\tr Y/N)$, whose radical contains the
central line $\C I_N$. This line commutes with $\slg_N$ and is
orthogonal to it. The BRST complex therefore splits off the
commutative level-zero affine vertex algebra of $\C I_N$, on which
the BRST differential is zero. Passing to $\slg_N$ means taking the
quotient in which this central current vanishes, giving $\sum_i h_i=0$.
Here $\rho-\rho_0=px$, so $-e_2=\frac12\sum h_i^2$ is the standard
conformal vector. PCT sends $e_d$ to $(-1)^de_d$. The image of this Miura homomorphism is
therefore simple by restriction of the positive Hermitian form. The class of its conformal vector in $R_{\mathcal H}$, where $\mathcal H$ is generated by the $h_i$
subject to $\sum_i h_i=0$, is nonzero, for example at $(1,-1,0,\ldots)$, and therefore is not
nilpotent. The rank of $\mathcal H$ is $a-1>0$ because $p<N$.
\end{proof}

\section{\texorpdfstring{Quasi-lisse simple affine vertex algebras of type $B$}{Quasi-lisse simple affine vertex algebras of type B}}\label{sec:type-B}
Wang and Xie propose boundary admissible affine vertex algebras
$$
 L_{-(2n-1)+(2n-1)/q}(B_n),\qquad
 q>1\text{ odd},\quad\gcd(q,2n-1)=1,
$$
from twisted $A_{2n-1}$ theories with a full regular puncture
\cite[Section~5]{WangXie2019}.
The rank-one list of Li, Xie and Yan includes $L_{-2}(B_3)$,
from a twisted $A_5$ construction, with predicted Higgs branch
$\overline{\OO_{(3,1^4)}}$
\cite[Table~4 and Appendix~B.1]{LiXieYan2023}.
Its associated variety is also discussed in \cite{ArakawaEtAl2025}.

Let $\g=\sog_{2n+1}$, $n\ge2$, and set $t=k+2n-1$.

\begin{theorem}\label{thm:B-classification}
For $k\ne-2n+1$, the vertex algebra $L_k(B_n)$ is quasi-lisse if and
only if $t=p/q\in\Q_{>0}$, with $p,q\in\Z_{>0}$ and
$\gcd(p,q)=1$, and
\begin{equation}\label{eq:B-classification}
 \begin{cases}
  2n\le3(p-1),&p,q\text{ both odd},\\
  p>n,\quad p\ne2n-1,&\text{otherwise}.
 \end{cases}
\end{equation}
\end{theorem}
\begin{proof}
By \cref{thm:denominator} and the exclusions in the introduction, it suffices
to treat $t=p$ for odd $q$ and $t=p/2$ for even $q$ (hence odd $p$).
Admissibility is covered by \cite{AAV}, and $B_2=C_2$ by
\cref{thm:C-seeds}, with $t=k+3$. Below we treat the nonadmissible
families for $n\ge3$.
\end{proof}

For $1\le j\le n$, write
$$
 H_j=(\underbrace{1,\ldots,1}_{j},0,\ldots,0)
$$
in the standard Cartan coordinates. Its derived centralizer has type
$A_{j-1}+B_{n-j}$, with absent factors omitted and rank-one factors
normalized as in the introduction. By \cref{com:finite-centers},
$X_{L_k(B_n)}\subseteq\N(B_n)$ holds if and only if
$H_j\notin X_{L_k(B_n)}$ for every $1\le j\le n$.
Write $H_f=H_f^0(L_k(\g))$. For the levels not satisfying \eqref{eq:B-classification}, other than $t=1$,
we prove $H_f\ne0$ and exhibit a non-quasi-lisse simple quotient.

The proofs of quasi-lisseness use singular vectors and their restrictions
in Sections~\ref{B:relations}--\ref{B:half-integral} for the ranges
$t=p\ge n$ and $t=p/2$. The remaining odd integers $p<n$ satisfying
$2n\le3(p-1)$ are treated by reductions and relations in
Section~\ref{B:lower-integral}. Sections~\ref{B:shift-one}--\ref{seed:B-coset}
establish failure of quasi-lisseness in the complementary ranges.

\subsection{\texorpdfstring{Non-quasi-lisseness of $L_{-2n+2}(\sog_{2n+1})$}{Non-quasi-lisseness of L(-2n+2)(so(2n+1))}}\label{B:shift-one}
The symplectic-boson vertex algebra associated with
$\C^{2n+1}\otimes\C^{2n-2}$ admits homomorphisms from
$V^{-2n+2}(\sog_{2n+1})$ and $V^{-(2n+1)/2}(\spg_{2n-2})$
whose images commute. For the second algebra, $k+h^\vee=-1/2$.
Every Weyl module induced from a finite-dimensional
$\spg_{2n-2}$-module at level $-(2n+1)/2$ is simple: for a positive root and
$A=\langle\lambda+\rho,\alpha^\vee\rangle$, its parabolic Jantzen
terms have multiplicities $A-s$, $1\le s<A$, and grades
$2s(A-s)$ for long roots or $s(A-s)$ for short roots.
The involution $s\leftrightarrow A-s$ preserves the grade and reverses
the finite Weyl sign. Its midpoint is singular. Negative roots give
no term. Nonnegativity of the characters of the nested Jantzen submodules
forces the first Jantzen submodule to vanish.

Using the oscillator contravariant form, split these simple Weyl
modules degree by degree, together with their nondegenerate multiplicity forms on their
finite-dimensional lowest spaces. Finite pieces make the splitting exhaustive.
The commutant theorem of Linshaw and Song
\cite[Theorem~6.1]{LinshawSong2024} identifies the image of
$V^{-2n+2}(\sog_{2n+1})$ with the multiplicity space of
$V^{-(2n+1)/2}(\spg_{2n-2})$. This identifies the image of $V^{-2n+2}(\sog_{2n+1})$.
Its simplicity follows from the preceding decomposition and the form.
Indeed, its nondegenerate form pulls back
to the normalized contravariant form on $V^{-2n+2}(\sog_{2n+1})$,
so the image is $L_{-2n+2}(B_n)$. The classes of the
$\sog_{2n+1}$ generators in the $C_2$ algebra of the symplectic-boson
vertex algebra are $J_{ab}=\omega(v_a,v_b)$. Set $v_1=e_1$, $v_2=f_1$ and
all other vectors zero. Only $J_{12}$ is nonzero, and the Sugawara
quadratic has nonzero value. Hence the class of the conformal vector
in the $C_2$ algebra of $L_{-2n+2}(B_n)$ is not nilpotent, so this
vertex algebra is not quasi-lisse.

\subsection{\texorpdfstring{Non-quasi-lisseness of $L_{p-2n+1}(\sog_{2n+1})$ for odd $p$ with $3p\le2n+1$}{Non-quasi-lisseness of L(p-2n+1)(so(2n+1)) for odd p with 3p at most 2n+1}}\label{B:odd-exclusions}
We prove that $L_{p-2n+1}(\sog_{2n+1})$ is not quasi-lisse for odd
$p\ge3$ with $2n+1\ge3p$. Set $t=p$ and $N=2n+1$. We choose the nilpotent element
according to the following four families. In each family, the
conformal-weight bound proves nonvanishing of the reduction, and the
Miura homomorphism produces a simple quotient that is not quasi-lisse.
The range $N\ge3p$ is exactly the following disjoint union:
$$
\begin{array}{c|l}
\mathrm I&N=(2j+1)p+2v+2,\ j\ge1,\ 0\le v\le(p-3)/2\\
\mathrm {II}&N=(2j+2)p+u,\ j\ge1,\ 1\le u\le p-4\text{ odd}\\
\mathrm {III}&N=(2j+1)p,\ j\ge1\\
\mathrm {IV}&N=(2j+1)p-2,\ j\ge2.
\end{array}
$$
In I set $E=2j+1$, $a=2v$, $b=p-a-1$, $c=p-1$ and choose
$f=((E+1)^a,E^b,1^{E+2})$. At the nonzero block endpoints of the grading,
the conformal-weight inequalities have the values
$$
 F(C_{\ell c})=\ell[(p-1)(\ell+2)-a],\qquad
 F(C_{\ell c+a})=(\ell+2)[\ell(p-1)+a]>0.
$$
Concavity on a constant-$x$ block and monotonicity on the zero block
prove the strict conformal-weight inequality, and hence $H_f\ne0$. The restriction of the invariant form to the generators corresponding
to $\sog_b$ is zero, including when $b=2$. The decomposition associated
with the Jordan blocks gives the tensor factor $L_{-2j}(B_{j+1})$,
for which $k+h^\vee=1$.
In II set $L=2j+2$, $b=p-u-2$, $c=p-2$ and
$f=((L+1)^u,L^b,1^{2L})$. At the block endpoints, the corresponding
conformal-weight differences are
$$
 F(C_{\ell c})=2\ell[(\ell-1)c+b],\qquad
 F(C_{\ell c+u})=2\ell[(\ell+1)c-b].
$$
Only $C_u$ has zero value, at degree $D=(j+1)u$.
Since $2D=Lu<Lp=N-u$, the bound for odd exterior powers proves nonvanishing of the reduction.
The vanishing of the generators corresponding to $\spg_b$ yields
a tensor factor $L_{4-2L}(\sog_{2L})$.
Here $L=2j+2$, and \cref{v6:scalar-partitions} identifies the simple
Miura image for the reduction at $((2j+1)^2,1^2)$ with
$\mathcal A_j$. The nonvanishing of this reduction is proved in
\cref{sec:D-exclusions}, so this tensor factor is not quasi-lisse.

For III/IV set $\epsilon=1,2$, $r=p-\epsilon$ and
$f=((2j+1)^r,(2j-1)^\epsilon,1^2)$. The conformal-weight inequalities use
$$
\begin{aligned}
 F(C_l)&=l(r+2-\epsilon-l),\quad1\le l\le r,\\
 F(C_{r+gp+l})&=(2-\epsilon)r+2gp+l(p+2-l).
\end{aligned}
$$
They are strict for $\epsilon=1$. For $\epsilon=2$ the only equality
is $C_r$ at degree $jr$, and $2jr=2j(p-2)<N-r=2jp$.
For $j\ge2$, the image of the Miura homomorphism is contained in
$$
 (L_\epsilon(\slg_r)\otimes\mathcal H)\otimes\mathcal A_{j-1},
 \qquad Q=\frac{1-\epsilon}{2}\sqrt{r/p}.
$$
Here $\mathcal H$ is generated by $h$, with
$[h_\lambda h]=\lambda\mathbf1$ and conformal vector
$\tfrac12:h^2:+Q\partial h$.
The calculation on the shorter Jordan chains following
\cref{seed:scalar-split} verifies its hypotheses in conformal
weights two and four, giving the tensor factor $\mathcal A_{j-1}$.
No invariant bilinear form on $\mathcal H$ with this conformal vector is required.

When $N=3p$, the tensor factor $\mathcal A_{j-1}$ is absent.
The Miura image for $f=(3^{p-1},1^3)$ is contained in the unitary vertex algebra
$L_1(\slg_r)\otimes\mathcal H$, $r=p-1$. It contains
$K=L_2(\sog_r)$ conformally in $L_1(\slg_r)$ and hence $T_{\mathcal H}$.
Its weight-three skew fields are generated by zero products of the
symmetric weight-two fields, since commutators of symmetric matrices
span the skew matrices. Consider the eigenspace decomposition of the weight-two subspace
for the operator $(T_{\mathcal H})_{(1)}$. Let $h$ generate
$\mathcal H$, with $[h_\lambda h]=\lambda\mathbf1$.
At eigenvalue one, the possible traceless subspace is
$h\otimes\Sym^2_0\C^r$. At eigenvalue two, the scalar subspace is $\C T_{\mathcal H}$.
At eigenvalue zero, the quasiprimary scalar and symmetric traceless subspaces
are in $K_2$. Indeed the complement's degree-one primary has Sugawara
conformal weight one and its only same-type degree-two descendant is its
derivative: the $\mathrm O(r)$-equivariant contraction
$\Lambda^2\C^r\otimes\Sym^2_0\C^r\to\Sym^2_0\C^r$ is unique.
This includes $r=2,4$, with both $\sog_4$ factors included.
All the projected spaces are PCT-stable. Thus the image is simple and
its conformal vector has nonnilpotent image $h^2/2$.

\subsection{\texorpdfstring{Non-quasi-lisseness of $L_{-2n+1+p/q}(\sog_{2n+1})$}{Non-quasi-lisseness of L(-2n+1+p/q)(so(2n+1))}}\label{seed:B-coset}
Write $k=-2n+1+p/q$, where $p,q\in\Z_{>0}$ are coprime and
$pq$ is even. We prove that $L_k(\sog_{2n+1})$ is not quasi-lisse
when $1\le p\le n$, and also when $p=2n-1$ and $q$ is even.
For $1\le p<n$, \cref{thm:denominator} allows us to use the level
$k=-2n+1+p/(p+1)$. Set $t=p/(p+1)$, $d=n-p$ and
$\kappa=t-1=-1/(p+1)$.
For the principal grading and $0<t<1$, one has
$$
 \Delta_t(\lambda)-\lambda(x_{\mathrm{prin}})
 =\frac1{2t}\sum_i a_i(a_i-1)
 +(1/t-1)\sum_i a_i(n-i+1)>0
$$
for every nonzero $B$ root-lattice dominant weight. Thus $H_{f_{\mathrm{prin}}}^0(L_{t-2n+1}(B_n))\ne0$.
Let $U_d(t)=\mathcal W^{t-2d+1}(B_d,f_{\mathrm{subreg}})$ for $d\ge2$,
with $d=1$ interpreted
as $V^{2\kappa}(\slg_2)$. It has a weight-one generator $J$
with $J_{(1)}J=\kappa\mathbf1$ and generators $G_d^\pm$ of weight $d$
satisfying $J_{(0)}G_d^\pm=\pm G_d^\pm$.
Deleting one grade-one root embeds it in
$\mathcal H_t\otimes U_{d-1}(t)$, where $A$ generates $\mathcal H_t$ and
$$
 T_d\longmapsto T_{d-1}+\frac{:AA:}{2t}
       +\frac{Q_d}{t}\partial A,\qquad Q_d=(d-1)\kappa-\tfrac12.
$$
The corresponding primary generators are
$\kappa\partial G^\pm\mp:JG^\pm:+:AG^\pm:$.
After removing the $J$ Fock factor their lowest conformal weight is
$h=d-1-1/(2\kappa)$ and their opposite-charge pairing is multiplied by
\begin{equation}\label{seed:B-norm}
 t-2\kappa^2h(2h+1)
 =2(2d-1)\kappa\bigl(1-(d-1)\kappa\bigr)\ne0.
\end{equation}
Starting with the $\slg_2$ pairing $2\kappa$ proves that the charged
fields survive in the simple quotient.

To apply \cref{com:heisenberg-quotient} to the commutant of
$\mathcal H_\kappa$, generated by $J$, note that for $d\ge2$ the
weight-one space is $\C J$.
The simultaneous determinant-minus-one transformations on the odd
$(2d-1)$ block and its two singleton blocks lie in $\mathrm{SO}_{2d+1}$,
centralize the triple, and give a conformal involution $J\mapsto-J$.
Hence $L_1J=0$. At $d=1$ the elements of $\slg_2$ are primary
with respect to the Sugawara conformal vector.
Thus $L_1V_1=0$, giving a normalized nondegenerate form by simplicity.
Let $\sigma$ denote the involution induced on the Heisenberg coset by
$J\mapsto-J$. For $d=1$, use the Weyl involution of $\slg_2$ with this
action on $J$. Since $\kappa\ne0$, \cref{com:heisenberg-quotient}
commutes with graded quotients and specialization. The
$\mathbb Z_2$-orbifold of the Heisenberg coset under $\sigma$ is
therefore the simple quotient of the corresponding universal orbifold.

We use the coincidence of simple vertex algebras proved by Creutzig
and Linshaw in \cite[Theorem~B.2(6)]{CL}:
\begin{equation}\label{seed:CL-full}
 \mathcal C_{\psi,1D}(a,m)=\mathcal W_s(\spg_{2R}),\qquad
 \psi=\frac{m+a-R}{m-R},\qquad
 s+R+1=\frac{R-m}{2(R-m-a)},
\end{equation}
for $a,m\ge0$, $R\ge1$, $R\ne m,m+a$.
Set $(a,m,R)=(1,d-1,n)$. Then $\psi=p/(p+1)$,
$s+n+1=(p+1)/(2p)$, and the required nonzero differences $R-m$ and $R-m-a$ are $p+1$ and $p$.
Compatibility with quotients and specialization identifies the left side
with this $\mathbb Z_2$-orbifold of the Heisenberg coset. Principal Feigin--Frenkel duality identifies the
right side with $\mathcal W_{t-2n+1}(B_n,f_{\mathrm{prin}})$.
This identification is a theorem about the simple vertex algebras.

Decompose the simple quotient of $U_d(t)$ as a module for
$\mathcal H_\kappa$:
$\bigoplus_jF_j\otimes M_j$, where $F_j$ is the Fock module on
which $J_{(0)}$ acts by $j$.
Simplicity and $M_1\ne0$ imply $M_j\ne0$ for every $j>0$. The coset lowest weights satisfy
$$
 \min\wt M_j\ge-\frac{j^2}{2\kappa}\longrightarrow\infty.
$$
Each $M_j$ is ordinary also when restricted to this orbifold.
A simple ordinary subquotient has the same lowest weight. There are
therefore infinitely many simple ordinary modules. The quasi-lisse
finiteness theorem (\cite{AK}) shows that this orbifold and, by
\cref{seed:CL-full}, the simple principal $W$-algebra
$\mathcal W_{t-2n+1}(B_n,f_{\mathrm{prin}})$ are not quasi-lisse.

At $p=n$ use $t=n/(n+1)$ and principal Levi $\glg_n$.
The central charge of the principal $W$-algebra for $\slg_n$ is
$(n-1)-n(n^2-1)(t-1)^2/t=0$, while the coefficient of the derivative term in the conformal
vector of the rank-one Heisenberg factor is
$t n(n+1)/2-n^2/2=0$. The simple principal $W$-algebra of central charge zero is
$\C$, since its only weight-two line is the conformal vector in the radical.
The resulting image is a conformal vertex subalgebra of the rank-one
Heisenberg vertex algebra with its standard conformal vector. Such an image is simple: the multiplicity-free unitary
Virasoro decomposition is $\bigoplus_{m\ge0}L(1,m^2)$, so every
conformal subalgebra is PCT-stable. Its conformal vector has polynomial image
$h^2/2$. The exceptional even-denominator numerator $p=2n-1$ is the
same argument at $t=(2n-1)/(2n)$ with Levi $\glg_1+B_{n-1}$:
the coefficient $tn-(n-1/2)$ of the central derivative term
in the Miura conformal vector and the central charge
of the principal $W$-algebra for $B_{n-1}$ are zero.
At $n=2$, the level $\ell$ of the short-root $\slg_2$ satisfies
$\ell+2=2t=3/2$.
These cases cover all the stated levels with $pq$ even at which
quasi-lisseness fails.

\subsection{\texorpdfstring{Singular vectors in $V^k(\sog_{2n+1})$ and their parabolic Wakimoto images}{Singular vectors in V(k)(so(2n+1)) and their parabolic Wakimoto images}}\label{B:relations}
For the cases in which we prove quasi-lisseness, we must distinguish
two uses of character
multiplicities. Comparing a constituent of the maximal ideal with
the full parabolic kernel can show that a semisimple point does not
belong to the associated variety. To
obtain a polynomial relation, however, we need an actual singular
subspace of the ideal. The next two lemmas make this distinction precise.

\begin{lemma}\label{com:whole-kernel}\label{seed:B-constituent}
Assume $k+h^\vee\in\Q_{>0}$, and let
$\Phi:V^k(\g)\to U$ be the universal maximal-parabolic map and
$T$ its quotient by the full derived-Levi maximal ideals, with kernel
$K\subset U$. If a weight space or a simple constituent of the maximal
proper graded ideal $I_k\subset V^k(\g)$ has positive multiplicity and its multiplicity in
$K$ is zero, then the corresponding center point is absent from $X_{L_k}$.
Composition multiplicities are taken in a finite highest-weight
truncation containing that constituent. In particular, the conclusion
holds for a simple module $S\not\cong L_k(\g)$ with $[V^k(\g):S]>0$ and $[K:S]=0$.
For the polynomial relations used below, we also need the relevant
finite-dimensional $\g$-module to occur as a subspace of $I_k$
annihilated by $x(r)$ for all $x\in\g$, $r>0$.
\end{lemma}
\begin{proof}
Injectivity (\cref{com:universal-parabolic-injective}) and the multiplicity
comparison give $\Phi(I_k)\not\subset K$, so \cref{com:point-map}
excludes the center. Every simple constituent of $V^k(\g)$ other than $L_k(\g)$ lies in $I_k$.
\end{proof}

\begin{lemma}\label{B:singular-lift}
Let $S$ be a simple constituent of $I_k\subset V^k(\g)$, of
multiplicity one and lowest conformal degree $D$. Work in a finite
highest-weight Serre quotient containing all simple labels through degree $D$.
Suppose that its standard-module multiplicity column has no label
other than $L_k(\g)$ and $S$, except possibly labels which do not occur
in $V^k(\g)$. Then the degree-$D$ highest-weight space of $S$ lifts
to a singular subspace of $I_k$.
\end{lemma}
\begin{proof}
We use the highest-weight order in which labels below $S$ have larger
Sugawara degree. These labels have been removed from the truncation.
Let $L_\lambda\ne S$ be a simple constituent of $I_k$ in this truncation.
For $\lambda>S$, the hypothesis gives
$[\Delta_\lambda:S]=0$. Let $P_\lambda$ be its projective cover.
The kernel of $P_\lambda\twoheadrightarrow\Delta_\lambda$ has a
standard filtration with labels strictly above $\lambda$. Each such
standard has a different simple head from $S$. Hence
$$
 \Hom(\ker(P_\lambda\to\Delta_\lambda),S)=0,
 \qquad
 \Ext^1(\Delta_\lambda,S)=0.
$$
Since $\operatorname{rad}\Delta_\lambda$ has no constituent $S$,
applying $\Hom(-,S)$ to
$0\to\operatorname{rad}\Delta_\lambda\to\Delta_\lambda
\to L_\lambda\to0$ gives $\Ext^1(L_\lambda,S)=0$.
Highest-weight duality gives $\Ext^1(S,L_\lambda)=0$ as well.
For incomparable labels the same projective argument applies:
$S$ cannot label a standard section of the projective kernel and
cannot occur in $\Delta_\lambda$. Hence both groups vanish. The simple module $L_k(\g)$ does not occur in the
proper graded ideal $I_k$. Thus the unique constituent $S$ splits from the
ideal in this finite truncation.

The simple modules removed by the truncation have no conformal-degree pieces through $D$.
Degree evaluation in this range is exact and factors through the
truncation, as do the maps $x(r)$, $x\in\g$, $r>0$, between these pieces.
The split degree-$D$ top is consequently annihilated by $x(r)$ for every $x\in\g$ and $r>0$ in the original ideal. Its finite-dimensional
$\g$-span is the required singular space.
\end{proof}

\subsubsection{Composition multiplicities}\label{com:hecke-multiplicities}
We use Kashiwara--Tanisaki's affine character formula
\cite[Theorem~1.4]{KT}. At positive rational $k+h^\vee$, choose the
closed dominant representative $\lambda$ for the integral root system.
Its real-coroot pairings lie in $\mathbb Q\setminus\mathbb Z_{<0}$,
and its stabilizer $W_0$ is finite. The theorem uses longest representatives
of the right cosets $wW_0$ and includes singular $\lambda$.
In each finite highest-weight truncation, inversion of its character
matrix gives ordinary Verma composition multiplicities. The finite BGG
resolution on the left then gives the parabolic multiplicities, which
are the antispherical coefficients at $v=1$. The full right
stabilizer is used in these computations. This interpretation gives composition
multiplicities, not positions in the Jantzen filtration, see also
\cite[\S6.1]{JiangSong2025} for the simply-laced conventions.

\begin{proposition}\label{seed:B-data}
The maximal ideal $I_k\subset V^k(B_n)$ contains the following
singular subspaces.
\begin{enumerate}
\item For $q=1$, $n\le p\le2n-4$, there is a singular subspace
$S$ of highest weight $r\theta$, $r=2p-2n+2$, and degree $2r$,
whose image under the parabolic homomorphism associated with $H_1$
is nonzero after replacing the affine Levi factors by their simple quotients.
\item For $q=2$, odd $p=2n-2s+1\ge n+2$, $s\ge2$, there is a singular subspace
$R$ of highest weight $(2^{2s+1})$ of degree $(2s+1)q$, nonzero at every
$H_j$, $j\ge s+1$. For $s\ge3$ there is also
$S=(p-2s+3)\theta$ of degree $2(p-2s+3)q$, nonzero at $H_1$.
For $s=2$ the replacement is $S=(2n-4,2n-4)$, degree $4(n-2)q$.
\item At $p=n+1$, $q=2$, for even $n\ge6$, there is a singular subspace
$R$ of highest weight $(3,2^{n-1})$ and degree $(n+2)q$ which has nonzero image under the parabolic homomorphism associated with $H_n$.
There is also a singular subspace $S$ of highest weight $4\theta$
and conformal degree $16$, whose image under the parabolic homomorphism
associated with $H_1$ is nonzero after replacing the affine Levi factors
by their simple quotients.
The remaining even-denominator endpoint $B_4,p=5$ is proved below
using the Jantzen sum and a bound for orthogonal representations
corresponding to partitions with two rows.
\end{enumerate}
\end{proposition}
\begin{proof}
Our original proof of the existence of these singular subspaces
used character identities and composition-multiplicity calculations,
together with \cref{B:singular-lift}. For the cases covered by
\cite[Theorem~4.8(5) and Corollary~5.11]{JiangWu2026},
we now use their results to shorten the proof.

At the half-integral levels in (2), Jiang--Wu's
\cite[Theorem~4.8(5) and Corollary~5.11]{JiangWu2026} gives
whichever of $R,S$ has the smaller listed conformal degree, and both
when the degrees agree. Their parameters satisfy $2(k+h^\vee)=p$.
In (3), $2(k+h^\vee)=n+1$ and their least degree is
$\min\{16,2n+4\}=16$. Hence the same result gives $S$ for every
even $n\ge6$, and also $R$ when $n=6$.

There are two separate questions here: whether the singular space
exists in the maximal ideal, and whether its image survives after
taking the quotient by the maximal proper graded ideals of the
affine vertex algebras associated with the derived Levi factors. The cited existence theorem answers the first
question for the lower-degree spaces just listed. It does not answer
the second question. In (1), and for the higher-degree spaces in (2)
and (3), we use the character calculations and the lifting lemma.

For these remaining singular spaces, use the character identities in
\hyperref[calc:B-FOLDS]{Calculation~\ref*{calc:B-FOLDS}} and the
multiplicity calculations in
\hyperref[calc:BCD]{Calculation~\ref*{calc:BCD}}.
The interpretation in \cref{com:hecke-multiplicities}
and \cref{B:singular-lift} give the required subspaces of $I_k$.
The full finite-parabolic columns are computed with the complete stabilizer.
For a column supported only on $L_k(\g)$ and the required singular
module, the lifting lemma applies directly. For $R$ in (3), the column
also contains $3\omega_1$. To apply the lifting lemma, we must show
that this additional label does not occur in the vacuum module.
Its usual orthogonal extension is odd under $-I$, whereas $-I$ acts
trivially on $\g=\Lambda^2V$ and on $V^k(\g)$.
The other extension, $\det\otimes\Sym^3_0V$, requires defining-representation
tensor degree at least $2n+2$, exceeding the bounds $6$ and $12$
in conformal degrees three and six. Thus this intervening label is absent.
The small ranks are treated below.

For the nonvanishing of $S$ in (3), the composition-multiplicity
column for $H_1$ has first-coordinate offsets
$4,-(n+2),0,-(n-2)$. Only the vacuum has zero offset, in degree zero.
Thus none of the nonvacuum labels can contribute the required
constituent to the full kernel of the quotient by the maximal proper
graded ideals of the affine Levi factors. That kernel has zero
multiplicity of the simple highest-weight module with top $S$.
Since $S$ occurs in the maximal ideal, \cref{com:whole-kernel}
gives the asserted nonvanishing.

For $S$ in (2) with $s=2$, set $p=2n-3$. For a finite-dominant
label $\lambda=(\lambda_1,\ldots,\lambda_n)$, define
$$
 c_0(\lambda)=\sum_{i=1}^n
 \left\lfloor\frac{\lambda_i+\rho_i}{p}+\frac12\right\rfloor,
 \qquad \rho_i=n-i+\frac12.
$$
The complete Jantzen cancellation through the degree $4(p-1)$ of $S$
gives multiplicity one for $S$. Its highest weight has $c_0=3$,
while every other label in the resulting truncation of the ideal
has $c_0\ge4$, as computed in
\hyperref[calc:B-FOLDS]{Calculation~\ref*{calc:B-FOLDS}}.
Monotonicity of $c_0$ in Bruhat order and strict increase of Sugawara
degree make these labels incomparable with $S$, so the projective and
degree-evaluation argument of \cref{B:singular-lift} gives its lift.
The calculations also compare each space with the full kernel of
the quotient by the maximal proper graded ideals of the affine Levi
factors. The following subsection extends nonvanishing to the remaining $H_j$.
\end{proof}

\subsection{\texorpdfstring{Evaluating $C_2$ relations at the semisimple points $H_m$}{Evaluating C2 relations at the semisimple points H(m)}}\label{B:covariant-restriction}
Let $E\subset I_k$ be a singular subspace from the preceding proposition.
For $v\in E$, its class $\bar v\in R_{V^k(\g)}\cong S(\g)$ is a
polynomial on $\g$, using the invariant form to identify $\g$ with
$\g^*$. These polynomials define the equivariant map
$F_E:\g\to E^*$ by $F_E(x)(v)=\bar v(x)$.
Thus $F_E(H_m)\ne0$ means that a polynomial relation defining
$X_{L_k(\g)}$ is nonzero at $H_m$, so $H_m\notin X_{L_k(\g)}$.
We first relate this evaluation to the image of $E$ under the parabolic
Wakimoto homomorphism. We then prove degree bounds under which
$F_E(H_{m-1})\ne0$ implies $F_E(H_m)\ne0$.
We use the self-duality of the finite type-$B$ representations to
identify $E^*$ with $E$.
\begin{lemma}\label{B:singular-detection}
Let $E\subset I_k$ be a homogeneous finite-dimensional affine singular
space of positive degree. For the parabolic homomorphism with simple
Levi factors one has
$$
 F_E(H_m)=0\quad\Longleftrightarrow\quad
 \Phi_{H_m}(t)(E)=0.
$$
\end{lemma}
\begin{proof}
If the image of $E$ is nonzero, choose a finite highest-weight vector
in that image. It is annihilated by $x(r)$ for every $x\in\g$ and $r>0$.
By \cref{com:leading-fiber}, its $C_2$ symbol is nonzero on the fiber
obtained by setting the semisimple Levi coordinates to zero. Its moment-map values have the form
$aH_m+y$, where $a\ne0$ and $y$ belongs to the opposite nilradical.
Every such point is conjugate to $aH_m$: the action of
$\ad H_m$ on the nilradical has nonzero degrees, so the
conjugation equations can be solved successively in those degrees.
Finite $\g$-equivariance and homogeneity would make all these symbols
zero if $F_E(H_m)=0$, a contradiction. Conversely, evaluation at zero
$\beta\gamma$ and semisimple Levi coordinates gives $F_E(H_m)$ from the
image of $E$. This evaluation descends to the simple affine quotients
of the Levi factors, since their proper graded ideals have zero constant term.
\end{proof}

The codomain of the universal parabolic homomorphism is
$\mathcal F_m\otimes V_A\otimes V_B$, where $\mathcal F_m$ is the
tensor product of the $\beta\gamma$ vertex algebra for the
nilradical and the rank-one Heisenberg vertex algebra associated
with $Z(\lev)$. Here $A=A_{m-1}$ and $B=B_{n-m}$.
Let $I_A\subset V_A$ and $I_B\subset V_B$ be the maximal proper graded
ideals. The kernel of the quotient map onto
$\mathcal F_m\otimes(V_A/I_A)\otimes(V_B/I_B)$ is exactly
$$
 \mathcal F_m\otimes
 (I_A\otimes V_B+V_A\otimes I_B).
$$
Let $J_A,J_B$ be the images of these affine ideals in
$S(A),S(B)$. Projecting to $C_2$ therefore sends every kernel symbol into
$$
 R_{\mathcal F_m}\otimes
 (J_A\otimes S(B)+S(A)\otimes J_B).
$$
This uses the displayed tensor kernel, without an assertion that the
$C_2$ functor is left exact. In an affine $C_2$ algebra, polynomial
degree equals conformal degree: the quotient is spanned by the
classes of $a_1(-1)\cdots a_d(-1)\mathbf1$, with $a_i\in\g$.
If $I_A$ starts in degree $M$, every polynomial in $J_A$ has degree
at least $M$, while $J_B(0)=0$. Setting the $B$ coordinates to zero
thus gives order at least $M$ in the $A$ coordinates.

Suppose $p\ge m$ and $F_E(H_m)=0$. By
\cref{B:singular-detection}, the singular space maps into this full
kernel. The $A_{m-1}$ level is admissible and its first ideal degree is
$M=(p-m+1)q$. Set the $\beta\gamma$ coordinates to zero. The image of each $x\in\g$ is then the sum of its central and
semisimple Levi components, and
quantum derivative terms vanish in $C_2$. Consequently
$F_E(aH_m+X_A)$ has order at least $M$ in $X_A$, with coefficients
polynomial in $a$.

On the line $h(z)=(1^{m-1},z,0,\ldots)$, write
$a=(m-1+z)/m$ and $X_A=h(z)-aH_m$. Every coordinate of $X_A$ is
divisible by $z-1$. Weyl sign equivariance gives the same zero order
at $z=-1$. If $F_E(H_{m-1})\ne0$, this polynomial is nonzero and has
degree at most $D=\deg E$. We obtain
\begin{equation}\label{B:forward-propagation}
 D<2(p-m+1)q,\quad F_E(H_{m-1})\ne0
 \quad\Longrightarrow\quad F_E(H_m)\ne0.
\end{equation}

There is also a useful equality case. Let $E$ have finite highest
weight $r\theta$, $D=2rq$, and $r=p-m+1$. If $F_E(H_m)=0$, degree
saturation gives $F_E(h(z))=(1-z^2)^{rq}v$, where
$v=F_E(H_{m-1})\ne0$. The radial identity \cref{com:radial} gives
half-Casimir eigenvalue $pr$ for the mixed-root sum at $z=1$.
The vector $v$ is fixed by $B_{n-m+1}$ on the remaining coordinates.
To see the resulting bound explicitly, write $\varepsilon_i$ for the
standard orthogonal coordinates. On the zero-weight space the mixed-root
operators are
$$
 A_{j,\epsilon}=\sum_{i=1}^{m-1}
 f_{\varepsilon_i-\epsilon\varepsilon_j}(0)
 e_{\varepsilon_i-\epsilon\varepsilon_j}(0),
 \qquad m\le j\le n,\quad \epsilon\in\{1,-1\}.
$$
These $2(n-m+1)$ operators are Weyl conjugates of the residue operator
at $z=1$, so each
has eigenvalue $pr$ on $v$. Their root sets are disjoint. For a compact
invariant Hermitian form each summand is positive semidefinite, and
the omitted positive-root terms are positive semidefinite as well.
Their sum is therefore bounded by the total half-Casimir, whose
eigenvalue is $2pr$. Hence
$$
 2pr\|v\|^2\ge2(n-m+1)pr\|v\|^2.
$$
This is impossible when $m<n$. Thus
\eqref{B:forward-propagation} also holds with equality for these
specified $r\theta$ spaces.

For the rectangular weight $(2^{2s+1})$ in
Proposition~\ref{seed:B-data}(2), set $m_0=\lceil2n/3\rceil$.
Order means the least total degree in the transverse semisimple
Levi coordinates at $H_m$. Here $e_A,e_B$ denote the half-Casimir
bounds on the corresponding finite Levi types. At $m=m_0$, $l=n-m$,
the coordinate bound two gives
$$
 e_A\le m(m+4)/2,\quad e_B\le l(l+2),\quad
 M_A=(p-m+1)q,\quad M_B=(p-2l+1)q/2.
$$
For $p\ge n+2$, the required strict inequalities are
$$
 \frac{m(m+4)q}{2p}<(p-m+1)q,
 \qquad
 \frac{l(l+2)q}{p}<\frac{(p-2l+1)q}{2}.
$$
They hold already at $p=n+2$, as follows by writing $n$ in its three
residue classes modulo three. After clearing the positive denominator,
the right sides increase with $p$. The radial equation bounds the
orders of the first nonzero transverse terms by $e_Aq/p$ and $e_Bq/p$.
These orders are smaller than the first possible degrees $M_A,M_B$
of the Levi ideals. If $F_E(H_{m_0})=0$, the image of $E$ under the universal
parabolic homomorphism would lie in the full kernel described
above. This would impose lower bounds on the vanishing orders
that contradict the bounds from the radial equation. 
Thus $F_E(H_{m_0})\ne0$.
Propagate forward from $H_j$ to $H_{j+1}$ for $m_0\le j<n$.
The spherical $A_j$ types are $a\theta$, $a\le2$, with order
at most $2(j+2)q/p<(p-j)q$, using $(p+2)(p-2-j)>0$.
For the reverse steps $H_j$ to $H_{j-1}$, $s+2\le j\le m_0$,
set $l=n-j+1\le n-s-1=(p-3)/2$.
The spherical $B_l$ types $(1),(2),(11),(21),(22)$ have half-Casimirs
$l,2l+1,2l-1,3l,4l$.
The first three are less than $p$, and $4l<2p$ cannot be a
multiple of the odd integer $p$. Thus only $3l$ can resonate,
but $(21)$ has no $\glg_1+B_{l-1}$ invariant:
its weight-zero tensors have degree one or three in the defining
representation, below the degree of the alternating invariant tensor. The ranks one and two are direct. The rank-one spectrum is
$1,3$. Thus $F_E(H_j)\ne0$ for every $s+1\le j\le n$.

\subsection{\texorpdfstring{Quasi-lisseness of $L_{p-2n+1}(\sog_{2n+1})$ for $p>n$ or odd $p=n$}{Quasi-lisseness of L(p-2n+1)(so(2n+1)) for p greater than n or odd p=n}}\label{B:upper-integral}
Let $n\ge3$ and $q=1$, with $p>n$ or $p=n$ odd, below the
admissible range. Write
$$
 p=2n-2s+1+\epsilon,\quad k=2-2s+\epsilon,
 \qquad\epsilon\in\{0,1\},\quad s\ge2.
$$
Then $n\ge2s-1$. At $k_0=2-2s$ the natural Pfaffian gives
$\rank X\le2s-2$. For $\epsilon=0$, this bound fails at $H_j$ for $j\ge s$.
For $\epsilon=1$, \cref{com:translation} gives the same bound on
$X\cap\N$. The following explicit degeneration rules out the points
$H_j$ with $j\ge s$ in this case as well.

Write $V=U\oplus W\oplus U^*$, $\dim U=j$, $\dim W=2n+1-2j=m$,
with paired first and last summands. The nilradical contains
$$
 N(A,B)=\begin{pmatrix}0&A&B\\0&0&-A^t\\0&0&0\end{pmatrix},
 \qquad B^t=-B.
$$
Set $d=\min(m,s)$. Take $A$ to have an identity $d$-block and $B$
to have $s-d$ alternating two-by-two blocks on unused coordinates of $U$.
These fit because, if $d=m<s$, $j+m\ge n+1\ge2s$.
The supports are disjoint, so $N^3=0$ and $\rank N=2s$.
Writing $N=N_1+N_2$ by center degree, for $a\ne0$ the element
$Z=-N_1/a-N_2/(2a)$ satisfies
$[Z,aH_j]=N$, $[Z,N]=0$. Hence $aH_j+N$ is conjugate to $aH_j$.
A closed conical support containing $H_j$ would contain $N$, contradicting
the translated nilpotent rank bound. Only the nilpotent part was translated.

For $s\ge3$, the singular space of highest weight $r\theta$ in \cref{seed:B-data} has degree
$D=2r$, $r=2p-2n+2$, and is nonzero at $H_1$. For
$2\le m\le s-1$,
$$
 D<2(p-m+1)\quad\Longleftrightarrow\quad m<2s-2-\epsilon.
$$
This is strict throughout the required interval. The order-of-vanishing
argument above gives $F_E(H_m)\ne0$ successively for $m\le s-1$. For $s=2$
only $H_1$ remains: use the six-term $k=-2$ or ten-term $k=-1$
formula for the degree-zero coefficient of the reduced Jantzen sum, both valid for every $n\ge3$.
This includes $B_3,p=3,4$, $B_5,p=6$, and $B_7,p=8$.

\subsection{\texorpdfstring{Quasi-lisseness of $L_{p/2-2n+1}(\sog_{2n+1})$}{Quasi-lisseness of L(p/2-2n+1)(so(2n+1))}}\label{B:half-integral}
We prove quasi-lisseness at $k=-2n+1+p/2$ for odd $p>n$,
$p\ne2n-1$. The admissible levels have already been treated, so it
remains to consider $p\le2n-3$.
For $p=2n-2s+1\ge n+2$, the singular subspaces of highest weight
a multiple of $\theta$ exclude $H_1,\ldots,H_{2s-2}$, while the
subspace of rectangular highest weight excludes $H_{s+1},\ldots,H_n$.
At $s=2$, the equality case of the nonvanishing argument in
\cref{B:covariant-restriction} shows that $H_2$ is absent as well.
For $p=n+1$ with even $n\ge6$, the composition-multiplicity column for $H_n$ has only
the labels of $L_k(B_n)$, the module of highest weight $3\omega_1$,
and $R$. The standard modules in the full Levi-ideal kernel have fixed
coordinate sum, whereas the two nonvacuum labels change it by $3$ and
$2n+1$. Thus $R$ shows that $H_n$ does not belong to the associated variety. The comparison with $H_{n-1}$ uses
the short-$B_1$ residues $1,3,6<p$ and rules out that point as well.
The singular subspace of highest weight $4\theta$ and conformal degree
$16$ in \cref{seed:B-data}(3) shows that $H_1,\ldots,H_{n-2}$ do not belong to the associated variety.
Indeed, the forward inequality is
$16<2(p-m+1)q=4(n+2-m)$ for $m<n-2$. Equality at $m=n-2$
is covered by the equality case for highest weight $r\theta$ in
\cref{B:covariant-restriction}.
This covers all maximal parabolics.

At $B_4,p=5,q=2$, the Jantzen calculation in
\hyperref[calc:B4-LIFT]{Calculation~\ref*{calc:B4-LIFT}}
gives a singular subspace $R=(3,2,2,2)$ of conformal degree $12$.
The existence of this subspace also follows from
\cite[Theorem~4.8(5) and Corollary~5.11]{JiangWu2026}.

For the additional space $S=(4,4,0,0)$ in degree $16$, use the
complete twenty-six-term Jantzen sum in
\hyperref[calc:B4-LIFT]{Calculation~\ref*{calc:B4-LIFT}}.
Its simple multiplicities are at most one through degree $16$,
so nestedness gives $J^2=0$ in this range,
so the first derived contravariant form $b$ is nondegenerate on $J^1$.
Let $K$ be the submodule generated by the first singular space $R$.
Then $K=J^1$ below degree $16$, while $J^1_{16}/K_{16}$ contains $S$
once and the standard module with top $R$ has no constituent $S$.
The orthogonal complement of $K_{16}$ in $J^1_{16}$ is dual to this
quotient and contains a copy of $S$. For a vector $u$ in that copy,
contravariance gives
$$
 b(x(a)u,v)=b(u,\sigma(x)(-a)v)=0
 \qquad(a>0,\ v\in J^1_{16-a}=K_{16-a}),
$$
where $\sigma$ is the standard anti-involution. Nondegeneracy on
$J^1_{16-a}$ gives $x(a)u=0$. This produces the singular $S$-space
without requiring the restriction of $b$ to $K_{16}$ to be nondegenerate.
To prove $F_S(H_1)\ne0$, it suffices to bound the finite-dimensional
Levi representations in the restriction of $S$. In general, for
$\dim W=N\ge5$ odd, the $\mathrm{SO}_2$-zero part of the harmonic
$(r,r)$ tensor on $\C_+\oplus W\oplus\C_-$ has only $\mathrm{SO}(W)$
types $(a,b,0,\ldots)$ with $b\le a\le r$ and $a+b$ even.
Indeed Young branching deletes equally many $+$ and $-$ boxes, and
orthogonal contractions delete pairs of boxes. All surviving shapes have
two rows. The first two columns have total length at most $4<N$, so
no unstable Hodge identification changes the parity. This fails for
$N=3$, where $\Lambda^2W\cong W$, which is outside the application.

For $S$ this bounds the possible $B_3$ half-Casimirs by
$
 0,5,7,12,14,18,21,23,32.
$
The only positive multiple of $5$ is the adjoint value $5$. Its local
order would be $q=2$, whereas its zero-weight longest-Weyl action is
$-1$. Even polynomial degree has action $+1$. Hence $F_S(H_1)\ne0$.
The weight-polytope and interlacing bounds in \hyperref[calc:B4-LIFT]{Calculation~\ref*{calc:B4-LIFT}}
give $F_R(H_3),F_R(H_4)\ne0$. If $F_S(H_2)=0$, its restriction to
$(1,z,0,0)$ would be $(z^2-1)^{4q}v$ with $v\ne0$.
The endpoint residue has value $20$. $B_3$ invariance makes all six
mixed-root expectations equal to $20\|v\|^2$. Their sum is bounded
by total half-Casimir $40$, contradicting $40\ge6\cdot20$.
None of the four semisimple points $H_1,H_2,H_3,H_4$ belongs to the
associated variety.

The remaining low-rank constituent needed for odd $p<n$ is
$(n,p,q)=(6,5,1)$, at $H_3$, with finite label $(4,2,2,0,0,0)$ and
$[V:S]=[U:S]=[T:S]=1$. Equality of the last two multiplicities means
that the kernel of the parabolic homomorphism with simple affine
quotients of the Levi factors in its target has multiplicity zero, the exact implication
required by \cref{seed:B-constituent}.

\subsection{\texorpdfstring{Quasi-lisseness of $L_{p-2n+1}(\sog_{2n+1})$ for odd $p<n$ with $2n\le3(p-1)$}{Quasi-lisseness of L(p-2n+1)(so(2n+1)) for odd p below n with 2n at most 3(p-1)}}\label{B:lower-integral}
\label{B:uniform-odd-sums}

Here $k=p-2n+1$, with $p<n$ odd and $2n\le3(p-1)$.
The proof of quasi-lisseness has two parts. We first determine the few
low-conformal-weight relations needed in the reduced ideal and use
them to prove that certain reductions vanish and that the reduction
at the boundary parameters is $C_2$-cofinite. We then use the resulting
restrictions on nilpotent orbits to show that no nonzero semisimple point can occur. The character
identities below are uniform in the rank, and the degree bounds explain
why no omitted terms can affect these relations.

Let
$$
 a\geq3\ \text{be odd},\qquad b\geq0,\qquad
 p=a+2b+2,\qquad n=\frac{3a+4b+3}{2},\qquad k=p-2n+1.
$$
Thus $p<n$, $2n\leq3(p-1)$ and $k=-2(a+b)$.
Finite specializations of the character identities proved below are recorded in
\hyperref[calc:B-ODD-STRIP]{Calculation~\ref*{calc:B-ODD-STRIP}}.
Write $\rho_i=n-i+\tfrac12$ and $\rho=(\rho_1,\ldots,\rho_n)$.
For a partition $\lambda$ of length at most $n$, write $E_\lambda$
for the irreducible $\mathfrak{so}_{2n+1}$-module with highest weight
$\lambda$, and set $\chi_\lambda=\operatorname{ch}E_\lambda$.

For a finite root $\alpha$ and $m\geq1$, set
$$
 c=\frac{2((\rho,\alpha)+mp)}{(\alpha,\alpha)}.
$$
The root contributes only if $c$ is a positive integer and
$\rho-c\alpha$ is regular. If
$\lambda+\rho=w(\rho-c\alpha)$ is dominant, its contribution to the
Jantzen numerator for $V^k(B_n)$ is
$\det(w)Q^{mc}\chi_\lambda$. After reduction at grading element $x$,
it becomes $\det(w)Q^{mc}\chi_\lambda(Q^{-x})$.
The least exponent is $mc-(\lambda,x)$.

\begin{lemma}\label{B:uniform-mode-bound}
Let $x$ be dominant. Every regular term above satisfies
$$
 mc-(\lambda,x)\geq
 \begin{cases}
 c(m-x_1-x_2),&\alpha\text{ long},\\
 c(m-x_1),&\alpha\text{ short}.
 \end{cases}
$$
\end{lemma}
\begin{proof}
Since $(w\rho,x)\leq(\rho,x)$,
$$
 mc-(\lambda,x)
 =c\bigl(m+(w\alpha,x)\bigr)+(\rho-w\rho,x).
$$
The minimum of $(\beta,x)$ on a long-root orbit is $-x_1-x_2$, and
on a short-root orbit it is $-x_1$.
\end{proof}

\begin{lemma}\label{B:uniform-mode-one}
The complete contribution with $m=1$ is
\begin{align*}
 &Q^{a+b}\sum_{r=1}^{p}(-1)^{r+1}Q^r\chi_{(r,1^{p-r})}+\sum_{1\leq i<j\leq p}(-1)^{i+j}
 Q^{p+2n-i-j+1}
 \chi_{(p-i+1,p-j+2,2^{i-1},1^{j-i-1})}\\
 &\qquad\qquad\qquad+\sum_{i=1}^{p}(-1)^iQ^{2p+2n-2i+1}
 \chi_{(2p-i+1,1^{i-1})}.
\end{align*}
\end{lemma}
\begin{proof}
Reflection in $\alpha$ sends $\rho-c\alpha$ to
$\rho+2mp\alpha/(\alpha,\alpha)$ and changes the alternant sign.
For $m=1$, the latter operation adds $p$ to a changed coordinate for
a long root and $2p$ for a short root, with the signs of $\alpha$.
All other coordinates have absolute values in
$\{\tfrac12,\tfrac32,\ldots,n-\tfrac12\}$. Thus a changed absolute
value in this set must equal $\rho_i$ for one of the modified
positions $i$. All resulting absolute values must be distinct.
Set $L=n-(p-1)/2=a+b+1$. Then $2L=2n+1-p$.

Consider $\alpha=-\varepsilon_i+\varepsilon_j$, $i\ne j$.
The changed value $|\rho_i-p|$ is in the original set because $p<n$.
A positive coefficient and regularity leave exactly the following
possibilities. If $i=L$ and $j\leq p$, $j\ne L$, the negative
coordinate returns to its own absolute value. The highest weight is
$(p-j+1,1^{j-1})$, its degree is $p+L-j$, and its sign is
$(-1)^{j+1}$. Otherwise $i+j=2L$, $i\geq a$, $i\ne L$. The highest
weight is $(i-a+1,1^{i-1})$, its degree is $2i-a$, and its sign is
$(-1)^i$. The only other possible matching, $j=i+p$, gives $c=0$.

The negative short root $-\varepsilon_i$ is regular exactly for
$i\geq a$. It has highest weight $(i-a+1,1^{i-1})$, degree $2i-a$
and sign $(-1)^{i-1}$. It cancels the second mixed-root family for
$i\ne L$. At $i=L$ it is the missing term of the first
family. Replacing $j$ by $p-r+1$ gives the first displayed sum.
For $-\varepsilon_i-\varepsilon_j$, the two modified absolute values
must be distinct elements of $\{\rho_i,\rho_j\}$. The resulting matching has
$i+j=2L$ and $c=0$, so there is no regular term with $c>0$.

For $\varepsilon_i+\varepsilon_j$, $i<j$, regularity requires either
$j\leq p$ or $j=i+p\leq n$. In the first case sorting the two new
largest coordinates gives the second displayed sum. In the second
case the highest weight is $(p-i+1,1^{j-1})$, the degree is
$2n-2i+1$, and the sign is opposite to that of the positive short
root $\varepsilon_j$. These terms cancel. The positive short roots
with $i\leq p$ give precisely the third sum.
These cases exhaust the possibilities in which the changed absolute
coordinates are distinct from each other and from all unmodified
absolute coordinates, so the list exhausts the regular roots.
\end{proof}

Let $V$ be the standard orthogonal module, and set
$$
 E_V(z)=\sum_{j\geq0}\operatorname{ch}(\Lambda^jV)z^j,
 \qquad
 H_V(z)=\sum_{j\geq0}\operatorname{ch}(\operatorname{Sym}^jV)z^j.
$$
The characters are evaluated at the grading under consideration when
these expressions are used below.

\begin{lemma}\label{B:uniform-hook-identity}
For odd $p<n$,
$$
 \sum_{r=1}^p(-1)^{r+1}Q^r\chi_{(r,1^{p-r})}
 =[z^p]\frac{Q}{1-Q}E_V(z)H_V(-Qz)
                   \frac{1-Q^2z^2}{1-Qz^2}.
$$
\end{lemma}
\begin{proof}
Let $G_{\rm GL}(u,v)=\sum_{r\geq1,l\geq0}
\operatorname{ch}S_{(r,1^l)}V\,u^rv^l$ and let $G_{\rm O}(u,v)$
denote the corresponding sum of harmonic orthogonal characters.
The hook Schur identity and removal of orthogonal traces give
$$
 G_{\rm GL}(u,v)=\frac{u}{u+v}\bigl(E_V(v)H_V(u)-1\bigr),
 \qquad
 G_{\rm GL}(u,v)=\frac{1+uv}{1-u^2}G_{\rm O}(u,v)
                         +\frac{u^2}{1-u^2}.
$$
For the second identity, traces remove $2j$ boxes from the first row,
or $2j-1$ from that row and one from the first column. The empty
partition occurs for an even row. These are all possibilities for
a hook. Only total degrees at most $p<n$ are used, so no orthogonal
modification relation occurs. Substitute $(u,v)=(-Qz,z)$ and take
the coefficient of $z^p$. The terms not containing
$E_V(z)H_V(-Qz)$ are rational functions of $z^2$ and have zero
coefficient because $p$ is odd.
\end{proof}

\begin{proposition}\label{B:uniform-five-hook}
For $f=(5,1^{2n-4})$, the conformal-degree-zero coefficient of the
reduced Jantzen sum is one. Consequently $H_f^0(L_k(B_n))=0$.
\end{proposition}
\begin{proof}
Here $x=(2,1,0,\ldots,0)$.
Lemma~\ref{B:uniform-mode-bound} leaves only long roots with
$m\leq3$ and short roots with $m\leq2$ in degrees at most zero.
In the first sum of Lemma~\ref{B:uniform-mode-one}, the standard
module has grading weights $Q^{-2},Q^{-1},1^{2n-3},Q,Q^2$, whence
$$
 E_V(z)H_V(-Qz)
 =\frac{1+Q^{-2}z}{1+Q^3z}
       \left(\frac{1+z}{1+Qz}\right)^{2n-4}.
$$
Lemma~\ref{B:uniform-hook-identity}, including the factor $Q^{a+b}$,
shows that this entire sum starts in degree at least $a+b-1>0$.
This also proves cancellation of any negative-degree terms before
multiplication by the character of $\mathcal W^k(B_n,f)$.

For a term indexed by $i<j$ in the second sum, the least reduced
degree is $a+i-4$. For the third sum it is $a-2$ if $i=1$ and
$a-3$ if $i\geq2$. These are positive unless $a=3$.
When $a=3$, the possible degree-zero terms are the second-sum terms
with $i=1$, $2\leq j\leq p$, and the third-sum terms with index $j$.
Their signs are opposite. Their highest $x$-weight spaces are the
same $\mathfrak{so}_{2n-3}$-module $\Lambda^{j-2}\mathbb C^{2n-3}$,
so they cancel, including their dimensions.

For $m=2$, set $a=2n-2p+1$ as above.
The regular negative-long terms $-\varepsilon_i-\varepsilon_j$
have $a\leq i<j$ and highest weight
$$
 (j-a+1,i-a+2,2^{i-1},1^{j-i-1}),\qquad D=2(i+j-a).
$$
Their least reduced degree is $i+a-4\geq2$.
For a mixed root $\varepsilon_i-\varepsilon_j$, if $j<a$ regularity
forces $i=a-j$. Its highest weight and degree agree with those of
the short root $-\varepsilon_j$, and their signs are opposite.
If $j\geq a$, both new absolute values lie above $\rho_1$. The
first two highest-weight coordinates are $2p-i+1,j-a+2$ and
$D=2(2p+j-i)$. Its least reduced degree is $j+a-4\geq2$.
After the indicated cancellations, the negative-short terms have
$j\geq a$ and least reduced degree $2j-3>0$.
Positive-long terms have least reduced degree
$4n-2p-j-2>0$, and positive-short terms have least reduced degree
$4n-2i-\mathbf1_{i\geq2}>0$. Thus $m=2$ contributes nothing in
degrees at most zero.

For $m=3$, equality in Lemma~\ref{B:uniform-mode-bound} requires
$(\rho-w\rho,x)=0$ and $(w\alpha,x)=-3$.
The first equality makes $w$ fix the first two coordinates, and the
second forces $\alpha=-\varepsilon_1-\varepsilon_2$.
Its coefficient is $A=3p-2n+2=2b+5$, its highest weight is $(A,A)$,
and its Sugawara degree is $3A$. Its maximal $x$-weight space is
one dimensional and its sign is positive. Every other term has
strictly positive reduced degree. Hence the complete numerator
has no negative powers and constant coefficient one.
The character of $\mathcal W^k(B_n,f)$ has constant coefficient
one. The completed Jantzen identity and the criterion for the vanishing of $H_f^0(L_k(B_n))$ now give
the assertion.
\end{proof}

\begin{proposition}\label{B:uniform-triple-hook}
For $f=(3^{a+1},1^{4b+1})$, the conformal-degree-zero coefficient
of the reduced Jantzen sum is $1+\mathbf1_{b=0}$. Consequently
$H_f^0(L_k(B_n))=0$.
\end{proposition}
\begin{proof}
Set $s=a+1$, so $x=(1^s,0^{n-s})$.
The mode bound leaves long roots with $m\leq2$ and short roots
with $m=1$. The standard module has weights
$Q^{-1}$ and $Q$ with multiplicity $s$, and weight one with
multiplicity $2n+1-2s$. Thus
$$
 E_V(z)H_V(-Qz)
 =\frac{(1+Q^{-1}z)^s(1+z)^{4b+1}}
 {(1+Qz)^{4b+1}(1+Q^2z)^s}.
$$
The first sum of Lemma~\ref{B:uniform-mode-one} starts in degree
$a+b+1-s=b$. If $b=0$, its constant coefficient is
$[z^{p-s}](1+z)=1$, since $p-s=1$.
For the second sum, taking the first $s$ parts of its partition
shows that its least reduced degree is at least
$2n-p-2s+2=2b+1>0$. The least degree of the third sum is
$2n-i-\min(s-1,i-1)>0$ for every $i\leq p$.

For $m=2$, equality in the mode bound requires that $w$ preserve
the first $s$ coordinates and that both nonzero coordinates of
$\alpha$ be negative and lie among them. As in the preceding
proof, a regular negative-long root has $a\leq i<j$.
Since $i,j\leq s=a+1$, only $(i,j)=(a,a+1)$ remains.
Its highest weight is $(2^{a+1})$, its degree is $2(a+1)$, and
its sign is positive. Its maximal $x$-weight space is the
one-dimensional determinant-square representation of
$\mathfrak{gl}_{a+1}$. It therefore contributes one.
There are no negative-degree terms, and the completed identity
gives the stated degree-zero coefficient and the vanishing of
$H_f^0(L_k(B_n))$.
\end{proof}

\begin{proposition}\label{B:uniform-tail-survival}
For odd $p\geq1$ and $p\leq2r-3$,
$$
 H^0_{(3,1^{2r-2})}\bigl(L_{p-2r+1}(B_r)\bigr)\ne0.
$$
\end{proposition}
\begin{proof}
Now $x=(1,0,\ldots,0)$ and $\rho_i=r-i+\tfrac12$.
By the mode bound, every $m\geq2$ contributes only in positive
degrees, and no term has negative degree. Equality for $m=1$ requires
$w$ to fix the first coordinate and $\alpha$ to contain
$-\varepsilon_1$.
For the short root $-\varepsilon_1$, a positive coefficient occurs
exactly when $p\geq r$, and it equals $C=2p-2r+1$.
The corresponding highest weight is $C\varepsilon_1$, the
degree is $C$, and the sign is positive.

For $-\varepsilon_1+\varepsilon_j$, regularity requires
$|\rho_1-p|=\rho_j$. If $p<r$, this gives $j=p+1$ and coefficient
zero. If $p\geq r$, it gives $j=2r-p$, which is at least three
under the hypothesis $p\leq2r-3$. The resulting regular term has
the same highest weight and degree as the short-root term, and
opposite sign. For $-\varepsilon_1-\varepsilon_j$, the matching
forces the coefficient to be nonpositive. These exhaust all
possibilities.
Thus for $p<r$ there is no degree-zero term, and for $p\geq r$ the two
terms cancel as characters of finite-dimensional modules. The degree-zero coefficient of the reduced Jantzen sum is zero,
so $H^0_{(3,1^{2r-2})}(L_{p-2r+1}(B_r))$ is nonzero.
In particular this applies to every $B_r$ reduction used in this subsection,
where $p\geq5$ and $r\geq(p+3)/2$.
\end{proof}

\medskip
Suppose now that $b=0$ and $a\geq5$. Thus $n=3(a+1)/2$,
$p=a+2$ and $k=-2a$.
Let $f\in\mathfrak{so}_{2n+1}$ have partition $(3^a,1^4)$, and set
$$
 \mathcal W=\mathcal W^k(\mathfrak{so}_{2n+1},f),\qquad
 K_i=H_f^0(J^i)\subset\mathcal W,\qquad H=\mathcal W/K_1,
$$
where $J^i$ is the Jantzen filtration of $V^k(\mathfrak{so}_{2n+1})$. All characters below record the action of
$\mathfrak{so}_a\oplus\mathfrak{so}_4$. We use $Q$ for conformal
degree. The defining representation decomposes under the triple as
$S_2\otimes U\oplus B$, where $\dim U=a$ and $\dim B=4$.

\begin{proposition}\label{B:boundary-relations}
The reduction $H$ is nonzero. The image of $\mathfrak{so}_4$ in $H$
is zero, and the homomorphism $V^4(\mathfrak{so}_a)\to H$ factors through
$L_4(\mathfrak{so}_a)$. There is a relation in $K_1$ of weight two
whose coefficient on the scalar strong generator is nonzero.
Consequently $H$ is $C_2$-cofinite.
\end{proposition}

\begin{proof}
Set $x=(1^a,0^{n-a})$, and let $s\in\Z_{>0}$ be the integer
index of a real affine root $\alpha+s\delta$ in the Jantzen sum.
We use the notation $E_\lambda$, $\chi_\lambda$, $E_V(z)$ and
$H_V(z)$ introduced above. In Lemma~\ref{B:uniform-mode-bound},
the two root bounds are $x_1+x_2=2$ and $x_1=1$.
Long roots with $s\geq3$ have $c\geq5$, and short
roots with $s\geq2$ have $c\geq a+6$. They cannot contribute in
degrees at most two.

By Lemma~\ref{B:uniform-mode-one}, the contribution from the mixed
and negative short roots in mode one is
\begin{equation}\label{B:boundary-hooks}
 \sum_{r=1}^{p}(-1)^{r+1}Q^{a+r}\chi_{(r,1^{p-r})}.
\end{equation}
The other mode-one roots have reduced degree greater than two.
A remaining positive short root has
$$
 \lambda=(2p-i+1,1^{i-1}),\quad D=2p+2n-2i+1,
 \quad D-(\lambda,x)\geq(a+5)/2.
$$
A remaining positive sum root has $1\leq i<j\leq p$, with
$$
 \lambda=(p-i+1,p-j+2,2^{i-1},1^{j-i-1}),\quad
 D=p+2n-i-j+1,
$$
and $D-(\lambda,x)\geq3$.

In mode two, the regular negative sum roots are exactly $a\leq i<j$.
They give
\begin{equation}\label{B:mode-two}
 \lambda=(j-a+1,i-a+2,2^{i-1},1^{j-i-1}),\quad
 D=2(i+j-a),\quad D-(\lambda,x)=i+j-2a+1.
\end{equation}
Only $i=a,j=a+1$ contributes through degree two, with positive sign:
\begin{equation}\label{B:degree-two-rectangle}
                    Q^{2a+2}\chi_{(2^{a+1})}.
\end{equation}
Indeed, the reflected positions for a negative sum
root in mode two are $a-j$ and $a-i$. If one is positive, its absolute
value equals an unmodified coordinate, violating regularity. The possible self-collision would require
$2i=a$, which is impossible. For a mixed root the regular alternatives
are $i<a,j=a-i$, or $i\geq a$. In the first alternative $D=2|\lambda|$.
In the second, setting $l=\min(i,j)$ and $h=\max(i,j)$ gives
$$
 \lambda=(2p-j+1,i-a+2,2^{l-1},1^{h-l-1}),
 \quad D=2(2p+i-j),\quad
 D-(\lambda,x)\geq(a+7)/2.
$$
A positive sum root has
$\lambda=(2p-i+1,2p-j+2,2^{i-1},1^{j-i-1})$ and
$D=2(2p+2n-i-j+1)$, giving $D-(\lambda,x)\geq a+7$.
This proves that \eqref{B:boundary-hooks} and
\eqref{B:degree-two-rectangle} are the full numerator through
degree two.

\medskip
In the graded representation ring,
$$
 V=Q^{-1}U+(U+B)+QU,
 \qquad
 E_V(z)H_V(-Qz)
 =\frac{E_U(Q^{-1}z)E_B(z)}{E_B(Qz)E_U(Q^2z)}.
$$
By Lemma~\ref{B:uniform-hook-identity}, the hook contribution is
$$
 [z^{a+2}]\frac{Q^{a+1}(1-Q^2z^2)}{(1-Q)(1-Qz^2)}
       \frac{E_U(Q^{-1}z)E_B(z)}{E_B(Qz)E_U(Q^2z)}.
$$
Its coefficient of $Q$ is $\Lambda^2B$. Its coefficient of $Q^2$ is
$U\otimes B-\operatorname{Sym}^2_0B$: the terms
$\Lambda^aU=\mathbb C$ and $\Lambda^{a-1}U=U$ give respectively
$$
 \Lambda^2B+\mathbb C-B\otimes B
       =-\operatorname{Sym}^2_0B,
 \qquad U\otimes\Lambda^3B=U\otimes B.
$$
The highest $x$-eigenspace of $E_{(2^{a+1})}$, restricted from
$\mathfrak{gl}_a\oplus\mathfrak{so}_{a+4}$, is
$\det(U)^2\otimes\operatorname{Sym}^2_0(U+B)$.
On $\mathfrak{so}_a\oplus\mathfrak{so}_4$ its character is
$$
 \operatorname{Sym}^2_0U+U\otimes B+
       \operatorname{Sym}^2_0B+\mathbb C.
$$
Consequently the complete reduced numerator begins
\begin{equation}\label{B:two-window}
 Q\operatorname{ch}\Lambda^2B+
 Q^2\operatorname{ch}
   \bigl(\operatorname{Sym}^2_0U+2(U\otimes B)+\mathbb C\bigr).
\end{equation}
The completed Jantzen formula multiplies this by
$\operatorname{ch}\mathcal W$, whose weight-one space is
$\Lambda^2U\oplus\Lambda^2B$. It follows that
\begin{equation}\label{B:actual-small-degrees}
 \sum_{i\geq1}\operatorname{ch}(K_i)_0=0,\qquad
 \sum_{i\geq1}\operatorname{ch}(K_i)_1=\operatorname{ch}\mathfrak{so}_4,
 \qquad
 \sum_{i\geq1}\dim (K_i)_2^{\mathfrak{so}_a\oplus\mathfrak{so}_4}=3.
\end{equation}
The last number is one from \eqref{B:two-window} and two
from the invariant pairings of the two simple summands of
$\mathfrak{so}_4$.

\medskip
The first equality in \eqref{B:actual-small-degrees} makes $K_1$
proper and proves $H\ne0$. The second gives
$(K_1)_1=\mathfrak{so}_4$ and $(K_i)_1=0$ for $i\geq2$.
The levels of the homomorphisms from the universal affine vertex algebras are
$$
 3p-(3a+2)=4\quad\text{on }\mathfrak{so}_a,
 \qquad p-(a+2)=0\quad\text{on }\mathfrak{so}_4.
$$
These follow from the invariant bilinear form on the weight-one subspace, using
$\mathfrak g_0=\mathfrak{gl}_a\oplus\mathfrak{so}_{a+4}$.
Let $C_+,C_-$ be the quadratic Casimir vectors of the two
$\mathfrak{sl}_2$ summands of $\mathfrak{so}_4$. Both belong to $K_1$.
For $u$ in the corresponding summand,
$u(1)C_\pm=4u(-1)\mathbf1$, and the other summand acts by zero.
The resulting map is injective on $\mathbb CC_+\oplus\mathbb CC_-$.
Each mode $u(1)$ vanishes on $(K_i)_2$ for $i\geq2$, since
$(K_i)_1=0$. Nestedness and the last equality in
\eqref{B:actual-small-degrees} therefore give
$$
 \dim(K_1)_2^{\mathfrak{so}_a\oplus\mathfrak{so}_4}=3,
 \qquad (K_i)_2^{\mathfrak{so}_a\oplus\mathfrak{so}_4}=0\quad(i\geq2).
$$
Indeed a nonzero invariant in $K_2$ would be independent of the two
Casimirs in $K_1$, contributing at least four copies to the total sum.

Let $\theta$ be the highest long root of $\mathfrak{so}_a$,
let $Z$ be the scalar strong generator of weight two, and set
$C_0=\sum_jv_j(-1)v^j(-1)\mathbf1$, where $v_j,v^j$ are dual
bases of $\mathfrak{so}_a$ for the normalized invariant form.
The invariant weight-two subspace of $\mathcal W$ is
$\C Z+\C C_0+\C C_++\C C_-$.
The canonical contravariant form on $V^4(\mathfrak{so}_a)$ is
nondegenerate through degree two: its maximal proper graded ideal is
generated by $e_\theta(-1)^5\mathbf1$, of degree five.
In particular, a nonzero vector in its degree-two subspace is sent
to a nonzero multiple of $\mathbf1$ by a composition of positive
modes. Its image under $V^4(\mathfrak{so}_a)\to\mathcal W$
therefore cannot belong to the proper ideal $K_1$.
For $C_0$ this is explicit: for $u,v\in\mathfrak{so}_a$,
$$
 u(1)C_0=2(a+2)u(-1)\mathbf1,
 \qquad
 v(1)u(1)C_0=8(a+2)(v,u)\mathbf1.
$$
Thus $C_0\notin K_1$. If every invariant vector in $(K_1)_2$
had zero coefficient on $Z$, subtracting a linear combination of
$C_+$ and $C_-$ would give a multiple of $C_0$ in $K_1$.
That multiple must be zero, so the invariant subspace of $(K_1)_2$
would have dimension at most two, contradicting
\eqref{B:actual-small-degrees}. Therefore a relation in $K_1$
has a nonzero coefficient on $Z$.

\medskip
Set $y=x+\theta^\vee/2$. Both gradings can be chosen dominant, with
$$
 y=((3/2)^2,1^{a-4},(1/2)^4,0^{(a-1)/2}).
$$
The character decomposition for $y$ is taken coefficientwise in the
original conformal degree. For $y$, Lemma~\ref{B:uniform-mode-bound} leaves only modes
$s\leq3$ for long roots and $s=1$ for short roots at reduced
$y$-degree zero. Set $S_l=\sum_{i\leq l}y_i$ and
$F_l=S_l-l/2$. Then $F_l\leq a/2$ for every $l$.
In mode one the positive sum roots with $j\leq p$ have
$D-(\lambda,y)=a+1-F_{i+1}-F_j\geq1$.
The positive short roots have
$D-(\lambda,y)=2a+2-i/2-S_i\geq(a+1)/4>0$.
Thus only the hooks in \eqref{B:boundary-hooks} remain.

All mode-two terms also have strictly positive $y$-degree.
For \eqref{B:mode-two} it is
$a+1-F_{i+1}-F_j\geq1$. The mixed-root case $i<a$ has
$D=2|\lambda|$ and $\lambda(y)\leq3|\lambda|/2$.
For $i\geq a$ its degree is at least $(a+3)/4$. A positive sum
root has degree at least $(a+5)/2$. These bounds follow by inserting
the displayed weights and using $S_l\leq a+1$ and $i+j\leq2n-1$.
In mode three equality in Lemma~\ref{B:uniform-mode-bound} requires
$w\in W_y$ and $(\alpha,y)=-3$, hence
$\alpha=-\epsilon_1-\epsilon_2$. Its contribution is precisely
$Q^{15}\chi_{(5,5)}$, with positive sign.

For the hook $(r,1^{p-r})$ the least $y$-degree is
$a/2-F_{a+3-r}\geq0$, with equality exactly when $1\leq r\leq5$.
Their highest $y$-eigenspaces have the following decomposition.
The two coordinates of $y$-weight $3/2$ give
$S_{(r,1)}\mathbb C^2$, of dimension $r$. The coordinates of
$y$-weight one give a determinant, and the remaining four give an
exterior power of degree $5-r$. Among these four coordinates two
have $x$-weight one and two have $x$-weight zero. Therefore their
combined contribution, with conformal degree recorded, is
$$
 Q^3\sum_{r=1}^{5}(-1)^{r+1}r[z^{5-r}]P(z)
 =Q^3\bigl(5P(-1)+P'(-1)\bigr)=0,
 \qquad P(z)=(1+Q^{-1}z)^2(1+z)^2.
$$
The highest $y$-eigenspace of $E_{(5,5)}$ is the one-dimensional
determinant power on the first two coordinates. Its $x$-weight is
ten, so the mode-three term has conformal degree five.

The only strong generator of $\mathcal W$ whose conformal weight
equals its $\theta^\vee/2$-weight is $e_\theta(-1)\mathbf1$.
Indeed the complete list of strong generators is
$$
 \begin{array}{c|c}
 1&\Lambda^2U\oplus\Lambda^2B\\
 2&\operatorname{Sym}^2U\oplus(U\otimes B)\\
 3&\Lambda^2U.
 \end{array}
$$
Every derivative has conformal weight strictly greater than its
$\theta^\vee/2$-weight. Thus the corresponding part of the universal character is
$(1-Q)^{-1}$, generated by powers of $e_\theta(-1)\mathbf1$.
The corresponding part of the total Jantzen character is consequently
$Q^5/(1-Q)$. The subspace of $\mathcal W$ of conformal weight five
and $\theta^\vee$-weight ten is
$\C e_\theta(-1)^5\mathbf1$. Its total multiplicity in the nested
ideals is one, so
$$
                  e_\theta(-1)^5\mathbf1\in K_1.
$$
The homomorphism $V^4(\mathfrak{so}_a)\to H$ therefore factors
through $L_4(\mathfrak{so}_a)$. The classes in $R_H$ of all weight-one
generators are nilpotent: the image of $\mathfrak{so}_4$ is zero,
and $L_4(\mathfrak{so}_a)$ is integrable. By Proposition~\ref{com:orthogonal-point} with $m=0$, the reduced
Poisson algebra $(R_H)_{\mathrm{red}}$ is generated by the class of the
scalar generator of weight two. The relation above sets this class to zero. Hence $H$ is $C_2$-cofinite, as asserted.
\end{proof}

\medskip
For general $a,b$, we now show that the remaining points $H_j$ do not belong to the associated variety.
Propositions
\ref{B:uniform-five-hook} and \ref{B:uniform-triple-hook} give
\begin{equation}\label{seed:B-strip-vacuum}
 H_{(5,1^{2n-4})}=0,\qquad H_{(3^{a+1},1^{4b+1})}=0.
\end{equation}
The Pfaffian of size $2(a+b+1)$ imposes rank at most $2(a+b)$,
so $H_j$ cannot belong to the associated variety for $j\geq a+b+1$.
Suppose that $H_j\in X_{L_k(B_n)}$ for some $j\leq a+b-1$, and set
$r=n-j$. Then
$r\geq(p+3)/2$ and Proposition~\ref{B:uniform-tail-survival} applies
to $L_{p-2r+1}(B_r)$. The support inclusion associated with the parabolic homomorphism induces the orbit
$(5,3^{j-1},1^{2r-j-1})$. Its closure contains
$(5,1^{2n-4})$, contradicting \eqref{seed:B-strip-vacuum}.

The remaining point is $H_{a+b}$, whose Richardson partition is
$(3^{a+b},1^{b+4})$. If $b\geq1$, this partition dominates
$(3^{a+1},1^{4b+1})$, again contradicting
\eqref{seed:B-strip-vacuum}.
For $b=0$ and $a\geq5$, Proposition~\ref{B:boundary-relations}
shows that $H_f^0(L_k(B_n))$ is $C_2$-cofinite for $f=(3^a,1^4)$.
If $H_a$ belonged to the associated variety, its conical orbit closure
would have dimension
$\dim\OO_f+1$ and hence a positive-dimensional $f$-slice,
contrary to $X_{H_f^0(L_k(B_n))}=X_{L_k(B_n)}\cap\mathcal S_f$.
The case $a=3,b=0$ is precisely the $B_6$ constituent case already
proved above. This completes the remaining cases of
\cref{thm:B-classification}.

\section{\texorpdfstring{Quasi-lisse simple affine vertex algebras of type $C$}{Quasi-lisse simple affine vertex algebras of type C}}\label{sec:type-C}

The representation theory of symplectic affine vertex algebras at
half-integral levels was studied in \cite{AdamovicThesis1996}.
The determinantal singular vectors constructed by Adamovi\'c in
\cite[Theorem~4.1]{Adamovic2003} provide the polynomial relations
used below. We combine these known relations with vanishing and
nonvanishing results for Drinfeld--Sokolov reductions.

The non-anomalous twisted $D_{n+1}$ construction gives the proposed
vertex algebras $L_{-n-1+(n+1)/q}(C_n)$ for odd $q>1$ coprime
to $n+1$ \cite[Section~5]{WangXie2019}.
For anomalous twisted $A_{2n}$ theories, Beem and Kulkarni instead
propose a finite extension of $L_\ell(C_n)$, where
$$
 \ell=-n-1+\frac{2n+1}{2s},\qquad
 s=2m+3,\quad m\ge0,\quad\gcd(s,2n+1)=1.
$$
In this description the $\mathbb Z_2$-gauged theory has vertex algebra
$L_\ell(C_n)$ itself. The ungauged Higgs branch can be a ramified
double cover of the affine associated variety
\cite[Sections~4--6]{BeemKulkarni2026}.
At level $-1/2$, $L_{-1/2}(C_n)$ is the even symplectic-boson
algebra, associated with the simultaneous sign gauging of $n$ free
hypermultiplets. Its associated variety is
$\overline{\OO_{\min}(C_n)}\cong\C^{2n}/\{\pm1\}$
\cite[Section~1.1]{BeemKulkarni2026}.

Let $\g=\spg_{2n}$, $n\ge2$. The denominator theorem reduces the
classification to levels $k=-n-1+P/2$ with $P\in\Z_{>0}$.

\begin{theorem}\label{thm:C-seeds}
The vertex algebra $L_{-n-1+P/2}(C_n)$ is quasi-lisse if and only if
$P\ge n+2$ and $P\ne2n$.
\end{theorem}
\begin{proof}
Admissibility covers $P\ge2n+1$. The range $n+2\le P\le2n-1$,
which is $k=-m/2$, $3\le m\le n$, follows from
\cref{seed:C-positive}, and the exclusions from \cref{sec:C-exclusions}.
For $n=2$ the affirmative nonadmissible range is empty.
\end{proof}

\subsection{\texorpdfstring{The associated variety of $L_{-m/2}(\spg_{2n})$}{The associated variety of L(-m/2)(sp(2n))}}\label{seed:C-positive}
We prove that $L_{-m/2}(\spg_{2n})$ is quasi-lisse for $3\le m\le n$
and identify its associated variety. The determinant relation bounds
the possible nilpotent and semisimple points. We combine it with
vanishing reductions to rule out the semisimple points, and with a
nonzero reduction to establish the required nilpotent orbit. Set $k=-m/2$, so $P=2n+2-m$, and
$f_m=(2^m,1^{2n-2m})$. The conformal-weight inequalities use
$$
 F(\omega_i)=i(m-i)\quad(i\le m),\qquad
 F(\omega_i)=(i-m)(2n+2-i-m)\quad(i>m).
$$
Only $\omega_m$ has zero value. If $m$ is odd it is not in the root
lattice. If $m$ is even its degree is $m/2$ and $f_m^2=0$, so
\cref{seed:covariants} shows that all $C_2$ polynomials from the
possible equality subspace vanish at $f_m$. Thus $H_{f_m}\ne0$.

\begin{proposition}\label{seed:C-vacuum}
Let $3\le m\le n$ and $k=-m/2$. Then
$H^0_{(4,1^{2n-4})}(L_k(C_n))=0$. If $m\ge4$, then
$H^0_{(3^2,1^{2n-6})}(L_k(C_n))=0$ as well.
\end{proposition}
\begin{proof}
For a root term in Euclidean epsilon coordinates set
$b=2(\rho\cdot\alpha+sP)/(\alpha\cdot\alpha)>0$ integral,
$d=sb$, and use the signed finite alternant of $\rho-b\alpha$.
For $(3^2,1^{2n-6})$ only $\lambda_1+\lambda_2\ge d$ contributes.
Since $|\lambda|\le2b$, the exhaustive bound is $s\le2$.
With $A=2(n-m+2)+1$, the complete sum is
$+\Delta((A,A,0,\ldots))[2A]$. Its maximal-$x$ eigenspace is the one-dimensional
$\det^A$ of $\glg_2$, so the corresponding degree-zero coefficient after reduction is one.

For the hook the confluent numerator is
$\det[z^{3u_i}-z^{-3u_i},z^{u_i}-z^{-u_i},u_i,u_i^3,\ldots]$.
Its denominator has top exponent $4n-1$.
A selected signed pair $\xi=3\epsilon e_R+\tau e_S$ can contribute
only when $b(2s+\xi\cdot\alpha)=\xi\cdot\rho-(4n-1)\le0$,
so $s\le3$. The exact signed sum is one for $m=3$ and
$1+\mathbf1_{m\text{ even}}$ for $m\ge4$.
The symbolic alternant cancellations, including singular raw terms
and the missing zero-$b$ partner at $m=3$, are
\hyperref[calc:BCD]{Calculation~\ref*{calc:BCD}}. With the displayed uniform mode bounds, completion and positivity
place $\mathbf1$ in the first reduced Jantzen ideal for every $3\le m\le n$.
\end{proof}

Suppose first $n>m$. Recall the size-$(m+1)$ determinant singular
vector at level $-m/2$ from \cite[Theorem~4.1]{Adamovic2003}. Its
$C_2$ image imposes
$\det(Q_X|_E)=0$ on every isotropic $(m+1)$-plane, where
$Q_X(v,w)=\omega(Xv,w)$. For a semisimple matrix $s$ of rank $2j$,
\begin{equation}\label{seed:C-rank}
 \max_{L\text{ Lagrangian}}\rank(Q_s|_L)=\min(n,2j).
\end{equation}
For the lower bound pair $\min(j,n-j)$ nonzero eigenpairs with zero
pairs and take $e_i+e_z,f_i-f_z$. Their $Q_s$ block has off-diagonal
entry $a_i\ne0$. On remaining nonzero pairs use $e_i+f_i$, of value
$2a_i$. The upper bound is immediate. Thus $\rank s\le m$.
A remaining rank at least four has two nonzero eigenpairs and a zero
pair. Companion-matrix degeneration on the isotropic three-space
$\langle e_1,e_0,f_2\rangle$, choosing distinct eigenvalues $a,0,-b$,
and its dual gives the forbidden $(3^2,1^{2n-6})$. For $m=3$ the
rank bound already leaves only rank two.

A rank-two point is, up to scale, $h_1$. The parabolic homomorphism associated with $h_1$
would give a parabolic homomorphism whose affine Levi factor is the simple
vertex algebra $L_{1-m/2}(C_{n-1})$. Its negative,
nonintegrable level has the minimal nilpotent in its support by the
minimal-DS criterion. Thus $ah_1+y\in X$. On a symplectic four-space
set
$$
 h_1e_0=e_0,\ h_1f_0=-f_0,\ yf_1=e_1,
 \qquad ze_0=f_1,\ ze_1=f_0.
$$
For $a\ne0$, $ah_1+y+z$ is conjugate to $ah_1+y$. At $a=0$ it has
the chain $e_0\to f_1\to e_1\to f_0\to0$, the forbidden hook.
None of the candidate semisimple points belongs to the associated
variety, proving quasi-lisseness. The same determinant rules out
nilpotent rank $>m$ by degeneration to the square-zero rank-$(m+1)$
orbit. The two vanishing reductions leave only square-zero possibilities
(only the hook is needed at $m=3$).
Consequently $X_{L_{-m/2}(C_n)}=\clO{2^m,1^{2n-2m}}$.

At $m=n$ there is no isotropic $(n+1)$-plane, so the preceding
determinant relation is unavailable. Instead, we identify the
reduction at $(2^n)$ as a simple integrable affine vertex algebra.
Its associated variety is a point, which gives the orbit
restriction needed at this endpoint.

\begin{proposition}\label{seed:C-full}
For $n=m\ge3$, $H^0_{(2^n)}(L_{-n/2}(C_n))$ is
$L_4(\sog_n)$, interpreted as $L_8(\slg_2)$ for $n=3$ and two
$L_4(\slg_2)$ factors for $n=4$. The vertex algebra $L_{-n/2}(C_n)$ is quasi-lisse, and its associated
variety is $\clO{2^n}$.
\end{proposition}
\begin{proof}
The depth-one grading $x=\frac12\sum h_i$ has
$\g_0=\glg_n$ and abelian $\g_{\pm1}$. Twisting by $w_0t_{-x}$
sends positive Levi roots to negative Levi zero modes and positive
degree-one roots to negative-root mode one. All annihilate $\mathbf1$.
The resulting simple highest-weight module is $L(k\Lambda_n)$, and
its highest-weight space for the Siegel parabolic is one-dimensional. It belongs to Arakawa's $\mathcal O_{0,k}$ because
$k\omega_n$ is integral dominant on the Levi $\slg_n$.
This argument does not require $L(k\Lambda_n)$ to belong to
the Kazhdan--Lusztig category or to have a finite-dimensional
$\g$-module as its lowest conformal-weight space. 
The ghost shifts $\psi_\alpha(r)\mapsto\psi_{-\alpha}(r+1)$ and
its inverse preserve both annihilator thresholds $r\ge0,r\ge1$.
Thus the plus/minus identification of \cite[Section~4.3]{Arakawa2011}
applies to $L_{-n/2}(\spg_{2n})$.
The corresponding generalized Verma $\g$-module is induced from
this one-dimensional representation of the Siegel parabolic. It is
free of rank one over $U(\g_{-1})$, and its Whittaker coinvariants
are one-dimensional. Actual nonvanishing
implies maximal GK dimension by
\cite[Theorem~5.5.4(ii)]{Arakawa2011}, hence a nonzero
Whittaker top by \cite[Theorem~5.1.1]{Arakawa2011}.
By \cite[Proposition~5.4.2]{Arakawa2011}, this top is
a quotient of the preceding one-dimensional space. It is simple.
Thus \cite[Corollary~5.5.5]{Arakawa2011} proves that
$H^0_{(2^n)}(L_{-n/2}(C_n))$ is simple. These hypotheses
are specific to the depth-one polarization.

The universal $W$-algebra is strongly generated by the image of
$\sog_n$ under $V^4(\sog_n)\to\mathcal W^{-n/2}(\spg_{2n},f)$
and a copy of $\Sym^2\C^n$ of weight two. Its central charge is
$2n(n-1)/(n+2)$, equal to the Sugawara central charge of $V^4(\sog_n)$.
Modulo the affine vertex subalgebra, the invariant weight-two subspace
has dimension at most one.
The invariant-form criterion gives $T=S$. By \cref{v3:resonance},
the weight-two strong generators in the traceless symmetric-square
representation vanish: their Sugawara conformal weight is
$n/(n+2)<1$, whereas their conformal weight is two. Thus $H^0_{(2^n)}(L_{-n/2}(C_n))$ is the asserted simple
integrable affine vertex algebra.

To exclude a semisimple $s$ of rank $2j\ge n$, choose a nonsingular
symmetric matrix $B$ supported on nonzero eigenvalue diagonals and
nonzero/zero pairs. Solve $DZ+ZD=-B$. The only denominators are
$2a_i$ or $a_i\ne0$. Then
$\varepsilon\operatorname{Ad}\exp(\left(\begin{smallmatrix}0&Z\\0&0\end{smallmatrix}\right)/\varepsilon)(s)$
tends to $\left(\begin{smallmatrix}0&B\\0&0\end{smallmatrix}\right)$,
of type $(2^n)$. The zero-dimensional Slodowy intersection forces every component
through this orbit to be its nilpotent closure, contradicting the
irreducible semisimple cone. The remaining rank at least four and
rank-two cases cannot occur, by the preceding degenerations to
$(3^2,1^{2n-6})$ and $(4,1^{2n-4})$, respectively.
For $n\ge4$, the vanishing reductions leave only square-zero nilpotents.
For $n=3$ every symplectic partition of six outside $\clO{2^3}$
dominates $(2^3)$, so its orbit cannot occur when the Slodowy intersection is zero-dimensional.
In particular $\dim X_{L_{-3/2}(C_3)}=3(6-3+1)=12$.
\end{proof}

\subsection{\texorpdfstring{Non-quasi-lisseness of $L_{P/2-n-1}(\spg_{2n})$}{Non-quasi-lisseness of L(P/2-n-1)(sp(2n))}}\label{sec:C-exclusions}
We prove that $L_k(\spg_{2n})$, with $k=P/2-n-1$, is not quasi-lisse
for $1\le P\le n+1$ and for $P=2n$.
For $2\le P\le n$, write $n=jP+r$, $0\le r<P$. The following
partitions specify the reductions used in each of the disjoint cases:
\begin{equation}\label{seed:C-negative-list}
\begin{array}{c|l}
r=0&((2j)^P)\\
r=P-c,\ c=1,2&((2j+2)^{P-c},(2j)^c)\\
P\text{ even},\ r=P/2-1&((2j+1)^{P-2},1^{4j})\\
1\le r\le\lfloor(P-3)/2\rfloor&
 ((2j+1)^{2r},(2j)^2,1^{2jp'}),\ p'=P-2r-2\\
\lceil(P-1)/2\rceil\le r\le P-3&
 ((2j+2)^v,(2j+1)^2,1^{2M}),\\
 &v=2r-P+2,\ b=P-r-2,\ M=(2j+1)b-1.
\end{array}
\end{equation}
The $c=2$ row requires $P\ge3$ and the third row $P\ge4$.
$P=2,3$ are exhausted by the remaining rows.
For the second row,
$A_{r+gP}=(2-c)r+2Pg$. Only $c=2,l=r$ is zero.
For the third row, set $n=(2a+1)b-1$, $T=2b$,
$f=((2a+1)^{2u},1^{2M})$, $M=(2a+1)(b-u)-1$, $1\le u<b$.
Then $A_{2ug}=4u(b-u)g^2>0$. The listed row has $a=j,u=b-1$.
For the fourth row, set $h=2r+2$ and
$B_g=g[2h+p'(h(g-1)+2)]$. The values at the block endpoints are
$A_{gh}=B_g$ and $A_{gh+2r}=B_g+4r(1+gp')$, all positive.
For the fifth row, set $B_g=4bg^2+2bvg(g+1)$. The values at the block endpoints are
$A_{v+g(v+2)}=B_g$ and
$A_{v+g(v+2)+2}=B_g+4b(2g+1)$. Only $l=v$ is zero.
These formulas follow by summing each arithmetic progression of chain
eigenvalues. Interpolation determines every interior value.
In the two rows where the conformal-weight inequality is not strict the equality degree is $r_0(2j+1)/2$,
$r_0=r$ or $v$, and $f^{2j+2}=0$. Odd $r_0$ is outside the root
lattice. For even $r_0$, \cref{seed:covariants} rules out the equality case.
Thus every listed reduction is nonzero.

For $r=0$, the image is $\mathcal A_j$ by Proposition~\ref{v6:scalar-partitions}.
For $r=P-c$, parabolic induction gives the image in
$(L_c(\slg_r)\otimes\mathcal H)\otimes\mathcal A_j$.
The even-chain calculation \eqref{seed:telescope}, with $s=j$ and
$\epsilon=c$, produces $1\otimes T_j$ and the independent quartic
symbol. Proposition~\ref{seed:scalar-split} then gives the tensor
factor $\mathcal A_j$.

In the third row, the homomorphisms from $V^0(\spg_{2u})$ and
$V^\ell(\spg_{2M})$, where $\ell=-2a(b-u)$ and $\ell+M+1=b-u$,
have commuting images. The image of $\spg_{2u}$ and all mixed
generators vanish in the simple quotient. The remaining strong
generators are the image of $\spg_{2M}$ and possible
$\spg_{2M}$-invariant generators of weights $2,4,\ldots,2a$.
The central charge is
$$
 c=-2aM(2M+1)=\frac{\ell\dim C_M}{\ell+M+1}.
$$
Modulo the affine vertex subalgebra generated by $\spg_{2M}$,
the invariant weight-two subspace has dimension at most one.
The invariant-form criterion and \cref{v3:resonance} give $T=S$ and
show that every remaining $\spg_{2M}$-invariant generator of positive
weight vanishes, since its Sugawara conformal weight is zero.
The simple quotient of the reduction is therefore $L_\ell(C_M)$.
Taking $u=b-1$ reduces the problem to the $P=2$ rectangle of rank $2a$, already treated above.

In the fourth and fifth rows, the homomorphisms associated with
$\sog_2$ and $\spg_2$, respectively, have level zero. The decomposition
associated with the Jordan blocks and its parity give tensor factors
$L_{p'/2-jp'-1}(\spg_{2jp'})$ and
$L_{-2jb}(\spg_{2((2j+1)b-1)})$, respectively. The first is covered
by the rectangular case, with its universal affine vertex algebra already
simple when $p'=1$. The second is covered by the third row, or by the
$P=2$ row when $b=1$.

It remains to prove non-quasi-lisseness for $P=1,n+1,2n$.
For $P=1$, $V^{-n-1/2}(\spg_{2n})$ is simple. For $P=2n$ use $f=(2^2,1^{2n-4})$.
The only equality is $\omega_2$ in degree one and the covariant lemma
proves nonvanishing of the reduction. The homomorphism associated with
$\spg_{2n-4}$ has level zero. The weight-one generator $J$ satisfies
$J_{(1)}J=2n\mathbf1$. The mixed generators vanish. The kernel of
$J_{(0)}$ in the weight-two subspace is exactly
$\langle:JJ:,\partial J,T\rangle$, with $c=1$.
The form criterion gives $T=:JJ:/(4n)$. \Cref{v3:resonance} shows that the weight-two generators on which $J_{(0)}$ acts by
$\pm2$ vanish, since their conformal weight for $T=:JJ:/(4n)$ is $1/n$.
The simple quotient of the reduction is the rank-one Heisenberg vertex
algebra generated by $J$.

For $P=n+1$, the conformal-weight inequalities for the partition
$(2^n)$ have the positive values $i(n+1-i)>0$.
For $n\ge3$, the Miura homomorphism takes values in
$L_1(\slg_n)\otimes\mathcal H$, where $\mathcal H$ has its standard conformal vector $\tfrac12:hh:$,
with $[h_\lambda h]=\lambda\mathbf1$. Its
weight-one subspace corresponding to $\sog_n$ generates $L_2(\sog_n)$ (use $L_4(\slg_2)$ at
$n=3$). The projected symmetric traceless weight-two fields have
Sugawara weight one and lie in the unique
$\Sym^2_0\C^n\otimes h$ space. This space, the weight-one subspace and the conformal vector
are PCT-stable and generate the image. The image is simple and its
conformal vector maps to $h^2/2$ after setting the classes from
$L_1(\slg_n)$ equal to zero in the reduced $C_2$ algebra.
For $n=2$, the exact realization in $V_{\Z\alpha}\otimes\mathcal H$,
$(\alpha,\alpha)=2$, is
$$
 J=\alpha/2,\quad T=:JJ:+:hh:/2,\quad
 G^\pm=\sqrt{3/8}\,e^{\pm\alpha}h.
$$
With lattice cocycle $\varepsilon(m\alpha,n\alpha)=1$, the operator
product expansion in \hyperref[calc:C2-IMAGE]{Calculation~\ref*{calc:C2-IMAGE}} is
as follows. Writing $s=z-w$ and evaluating the fields on the right at $w$, we have
$$
G^+(z)G^-(w)\sim\frac{3}{8s^4}+\frac{3J}{4s^3}
+\frac{3T/4+3\partial J/8}{s^2}
+\frac{3:TJ:/2+3\partial T/8-:J^3:-5\partial^2J/8}{s}.
$$
The equal-sign products are regular and $J,G^\pm$ primary. These are
the specialized relations of \cite[Section~5]{Fasquel2022}, defining the map. The involution
$\alpha,h\mapsto-\alpha,-h$, $e^{m\alpha}\mapsto(-1)^{m^2}e^{-m\alpha}$
preserves the image. Its restricted positive form proves simplicity.
Setting the classes from $V_{\Z\alpha}$ equal to zero in the reduced
$C_2$ algebra sends its conformal vector to $h^2/2$.

\section{\texorpdfstring{Quasi-lisse simple affine vertex algebras of type $D$}{Quasi-lisse simple affine vertex algebras of type D}}\label{sec:type-D}

Per\v se constructed the harmonic singular-vector powers used in
this section in \cite[Theorem~3.1]{PerseD2013}. The associated
varieties at $k=-2,-1$ are the minimal nilpotent orbit closure by
\cite[Theorem~1.1]{ArakawaMoreau2018}. Our argument treats the full
range of integer numerators, including these previously known cases.
For $D_4$, \cite{AVD42026} proves quasi-lisseness at $k=-14/3$ and identifies the associated variety with the subregular orbit closure.

The algebra $L_{-2}(D_4)$ is the proposed vertex algebra of $SU(2)$
gauge theory with four fundamental hypermultiplets
\cite[Section~5.1.3]{BeemEtAl2015}. Its minimal nilpotent associated
variety agrees with the Higgs branch \cite[Section~4]{BeemRastelli2018}.
The rank-two candidate $L_{-4}(D_6)$, with Coulomb-branch dimensions
$(4,2)$ and predicted Higgs branch $\overline{\OO_{(2^4,1^4)}}$,
appears in \cite[Table~5]{LiXieYan2023} and agrees with our orbit table.

Let $\g=\sog_{2n}$, $n\ge4$, and write $t=k+2n-2$.
By the denominator theorem, it suffices to consider $t=p\in\Z_{>0}$.

\begin{theorem}\label{thm:D-integer-classification}
For $p\in\Z_{>0}$, $L_{p-2n+2}(D_n)$ is quasi-lisse if and only if
$$
 p\ge2\lceil n/2\rceil,
 \qquad\text{or}\qquad
 p\text{ is even and }2n/3<p\le n-2.
$$
\end{theorem}
\begin{proof}
The first range is proved in \cref{sec:D-upper}, and the second in
\cref{seed:D-positive}. The remaining values are excluded by
\cref{sec:D-exclusions}, \cref{seed:D-growth}. For $p=1$, the assertion follows from the simplicity of
$V^{3-2n}(\sog_{2n})$.
\end{proof}

\subsection{\texorpdfstring{Quasi-lisseness of $L_{p-2n+2}(\sog_{2n})$ for $p\ge2\lceil n/2\rceil$}{Quasi-lisseness of L(p-2n+2)(so(2n)) for p at least 2ceil(n/2)}}\label{sec:D-upper}

We prove quasi-lisseness at $k=p-2n+2$ for
$p\ge2\lceil n/2\rceil$. When $2\lceil n/2\rceil\le p<2n-2$, the
harmonic singular-vector power of \cite[Theorem~3.1]{PerseD2013}
(see also \cite[Lemma~4.1 and Theorem~4.2]{ArakawaMoreau2018}) imposes $X^2=(\tr X^2/2n)I$. The natural Pfaffian for
even $p$ or the square for odd $p$ gives
$$
 \rank X\le r,\qquad
 r=2n-p-2\ (p\text{ even}),\quad r=2n-p-1\ (p\text{ odd}).
$$
Since $r<2n$, $X^2=0$. The strict conformal-weight inequality at $(2^r,1^{2n-2r})$ proves
nonvanishing of the reduction. The sole positive value at a block
endpoint is $A_r=2r$ or $r$. Thus the associated variety of $L_{p-2n+2}(\sog_{2n})$ is the closure of this orbit,
and in particular the affine vertex algebra is quasi-lisse. For $p\ge2n-2$ use
integrability.

For the endpoint $D_4,p=5$, \cite{ArakawaMoreau2018}, Theorem~1.1(2),
explicitly gives $X_{L_{-1}(D_4)}=\overline{\OO_{\min}}$. This proves
quasi-lisseness at $k=-1$. At $t=5/6$,
$$
 \Delta_{5/6}(\lambda)-\lambda(\rho)
 =\tfrac35\|\lambda\|^2+\tfrac15(\rho,\lambda)>0
 \quad(\lambda\ne0).
$$
Nonvanishing of the principal reduction and the denominator theorem therefore imply
$X_{L_{5/6-6}(D_4)}=\N(D_4)$. These distinct varieties illustrate that Theorem~\ref{thm:denominator} preserves quasi-lisseness,
not the associated variety.

\subsection{\texorpdfstring{Quasi-lisseness of $L_{p-2n+2}(\sog_{2n})$ for even $p$ with $2n/3<p\le n-2$}{Quasi-lisseness of L(p-2n+2)(so(2n)) for even p with 2n/3 less than p at most n-2}}\label{seed:D-positive}
We prove that $L_{p-2n+2}(\sog_{2n})$ is quasi-lisse for even $p$
with $2n/3<p\le n-2$. The required orbit bound comes from relations
in the Drinfeld--Sokolov reduction. We first show which weight-one
generators vanish, then identify the surviving integrable affine
vertex subalgebra. The character sums below establish these relations in the
reduction of the maximal ideal itself.

Write $n=3j+r$, $p=2j+r$, $j\ge2$, $r\ge2$ even, and
$f=(3^{2j},2^{r-2},1^4)$. The conformal-weight differences at the endpoint weights are
$A_{2j}=0$ and $A_{2j+r-2}=2r-4$ (the latter only when $r>2$).
The sole equality is $(1^{2j})$ at degree $2j$, since
$j<2j<n-j=2j+r$ and $f^3=0$, \cref{seed:covariants} proves
nonvanishing of the reduction, even for $j=1$ when needed in degenerations.

The weight-one subspace defines homomorphisms from
$V^{r+2}(\sog_{2j})$, $V^0(\spg_{r-2})$ and $V^0(\sog_4)$
to $\WW^k(\sog_{2n},f)$. The image of $\sog_4$ vanishes in
$H_f^0(L_k(\sog_{2n}))$, as follows from a linear equation on the Slodowy slice. Set
$R=\rank f=4j+r-2=2n-p-2$. Choose an $x$-graded $R$-space $E$
with $\sum_{E}x=R/2$ on which $(fu,v)$ is nondegenerate. A size
$R+2$ Pfaffian belonging to the defining ideal, evaluated on $E$ and two one-chain vectors,
has Kazhdan weight one and linear term $\operatorname{Pf}(f|E)$
times the corresponding $\sog_4$ coordinate. All slice generators have
weight at least one, so no nonlinear term is possible. Equivariance
excludes contributions from $\sog_{2j}$ and $\spg_{r-2}$.
The scheme-theoretic Drinfeld--Sokolov identity places this linear
coordinate in the defining ideal of $R_{H_f^0(L_k(\sog_{2n}))}$.
Since $C_2(H_f^0(L_k(\sog_{2n})))_1=0$, the corresponding
weight-one vector vanishes in the reduction.

\begin{proposition}\label{seed:D-windows}
Set $\ell=2j$. Write $(\mu)[D]$ for a character term with finite
highest weight $\mu$ and original Sugawara degree $D$.
For a Cartan element $y$, let $C$ be twice its dominant Cartan
coordinates. The least eigenvalue of $D-y(0)$ on the finite-dimensional
$\g$-module of highest weight $\mu$ is $D-C\cdot\mu/2$.
Only terms for which this value is at most $B$ enter the calculation.
In (1) and (2), take $B=0$ and $y=x+\beta^\vee/2$, where
$\beta^\vee$ is the coroot of a highest root of $\spg_{r-2}$ and
$\sog_\ell$, respectively (either simple factor when $\ell=4$).
This subtracts half the corresponding Cartan eigenvalue from the
reduced conformal degree. In (3) and (4), use the original
Drinfeld--Sokolov grading $y=x$ and $B=2$.
The complete character sums are as follows:
\begin{enumerate}
\item For $r\ge4$, $C=(2^{\ell+1},1^{r-4},0^{j+3})$ gives
$+(2^{\ell+1})[2\ell+2]$.
\item For $r\ge2$, $C=(3,3,2^{\ell-4},1^{r+2},0^j)$ gives
$$
 \sum_{u=1}^{r+3}(-1)^{u-1}(u,1^{p-u})
       [\ell+r/2+u-1]+(r+3,r+3)[3(r+3)].
$$
\item At $r=2$, through reduced degree two the list is
$$
 (1^{\ell+2})[\ell+1]-(2,1^\ell)[\ell+2]
 +(3,1^{\ell-1})[\ell+3]-(4,1^{\ell-2})[\ell+4]
 +(2^{\ell+1})[2\ell+2].
$$
\item For $p\ge5$ odd, $n=(m+1)p-1$, $H=2m+1$, and
$f=((2m+2)^{p-3},(2m+1)^3,1)$, the degree-two list is
$$
 (4)[4(m+1)]+(2^{p-1})[H(p-1)]-(2^p)[Hp-1].
$$
All three terms in the last list start in reduced degree two.
\end{enumerate}
\end{proposition}
These identities hold for all parameters in the stated ranges.
For a root contribution write $c=(\rho,\alpha)+hp>0$, $D=hc$,
and let $\rho+\mu$ be the dominant Weyl conjugate of $\rho-c\alpha$.
The passage to this conjugate is nonexpansive, hence
$\sum_i|\mu_i|\le2c$, $\max|\mu_i|\le c$, and
$(\mu,\mu+2\rho)=2pD$.
For $0\le C_i\le c_0$, we have
$D-C\cdot\mu/2\ge c(h-c_0)$. Since $c\ge1$, the cutoff
$D-C\cdot\mu/2\le B$ requires $h\le c_0+B$.
If $\rho-c\alpha$ is regular, its absolute coordinates have one or
two entries outside $\{0,\ldots,n-1\}$. Otherwise its quadratic
Casimir eigenvalue would force $D=0$.
The four sign choices with two coordinates outside this interval,
and the paired terms with one such coordinate, give the stated lists
by the polynomial inequalities in
\hyperref[calc:BCD]{Calculation~\ref*{calc:BCD}}.
When the removed $\rho$-coordinate is zero, the type-$D$ Weyl
alternants with negative and positive last coordinates cancel separately. The inequalities and Weyl alternants, including the terms whose removed
$\rho$-coordinate is zero, are recorded in \hyperref[calc:BCD]{Calculation~\ref*{calc:BCD}}. These are identities of full finite-dimensional Lie algebra characters.
By \cref{com:completion}, they determine the sum of the characters of
the nested ideals obtained by Drinfeld--Sokolov reduction.

We now determine the relations in the reduced affine ideal. The first character sum has a single
line of conformal weight one and charge two, generated by a long-root
vector of $\spg_{r-2}$. Positivity of multiplicities and nestedness place
this vector in the first reduced ideal.
For the sum associated with $\sog_{2j}$ let $U_+$ have dimension two and
$Z=Q^{-1}U_+^*+Q^{-1/2}W+U_+$, $\dim W=r-2$.
The maximal-$x$ eigenspaces of the hook representations, with their
original reduced conformal weights, sum to a common
factor times
$$
 [z^{r+2}]\frac{\Lambda_z(Z)}{\Lambda_z(U_+)}
 =[z^{r+2}]\Lambda_z(Q^{-1}U_+^*+Q^{-1/2}W)=0.
$$
The polynomial on the right has degree $r$. The square leaves the
unique one-dimensional subspace spanned by $e_\theta(-1)^{r+3}$,
of degree $r+3$ and charge
$2r+6$. Only the corresponding root vector of $\sog_{2j}$ has charge equal to twice its conformal weight.
No universal descendants alter this multiplicity. Nestedness implies that this integrability relation vanishes in
$H_f^0(L_k(\sog_{2n}))$, separately in both
$\slg_2$ factors when $\ell=4$.

The homomorphisms from the three universal affine vertex algebras
therefore factor through $L_{r+2}(\sog_{2j})$,
$L_0(\spg_{r-2})=\C$ and $L_0(\sog_4)=\C$.
Thus the classes of all weight-one vectors in the $C_2$ algebra of
the reduction are nilpotent. For $r\ge4$, the two-chain
multiplicity is $m=r-2>0$, and \cref{com:orthogonal-point} implies that
$H_f^0(L_k(\sog_{2n}))$ is $C_2$-cofinite.

For $r=2$, one invariant strong generator remains. The contributions of
the two largest $x$-eigenvalues are as follows. With $Z=U+B$,
$\dim B=4$, the five terms reduce to
$$
 Q\sog(B)+Q^2(\Sym^2U+2UB)+O(Q^3).
$$
The identity follows from $\Lambda^3B=B$ and
$\Sym^3U-U\Sym^2U+U\Lambda^2U-\Lambda^3U=0$.
There are three copies of the trivial $\sog_{2j}\oplus\sog_4$-module
at conformal weight two. Two are the quadratic Casimir vectors of
the two simple summands of $\sog_4$. Neither belongs to the second
reduced ideal: a positive mode of $\sog_4$ maps it to a nonzero
weight-one vector. The third invariant vector therefore belongs to
the first reduced ideal. It cannot belong to the vertex subalgebra
generated by $\sog_{2j}$, whose invariant form in conformal weight
two at level four is nondegenerate. Its coefficient on the remaining
invariant strong generator is consequently nonzero. Applying
\cref{com:orthogonal-point} also in this case gives
\begin{equation}\label{seed:D-actual}
 X_{H_f^0(L_{p-2n+2}(D_n))}=\{f\}.
\end{equation}

By the associated-variety theorem for reduction,
\eqref{seed:D-actual} is equivalent to
$X_{L_{p-2n+2}(D_n)}\cap\mathcal S_f=\{f\}$.
It remains to use this identity to show that no nonzero semisimple
point belongs to the associated variety. Starting from such a point,
we use the parabolic homomorphism to obtain a family in the associated
variety and degenerate it to $f$. The nonconstant quadratic invariant
of this family will contradict the one-point intersection with
$\mathcal S_f$.
Pfaffians leave only $h_b$, $1\le b\le n-p/2-1$. The parabolic homomorphism gives
$zh_b+y$ for every $y$ in the associated variety of the corresponding
simple affine vertex algebra of type $D$. The required degenerations
follow from the following three orthogonal families:
$$
 (\text{pair})+2^2\to3^2,\qquad
 (\text{two equal pairs})+1^2\to3^2,\qquad
 (\text{pair})+1^2\to2^2.
$$
The first two are the exterior-square images of
$\operatorname{diag}(z/2,z/2,-z/2,-z/2)+E_{12}+E_{23}$ and
$\operatorname{diag}(z,0,-z,0)+E_{12}+E_{23}$ in $\slg_4\simeq\sog_6$.
The third is $\left(\begin{smallmatrix}zD&B\\0&-zD\end{smallmatrix}\right)$,
$D=\operatorname{diag}(1,0)$, $B=\left(\begin{smallmatrix}0&1\\-1&0\end{smallmatrix}\right)$.
If $b<j$, use the nonzero reduction for the smaller orthogonal Lie algebra with $j'=j-b,r'=r+2b$ and the first
move $b$ times. The nonvanishing assertion includes $j'=1$. If $b\ge j$, set $u=b-j$,
$c=p-2u-2$, and use the orbit $(2^c,1^{2u+4})$.
With $v=\max(0,u-(r-2)/2)$ use the three moves $j-v,v,u-v$ times.
All counts are nonnegative and the limit is $f$. Along each curve
$\tr X^2/2=bz^2\ne0$. Smoothness of $G\times(f+\g^e)\to\g$
near $(1,f)$ places generic nearby conjugates in the slice, contradicting
the zero-dimensional Slodowy intersection in \eqref{seed:D-actual}.
This proves quasi-lisseness for $2n/3<p\le n-2$ with $p$ even.
No maximality assumption on the orbit of $f$ is needed.

\subsection{\texorpdfstring{Non-quasi-lisseness of $L_{p-2n+2}(\sog_{2n})$}{Non-quasi-lisseness of L(p-2n+2)(so(2n))}}\label{sec:D-exclusions}\label{D:negative-factors}\label{D:negative-assembly}
We consider $k=p-2n+2$ with $p\in\Z_{>0}$ outside the two ranges
in \cref{thm:D-integer-classification}: $p\ge2\lceil n/2\rceil$,
and even $p$ with $2n/3<p\le n-2$.
We prove that $L_k(\sog_{2n})$ is not quasi-lisse at these remaining levels.
Write $\mathsf A(n,t)=L_{t-2n+2}(\sog_{2n})$. At small ranks,
$\mathsf A(1,t)$ denotes the rank-one Heisenberg vertex algebra
associated with $\sog_2$ at nonzero level $t$,
$\mathsf A(2,1)=L_{-1}(\slg_2)^{\otimes2}$ and
$\mathsf A(3,t)=L_{t-4}(\slg_4)$.
Apart from $p=1$ and the affirmative ranges, write $n=mp+j$, $0\le j<p$.
For odd $p\ge3$ the rows are $j<p/2$, $j=(p+1)/2$ and
$j\ge(p+3)/2$. For even $p\ge4$, $n>p$. The $m=1$, $2\le j<p/2$
case was proved quasi-lisse above, whereas $j=1$ gives a level at
which quasi-lisseness fails. For $m\ge2$, quasi-lisseness fails
throughout $j<p/2$. The remaining cases are $j=p/2,p/2+1$ and
$j\ge p/2+2$. The following table gives the reductions used to prove
failure of quasi-lisseness, apart from the last odd case, which is
treated in \cref{seed:D-growth}.
\begin{longtable}{@{}>{\raggedright\arraybackslash}p{.26\textwidth}>{\raggedright\arraybackslash}p{.35\textwidth}>{\raggedright\arraybackslash}p{.33\textwidth}@{}}
\toprule
Parameters & Partition of $2n$ & Non-quasi-lisse quotient or tensor factor\\\midrule\endhead
$p=2$, $n=2s+2$ & $((2s+1)^2,1^2)$ & $\mathcal A_s$\\ \midrule
$p=2$, $n=2s+3$ & $(2s+3,2s+1,1^2)$ & tensor factor $\mathcal A_s$\\ \midrule
$n=mp+j$, $0<2j<p$ & $((2m+1)^{2j},1^{2m(p-2j)})$
 & $\mathsf A(j,p)\otimes\mathsf A(m(p-2j),p-2j)$\\ \midrule
$n=ma$, $p=a$, even $0<\ell<a$, $m\ge1$ for odd $a$, $m\ge2$ for even $a$ & $((2m)^\ell,1^{2m(a-\ell)})$
 & $\mathsf A(m(a-\ell),a-\ell)$\\ \midrule
$n=da+1$, $p=2a$, $a\ge2$, $d\ge3$ odd & $(d^{2a-2},1^{2d+2})$
 & $\mathsf A(d+1,2)$\\ \midrule
$p$ even, $j\ge p/2+2$ & $((2m+2)^v,(2m+1)^2,1^{2M})$,
$v=2j-p-2$, $b=p-j$, $M=(2m+1)b+1$ & $\mathsf A(M,2b)$\\ \midrule
$n=(2m+1)p/2$ & $((2m+1)^{p-1},2m-1,1^2)$
 & $\mathcal A_{m-1}$ if $m\ge2$, unitary image if $m=1$\\ \midrule
$p$ odd, $j=(p+1)/2$ & $((2m+1)^p,1)$ & simple Miura image in $\mathcal H^{\otimes m}\otimes F(p+1)_{\rm even}$\\
\bottomrule
\end{longtable}
Repeated parts equal to one are combined in these partitions.
The conformal-weight differences at block endpoints are:
$$
\begin{array}{c|c}
\text{row}&\text{positive-block endpoints}\\\hline
\text{odd-rank numerator two}&A_{1+2g}=4g+1\\
\text{divisible}&A_{g\ell}=\ell(a-\ell)g^2\\
n=da+1,\ h=2a-2&A_{gh}=2hg(g+1)\\
n=(2m+1)p/2&A_{p-1+gp}=(2g+1)p-1\\
p\text{ odd},\ j=(p+1)/2&A_{gp}=gp
\end{array}
$$
Here $g$ runs over the block endpoints of the displayed partition.
For the rectangular partition at $p=2$, the endpoint values follow
from the formulas of \cref{seed:scalar}.
For the row with even $p$ and $j\ge p/2+2$, its endpoint values equal
those of the fifth type-$C$ row above plus twice the index. They are strictly positive.
The last odd family has $A_l=l$ at every positive-block endpoint.
In the row $n=mp+j$ with $0<2j<p$, one has
$A_{2jg}=2jg(p-2j)(g-1)$, so the only equality is $(1^{2j})$ at degree $2jm$. Since
$$
 j(2m-1)<2jm<mp=n-j,
$$
\cref{seed:covariants} proves nonvanishing of the reduction there too.

The simple quotients for $0<2j<p$ and for $p\mid n$ are the unconditional entries of
\cite[arXiv v3, Propositions~4.10--4.11, Table~1]{AdamovicEtAl2026}, with parameter values
$(2n-2j)/(2m)=p$ and $(h^\vee+2)/(2m)=a$.
For odd $a$ take $\ell=a-1$. For even $a>2$ take $\ell=a-2$.
They reach numerator one or even-rank numerator two. Rank three is
type $A$, and when $m=j=1$, the rank-one vertex algebra $\mathsf A(1,p)$ occurs as a tensor factor.
In particular the divisible row includes $n=p=a\ge5$ odd.
Here $f=(2^{a-1},1^2)$, and the endpoint values
$A_{a-1}=a-1$, $A_a=a$ give the strict conformal-weight inequality.
In the cited even-block formula, take block size $2$, multiplicity
$a-1$ and multiplicity $2$ for the one-dimensional blocks: this is the unconditional $k'_1$
column, including block size two and two one-dimensional blocks. At $k=2-a$
the induced levels are $k+a-2=0$ on $\spg_{a-1}$
and $k+a-1=1$ on $\sog_2$. Thus its simple quotient is the
rank-one Heisenberg vertex algebra associated with $\sog_2$
at level one. This also directly checks the first
such case, $D_5$, $p=5$.
For $n=da+1$, the image of $\sog_{2a-2}$ at level zero
vanishes. The Lie algebra $\sog_{2a-2}$ is one-dimensional when $a=2$. Every nontrivial
$\sog_{2a-2}$-module of strong generators therefore vanishes.
The tensor product decomposition theorem yields the displayed
vertex subalgebra $\mathsf A(d+1,2)$. For even $p$ and $j\ge p/2+2$, the image of $\sog_2$
belongs to the radical of the invariant bilinear form. The mixed
strong generators for adjacent blocks have weights $3/2,m+1$,
and their zeroth product has weight $m+3/2$. The only generator of
smaller weight that is odd under the reflection of the
one-dimensional blocks has weight $m+1$. Hence no derivative or
normal product of smaller-weight generators can contribute.
The resulting orthogonal vertex subalgebra has parameter $t=2$
if $b=1$. Otherwise it belongs to the family $n=da+1$.

The odd-rank numerator-two row and the row $n=(2m+1)p/2$ with $m\ge2$ use
\cref{seed:telescope}, \cref{seed:scalar-split}. The vertex algebras $U$
in these Miura realizations are, respectively, the rank-one Heisenberg vertex algebra associated with the outer
$\glg_1$ block, and $L_1(\slg_{p-1})\otimes\mathcal H$, where
$\mathcal H$ is associated with the center of the outer $\glg_{p-1}$ block. The proposition
gives the stated $\mathcal A_s$ tensor factor in the simple quotient
of the Miura image.
In the latter case, let $a_0$ generate $\mathcal H$, normalized by
$[a_0{}_{\lambda}a_0]=\frac{p}{p-1}\lambda\mathbf1$. Its conformal
vector is $\frac{p-1}{2p}:a_0^2:$, with no derivative term.
For $n=(2m+1)p/2$ with $m=1$, the exact image is generated in weights
one and two in $L_1(\slg_{p-1})\otimes\mathcal H$.
The induced level of $\sog_{p-1}$ is two. The only derivative
from the complement of $\sog_{p-1}$ in the weight-one space
is that of the symmetric traceless primary vector. The condition $L_1B_2\subset B_1$ rules out
both this derivative and $\partial a_0$. The remaining real quadratics
are stable under $M\mapsto-M^t$, so the image is simple. Its conformal vector
maps to $(p-1)a_0^2/(2p)$ under the map
$R_{L_1(\slg_{p-1})}\otimes R_{\mathcal H}\to R_{\mathcal H}$
that sets the positive-degree elements of $R_{L_1(\slg_{p-1})}$ to zero.

The image for $j=(p+1)/2$ lies in $\mathcal H^m\otimes F(p+1)_{\rm even}$
with its standard conformal vector. For $m=1$, rational branching over
$U=F(p)_{\rm even}$ leaves in weight two exactly
$U_2+\C(T_{\mathcal H}+T_{\psi_0})+\langle h\psi_i\psi_0\rangle$,
including $p=3,5$. These PCT-stable spaces generate weight three.
Write $B_m$ for this image and $B_0=F(p+1)_{\rm even}$.
Deleting an outer pair of grading columns means parabolic induction
for $\glg_p\oplus\sog_{(2m-1)p+1}$, followed by
$V^0(\slg_p)\twoheadrightarrow\C$ \cite[Theorem~5.10]{Genra}:
$$
 B_m\lhook\joinrel\longrightarrow\mathcal H\otimes B_{m-1}
 \lhook\joinrel\longrightarrow\mathcal H^{\otimes m}\otimes B_0.
$$
The remaining orthogonal level is $1-2(m-1)p$. No Heisenberg
specialization is made. For adjacent equal columns the identifications
act by $\sigma_i(h_i)=h_{i+1}$, $\sigma_i(h_{i+1})=h_i$ and fix
the other generators. The final column reflection acts by
$\tau(h_m)=-h_m$, $\tau(\psi_0)=-\psi_0$, fixing the other generators.
Together with the one-chain reflection these maps give the sign changes
in every column. The paired chain forms differ by $(-1)^m$ in the
middle block, and the same fermion reflection identifies their conjugations.
Thus the one-column PCT calculation applies to each parabolic map.
The real Miura coefficients and these identifications make $B_m$
PCT-stable, hence simple, with conformal class $\frac12\sum h_i^2$.

\subsection{\texorpdfstring{Non-quasi-lisseness of $L_{p-2n+2}(\sog_{2n})$ by character growth}{Non-quasi-lisseness of L(p-2n+2)(so(2n)) by character growth}}\label{seed:D-growth}
\begin{proposition}\label{seed:D-growth-result}
The vertex algebra $L_{p-2n+2}(\sog_{2n})$ is not quasi-lisse for odd
$p\ge5$, $n=mp+j$, $m\ge1$, and $(p+3)/2\le j<p$.
\end{proposition}

For the remainder of this subsection, fix the parameters in the
proposition and set
$$
 \begin{gathered}
 v=2j-p-1\ge2,\quad w=2(p-j)+1\ge3,\quad
 b=(w-1)/2,\quad p=v+w,\\
 f=((2m+2)^v,(2m+1)^w,1).
 \end{gathered}
$$
Every positive-block endpoint has $A_l=l>0$. The Drinfeld--Sokolov reduction is therefore nonzero. Let $S$ be its simple quotient.
Let $V=\mathbb C^v$ be symplectic and let $W=\mathbb C^w$ be
orthogonal. Write $F(w)$ for the vertex superalgebra of $w$ free fermions.
The main construction will separate two known factors from an unknown
commutant:
$$
 \begin{gathered}
 S=(B\otimes F(w))_{\rm even},\qquad
 B=A\otimes C,\qquad A=L_b(\mathfrak{osp}(1|v)),\\
 T_B=T_S-T_{F(w)},\qquad T_C=T_B-T_A.
 \end{gathered}
$$
We first separate the free-fermion factor and then the rational factor $A$. Quasi-lisseness of
the original affine vertex algebra would make $C_{\rm even}$ lisse.
Its strong generators and central charge then give incompatible
character-growth bounds.

The following results isolate the four parts of the argument: the
current algebras and generators of $S$, the tensor decomposition, the
properties of $C$, and the general growth estimate. We return to the
proof of the proposition after establishing them.

\begin{lemma}\label{seed:D-centralizer-generators}
The simple quotient $S$ is generated as a vertex algebra
by its weight-one subspace, the fields $G$ in $V\otimes W$ of weight $3/2$,
$\mathfrak{sp}(V)\oplus\mathfrak{so}(W)$-invariant fields $Q_{2r}$ with linear terms $f^{2r-1}$, $1\leq r\leq m$,
and a field space $P$ of type $W$ and weight $m+1$. More precisely,
the image in $S$ of every universal strong generator of weight less than $a$ belongs to
the subalgebra generated by these fields of weight less than $a$.
\end{lemma}
\begin{proof}
Set $d=2m+1$. For $r\geq1$, write
$R_r=\mathbb C[z]/(z^r)$ and equip it with the form
$$
 (z^i,z^j)_r=(-1)^i\delta_{i+j,r-1}.
$$
Multiplication by $z$ is skew-adjoint. The defining orthogonal
representation is
$$
 (R_{d+1}\otimes V)\oplus(R_d\otimes W)\oplus R_1.
$$
The form on the first summand is the product of two alternating forms.
Consequently the $\mathfrak{sp}(V)\oplus\mathfrak{so}(W)$ modules of
strong generators of the universal $W$-algebra are as follows:
$$
\begin{array}{c|c}
\text{module}&\text{conformal weights}\\\hline
\operatorname{Sym}^2V&1,3,\ldots,2m+1\\
\Lambda^2V&2,4,\ldots,2m+2\\
\Lambda^2W&1,3,\ldots,2m+1\\
\operatorname{Sym}^2W&2,4,\ldots,2m\\
V\otimes W&\frac32,\frac52,\ldots,2m+\frac32\\
W&m+1\\
V&m+\frac32.
\end{array}
$$
Every entry has multiplicity one. This follows either by decomposing
the exterior square of the defining representation under the
$\mathfrak{sl}_2$-triple, or directly by taking the skew-adjoint
endomorphisms commuting with $z$.

For $A\in\operatorname{Hom}(W,V)$, define $A^\dagger$ by
$(Aw,v)_V=(w,A^\dagger v)_W$. The mixed elements of the centralizer are
$$
 K_s(A)=
 \begin{pmatrix}
 0&z^sA&0\\
 (-1)^{s+1}z^{s-1}A^\dagger&0&0\\
 0&0&0
 \end{pmatrix},\qquad 1\leq s\leq d.
$$
They correspond to generators of weight $s+\frac12$. Truncation in
$R_{d+1}$ and $R_d$ is understood. Matrix multiplication gives
\begin{equation}\label{seed:D-block-bracket}
 [K_1(A),K_s(B)]
 =z^s\begin{pmatrix}
 (-1)^{s+1}AB^\dagger-BA^\dagger&0&0\\
 0&A^\dagger B-(-1)^{s+1}B^\dagger A&0\\
 0&0&0
 \end{pmatrix}.
\end{equation}
For even $2\le s\le2m$, these brackets span
$\mathfrak{sp}(V)\oplus\mathfrak{so}(W)$.
For odd $1\le s\le2m-1$, they span the self-adjoint coefficient
pairs $(C,D)$ satisfying $\operatorname{tr}C+\operatorname{tr}D=0$.
Indeed, rank-one maps $A=u\otimes(w,\cdot)_W$ and
$B=u'\otimes(w',\cdot)_W$ give the two traceless subspaces
separately and a pair with nonzero opposite traces. The element
$f^s$ has coefficient pair $(1_V,1_W)$, whose trace sum is
$v+w=p\ne0$.
At $s=d=2m+1$, multiplication by $z^d$ vanishes on $R_d$,
so the $W$ block is zero. The remaining brackets span the entire
last $\Lambda^2V$ coefficient space, including its trace.
For example, when $w'=w$ and $(w,w)_W=1$, the coefficient on $V$ is
$$
 C(v)=u(u',v)_V-u'(u,v)_V,
 \qquad \operatorname{tr}C=-2(u,u')_V.
$$
It has nonzero trace whenever $(u,u')_V\ne0$. Such tensors,
together with those for $(u,u')_V=0$, span $\Lambda^2V$.
This argument also applies when $v=2$, where the traceless
self-adjoint subspace is zero.

The non-scalar coefficient spaces in $z^s\operatorname{End}(W)$ act
nontrivially on $W$. Their brackets with $K_1$ span $K_{s+1}$ whenever
$s+1\leq d$. For odd $s$ one may use the traceless self-adjoint
coefficients, and for even $s$ the skew-adjoint coefficients. This
holds also for $w=3$. The map from $R_1$ to $R_d\otimes W$ is
$1\mapsto z^{d-1}w$. Its skew-adjoint completion, denoted $P(w)$,
satisfies that $[K_1(A),P(w)]$ spans the maps from $R_1$ to
$R_{d+1}\otimes V$. Thus the centralizer is generated by its degree-zero
part, the $K_1(A)$, the elements $f,f^3,\ldots,f^{2m-1}$, and the
$P(w)$. These statements include $v=2$. In that case the traceless
part of the self-adjoint $V$ coefficient space is simply absent.

These brackets prove the claimed generation statement by induction from \eqref{seed:D-block-bracket}: a zeroth product, by \cref{com:linear-bracket},
has the displayed centralizer bracket as its linear term, and its
remaining terms use derivatives and products of fields of smaller
weight.

\end{proof}

\begin{lemma}\label{seed:D-current-pairing}
In $S$, the symplectic currents satisfy the integrability relation
$e_\theta(-1)^{b+1}\mathbf1=0$, and the orthogonal currents generate
$L_1(\sog_w)$, interpreted as $L_2(\slg_2)$ when $w=3$.
The scalar pairing of the mixed generators of conformal weight $3/2$
is nondegenerate.
\end{lemma}
\begin{proof}
The weight-one subspace defines homomorphisms
$V^b(\spg_v)\to S$ and $V^1(\sog_w)\to S$.
The following calculations of the radical take place at integral conformal weights,
using the full conformal-weight spaces and the form of \cref{seed:integral-form}.
For the first symplectic coordinate, charge is at most twice weight,
with equality only for a long-root vector of $\spg_v$. Thus the type
$(2b+2)\omega_1$ at weight $b+1$ is exactly the affine integrability
line. Its zero norm and multiplicity one in the \emph{full} weight
space make it radical. For $w\ge5$ the orthogonal null type $2\theta$
at weight two likewise occurs only among quadratic products of $\sog_w$, not in
$\Lambda^2V$, $\Sym^2W$ or the extra $W$ at $m=1$.
For $w=3$, \cref{seed:D-windows}(4) has maximal-$x$ eigenspace
$$
 \Sym^4V+S_{(2,2)}W-1=\Sym^4V+\Sym^2_0W.
$$
Its multiplicity-one spin-two relation belongs to the first reduced ideal.
The invariant form on the spin-two subspace of $V^2(\slg_2)$ is
nondegenerate, so the additional spin-two generator belongs to
the vertex subalgebra generated by $\slg_2$. At weight three,
the $\operatorname{Sp}(V)$-invariant contraction of two mixed generators has
$\Lambda^2W$ type, not spin three. The spin-three null line is therefore
unique in the full space and radical. This gives
$U=L_1(\sog_w)$, read as $L_2(\slg_2)$ for $w=3$.

Let $G^u$, $u\in\g^f_{-1/2}\simeq V\otimes W$, be the weight-$3/2$
strong generators corresponding to $\g^f_{-1/2}$. With the partition fixed,
the two induced affine levels as functions of $t=k+2n-2$ are
$$
 a(t)=(m+1)t-n,\qquad d(t)=(2m+1)(t-p)+1.
$$
Write $G^u_{(2)}G^{u'}=P(t;u,u')\mathbf1$.
The formula of \cite[Corollary~3.1]{KacWakimoto2022} expresses the image
of $G^u$ in the free-field realization as a quadratic expression
in the current fields and neutral free fields, a cubic expression
in the neutral free fields,
and
$$
 -\sum_{i\in S_{1/2}}
 \bigl(k(u\mid u_i)+\kappa_{>0}(u,u_i)\bigr)\partial\Phi^i.
$$
The neutral pairing is $\Omega(q,q')=(f\mid[q,q'])$.
All scalar contractions have degree at most two in $k$, and hence in $t$. The unique contribution of degree two is the contraction
of the two displayed derivative terms. Identifying their neutral
labels by $q_u=(\operatorname{ad}f)^{-1}u$, it is a fixed nonzero
scalar multiple of $k^2\Omega(q_u,q_{u'})$.
Indeed $\g^f_{-1/2}$ consists of the lowest vectors of the
spin-$1/2$ summands for the defining $\mathfrak{sl}_2$-triple.
The invariant form and $\operatorname{ad}f$ give a nondegenerate
pairing on their highest vectors. Thus the leading coefficient
is nonzero and nondegenerate on $V\otimes W$.

The zeros of the two induced affine levels are zeros of this scalar pairing
already in the universal $\mathcal W$-algebra. At a zero of one level,
let $J$ belong to the corresponding simple Lie algebra in the weight-one space. There are no states of
weight $1/2$, so $J_{(1)}G^u=J_{(1)}G^{u'}=0$.
Since $G^u_{(1)}G^{u'}$ has weight one and the chosen factor pairs
trivially with the entire weight-one space, the mode commutator gives
$$
 0=J_{(1)}(G^u_{(1)}G^{u'})
   =(J_{(0)}G^u)_{(2)}G^{u'}.
$$
The chosen simple factor acts nontrivially on $V\otimes W$, and
its action spans that space. Hence $P$ vanishes at the zero of
either $a$ or $d$, without passing to a simple quotient.
Both zeros are noncritical: they are $n/(m+1)>0$ and
$p-1/(2m+1)>0$. They are distinct because
$b(2m+1)>m+1$, with $m\ge1$ and $b\ge1$.
Invariance and irreducibility now give
$$
 P(t;u,u')=c\,a(t)d(t)\,\beta(u,u'),\qquad c\ne0,
$$
where $\beta$ is a fixed nondegenerate invariant pairing on
$V\otimes W$. At $t=p$ the levels are $b$ and $1$, so this
pairing remains nonzero in every nonzero quotient, including $S$.

\end{proof}

\begin{proposition}\label{seed:D-tensor-realization}
There are simple vertex superalgebras $B$ and $C$ with
$A=L_b(\mathfrak{osp}(1|v))$ rational and lisse, such that
\begin{equation}\label{seed:D-tensor}
 B=A\otimes C,\qquad
 S=(A\otimes C\otimes F(w))_{\overline0},\qquad
 C_0=\mathbb C,\quad C_{1/2}=C_1=0.
\end{equation}
The conformal vector of $C$ is $T_B-T_A$. The central charges are
$$
 c(S)=m+(w+1)/2+bv(v-1)/p,\quad c(A)=bv(v-1)/p,
 \quad c(F(w))=w/2,\quad c(C)=m+\tfrac12.
$$

\end{proposition}
\begin{proof}
Let $U=L_1(\mathfrak{so}_w)$, with
$U=L_2(\mathfrak{sl}_2)$ when $w=3$, and let $M$ be its vector module.
The integrability of the $\widehat{\sog}_w$-action proved above makes each half-integral grading
coset of $S$ an admissible $U$-module. Rationality of $U$ therefore
implies that its restriction is a direct sum of ordinary modules.
Only $U$ and $M$ occur: every finite orthogonal type in $S$ is a tensor
type, whereas a spinor module has only spinorial finite weights.
There are therefore multiplicity spaces $B_{\overline0},B_{\overline1}$
such that
$$
 S=(U\otimes B_{\overline0})\oplus(M\otimes B_{\overline1}).
$$
Write $F_{\overline0}=U$ and $F_{\overline1}=M$. Their fusion rules are
$$
 F_{\bar\epsilon}\boxtimes_U F_{\bar\delta}
   =F_{\overline{\epsilon+\delta}},
 \qquad \epsilon,\delta\in\{0,1\},
$$
with one-dimensional spaces of intertwining operators. These are
the level-one orthogonal fusion rules. For $w=3$ they are the
level-two $\slg_2$ rule
$L_2(2\omega_1)\boxtimes L_2(2\omega_1)=L_2(0)$.
Fix the intertwining operators $I_{\epsilon,\delta}$ obtained by
restricting the vertex operator of $F(w)$, with
$\psi(u)_{(0)}\psi(v)=(u,v)_W\mathbf1$.
Their powers are integral, their associativity scalar is one,
and their skew-symmetry scalar is $(-1)^{\epsilon\delta}$.
These identities are restrictions of the corresponding identities
in $F(w)$.

Identify the summand $F_{\bar\epsilon}\otimes B_{\bar\epsilon}$
with $B_{\bar\epsilon}\otimes F_{\bar\epsilon}$ as a vector space.
The one-dimensional intertwining-operator spaces uniquely define
operators $Y_B$ by the factorization
$$
 Y_S(b\otimes u,z)(b'\otimes u')
 =(-1)^{\epsilon\delta}
    Y_B(b,z)b'\otimes I_{\epsilon,\delta}(u,z)u',
$$
where $b,u$ have parity $\epsilon$ and $b',u'$ have parity $\delta$.
The displayed sign is the tensor-product sign. It satisfies
$$
 (-1)^{\epsilon\delta+(\epsilon+\delta)\eta}
 =(-1)^{\delta\eta+\epsilon(\delta+\eta)}.
$$
Consequently weak associativity for $S$, after factorization through
the fixed operators $I_{\epsilon,\delta}$, gives weak associativity
for $Y_B$. Locality gives superlocality with sign
$(-1)^{\epsilon\delta}$. The unit components give the vacuum and
creation properties, and the translation operator on $B$ is
induced by $L^S_{-1}-L^U_{-1}$. Hence
$B=B_{\overline0}\oplus B_{\overline1}$ is a vertex superalgebra and
$$
 S\cong(B\otimes F(w))_{\overline0}.
$$
The half-integral grading on $B$ is the difference of the conformal
gradings of $S$ and $U$. Simplicity of $S$ implies simplicity of $B$:
a nonzero proper graded superideal would give a nonzero proper ideal
in the right side of the displayed tensor realization. The low-weight spaces in the table give
$$
 B_0=\mathbb C,\qquad B_{1/2}=0,\qquad
 B_1=\mathfrak{sp}(V)\oplus\Pi V.
$$
Here the weight-one copy of $\sog_w$ belongs to the $U$-summand,
including when $w=3$. It contributes no affine highest-weight vector for $M$
at weight one.
The nonzero pairing of the mixed generators and the $\mathfrak{sp}(V)$-action identify $B_1$ with
$\mathfrak{osp}(1|v)$ at level $b$. Its even highest-root relation is
the relation $e_\theta(-1)^{b+1}\mathbf1=0$ already proved in $S$.
Since this generates the maximal ideal of
$V^b(\mathfrak{osp}(1|v))$ (see \cite[Section~4]{JiangWang2022}), it gives
$$
 A=L_b(\mathfrak{osp}(1|v))\subset B.
$$
Since $b\ge1$ and $v/2\ge1$ are integers,
\cite[Theorem~7.1]{CL} implies that $A$ is rational and lisse.

We choose the generators $Q_{2r}$ and $P$ as follows.
Let $L_0^U$ be the Sugawara operator for $U$. On a $\mathfrak{so}(W)$-invariant vector $Q_a$, $L_0^UQ_a$ is a sum of states $u(-n)u'(n)Q_a$, $n>0$.
On the copy $P$ of the vector representation of $\mathfrak{so}(W)$, the same assertion holds for
$(L_0^U-\frac12)P$. In both cases $u'(n)$ lowers total conformal
weight, so these states belong to the algebra already generated in
smaller weights. On each finite-dimensional conformal-weight space
$L_0^U$ is diagonalizable. Projecting to eigenvalue zero for $Q_a$,
and to eigenvalue $1/2$ for $P$, subtracts only $U$ descendants of
smaller-weight states. Up to derivatives, the resulting fields have
the same linear terms corresponding to elements of $\g^f$, and they
have the form
$$
 Q_a=q_a\otimes\mathbf1,
 \qquad P(w)=y\otimes\psi(w),
$$
where $q_a\in B_{\overline0}$ has weight $a$ and
$y\in B_{\overline1}$ has weight $m+1/2$.
This also shows that all multiplicity coefficients of the original
generating fields are generated by $A$, the $q_a$ and $y$.
The argument treats possible finite-type coincidences at $w=3$
while distinguishing the two irreducible $U$-modules $U$ and $M$.

The generation assertion respects $B$-weight. To see this, an even
state of $B$ of weight $c$ is recovered from its lift in $\mathbf1_U\otimes B$ of weight $c$, and an odd state from its lowest-vector lift of weight
$c+1/2$. Expand these lifts using the full universal strong generators.
For $c<a$, all source weights are less than $a$, except possibly the
odd case $c=a-1/2$. At source weight $a$, the only new invariant generator
could be $Q_a=q_a\otimes\mathbf1$. Its vector coefficient is zero.
A generator $P$ of source weight $a$ already has $B$-weight $a-1/2$.
All other source generators are expressions in smaller-weight fields
and these generators, by \eqref{seed:D-block-bracket}. Thus every state of $B$ of weight less
than $a$ belongs to the algebra generated by $A$ and the coefficients
of $B$-weight less than $a$. This statement also applies to the
coefficients of derivatives and products in the induction above.

We next construct commuting lifts in increasing $B$-weight. Suppose
the earlier lifts and $A$ generate $D\subset B$. The vector $q_a$
is $\mathfrak{sp}(V)$-invariant. Write $G=X\otimes\psi$, where
$X\in\Pi V\subset A_1$. Since $[K_1,f^{a-1}]=0$, the zeroth
product $G_{(0)}Q_a$ has no new linear term. Its product and derivative
terms involve universal strong generators of weight at most
$a-\frac12$. Their multiplicity coefficients are already in $D$.
The coefficient of the lowest fermion in this zeroth product is
$X_{(0)}q_a$, hence $X_{(0)}q_a\in D$.

For $y$, project $G_{(1)}P$ onto the $U$-summand and then onto
$\C\mathbf1_U$. Its coefficient is a nonzero scalar multiple of
$X_{(0)}y$, of weight $m+1/2$. The first strong
generator of type $V\otimes\mathbf1$ has weight $m+3/2$, by the
table. Thus there is no linear term of this type at weight $m+1/2$.
The remaining terms involve smaller-weight generators, and their
coefficients belong to $D$. Therefore $X_{(0)}y\in D$.
Positive modes of every element of $A_1$ lower $B$-weight and likewise
have image in $D$, and the even zero modes annihilate $y$ because
it is $\mathfrak{sp}(V)$-invariant. Each new class in $B/D$ is
therefore annihilated by
$\mathfrak{osp}(1|v)\otimes\C[z]$. If nonzero, it generates a copy
of $A$ for an even class or of $\Pi A$ for an odd class, with the
corresponding translation of the degree grading.
The $B$-grading is bounded below and every element of $A_1$ has
weight one. Separating its two half-integral grading cosets and
translating their lowest degrees gives admissible untwisted
$A$-modules. Rationality of $A$ splits the preimage of the new
simple submodule. The splitting may be chosen homogeneous and
parity preserving: its degree-zero, parity-even component is still
an $A$-module map and a section. The image of this section contains a lift annihilated
by all nonnegative modes of $A_1$.

Induction shows that $B$ is generated by $A$ and commuting lifts.
For $C=\operatorname{Com}(A,B)$, multiplication therefore gives a
surjective homomorphism $A\otimes C\to B$. Density on the simple
adjoint $A$-module makes it injective: any kernel is of the form
$A\otimes I$, and its intersection with $\mathbf1_A\otimes C$ is zero.
Simplicity of $B$ also implies simplicity of $C$. The conformal vectors and central charges are obtained by subtracting
those of the two separated factors, which proves the displayed formulas.
\end{proof}

\begin{lemma}\label{seed:D-commutant-generators}
The commutant $C$ has a strong generating set with at most one even
generator at each weight $2,4,\ldots,2m$ and at most one odd generator
of weight $m+1/2$. The weight-two generator may be chosen to be $T_C$.
\end{lemma}
\begin{proof}
Set
$E=A\otimes F(w)$ and let $\pi_0,\pi_\psi$ denote coefficient
projections onto $\C\mathbf1_E$ and a fixed one-fermion vector, respectively.
They recover respectively the even and odd states of $C$. Set
$$
 \mathcal D(C)=\partial C+
 \operatorname{span}\{u_{(-1)}v:\operatorname{wt}u,
 \operatorname{wt}v>0\}.
$$
Both projections take $\mathcal D(S)$ into $\mathcal D(C)$.
Indeed, for a normal product of vectors with $E$-weights
$\alpha,\beta$, its component of $E$-weight $t=0$ or $1/2$ has
$C$-product index $-1-(\alpha+\beta-t)$.
It is at most $-1$. Equality for $t=0$ requires both $E$ components to belong to $\C\mathbf1_E$.
Equality for $t=1/2$ requires one $E$ component to be $\mathbf1_E$ and the other to be a one-fermion vector. Parity and
positive total weights of the two vectors then force both $C$ factors to have
positive weight. Strict inequality gives a derivative or a normal product of
positive-weight vectors. If one $C$ component is $\mathbf1_C$,
the creation property gives the same conclusion. Derivatives project
to derivatives, since $E$ has no vectors of weight $-1$ or $-1/2$.

Apply these projections to the full table of strong generators. Mixed
$V\otimes W$ fields and the last $V$ field have zero projections
by $\mathfrak{sp}(V)$-equivariance. Formula \eqref{seed:D-block-bracket}, lifted by the
identity of \cref{com:linear-bracket}, removes both internal skew-adjoint coefficient
spaces and the space of self-adjoint pairs $(C,D)$ satisfying
$\operatorname{tr}C+\operatorname{tr}D=0$, modulo $\mathcal D(C)$. Indeed $G\in E\otimes\mathbf1$, so
$G_{(0)}$ cannot have a coefficient in $\C\mathbf1_E$. To have a
single-fermion coefficient it would have to act on $\mathbf1_E$,
where its zeroth mode vanishes. At $s=d$, \eqref{seed:D-block-bracket} removes the whole
last $\Lambda^2V$ space. This also removes the possible
$\Lambda^2W\cong W$ projections when $w=3$.
Only one even class at each weight $2,4,\ldots,2m$ and the odd
class contributed by $P$ at weight $m+1/2$ remain. These classes
span $C_+/\mathcal D(C)$, so their lifts strongly generate $C$ by
weight induction. Since $c(C)\ne0$, its weight-two generator can be
taken to be its conformal vector.
\end{proof}

\begin{example}\label{ex:D9-growth}
Take $m=1$, $p=5$ and $j=4$. Then $n=9$, $k=-11$, and
$$
 v=2,\qquad w=3,\qquad b=1,\qquad f=(4^2,3^3,1).
$$
The preceding results give
$$
 S=\bigl(L_1(\mathfrak{osp}(1|2))\otimes C\otimes F(3)\bigr)_{\overline0},
 \qquad c(C)=\frac32.
$$
Here $C$ is strongly generated by its conformal vector of weight two
and at most one odd field of weight $3/2$. The contradiction below
already appears in this case. If the affine vertex algebra were
quasi-lisse, $C_{\mathrm{even}}$ would be lisse. This would force a
relation on the conformal vector, lowering the growth bound for $C$
strictly below $3/2$, whereas modularity of $C_{\mathrm{even}}$
requires a bound of at least $3/2$.
\end{example}

\begin{lemma}\label{seed:D-commutant-lisse}
If $L_{p-2n+2}(\sog_{2n})$ is quasi-lisse, then $C_{\mathrm{even}}$
is lisse.
\end{lemma}
\begin{proof}
The images of $\spg_v\oplus\sog_w$ are nilpotent in $R_S$. Its nilradical
is a Poisson ideal, so their zero modes kill every nontrivial finite
$\operatorname{Sp}(V)\times\operatorname{SO}(W)$-module of strong
generators in the reduced $C_2$ algebra. Only invariant generators
of even conformal weight remain. Their Poisson brackets have odd conformal weight and are therefore
zero. The associated variety of $S$ has finitely many symplectic
leaves as a Poisson subvariety of the reduction of the quasi-lisse
affine vertex algebra. Since its Poisson structure is zero, it is
a finite set. Its nonnegative grading and one-dimensional degree-zero
part then imply that it is the origin.
Finite generation then makes $S$ $C_2$-cofinite.
Replace the conformal vector by its sum with $\partial J^h$, where $h$ has eigenvalues
$\pm1/2$ on $V$. Under $L_0-h_0$ the generators from the even part of $\mathfrak{osp}(1|v)$ have
weights $0,1,2$ and the odd ones $1/2,3/2$. The even subalgebra
$S=(A\otimes C\otimes F(w))_{\overline0}$, where parity is the sum
of the parities of the three factors, is nonnegatively integer graded. Its only zero-weight creation
modes are the finitely many commuting $e_\alpha(-1)$ with
$\alpha(h)=1$. No sum of two such roots is a root. Integrability
makes each locally nilpotent. Moving them to one side of each
product of positive-weight generators gives a finite joint span, so graded pieces
remain finite although degree zero may enlarge.
Miyamoto's solvable-orbifold theorem and the non-CFT-type version
noted in \cite[Theorem~1 and Remark~1]{Miyamoto2015}
apply to the two parity involutions. Their common fixed algebra is
$A_{\mathrm{even}}\otimes C_{\mathrm{even}}\otimes F(w)_{\mathrm{even}}$.
The tensor $C_2$ identity makes $C_{\mathrm{even}}$ lisse in its
original grading of CFT type.

\end{proof}

For a nonnegatively graded vector space $V$ with finite-dimensional
graded pieces, suppose that its Hilbert series
$H_V(u)=\sum_d\dim V_d\,e^{-ud}$ converges for $u>0$. Its exponential
growth is measured by
$$
 g(V)=\limsup_{u\to0^+}\frac{6u}{\pi^2}\log H_V(u).
$$
This definition also applies to half-integrally graded vertex
superalgebras, using ordinary dimensions of both parity components.
Finite strong generation ensures convergence for the algebras used here.

\begin{lemma}\label{seed:modular-growth-lower-bound}
Let $V$ be a lisse vertex operator algebra of CFT type with real central
charge. Then $g(V)\ge c(V)$.
\end{lemma}
\begin{proof}
The argument does not assume rationality. It combines
module spanning with the pseudotrace description of genus-one
conformal blocks. Let $P$ be a projective generator of the finite category of
generalized $V$-modules. Module spanning (\cite{Buhl2002}, Theorem~1)
gives finite-dimensional $P/C_2(P)$. Choose homogeneous representatives
of a basis, with generalized conformal weights $d_i$. The Li module
lemma (\cite{Arakawa2012}, Lemma~3.1.3) says that their symbols generate
$\operatorname{gr}P$ over $\operatorname{gr}V$. They therefore give a
surjection from a finite direct sum of graded translates of
$\operatorname{gr}V$ to $\operatorname{gr}P$.
For a generalized module, define $H_P(u)$ by using the real parts of
its generalized conformal weights. Comparing dimensions gives
$H_P(u)\le\sum_i e^{-u\operatorname{Re}d_i}H_V(u)$.
Each generalized weight space of $P$ is projective over
$\operatorname{End}(P)^{\mathrm{op}}$ by \cite[Section~3.2.1]{GuiZhang}. The nilpotent part of $L_0$
is a fixed central nilpotent endomorphism, with one nilpotence bound
on all of $P$. There are finitely many isomorphism classes of
indecomposable projectives over this finite-dimensional endomorphism
algebra. On each class, every power of that central endomorphism has
a fixed pseudotrace. Additivity under direct sums therefore bounds
each pseudotrace by a constant times the ordinary weight-space
dimension, with a constant independent of conformal degree. The pseudotrace theorem of Gui--Zhang \cite[Theorem~3.7]{GuiZhang} and modular invariance \cite[Theorem~5.2]{Miyamoto2004} therefore give  for $y\ge1$, constants $K>0$ and $r\ge0$ such that
$$
 \chi_V(iy)\le K(1+y)^r e^{\pi c(V)/(12y)}H_V(2\pi/y).
$$
Here $\chi_V(\tau)=\operatorname{tr}_V
 e^{2\pi i\tau(L_0-c(V)/24)}$. Since $V$ is of CFT type,
$$
 \log\chi_V(iy)=\frac{\pi c(V)y}{12}+o(1)
 \qquad(y\to+\infty).
$$
Taking logarithms and setting $u=2\pi/y$ gives $g(V)\ge c(V)$.
The cited pseudotrace theorem applies to an $\mathbb N$-graded $C_2$-cofinite vertex operator algebra without rationality or self-duality assumptions.

\end{proof}

\begin{proof}[Proof of \cref{seed:D-growth-result}]
Suppose that the affine vertex algebra were quasi-lisse. By
\cref{seed:D-commutant-lisse}, $C_{\mathrm{even}}$ would be lisse.
The lower bound of \cref{seed:modular-growth-lower-bound} therefore
applies. We now derive a strictly smaller upper bound from
\cref{seed:D-commutant-generators}.
Since $C_{\mathrm{even}}$ is lisse, $[T_C]^N=0$ in $R_C$ for
some finite $N\ge1$. In the Li graded algebra set $x_i=\partial^iT_C$.
The relation $\partial^{Ni}(T_C^N)$ has coefficient
$(Ni)!/(i!)^N\ne0$ at $x_i^N$. Every other monomial has a strictly larger sum of squared derivative
indices, by the equality case of Cauchy--Schwarz. On each fixed
conformal-weight space this sum is bounded. Descending induction on
that sum therefore shows that monomials in which every $x_i$ has
exponent at most $N-1$ span.
The resulting spanning set gives the Hilbert-series bound
$$
 H_C(u)\le
 \prod_{i\ge0}\frac{1-e^{-Nu(i+2)}}{1-e^{-u(i+2)}}
 \prod_{r=2}^{m}\prod_{i\ge0}\frac1{1-e^{-u(i+2r)}}
 \prod_{i\ge0}\bigl(1+e^{-u(i+m+1/2)}\bigr).
$$
The second product is empty when $m=1$. The first product counts
derivatives of $T_C$ with multiplicity at most $N-1$, the second
counts the other even generators, and the third counts the odd
generator. Taking the logarithmic asymptotics as $u\to0^+$ gives
contributions at most $1-1/N$, $m-1$ and $1/2$, respectively.
Thus
$$
 m+\tfrac12=c(C_{\mathrm{even}})
 \le g(C_{\mathrm{even}})\le g(C)
 \le m+\tfrac12-\frac1N,
$$
a contradiction.
\end{proof}
This also completes the proof of \cref{thm:D-integer-classification}.

\section{\texorpdfstring{Quasi-lisse simple affine vertex algebras of types $G_2$ and $F_4$}{Quasi-lisse simple affine vertex algebras of types G2 and F4}}\label{seed:FG}

The admissible one-third levels in type $G_2$ were studied by Axtell and Lee \cite{AxtellLee2011}.
Per\v se's conformal embedding $L_{-5/2}(B_4)\subset L_{-5/2}(F_4)$ \cite{PerseBF2007} and the embeddings
$L_{-2}(G_2)\subset L_{-2}(B_3)\subset L_{-2}(D_4)$ of Adamovi\'c and Per\v se \cite{APConformal2013}
are earlier uses of conformal embeddings to study these nonadmissible levels.

Li, Xie and Yan propose $L_{-2}(G_2)$ and $L_{-3}(F_4)$ for
rank-one Argyres-Douglas theories from twisted $D_4$ and $E_6$ constructions, with
Higgs branches $\overline{\OO_{G_2(a_1)}}$ and
$\overline{\OO_{\widetilde A_1}}$, respectively
\cite[Table~4 and Appendix~B.1]{LiXieYan2023}.
The associated-variety statements are proved in
\cite[Theorem~A]{ArakawaEtAl2025} and
\cite[Theorems~4.2 and~4.5]{ADS2026}, respectively.
The pair $L_{-6}(F_4)$ and $\overline{\OO_{F_4(a_3)}}$ was predicted
in \cite[Table~19]{XieYanLisse2021} and also appears in
\cite[Table~17]{LiXieYan2023}.
The minimal-orbit algebras $L_{-5/3}(G_2)$ and $L_{-5/2}(F_4)$
are discussed separately in \cite[Section~4]{BeemRastelli2018},
where their four-dimensional interpretation was left open.
Boundary admissible examples with $k+4=4/q$ and $k+9=9/q$, respectively, arise
in \cite[Section~5]{WangXie2019}, for $q>1$ with
$\gcd(q,6)=1$ and $\gcd(q,18)=1$, respectively.

\begin{theorem}\label{thm:FG-classification}
Let $k\ne-h^\vee$ and write $k+h^\vee=p/q>0$ in lowest terms.
For $F_4$, $L_k(F_4)$ is quasi-lisse exactly when $p\ge5$, or when
$p=3$ and $q$ is odd. For $G_2$, it is quasi-lisse exactly when
$p\ge2$ if $3\nmid q$, and $p\ge4$ if $3\mid q$.
\end{theorem}
\begin{proof}
Theorem~\ref{thm:denominator} and admissibility reduce to the representatives
below. We prove that the reductions of the simple affine vertex algebras
are nonzero and have simple quotients that are not quasi-lisse at the
levels where quasi-lisseness fails. At the remaining levels, relations
in the reduced maximal ideals and the absence of nonzero semisimple
points show that the associated varieties are contained in the nilpotent cone.
\end{proof}

\subsection{Drinfeld--Sokolov reductions of \texorpdfstring{$L_k(F_4)$ and $L_k(G_2)$}{Lk(F4) and Lk(G2)}}
The nonadmissible values of $t=k+h^\vee$ at which we prove quasi-lisseness
are $3,5,6,7,8,5/2,7/2,9/2,11/2$ for $F_4$,
and $2,3,4/3,5/3$ for $G_2$. The values for which quasi-lisseness
fails are $1,2,4,1/2,3/2$ and $1,1/3,2/3$, respectively.
At $t=1/2$ for $F_4$ and $t=1/3$ for $G_2$, the universal affine
vertex algebra is simple.

For $F_4$, positivity of the fundamental-weight Gram matrix reduces
the strict conformal-weight inequality to the following values on fundamental weights, where
$Q(\lambda)=(\lambda,\lambda+2\rho)/2$:
$$
\begin{array}{c|c|c|c}
t&f&x&Q(\omega_i)/t-\omega_i(x)\\\hline
1&F_4(a_1)&(5,2,1,0)&2,5,3,1\\
2&F_4(a_3)&(2,1,1,0)&3/2,3,2,1\\
4&A_2&(1,1,0,0)&1/4,3/2,1,1/2\\
3/2&B_3&(3,2,1,0)&1,3,2,1
\end{array}
$$
For the $G_2$ convention in which the first simple root is short,
$Q(a,b)=(a^2+3ab+3b^2+5a+9b)/3$ and
$\lambda(x_{\rm sub})=a+2b$. These values are $(1,2)$ and $(2,4)$
at $t=1,2/3$. Thus $H_f^0(L_k(\g))\ne0$ for every pair $(k,f)$
listed above.

\subsection{\texorpdfstring{Non-quasi-lisseness of $L_k(F_4)$ and $L_k(G_2)$}{Non-quasi-lisseness of Lk(F4) and Lk(G2)}}
We treat the following levels, with the nilpotent elements specified above:
$$
 \begin{array}{c|l}
 \g&k\\\hline
 F_4&-8,-7,-5,-15/2\\
 G_2&-3,-10/3
 \end{array}
$$
The nonvanishing established above allows us to use a simple quotient of each
reduction as an obstruction. We use symplectic-fermion realizations for
$F_4$ at $k=-8,-7$ and for $G_2$ at $k=-3$. In the other three cases,
the simple quotient is a Virasoro vertex algebra. In each realization,
the class of the conformal vector is not nilpotent in the $C_2$ algebra,
so \cref{com:conformal-obstruction} excludes quasi-lisseness of the original
affine vertex algebra. The generators below establish the required
quotient identifications.

For $f$ in the nilpotent orbit $F_4(a_1)$, the universal algebra
$\mathcal W^{-8}(F_4,f)$ is freely strongly generated by six vectors
of weights $2,3,4,6,6,8$.
The \hyperref[calc:V5-F-T1]{Calculation~\ref*{calc:V5-F-T1}} gives the canonical conformal vector
and two BRST cocycles representing strong generators $C,Q$ of weights three and four.
The linear symbols of $C_{(0)}Q$ and $Q_{(1)}Q$ are independent in the
two-dimensional weight-six space, and that of
$C_{(0)}(C_{(0)}Q)$ is nonzero in weight eight. The first product
includes the normal-ordering terms in the universal $W$-algebra.
By \cref{v5:triangular-generation}, the universal $W$-algebra is vertex-generated
by $T,C,Q$. Let $\mathcal H^{\otimes3}$ be generated by
$h_1,h_2,h_3$, with $[h_i{}_{\lambda}h_j]=\delta_{ij}\lambda\mathbf1$.
Their images under the homomorphism to this vertex algebra are
$$
 T_i=\tfrac12(:h_i^2:-\partial h_i),\qquad T\longmapsto\sum_iT_i,
 \quad C\longmapsto0,\quad
 Q\longmapsto\bigl(\sum_i:T_iT_i:-10:TT:+3\partial^2T\bigr)/17.
$$
The image of the weight-four vector modulo $C_2$ has nonzero
coefficient of $p_4$, so
\cref{com:sf-recognition} identifies this quotient of
$\mathcal W^{-8}(F_4,F_4(a_1))$ with $\mathcal A_3$.
At $k=-7$ ($t=2$), for $f\in\OO_{F_4(a_3)}$, \cref{com:quadratic-miura} and
\hyperref[calc:QUAD]{Calculation~\ref*{calc:QUAD}} give two complete
quadratic lifts with matrices
$D=\operatorname{diag}(0,-4,4)$ and
$M=-2(E_{12}+E_{21})+E_{23}+E_{32}$.
The first matrix has distinct eigenvalues, and the graph defined by the
nonzero off-diagonal entries of the second is connected. The universal
$W$-algebra is generated by its weight-two subspace, whose image has
dimension at most six. Hence \cref{v3:quadratic-generation} identifies
the image with $\SF(\C^3\otimes\C^2)^{\mathrm{SL}_2}$.
The invariant bilinear forms inherited from the symplectic-fermion
realizations prove simplicity. In their lattice realizations, restriction
to the polynomial $C_2$ algebra of the corresponding Heisenberg vertex
algebra sends each conformal vector to a nonzero quadratic polynomial.
Its class is therefore not nilpotent.
At $k=-5$ ($t=4$), the image of $V^0(\slg_3)$ in
$\mathcal W_{-5}(F_4,A_2)$ is trivial, leaving a quotient of the positive part of the centralizer
whose abelianization is one line of degree two, by the published
centralizer row in \cref{v3:abelianization}. That lemma gives
$\mathcal W_{-5}(F_4,A_2)\cong\Vir_1$.
At $k=-15/2$ ($t=3/2$), with $f\in F_4$ in the nilpotent orbit labelled $B_3$,
the universal algebra
$\mathcal W^{-15/2}(F_4,f)$ has three strong generators of weight one,
one of weight two, five of weight four, and one of weight six.
The dimensions of the reduction $(1,0,1,1,2,2,4)$ through degree six,
including all intervening half-degrees, are \hyperref[calc:F-CHAR]{Calculation~\ref*{calc:F-CHAR}}.
The universal Virasoro vertex algebra $\Vir_{-2}$ is simple.
Its image exhausts the reduction in every weight of a strong generator,
and therefore equals $H^0_{B_3}(L_{-15/2}(F_4))$.

For $L_{-3}(G_2)$ ($t=1$), define the following fields in $\SF(2)$:
$T_i=-:a_ib_i:$, $L=T_1+T_2$, $P=T_1-T_2$,
$M=-\sqrt{-1}(:a_1b_2:+:a_2b_1:)$ and $F=-P_{(0)}M/2$.
The operator-product \hyperref[calc:G-IMAGE]{Calculation~\ref*{calc:G-IMAGE}} matches
\cite{Fasquel2025} and identifies the simple quotient of
$H^0_{f_{\rm sub}}(L_{-3}(G_2))$ with
$\SF(2)^{\mathrm{SL}_2}$, with its nonnilpotent conformal class.
At $k=-10/3$ ($t=2/3$), the simple quotient
$\mathcal W_{-10/3}(G_2,f_{\rm sub})\cong\Vir_{-24}$ is already known
from Fasquel \cite[Corollary~4.7(3)]{Fasquel2025}. We recall the short
argument to connect this identification with the nonzero reduction used here.
The complete pairings of $P,M$, including with the conformal vector, vanish. The centralizer abelianization has three degree-two lines.
Killing the two represented by $P,M$ leaves only the conformal-vector line.
Thus \cref{v3:abelianization} gives $\Vir_{-24}$ as the simple
quotient. The universal Virasoro vertex algebra at central charge $-24$ is simple:
the Virasoro minimal-model central-charge formula
$c=1-6(p-q)^2/(pq)$ gives $p/q=6$ or $1/6$, incompatible with
coprime $p,q\ge2$. Each Virasoro vertex algebra obtained here has $C_2$ algebra $\C[T]$.

\subsection{\texorpdfstring{Quasi-lisseness of $L_k(F_4)$}{Quasi-lisseness of Lk(F4)}}
We prove $X_{L_k(F_4)}\subseteq\N(F_4)$ for
$$
 k\in\{-6,-13/2\}\cup\{-3+m,-9/2+m:m\in\Z_{\ge0}\}.
$$
At $t=6$, Adamovi\'c, Dai and Sun \cite[Theorems~4.2 and~4.5]{ADS2026} prove
$X_{L_{-3}(F_4)}=\overline{\OO}_{\widetilde A_1}$.
The argument below also treats $L_{-3+m}(F_4)$ for $m\in\Z_{\ge0}$.

At $t=3$, we first rule out nonzero points on a coordinate axis by identifying the two
conformal vectors in a reduction. Two sheet bounds then reduce every
possible nonzero semisimple point to this axis.
Reduce at $\widetilde A_2$, $x=(2,0,0,0)$.
The triple decomposition is $V_2+(7\otimes V_4)+(14\otimes V_0)$.
The weight-one subspace generates a quotient of $V^{-2}(G_2)$,
and both conformal vectors have central charge
$-14$. The complete fifteen-term truncated Jantzen sum \hyperref[calc:F-WINDOWS]{Calculation~\ref*{calc:F-WINDOWS}}
gives $I_{<2}=0$, $\dim I_2=1$ for the reduced maximal ideal $I$.
The $G_2$-invariant subspace of the conformal-weight-two space of
$\mathcal W^{-6}(F_4,\widetilde A_2)$ is
$\C T_W+\C T_{G_2}$. Applying $J^u_{(1)}$ to a relation
forces opposite coefficients. Hence $T_W=T_{G_2}$ in $H_f^0(L_{-6}(F_4))$.
The remaining generating $G_2$-module of dimension seven has
conformal weight three and Sugawara conformal weight one, so
\cref{v3:resonance} shows that the reduction is generated by $G_2$ and
\cref{vtwo:sugawara-axis} shows that $X_{L_{-6}(F_4)}$ contains no
nonzero point of $\C e_1$.
The long-$A_2$ and $B_2$ level-zero sheet bounds
\cite[Theorem~3.5]{ArakawaFutornyKrizka2024} have center planes $(a,b,b,b)$ and
$(a,b,0,0)$. Their common semisimple locus consists of Weyl-conjugate
coordinate axes: a long $A_2$ forces three absolute coordinates equal.
Consequently $X_{L_{-6}(F_4)}$ contains no nonzero semisimple point.

We use the following consequence of $H_f^0(L_k(\g))=0$ for $f\in\OO$.
For each maximal parabolic, find $f=n+v\in\OO$ with $v$ in its nilradical
and $n$ in the associated variety of the corresponding simple affine
vertex algebra of the Levi subalgebra. A semisimple point would give
$az+n\in X$. For $a\ne0$, invertibility of $\ad(az+n)$ on positive
center degrees conjugates $az+n+v$ to it. Closure then gives the forbidden $f$.

At $k=-13/2$ ($t=5/2$), the twenty-five-term truncated Jantzen sum gives
$$
 H^0_{B_3}(L_{-13/2}(F_4))=0.
$$
The four explicit root triples of \hyperref[calc:F-WINDOWS]{Calculation~\ref*{calc:F-WINDOWS}} have Levi
parts $C_3:(2^3)$, long-root $A_1$, long-root $A_2$, and
$B_3:(3^2,1)$, respectively. For the two outer Levi subalgebras, the parameter $k+h^\vee$ is $5/2$.
For $C_3$ the conformal-weight differences $(1/5,1/5,0)$ have an equality case with Sugawara conformal weight
$3/2$, which is not integral. For $B_3$ the differences $(1/5,0,1/20)$ leave
only the degree-two adjoint, excluded because $k=-5/2\ne0$ and
$\Sym^2(21)=1+27+35+168$ (\cref{vtwo:current-descendant}).
The long $A_1,A_2$ levels $1/2,-1/2$ give the minimal-orbit inclusions by
\cite[Proposition~4.7.2]{AAV}. Short factors contribute only zero.
Thus the required nilpotent points belong to the corresponding associated varieties.

At $k=-3,-9/2$ ($t=6,9/2$), the eight- and nine-term truncated Jantzen sums show that
$\OO_{A_1+\widetilde A_1}$ and $\OO_{A_2+\widetilde A_1}$, respectively,
are absent from the associated varieties.
The root placements in the same calculation place the first orbit in
every nilradical. For the second, only the two outer parabolics need
a Levi point, a long root at levels $1/2,-1/2$.
Translation by an integrable affine vertex algebra preserves the absence of these orbits, while the
levels on the simple ideals of the two outer Levi subalgebras remain nonintegral.
This proves quasi-lisseness of $L_{-3+m}(F_4)$ and $L_{-9/2+m}(F_4)$
for $m\ge0$. The argument uses only the stated orbit exclusions at
these levels. It does not require equality of their associated varieties.

\subsection{\texorpdfstring{Quasi-lisseness of $L_{-4}(F_4)$ and $L_{-11/2}(F_4)$}{Quasi-lisseness of L-4(F4) and L-11/2(F4)}}
At $t=5$, the minimal Drinfeld--Sokolov reduction of $L_{-4}(F_4)$ is a simple conformal
extension of $U=L_{-3/2}(C_3)$ by one simple module whose lowest
conformal-weight space has highest weight $\omega_3$. This decomposition
is due to Adamovi\'c, Kac, M\"oseneder Frajria, Papi and Per\v se
\cite[Theorem~6.8]{AdamovicEtAl2017}. The simplicity and nonvanishing of the reduction
follow from \cite[Theorem~6.2]{ArakawaMoreau2017}.
By \cref{seed:C-full}, $\dim X_U=12$. The additional simple
$U$-module is generated by its finite-dimensional lowest
conformal-weight space under the operators $a(-n)$, $a\in C_3$,
$n>0$. Modulo $C_2$ of the extension, terms with $n\ge2$ vanish,
and $a(-1)$ acts by multiplication by $[a]$. Consequently the
$C_2$ algebra of the extension is a finite module over $R_U$.
Its associated variety, identified with the corresponding Slodowy intersection, has dimension
at most twelve. The $B_3$ sheet has a single nonnilpotent stratum
$G\C^*e_1$, whose intersection with the minimal Slodowy slice has
dimension $31-16=15$. Thus this stratum cannot occur in $X_{L_{-4}(F_4)}$.
At $t=7/2$, the dimensions of the $B_2$ reduction
$(1,0,1,1,2,2,3)$ through six, with all half-degrees, are
\hyperref[calc:F-CHAR]{Calculation~\ref*{calc:F-CHAR}}. Strong weights are at most four, and the
maximal proper ideal of the universal Virasoro vertex algebra at central charge $c=-68/7$ is generated in degree six.
The graded-dimension criterion of \cref{com:finite-image} therefore gives
\begin{equation}\label{seed:F7actual}
H^0_{B_2}(L_{-11/2}(F_4))=\Vir(2,7).
\end{equation}
The two accompanying vanishing results for Drinfeld--Sokolov reduction exclude $C_3(a_1),\widetilde A_2$,
giving $X_{L_{-11/2}(F_4)}\cap\N(F_4)\subseteq\overline{\OO_{B_2}}$,
of dimension $36$. The maximal Levi subalgebras with more than one
simple factor give induced nilpotent orbits of dimension $40$.
A nonzero point of the center of the $B_3$ Levi subalgebra would force
$\OO_{\widetilde A_2}$ into the associated variety, contrary to the
vanishing reduction above. The associated variety of the remaining
affine vertex algebra $L_{-1/2}(C_3)$ is the minimal nilpotent orbit
closure of dimension six, giving an induced orbit of dimension $52-22+6=36$.
The closure of this semisimple family has dimension $37$, contrary
to the zero-dimensional Slodowy intersection.

\subsection{\texorpdfstring{Quasi-lisseness of $L_k(G_2)$ at $k=-2,-5/2,-8/3,-7/3$}{Quasi-lisseness of Lk(G2) at k=-2,-5/2,-8/3,-7/3}}
We show that the associated varieties at these four levels are contained
in $\N(G_2)$. For numerator three, we use $k=-5/2$ ($t=3/2$)
and then apply Theorem~\ref{thm:denominator}.
At $k=-2$, the associated variety is the subregular orbit closure by
\cite[Theorem~A]{ArakawaEtAl2025}.
For principal reduction, the bound $D\le26$ gives thirteen root contributions.
Five equal-degree pairs and two singular alternants cancel, leaving
$M((7,1),26)$ whose reduction has a one-dimensional conformal-weight-zero space
(\hyperref[calc:G-PRINC]{Calculation~\ref*{calc:G-PRINC}}). Hence the principal Drinfeld--Sokolov reduction of $L_{-5/2}(G_2)$ vanishes.
The simple affine vertex algebras for the long-root and short-root Levi subalgebras
are $L_{-1/2}(\slg_2)$ and $L_{5/2}(\slg_2)$, respectively, with associated varieties
equal to the nilpotent cone of $\slg_2$. A regular
Levi vector plus the other simple root is regular in $G_2$.
Neither maximal Levi center has a nonzero point in
$X_{L_{-5/2}(G_2)}$. Theorem~\ref{thm:denominator} then
applies because $\gcd(2,3)=1$.

The proof for $L_{-8/3}(G_2)$ and $L_{-7/3}(G_2)$ uses an elementary bound on polynomial
solutions of a regular singular differential equation. We record the
coefficient comparison explicitly: a putative polynomial solution is
determined by its constant term, and one further coefficient suffices
to rule it out.

\begin{lemma}\label{vtwo:one-return}
Consider $u y'=A(u)y$, $A(u)=\sum_{j\ge0}A_ju^j$ on a finite-dimensional
vector space $E$. Set $K=\ker A_0$, and assume $nI-A_0$ is invertible for
$1\le n\le N+1$. Define $T_0:K\hookrightarrow E$ and
$$
 T_n=(nI-A_0)^{-1}\sum_{j=1}^n A_jT_{n-j}.
$$
If, for a basis $k_j$ of $K$ and linear functionals $\ell_i$, the
square matrix $(\ell_i(T_{N+1}k_j))$ is invertible, there is no nonzero
polynomial solution of degree at most $N$ with initial value in $K$.
\end{lemma}
\begin{proof}
Coefficient comparison forces $y_n=T_ny_0$. A polynomial of degree
at most $N$ has $y_{N+1}=0$, so the matrix forces $y_0=0$. The invertibility of $nI-A_0$, $1\le n\le N$, then gives
$y_n=0$ for all $n$.
\end{proof}

At $t=p/3$, $p=4,5$, the ten- and thirteen-term initial Jantzen sums give
spaces of singular vectors $V((p-2)\rho)$ in degree $D_0=7(p-2)$.
Their nonzero symbols satisfy
$dF=t^{-1}\sum_{\alpha>0}A_\alpha\,d\log\alpha\,F$,
$A_\alpha=f_\alpha e_\alpha/d_\alpha$, with short/long $d=3,1$.
The decomposition under the root $\slg_2$-subalgebras and the Weyl
action imply that the first possible vanishing orders on a short-root
hyperplane are $3,4$, and on a long-root hyperplane are $9,12$,
respectively. Vanishing on all three long-root hyperplanes would
require degree at least $27,36>D_0$. Vanishing on a short-root
hyperplane would imply
$F=\Delta_{\rm short}^{p-1}G$, $\deg G=N=5,9$.
Here, for a short simple root $\alpha$ and a long simple root $\beta$,
$$
 \Delta_{\rm short}=\alpha(\alpha+\beta)(2\alpha+\beta).
$$
To apply \cref{vtwo:one-return}, take the chart $\alpha=1$, $\beta=u$.
Dividing out this factor replaces the residue $A_r$ by
$A'_r=A_r-p(p-1)I/3$ for short roots and does not alter it for
long roots. The equation for $G$ is
$$
 uG'(u)=A(u)G(u),\qquad
 A(u)=\frac3p\sum_{\substack{r=a\alpha+b\beta>0\\b\ne0}}
                   \frac{bu}{a+bu}A'_r,
 \qquad A(0)=\frac3p A_\beta.
$$
The apparent singularity of $A(u)$ at zero is removable.
Write $m=p-2$. After identifying $V(m\rho)$ with its dual by the
contravariant form, work in the opposition eigenspace of sign
$\epsilon=(-1)^{D_0}$ in its zero-weight space. The residue operators
preserve this eigenspace. Within it, $K=\ker A(0)$ is the zero-weight
part of $V(m\rho)^{\slg_{2,\beta}}$ with this opposition eigenvalue.
The sign is determined by $F$, before division
by $\Delta_{\rm short}^{p-1}$. The six-factor Weyl character product
shows that this initial-value space $K$ has dimension four for both
$m=2$ and $m=3$.

Let $Q$ be the orthogonal projector onto $K$ for the contravariant
form. In the normalization of
\hyperref[calc:G-RETURN]{Calculation~\ref*{calc:G-RETURN}}, the basis vectors are
$k_j=Qw_j$ for four specified weight-zero PBW vectors, and the
functionals are $\ell_i(v)=(w_i,Qv)$. For $p=4$, the first forbidden
coefficient is the coefficient of $u^6$. Its matrix is
$$
 (\ell_i(T_6k_j))\equiv
 \begin{pmatrix}
 335&238&682&32\\
 139&621&379&0\\
 782&757&540&226\\
 168&582&501&894
 \end{pmatrix}\pmod {1009},
 \qquad \det\equiv916\pmod {1009}.
$$
For $p=5$, the analogous coefficient is that of $u^{10}$ and its
determinant is congruent to $382$ modulo $1009$.
The Gram determinants of the four vectors are $225,980$ modulo
$1009$, respectively, so they do form bases of $K$.
The operators, PBW vectors, denominator bounds and invertibility
of $nI-A_0$ for $1\le n\le N+1$ are given in the same calculation.
All denominators are prime to $1009$, so these nonzero modular
determinants imply nonzero rational determinants.
By \cref{vtwo:one-return}, the resulting differential equation has no
nonzero polynomial solution of the required degree. Thus the associated
varieties of $L_{-8/3}(G_2)$ and $L_{-7/3}(G_2)$ contain no nonzero
point of either maximal Levi center.

This proves \cref{thm:FG-classification}. The associated varieties at
other denominators are determined separately in the appendix.

\section{\texorpdfstring{Quasi-lisse simple affine vertex algebras of types $E_6,E_7,E_8$}{Quasi-lisse simple affine vertex algebras of types E6,E7,E8}}
\label{sec:exceptional-e}

Uniform features of exceptional Lie algebras have long been studied, see \cite{Cvitanovic1981,Deligne1996}.
The minimal-orbit examples of Arakawa and Moreau
\cite[Theorem~1.1(2)]{ArakawaMoreau2018} form one part of the classification below.
For the remaining non-quasi-lisse levels we prove that the chosen
reduction is nonzero and identify a simple quotient that is not
quasi-lisse. The quasi-lisse cases in Section~\ref{exc:seeds} instead
require relations in $H_f^0(L_k(\g))$ itself.

The affine VOAs $L_{-3}(E_6)$, $L_{-4}(E_7)$ and $L_{-6}(E_8)$ occur as the vertex algebras of the rank-one Minahan--Nemeschansky theories, whose Higgs branches are minimal nilpotent orbit closures
\cite[Section~4]{BeemRastelli2018}. Free-field constructions from
these Higgs branches are given in
\cite[Section~3]{BeemMeneghelliRastelli2019}.
Further full-puncture predictions of Xie and Yan
\cite[Tables~20 and~22]{XieYanLisse2021} agree with our associated
varieties in the cases listed in \cref{tab:physics-E}. The orbit column
specifies the closure predicted as the Higgs branch, and the final column
is the Coulomb-branch rank of the proposed 4d theory. The cases $E_6$ at $k=-6$ and $E_8$ at $k=-10,-15$ also appeared in \cite[Tables~5, 7 and~17]{LiXieYan2023}.

\begin{table}[htbp]\centering\small
\begin{tabular}{@{}lrlr@{}}\toprule
$\g$&$k$&Orbit&4d rank\\\midrule
$E_6$&$-6$&$A_2$&3\\
$E_7$&$-12$&$D_4(a_1)$&7\\
$E_8$&$-10$&$2A_1$&2\\
$E_8$&$-15$&$A_2+A_1$&4\\
$E_8$&$-18$&$2A_2$&6\\
$E_8$&$-20$&$D_4(a_1)+A_1$&8\\
$E_8$&$-24$&$E_8(a_7)$&16\\\bottomrule
\end{tabular}
\par\smallskip
\caption{Exceptional integer-level examples predicted in \cite{XieYanLisse2021}.}
\label{tab:physics-E}
\end{table}

We use Bourbaki numbering and the normalization fixed in the introduction.
The classification is as follows.
\begin{theorem}\label{thm:exceptional-e}
Let $\g$ have type $E_6$, $E_7$ or $E_8$, and let $k\ne-h^\vee$.
Then $L_k(\g)$ is quasi-lisse if and only if
$t=k+h^\vee=p/q>0$, with $\gcd(p,q)=1$, and respectively
$$
 \begin{array}{c|l}
 E_6&p\in\{6,7\}\text{ or }p\ge9,\\
 E_7&p\in\{6,8\}\text{ or }p\ge10,\\
 E_8&p=6\text{ or }p\ge9.
 \end{array}
$$
\end{theorem}
\begin{proof}
By \cref{thm:denominator}, it suffices to consider $k=-h^\vee+p$.
The levels with $p\ge h^\vee$ are admissible. For $p=1$, the vertex
algebra $V^{-h^\vee+1}(\g)$ is simple and its associated variety is
$\g^*$. The remaining levels are treated in
Sections~\ref{exc:negative} and~\ref{exc:seeds}.
\end{proof}

\subsection{\texorpdfstring{Drinfeld--Sokolov reductions of $L_k(E_6)$, $L_k(E_7)$ and $L_k(E_8)$}{Drinfeld--Sokolov reductions of Lk(E6), Lk(E7) and Lk(E8)}}
\label{exc:negative}
Table~\ref{exc:negative-diagrams} specifies the reductions
$H_f^0(L_{-h^\vee+p}(\g))$ used to prove non-quasi-lisseness.
For most selected orbits, nonvanishing follows from conformal-weight
inequalities and the simple quotient admits an explicit description.
Most quotients are Virasoro vertex algebras,
symplectic-fermion invariant algebras, or algebras with a Heisenberg
tensor factor. These descriptions provide the obstruction to
quasi-lisseness. At $E_6$, $p=4$ ($k=-8$), nonvanishing is conditional
on quasi-lisseness, which the quotient argument then contradicts.
The equality cases in the conformal-weight inequalities follow the table.

In Bourbaki numbering, the chain is $1,3,4,5,6,7,8$, with node $2$
attached to node $4$. Write $p=k+h^\vee$ for the integer representative.
The conformal-weight formula used to prove nonvanishing is
\begin{equation}\label{exc:energy}
 \Delta_p(\lambda)-\lambda(x)
 =\sum_i a_i\bigl(\Delta_p(\omega_i)-\omega_i(x)\bigr)
 +\sum_i\frac{K_{ii}a_i(a_i-1)}{2p}
 +\sum_{i<j}\frac{K_{ij}a_i a_j}{p},
 \qquad\lambda=\sum_i a_i\omega_i,
\end{equation}
where $K$ is the inverse Cartan matrix, whose entries are positive. The right-hand side is positive for every nonzero dominant weight
whenever it is positive on all fundamental weights.
The completion of the character sums and the equality cases are treated in
\cref{com:completion}, \cref{com:energy}, \cref{com:equality}.

\Needspace{10\baselineskip}
\begin{longtable}{@{}ccll@{}}
\caption{Reductions $H_f^0(L_{-h^\vee+p}(\g))$ used to prove that the indicated affine vertex algebras are not quasi-lisse. The last column gives the weighted Dynkin diagram $2x$.}
\label{exc:negative-diagrams}\\\toprule
Type&$p$&orbit&diagram\\\midrule\endfirsthead
\toprule Type&$p$&orbit&diagram\\\midrule\endhead
$E_6$&2&$E_6(a_3)$&$200202$\\
&3&$D_4(a_1)$&$000200$\\
&4&$D_4(a_1)$&$000200$\\
&5&$A_2+2A_1$&$001010$\\ \midrule
$E_7$&2&$E_6(a_1)$&$2002020$\\
&3&$E_7(a_5)$&$0002002$\\
&4&$A_4+A_2$&$0002000$\\
&5&$A_3+A_2+A_1$&$0000200$\\
&7&$A_2+3A_1$&$0200000$\\
&9&$A_2+A_1$&$1000010$\\ \midrule
$E_8$&2&$E_8(a_4)$&$20020202$\\
&3&$E_8(a_6)$&$00020020$\\
&4&$D_7(a_2)$&$10010101$\\
&5&$E_8(a_7)$&$00002000$\\
&7&$A_4+A_2+A_1$&$00100100$\\
&8&$D_4(a_1)+A_2$&$02000000$\\\bottomrule
\end{longtable}
Except at $E_6,p=4$ and $E_7,p=9$, the conformal-weight inequalities are strict.
For numerator two the three lists of fundamental conformal-weight differences are, explicitly,
$$
\begin{array}{c|l}
E_6&1/3,1,4/3,2,4/3,1/3\\
E_7&1,17/8,3,5,29/8,2,9/8\\
E_8&2,4,6,11,8,5,3,1.
\end{array}
$$
Thus the same conformal-weight argument treats these three cases. At $E_7,p=9$,
equality can occur only for the adjoint highest weight $\omega_1$
in degree two of $V^{-9}(E_7)$.
Since $k=-9\ne0$ and $\Sym^2(133)=1+1539+7371$ has no adjoint,
\cref{vtwo:current-descendant} excludes it from the maximal ideal of
$V^{-9}(E_7)$.
At $E_6,p=4$, the possibilities $\omega_2$ at degree three and
$\omega_4$ at degree six remain. The sheet argument below proves
nonvanishing under the hypothesis of quasi-lisseness.

Write $\mathcal A_r=\SF(\C^{2r})^{\Sp_{2r}}$. Below we identify simple conformal
quotients of these nonzero reductions. Their realizations using
symplectic fermions and their tensor product decompositions show that
the classes of their conformal vectors in the corresponding $C_2$
algebras are not nilpotent. Hence
\cref{com:conformal-obstruction} proves that the corresponding
$L_k(\g)$ are not quasi-lisse.

\subsection{\texorpdfstring{Non-quasi-lisseness of $L_{-10}(E_6)$}{Non-quasi-lisseness of L-10(E6)}}
For the nilpotent element $f$ specified below, we construct a quotient
$B$ of $\mathcal W^{-10}(E_6,f)$ with a tensor factor $\Vir_{-2}$. The role of the quadratic calculation below is to show that
the entire factor belongs to the image, so that it survives in every
simple graded quotient.
Let $E_1,E_2,E_3$ be the natural two-dimensional modules of the
three simple factors of $\g_0=\slg_2^{\oplus3}\oplus\C^3$, with bases
$a_0,a_1$, $b_0,b_1$ and $c_0,c_1$, respectively. Write
$v_{ij\ell}=a_i\otimes b_j\otimes c_\ell$. As a module for the
semisimple part of $\g_0$, the degree $-1$ subspace is
$E_1\oplus(E_1\otimes E_2\otimes E_3)\oplus E_3$.
We use the central action and root-vector normalization specified in
\hyperref[calc:E6-LIFT]{Calculation~\ref*{calc:E6-LIFT}}, and set
$
 f=a_0+v_{000}+v_{011}+v_{101}+v_{110}+c_0.
$ 
The infinitesimal stabilizer is zero. The equations from $a_0,c_0$
force the lowering entries of the first and third $\slg_2$ matrices
to vanish. The four odd tensor coordinates force the remaining
off-diagonal entries to vanish, and the four even coordinates force
the diagonal and central entries to vanish. Thus the orbit is
open in the even Dynkin grading $x=\omega_1^\vee+\omega_4^\vee+
\omega_6^\vee$. Deleting the degree-one nodes $1,6$ gives a compatible
Levi subalgebra $D_4\oplus\C^2$. The restricted nilpotent element has
partition $(3^2,1^2)$, and the induced level of $D_4$ is $-4$.
The homomorphism of \cite[Theorem~5.10]{Genra} therefore gives
$$
 B=\im\mathcal W^{-10}(E_6,f)\ \subset\ \mathcal H^2\otimes\Vir_{-2}.
$$
The image of $\mathcal W^{-4}(D_4,(3^2,1^2))$ is $\Vir_{-2}$.
Indeed, write the natural $D_4$ module as
$V_2\otimes U_3\oplus V_0\otimes U_1$, where $V_j$ is the
$(j+1)$-dimensional module for the defining $\slg_2$-triple and
$U_3,U_1$ are its two-dimensional multiplicity spaces. The weight-one
generators associated with the torus vanish. Among the weight-two
subspaces $\Sym^2U_3\oplus(U_3\otimes U_1)$, only the trace remains.
The bracket $[zS_+,zS_-]=z^2[S_+,S_-]\ne0$ eliminates the weight-three
line by \cref{vtwo:cotangent}. 

Let $z_1,z_2,z_3$ be a basis of the center of $\g_0$, with Gram matrix
$$
 G=\begin{pmatrix}4/3&2&2/3\\2&6&2\\2/3&2&4/3\end{pmatrix}.
$$
At $\tau=k+12=2$, set $u=z_1,v=z_3,z=z_2-z_1-z_3$. The conformal vector is
$$
 T_++T_-,\qquad T_+=(u^2-uv+v^2)/4,\qquad
 T_-=z^2/8-\partial z/4.
$$
For the vector $-2v_{000}$, \cref{com:quadratic-miura} gives the
following expression. Write $e_i,h_i$ for the indicated elements of the
three $\slg_2$ summands of $\g_0$, and set
$$
 K=h_1+h_2+h_3+z_1-z_2+z_3,\qquad
 F_\tau=\tfrac18:K^2:-\tfrac12\sum_i:e_i^2:
       +\sum_{i<j}:e_ie_j:+(\tfrac54-\tfrac\tau2)\partial K.
$$
The differential at the identity of the orbit map $G_0\to\g_{-1}$,
$g\mapsto\operatorname{Ad}(g)f$, is $a\mapsto[a,f]$.
Its matrix is $12\times12$ with determinant $-32$
in the bases specified in
\hyperref[calc:E6-LIFT]{Calculation~\ref*{calc:E6-LIFT}}. Substitution in \eqref{com:quadratic-miura-formula} gives exactly
$F_\tau$ as the complete universal Miura image of $W_{-2v_{000}}$,
including the displayed derivative coefficient. The finite equations
and normalization are given in
\hyperref[calc:E6-LIFT]{Calculation~\ref*{calc:E6-LIFT}}.
At $\tau=2$, scalar projection sends $K$ to $-z$ and $F_2$ to $T_-$.

Hence $1\otimes\Vir_{-2}\subset B$. The tensor product lemma \cref{vtwo:whole-factor}
gives $B=C_B\otimes\Vir_{-2}$, where
$C_B=\{u\in\mathcal H^2:u\otimes\mathbf1\in B\}$.
Its simple graded quotient is $C_{B,\mathrm s}\otimes\Vir_{-2}$,
where $C_{B,\mathrm s}$ is the simple graded quotient of $C_B$.
Under augmentation of $R_{C_{B,\mathrm s}}$, the total conformal class
maps to the polynomial generator of $R_{\Vir_{-2}}$.
The nonvanishing proved above now implies that $L_{-10}(E_6)$ is not quasi-lisse.

\subsection{\texorpdfstring{Non-quasi-lisseness of $L_{-16}(E_7)$ and $L_{-28}(E_8)$}{Non-quasi-lisseness of L-16(E7) and L-28(E8)}}
We realize the simple quotients of
$H^0_{E_6(a_1)}(L_{-16}(E_7))$ and
$H^0_{E_8(a_4)}(L_{-28}(E_8))$ as the symplectic-fermion invariant
vertex algebras $\mathcal A_3$ and $\mathcal A_4$, respectively.
For $E_7$ use the regular symmetric-pair model
$E_7=\slg_8\oplus\Lambda^4\C^8$ with a principal $A_7$ triple.
Its natural $h$ is $(7,5,3,1,-1,-3,-5,-7)$ and its dominant diagram
is $2002020$. The three orthogonal roots defining the generators $h_1,h_2,h_3$
of the rank-three Heisenberg vertex algebra are given in
\hyperref[calc:NUM2]{Calculation~\ref*{calc:NUM2}}. Normalize them by $[h_i{}_{\lambda}h_j]=\delta_{ij}\lambda\mathbf1$.
The eight diagonal entries are
$(\epsilon_1h_1+\epsilon_2h_2+\epsilon_3h_3)/2$. Thus
$$
 \tr Z^4=\tfrac32p_2^2-p_4,
 \qquad p_j=\sum_i h_i^j.
$$
Restriction to the principal Slodowy slice gives this homogeneous
polynomial in $R_{\mathcal W^k(\g,f)}=\C[\mathcal S_f]$.
Lifting it through the graded quotient map gives a weight-four vector
in $\mathcal W^k(\g,f)$ with the stated image modulo $C_2$.
The filtered restriction uses the same integral triple in the regular
subalgebra.

The images of the eight centralizer generators on which the involution
acts by $-1$ are zero. The five explicit four-forms in \hyperref[calc:NUM2]{Calculation~\ref*{calc:NUM2}} give
wedge pairings $6,-4,-4,4$, eliminating the even generators of weights
$3,5,7,8$ by \cref{v5:equivariant-generation}, \cref{vtwo:cotangent}. The remaining strong generators have weights $2,4,6$. Deleting node $1$ gives an explicit
$\mathcal H\otimes \mathcal A_2$ containment: its two $D_6$ column axes are
$e_1+e_2,e_3+e_4$, and intrinsically the three unhalved axes are
$\omega_1+\omega_6-\omega_4,\omega_1-\omega_6,-\omega_1$.
They identify the three roots defining the Heisenberg generators,
up to signs and permutation.
The conformal vector, the symmetric quartic polynomial independent of
$p_2^2$, and this bound on the strong generators satisfy
\cref{vtwo:scalar-recognition}. Thus the image of
$\mathcal W^{-16}(E_7,E_6(a_1))$ in this realization is $\mathcal A_3$.

For $E_8$ use $E_8=D_8\oplus S^+$ and the principal $D_8$ grading
$x=(7,6,5,4,3,2,1,0)$. Its dominant diagram is $20020202$.
The four orthogonal roots defining the Heisenberg generators are the rows of the signed Walsh matrix
in \hyperref[calc:NUM2]{Calculation~\ref*{calc:NUM2}}. Eight fourth powers give the natural $D_8$ identity
$\tr_{16}Z^4=3p_2^2-2p_4$, and the same filtered slice argument
gives a weight-four vector whose leading term is independent of $p_2^2$.
The images of the half-spin generators vanish by the $D_8$ involution. Five explicit lowest spinors give the four
pairings $1,4,2,-8$ in \hyperref[calc:NUM2]{Calculation~\ref*{calc:NUM2}}.
\Cref{v5:equivariant-generation} eliminates one
weight-eight direction and the weights $10,12,14$, leaving strong generators of weights $2,4,6,8$. Deleting grade-one node $8$ gives an explicit
$\mathcal H\otimes \mathcal A_3$ containment. The additional Heisenberg generator corresponds to the root $\omega_8$, and
the three orthogonal roots from the $E_7$ realization give the other
three generators. Proposition~\ref{com:sf-recognition} therefore identifies
the image of $\mathcal W^{-28}(E_8,E_8(a_4))$ with $\mathcal A_4$.

\subsection{\texorpdfstring{Relations in $H^0_{E_8(a_6)}(L_{-27}(E_8))$}{Relations in the E8(a6) reduction of L-27(E8)}}
Let $J^1$ be the maximal proper ideal of $V^{-27}(E_8)$ and
$f\in\OO_{E_8(a_6)}$. We use the lowest nonzero conformal-weight
subspace of $H_f^0(J^1)\subset\mathcal W^{-27}(E_8,f)$ to identify a
simple quotient of $H_f^0(L_{-27}(E_8))$.
The lemma computes the action of quadratic generators on the lowest
conformal-weight subspace of a reduced standard module, using Genra's
Wakimoto reduction and Miura homomorphism
\cite[Lemma~3.1, Proposition~4.5 and Corollary~4.6]{Genra}.
We give the argument to obtain the required zero-mode formula,
including its derivative correction.

\begin{lemma}
\label{vtwo:miura-top}
Let $x$ define an integral dominant good grading, and let $\g_0$ be
its degree-zero Levi subalgebra. The lowest space of $H_f^0(\mathbb V_k(\lambda))$ has degree
$\delta_x(\lambda)$ and is the finite highest-$x$ space
$E_{\lambda,\max}$, with the action pulled back through the universal
Miura map. In particular, if a quadratic Miura field is
$$
 F=\sum_{a,b}q_{ab}:J^aJ^b:+\sum_a d_a\partial J^a,
 \qquad q_{ab}=q_{ba},\quad a,b\in\g_0,
$$
the action of $o(F)=F_{(1)}$ on this space is
\begin{equation}\label{vtwo:miura-action}
 o(F)=\sum_{a,b}q_{ab}\rho(a)\rho(b)-\sum_a d_a\rho(a).
\end{equation}
The highest weight for $\g_0$ is $\lambda$.
\end{lemma}
\begin{proof}
The finite big-cell realization embeds $E_\lambda$ in
$\C[N_+]\otimes|\lambda\rangle$, with highest vector $1$.
For $N_+=G_{>0}G_0^+$, the highest-$x$ subspace involves only the
$G_0^+$ coordinates. This embedding extends by affine induction to
$\mathbb V_k(\lambda)\to\operatorname{Wak}_{\g}^k(\lambda)$:
the operators $a(n)$, $a\in\g$, $n>0$, annihilate the polynomial
subspace of lowest degree, and $a(0)$ acts by the finite-dimensional
representation realized by differential operators. In the BRST
complex, vectors of $E_{\lambda,\max}$ tensored with the vacuum vector
of the charged fermion vertex superalgebra are cycles: $a(n)$ vanishes
on the lowest-degree subspace for $n>0$, $a(0)$ vanishes on
$E_{\lambda,\max}$ for $a\in\g_{>0}$, and the remaining terms of
the BRST differential annihilate the fermion vacuum vector. There are no neutral ghosts for this integral grading.

Wakimoto reduction contracts the $G_{>0}$ coordinates and ghosts and is
the identity on the $G_0^+$ beta-gamma and Heisenberg factors. Our cycles
therefore map to the same independent finite polynomials. Their
number equals the lowest coefficient of
\eqref{com:standard-ds-character}, so they exhaust the lowest
space and carry exactly the asserted Miura action. These maps and
contractions, valid at every level,
are given by Genra, \cite{Genra}, Lemma~3.1,
Proposition~4.5 and Corollary~4.6. The specialized Miura homomorphism
is covered by Corollary~4.13 of the same paper.
A generic Wakimoto resolution is not specialized here. In Zhu's algebra
$[\partial J^a]=-[J^a]$. Symmetric quadratic ordering cancels the
antisymmetric commutator corrections, giving
\eqref{vtwo:miura-action}. The background correction is already in
$d_a$ and the invariant bilinear form
$k(\ ,\ )+(\kappa_\g-\kappa_{\g_0})/2$.
\end{proof}

Set $U=\mathcal W^{-27}(E_8,E_8(a_6))$, with $x=\omega_4+\omega_7$,
and $I_i=H_f^0(J^i)\subset U$. The positive fundamental conformal-weight differences are
$(2,3,14/3,8,6,4,2,1)$, so the quotient $H=U/I_1$ is nonzero.
Its universal dimensions in weights two and four are $3,15$.
The centralizer abelianization has three generators each in $h$-degrees
$2,6$. \Cref{v3:abelianization} therefore gives vertex generation in
weights two and four. The sum of nineteen root contributions in
\hyperref[calc:E8-TOP]{Calculation~\ref*{calc:E8-TOP}} is complete for modes $a\le13$, since
$\max_\alpha\alpha(x)=9$ and only reduced degrees at most four are considered. It cancels
characters of standard modules of the same degree and leaves only
$M(\omega_7+\omega_8,31)$, whose subspace of maximal $x$-eigenvalue is a two-dimensional module
for the $A_1$ summand of $\g_0$.

The Miura action on the lowest conformal-weight space of \cref{vtwo:miura-top} evaluates
the three complete quadratic fields of \cref{com:quadratic-miura},
with projections in \hyperref[calc:CENTRALIZER]{Calculation~\ref*{calc:CENTRALIZER}} as
$$
 M_1=\begin{pmatrix}0&-2\\-24&0\end{pmatrix},\qquad
 M_2=\operatorname{diag}(-3/2,-1/2),\qquad
 M_3=\operatorname{diag}(-21/2,-23/2).
$$
The conformal vector acts by $-(M_2+M_3)/3=4I$, as required by the grading. Since $M_2$ has
distinct eigenvalues and both off-diagonal entries of $M_1$ are nonzero,
the lowest conformal-weight subspace $R$ is irreducible. By
\cref{vtwo:first-top}, the Jantzen identity in the completed Grothendieck
group and the exactness of the Drinfeld--Sokolov reduction functor
therefore imply
$$
 (I_1)_{<4}=0,\quad(I_1)_4=R,\quad(I_i)_4=0\ (i\ge2),
 \qquad(\dim H_2,\dim H_4)=(3,13).
$$
Since the reduced Jantzen ideals are nested, irreducibility prevents this
lowest subspace from contributing again in a later ideal.

The complete images in \hyperref[calc:CENTRALIZER]{Calculation~\ref*{calc:CENTRALIZER}}
give a conformal surjection from the subalgebra $C\subset U$
generated by these quadratic fields onto $A=\SF(\C^2\otimes\C^4)^{\Sp_4}$, simple with weight-two and
weight-four dimensions $3,13$. Since $H\ne0$, the proper graded ideal
$C\cap I_1$ lies in $\ker(C\to A)$. The image of $C$ in $H$ therefore
fills both generating weights. By \cref{vtwo:small-subalgebra}, it is
$H$, giving $H\twoheadrightarrow A$. By \cref{vtwo:torus-detector},
the conformal class of $A$ is not nilpotent. Therefore $L_{-27}(E_8)$
is not quasi-lisse.

\Needspace{6\baselineskip}
\subsection{\texorpdfstring{The remaining non-quasi-lisse algebras $L_k(E_6)$, $L_k(E_7)$ and $L_k(E_8)$}{The remaining non-quasi-lisse algebras Lk(E6), Lk(E7) and Lk(E8)}}
We treat
$$
 \begin{array}{c|l}
 \g&k\\\hline
 E_6&-9,-8,-7,-4\\
 E_7&-15,-14,-13,-11,-9\\
 E_8&-26,-25,-23,-22
 \end{array}
$$
For the reductions, we use the nilpotent elements in
Table~\ref{exc:negative-diagrams}. The last case for $E_6$, at $k=-4$,
follows directly from the known associated variety.
\subsubsection{\texorpdfstring{The levels $k=-9$ for $E_6$ and $k=-15,-14,-13$ for $E_7$}{The levels k=-9 for E6 and k=-15,-14,-13 for E7}}
For $L_{-9}(E_6)$, use the tensor basis $v_{ij\ell}$ of
$\C^3\otimes\C^3\otimes\C^2$, with indices beginning at zero, and set
$$
 f=v_{000}+v_{111}+v_{220}+v_{221},\qquad
 u_1=v_{221}-v_{001},\quad u_2=v_{220}-v_{110}.
$$
The vectors $u_1,u_2$ commute with $f$. Normalize the central element
$z\in\g_0$ to act by $-I$ on this tensor module, with $(z,z)=6$.
At $k=-9$, the two vertex algebras $V^0(\slg_3)$ associated with
$\g_0$ map to $L_0(\slg_3)=\C$. The remaining vertex algebra
$L_1(\slg_2)$ is realized by odd fields $b_i,c_i$, $i=1,2$, with
$b_i(z)c_j(w)\sim\delta_{ij}/(z-w)$. The center maps to
$-3I_2$ in this matrix realization, and the images of both weight-one
generators associated with the centralizer torus vanish. Write
$\mathcal Q(A)=\sum_{i,j}A_{ij}:\partial b_i c_j:$.
\Cref{com:quadratic-miura} and the finite projection in
\hyperref[calc:E6-LIFT]{Calculation~\ref*{calc:E6-LIFT}} give
$$
 \mu(W_{u_1})=\mathcal Q\!\begin{pmatrix}0&0\\1&-2\end{pmatrix},\qquad
 \mu(W_{u_2})=\mathcal Q\!\begin{pmatrix}-2&1\\0&0\end{pmatrix}.
$$
The image of the conformal vector is $\mathcal Q(I_2)$. These three matrices span
the self-adjoint matrices for
$G=\left(\begin{smallmatrix}1&-1/2\\-1/2&1\end{smallmatrix}\right)$.
The torus-fixed centralizer is the six-dimensional centralizer for
partition $(5,3)$ in $\sog_8$, generated by three degree-two directions.
Thus \cref{v3:abelianization,v3:quadratic-generation} identify the
image of $\mathcal W^{-9}(E_6,D_4(a_1))$ with the simple vertex algebra
$\SF(\C^2\otimes\C^2)^{\mathrm{SL}_2}$.
For $L_{-15}(E_7)$, \cref{com:quadratic-miura} gives the six complete
quadratic images in \hyperref[calc:CENTRALIZER]{Calculation~\ref*{calc:CENTRALIZER}}.
They show that the image of $\mathcal W^{-15}(E_7,E_7(a_5))$
contains the commuting tensor factor $\Vir_1$. The six degree-two
generators of the centralizer abelianization give vertex generation,
and \cref{vtwo:whole-factor} shows that this factor survives in the
simple quotient of $H^0_{E_7(a_5)}(L_{-15}(E_7))$.
For $L_{-14}(E_7)$ and $L_{-13}(E_7)$, the differentials of the
$G_0$-orbit maps at the chosen nilpotent elements are surjective by
the scalar and $2\times2$ blocks for the principal $\slg_2$ in
\hyperref[calc:TANGENT]{Calculation~\ref*{calc:TANGENT}}, of ranks $24,30$ and split diagonal
$\slg_2$ stabilizers. The even Dynkin grading, not openness alone,
identifies these with the required Richardson orbits. The bound on the generators and their parities given in
\hyperref[calc:E7-FORM]{Calculation~\ref*{calc:E7-FORM}} show that the images of
$\mathcal W^{-14}(E_7,A_4+A_2)$ and
$\mathcal W^{-13}(E_7,A_3+A_2+A_1)$ in these realizations are PCT invariant
and simple. In each case, the conformal vector has a nonzero quadratic
term in the generators of the Heisenberg vertex algebra and no
term linear in derivatives of those generators. Its class in the $C_2$ algebra is therefore not
nilpotent, so the corresponding $L_k(E_7)$ is not quasi-lisse.

\subsubsection{\texorpdfstring{Virasoro quotients of reductions of $L_{-11}(E_7)$ and $L_{-22}(E_8)$}{Virasoro quotients of reductions of L-11(E7) and L-22(E8)}}
At $k=-11$ for $E_7$, the homomorphism from the universal affine
vertex algebra of type $G_2$ to $\mathcal W^{-11}(E_7,A_2+3A_1)$
has level $2k+22=0$, and the
central charge is one, see the $A_2+3A_1$ row of
\cite[arXiv v3, Table 8]{AdamovicEtAl2026}. The full weight-one space lies
in the radical of the invariant bilinear form.
The centralizer abelianization has a single $h$-degree-two line, so
\cref{v3:abelianization} gives
$\mathcal W_{-11}(E_7,A_2+3A_1)\cong\Vir_1$.
The same argument applies to $L_{-22}(E_8)$: the homomorphism
$V^0(A_2)\to\mathcal W^{-22}(E_8,D_4(a_1)+A_2)$ has level zero, its abelianization again has one degree-two line, and its central
charge is $-2$, giving
$\mathcal W_{-22}(E_8,D_4(a_1)+A_2)\cong\Vir_{-2}$.
Since the corresponding reductions of $L_{-11}(E_7)$ and $L_{-22}(E_8)$
are nonzero, these simple quotients prove that both affine vertex
algebras are not quasi-lisse.

\subsubsection{\texorpdfstring{A Heisenberg quotient of $H^0_{A_2+A_1}(L_{-9}(E_7))$}{A Heisenberg quotient of the A2+A1 reduction of L-9(E7)}}
At the collapsing level $k=-9$, the simple $W$-algebra satisfies
$\mathcal W_{-9}(E_7,A_2)=L_{-1}(A_5)$
\cite[arXiv v3, Table 8, $A_2$ row]{AdamovicEtAl2026}.
We apply the reduction-by-stages theorem of Genra and Juillard
\cite[Theorem~5.3.4]{GenraJuillard2025}. The following grading check verifies its
hypotheses for the two nilpotent elements used here.
For the first triple,
$$
 E_7=(\ad A_5)\otimes V_0\oplus
       (\Lambda^2\C^6\oplus\Lambda^4\C^6\oplus\C)\otimes V_2
       \oplus V_4.
$$
Add a minimal triple in $A_5$ and write $H_2=H_1+H_0$ for the
doubled gradings, with $\g_{a,b}$ denoting the joint eigenspace of
$H_1,H_2$. On $\C^6$, $H_0$ has weights $1,-1,0,0,0,0$.
Its wedge-square and wedge-fourth weights are $0,\pm1$, and its
adjoint weights are $0,\pm1,\pm2$. Consequently
$\g^{(1)}_{\ge2}\subseteq\g^{(2)}_{\ge1}\subseteq\g^{(1)}_{\ge0}$,
$\g^{(1)}_1=0$, $\g^{(2)}_1\subseteq\g_{0,1}\oplus\g_{2,1}$,
and $f_0\in\g_{0,-2}$. The stage nilalgebra is
$\mathfrak a\oplus\C e_\theta$, where $\mathfrak a$ is a
four-dimensional Lagrangian in the eight-dimensional minimal
$A_5$ half-grade space. Thus \cite[Theorem 5.3.4]{GenraJuillard2025}
applies to the universal algebras and ordinary minimal $A_5$ reduction.
The natural affine module is locally finite and has finite-dimensional
homogeneous subspaces for the $W$-algebra grading, so exactness in $KL_{-1}(A_5)$ gives
$$
 \mathcal W^{-9}(E_7,A_2+A_1)\twoheadrightarrow
 H^0_{\min}(L_{-1}(A_5))\twoheadrightarrow\mathcal H.
$$
The middle algebra is nonzero: in simple-root coordinates
$\rho_{A_5}-5\theta^\vee/2=(0,3/2,2,3/2,0)$ gives strictly positive conformal-weight differences.
In its simple quotient, the images of the generators of $V^0(A_3)$
vanish, and their action forces every weight-$3/2$ generator to vanish.
The remaining weight-one generator, normalized by
$\operatorname{diag}(2,-1,-1,-1,-1,2)$, has level $24$. The remaining
coset conformal vector has central charge zero and trivial simple quotient.
Uniqueness of the simple graded quotient identifies the rank-one
Heisenberg vertex algebra $\mathcal H$ in the displayed surjection as
the simple quotient of $H^0_{A_2+A_1}(L_{-9}(E_7))$.

\subsubsection{\texorpdfstring{The reductions of $L_{-26}(E_8)$ and $L_{-25}(E_8)$}{The reductions of L-26(E8) and L-25(E8)}}
In $\mathcal W^{-26}(E_8,D_7(a_2))$, the unique weight-one generator has self-bracket
$(4k+104)\lambda\mathbf1=0$ and vanishes in the simple quotient. The positive abelianization degrees are $2,6$.
The degree-two generator is represented by the central element $f$:
its image is nonzero since the conformal vector has central charge $-4$ and cannot be
a normally ordered product of positive-weight vectors once the
weight-one subspace vanishes. The two Lie generators
therefore commute, and \cref{v3:abelianization} gives strong generators
of weights two and four. The Virasoro Gram form in weight four is
nondegenerate since $c(5c+22)=-8$ at $c=-4$. If the additional primary
vanishes, the simple quotient is $\Vir_{-4}$. Otherwise its norm is
nonzero and we normalize it by $W_{(7)}W=c/4$. The complete singular
products are then those of \cite[\S5.2.2, equations~(5.10)--(5.13)]{BSW}.
The denominator factors $c,c+24,5c+22,2c-1,7c+68$ are nonzero at $c=-4$,
and the squared structure constant is $(C_{44}^4)^2=-1350$. The simple algebra $\mathcal A_2$ has the
same normalized presentation, as follows from \cite[Theorem~3.9]{CLSF}.
Uniqueness of its graded simple quotient identifies the second alternative.
\hyperref[calc:W24]{Calculation~\ref*{calc:W24}} checks these coefficient
and denominator specializations. The complete presentation is the cited input.

For $\mathcal W^{-25}(E_8,E_8(a_7))$, the centralizer abelianization has ten degree-two lines.
The two complete fields obtained from \cref{com:quadratic-miura} and
\hyperref[calc:QUAD]{Calculation~\ref*{calc:QUAD}} have,
up to a common sign, self-adjoint matrices
$$
 D=\operatorname{diag}(3,0,2,5/2),\qquad
 M=-E_{13}+2E_{31}+2E_{34}+E_{43}-E_{12}+10E_{21}
$$
for $G=\operatorname{diag}(-1,1/10,1/2,1)$. The spectrum of $D$ is
simple and the graph of $M$ has edges $13,34,12$. Lemma~\ref{v3:quadratic-generation} identifies the image generated
by these quadratic fields with
$\SF(\C^4\otimes\C^2)^{\mathrm{SL}_2}$. The two fields include all derivative corrections.
Lemma~\ref{vtwo:torus-detector} proves nonnilpotence of the conformal class in these invariant
vertex algebras and in the quotients used for $L_{-9}(E_6)$ and $L_{-15}(E_7)$.

\subsubsection{\texorpdfstring{The quotient $\mathcal W_{-23}(E_8,A_4+A_2+A_1)$}{The quotient W-23(E8,A4+A2+A1)}}
For $f\in\OO_{A_4+A_2+A_1}$, write $K_j=(\g^f)_{-j}$ for the Dynkin $h$-grading
and $V_m$ for the natural centralizer $A_1$ module of highest weight $m$.
The complete graded decomposition is
$$
\begin{array}{c|ccccccccc}
j&0&1&2&3&4&5&6&7&8\\\hline
K_j&V_2&V_5&V_4+2V_0&V_5+V_1&V_6+V_2&V_3+V_1&V_4&V_3&V_2.
\end{array}
$$
The fourteen highest-weight vectors and their graded character give
this decomposition.
For a highest vector $u\in K_1$, the pairing of the complete highest- and lowest-weight fields is
$$
 G(u)_{(2)}G((\ad F)^5u)=-240(k+23)(3k+70)\mathbf1.
$$
The complete Miura formula has degree at most two in $k$: neutral
contractions are constant, while the invariant bilinear forms on the
weight-one subspaces and the coefficients of derivatives are affine in $k$. Three exact evaluations therefore determine the displayed
pairing. It is $A_1$-invariant. The full weight-one space consists
of primary vectors, and there is no weight $1/2$. Hence every
weight-$3/2$ field is quasiprimary. Irreducibility and
\cref{v5:half-pairing} imply that the generators corresponding to $V_5$
vanish in $\mathcal W_{-23}(E_8,f)$. This uses the complete field, including derivative and neutral
corrections.

The eleven equivariant bracket maps in \hyperref[calc:V5-E8-P7]{Calculation~\ref*{calc:V5-E8-P7}}, with
one argument fixed to $u$, generate $V_4+\C z$ in degree two and every component
of $K_j$ for $j\ge3$. The Clebsch--Gordan coefficients and the fourteen highest-weight
vectors are given there over the rationals.
Since $\Lambda^2V_5=V_8+V_4+V_0$, the degree-two bracket image is
$V_4+\C z$, with only one scalar direction. In the root-vector
normalization of \hyperref[calc:V5-E8-P7]{Calculation~\ref*{calc:V5-E8-P7}},
the matrix whose columns are the $e_{-\alpha_1},e_{-\alpha_2}$ coordinates
of $z,f$ is $\left(\begin{smallmatrix}1&-1\\-6&-1\end{smallmatrix}\right)$.
Its determinant is $-7$, so $z$ and $f$ are linearly independent. All degrees at least
three are then eliminated by \cref{v5:equivariant-generation}. Grading
shows the displayed subspace is an ideal. Thus
$K/\langle K_1\rangle=A_1\oplus\C f$.
The simple quotient is therefore generated by the image of
$V^5(A_1)$ and the conformal vector.
Their total and residual central charges are $22/7$ and $1$, so the whole simple
quotient is $L_5(A_1)\otimes\Vir_1$. The conformal-weight argument establishes nonvanishing independently.
The $\Vir_1$ tensor factor shows that $L_{-23}(E_8)$ is not quasi-lisse.

\subsubsection{\texorpdfstring{Non-quasi-lisseness of $L_{-8}(E_6)$, $L_{-7}(E_6)$ and $L_{-4}(E_6)$}{Non-quasi-lisseness of L-8(E6), L-7(E6) and L-4(E6)}}
We treat these levels using sheet closures and their reductions.
At $k=-8$, \cite[Theorem~3.5]{ArakawaFutornyKrizka2024} gives the $A_3$
sheet bound, whose nilpotent intersection is the Richardson orbit
closure labelled $A_4$. This bound, combined with the obstructions from affine vertex subalgebras at negative levels, leaves only $D_4(a_1)$ and $A_4$ under quasi-lisseness.
Both contain $D_4(a_1)$, so its reduction is nonzero. The rank-two Heisenberg vertex algebra is generated by $J_1,J_2$,
with $[J_i{}_{\lambda}J_j]=K_{ij}\lambda\mathbf1$, where
$K=\left(\begin{smallmatrix}4&-2\\-2&4\end{smallmatrix}\right)$.
The centralizer of the $\slg_2$-triple is $T^2\rtimes S_3$. Its weight-two invariants
are the two conformal vectors, whose difference lies in the radical
of the invariant bilinear form. The charged
generators have $J_i(0)$-weights $q\in\{\pm(1,1),\pm(2,-1),\pm(1,-2)\}$
and conformal weight $\tfrac12q^{\mathsf T}K^{-1}q=1/2$
for $T_{\mathcal H}=\tfrac12\sum_{i,j}(K^{-1})_{ij}:J_iJ_j:$,
incompatible with the integral grading. They vanish, leaving
$\mathcal H^2$, a contradiction.

At $k=-7$, the $A_4$ sheet bound from the same theorem leaves
$A_3$ or $A_2+2A_1$. If
$A_3$ occurs, its reduction is nonzero. The image of $V^0(B_2)$ is trivial, and its action forces the generators
in nontrivial $B_2$-modules to vanish and leave the rank-one Heisenberg vertex algebra associated with the
center of the reductive centralizer, a contradiction. Thus under quasi-lisseness the maximal nilpotent
boundary is $A_2+2A_1$. The intersection of its Slodowy slice with
$X_{L_{-7}(E_6)}$ is zero-dimensional. No irreducibility assumption is needed.

Let $B$ be the simple graded quotient of the nonzero reduction at
$f\in A_2+2A_1$. Under the hypothesis of quasi-lisseness, this reduction,
and hence $B$, is lisse. Write $V_m$ for the $(m+1)$-dimensional
$\slg_2$-module, and use the second subscript below for the eigenvalue
of $J_0$, where $J$ spans the centre of $\g^\natural$. The strong generators
indexed by $\g^f$ have the following types:
$$
\begin{array}{c|l}
\text{conformal weight}&\slg_2\oplus\C J\text{-module}\\\hline
1&V_{2,0}\oplus V_{0,0}\\
3/2&V_{3,3}\oplus V_{3,-3}\\
2&V_{4,0}\oplus V_{2,0}\oplus V_{0,0}\\
5/2&V_{1,3}\oplus V_{1,-3}\\
3&V_{2,0}.
\end{array}
$$
These finite-dimensional data can be obtained directly from the
root decomposition. Choose $f$ principal in the standard Levi subalgebra
on Bourbaki nodes $\{1,2,3,6\}$. In simple-coroot coordinates, the
semisimple element of its $\slg_2$-triple and the two Cartan elements
of $\g^\natural$ are
$$
 \begin{gathered}
 h=(2,1,2,0,0,1),\\
 H=(-2,-3,-4,-6,-6,-3),\qquad J=(2,3,4,6,2,1).
 \end{gathered}
$$
Here $H$ is the coroot of the natural $\slg_2$, and $J$ is central.
For the Dynkin grading defined by $h/2$, the $\g^\natural$-characters
satisfy
$$
 \ch(\g^f\cap\g_{-n/2})
   =\ch\g_{-n/2}-\ch\g_{-(n+2)/2}\qquad(n\ge0),
$$
which gives the table above. Moreover, $(H,H)=(J,J)=12$ and
$(H,J)=0$. The traces of $(\ad H)^2$ on $\g_0,\g_{1/2}$ are
$64,44$, and those of $(\ad J)^2$ are $0,108$. The form on $\g^\natural$ is
$$
 K(a,b)=(k+12)(a,b)
 -\tfrac12\tr_{\g_0}(\ad a\,\ad b)
 -\tfrac12\tr_{\g_{1/2}}(\ad a\,\ad b),
$$
so $K(H,H)=K(J,J)=6$ and $K(H,J)=0$ at $k=-7$.

The induced $\slg_2$ level is $6k+45=3$ at $k=-7$
(see also \cite[Table~13]{AVEMSlodowy2024}). In the normalization of the
displayed charges, $J_{(1)}J=6\mathbf1$. Since $B$ is lisse,
\cref{com:bounded-current} identifies the vertex subalgebra generated by $B_1$ with
$A=\mathcal H_6\otimes L_3(A_1)$.
The dominant weighted Dynkin diagram is $001010$, with
$\dim\g_0=16$ and $\dim\g_{1/2}=12$. The Drinfeld--Sokolov central-charge
formula gives $c(B)=14/5=1+9/5=c(A)$. The table shows that
the only $\g^\natural$-invariant generator of weight two is the one
indexed by $f$, represented by the conformal vector modulo products of weight-one generators. Thus $(B_2/A_2)^{\g^\natural}$ is spanned by the image of $\omega_B$,
and \cref{v3:resonance} gives
$\omega_B=\omega_A$. This verifies the conformal embedding needed
for the following fusion argument.

Suppose one of the two weight-$3/2$ generating spaces survives.
Its positive modes for $\g^\natural$ vanish for degree reasons, so it generates
an $A$-submodule
$M_\pm=F_{\pm3}\otimes L_3(3)$, with lowest conformal weight
$3/4+3/4=3/2$. Here $F_q$ denotes the $\mathcal H_6$ Fock module
of charge $q$. The strong generators are $J_0$-eigenvectors, so $J_0$ is semisimple on
the finite-dimensional graded pieces of $B$, and its positive modes act locally nilpotently.
The usual Fock decomposition therefore applies. On the vectors annihilated by $J(n)$ for $n>0$ and satisfying
$J(0)v=qv$, the $A_1$ Sugawara grading is
$L_0-q^2/12$, bounded below with finite-dimensional pieces. Complete
reducibility of integrable $\widehat{\slg}_2$-modules at level three
\cite[Proposition~5.6]{DLMSimpleCurrents1996} now makes $B$ a direct
sum of simple $A$-modules. The module $L_3(3)$ is a
self-dual simple current \cite[Proposition~5.1]{DLMSimpleCurrents1996}.
Together with the fusion rules for the Fock modules of $\mathcal H_6$
\cite[Section~2.2]{CKLRHeisenberg2019}, this gives
$$
 M_\pm\boxtimes_A M_\pm
   =F_{\pm6}\otimes L_3(0),
$$
with fusion multiplicity one. This module has lowest conformal weight
$3$, so the lowest coefficient of such an intertwining operator is
of degree $z^0$. It is a scalar multiple of the invariant pairing
$V_3\otimes V_3\to\C$.
Moreover, that lowest coefficient determines the intertwining operator:
the intertwining operator between the $\mathcal H_6$ Fock modules
has nonzero leading coefficient, and the intertwining operator of
type $\binom{L_3(\slg_2)}{L_3(3)\ L_3(3)}$ is obtained from the
nondegenerate contragredient pairing on $L_3(3)$. Their tensor product
spans the one-dimensional intertwining space.
The invariant pairing on $V_3$ is alternating, whereas bosonic skew
symmetry requires this degree-zero coefficient to be symmetric.
It therefore vanishes, and so does the entire intertwining operator.
This argument applies to the projection onto every copy of this
module. It does not require a multiplicity-one assertion for the
branching of $B$. Thus $Y(M_\pm,z)M_\pm=0$, contradicting the absence
of zero divisors in a simple vertex algebra
\cite[Corollary~4.3]{DAndrea2008}. Both weight-$3/2$ generating spaces
must vanish.

The only other charged strong generators have weight $5/2$ and type
$V_{1,\pm3}$. Their positive modes for $\g^\natural$ now vanish, since all
lower-weight spaces have charge zero. The Sugawara conformal weight
on their tops would be
$$
 \frac{1(1+2)}{4(3+2)}+\frac{3^2}{2\cdot6}
 =\frac9{10}\ne\frac52,
$$
so these generators vanish as well. All remaining strong generators
are annihilated by $J_0$. Hence $J_0=0$ on $B$, and
\cref{com:heisenberg} gives $\mathcal H_6$, generated by $J$, as a tensor factor of $B$. Its polynomial $C_2$ algebra contradicts lisseness.

Finally, $L_{-4}(E_6)$ is already treated by Arakawa, Futorny and K\v{r}i\v{z}ka \cite[Theorem~3.24 and Table~6]{ArakawaFutornyKrizka2024}. At $k=-4$,
the associated variety is the $D_5$ Dixmier-sheet closure of dimension
$33$. In particular, it is not contained in the nilpotent cone.

\subsection{The nonadmissible quasi-lisse levels}\label{exc:seeds}
The table lists the nonadmissible numerators giving quasi-lisse vertex
algebras. \Cref{exc:negative} and the subsequent subsections prove that
the other nonadmissible numerators do not give quasi-lisse algebras.

\begin{table}[htbp]\centering\small
\begin{tabular}{@{}lll@{}}\toprule
Type&quasi-lisse numerators below $h^\vee$&$p\ge2$ for which quasi-lisseness fails\\\midrule
$E_6$&$6,7,9,10,11$&$2,3,4,5,8$\\
$E_7$&$6,8,10,11,12,13,14,15,16,17$&$2,3,4,5,7,9$\\
$E_8$&$6,9,10,\ldots,29$&$2,3,4,5,7,8$\\\bottomrule
\end{tabular}
\par\smallskip
\caption{The nonadmissible numerators in types $E_6,E_7,E_8$.}
\label{exc:coverage}\end{table}

For the $37$ classes giving quasi-lisse algebras and the $17$
classes where quasi-lisseness fails, take $t=p$ except
at $(E_6,7),(E_7,11),(E_8,19)$, where $t=p/2$ shortens the proof.
Associated-variety bounds refer to these representatives.
We begin with known associated varieties and cases where an upper
bound on the nilpotent intersection rules out every semisimple family.
When an induced nilpotent orbit attains that bound, a dimension
comparison alone is insufficient. In those cases we prove that the
corresponding reduction has a point associated variety, and use its
Slodowy slice to show that the remaining family cannot occur.
Set $\delta_x(\lambda)=(\lambda,\lambda+2\rho)/(2t)-\lambda(x)$
for $x=h/2$. By \cref{com:completion,com:vacuum-test}, the constant terms of the
complete reduced character sums in
\hyperref[calc:EXCEPTIONAL]{Calculation~\ref*{calc:EXCEPTIONAL}}
determine whether the corresponding reductions vanish. The sums are
formed from characters of standard modules before specialization.

\subsubsection{Associated varieties from minimal reductions and Dixmier sheets}
\label{exc:positive-geometric}
Thirteen levels are already covered by Arakawa and Moreau's theorem on
minimal-orbit associated varieties \cite[Theorem~1.1(2)]{ArakawaMoreau2018}:
$L_k(E_6)$ at $k=-3,-2,-1$, $L_k(E_7)$ at $k=-4,-3,-2,-1$,
and $L_k(E_8)$ at $k=-6,-5,\ldots,-1$.
The minimal Drinfeld--Sokolov reductions
$$
 H_{f_{\min}}^0(L_{-6}(E_7))=L_{-2}(D_6),\qquad
 H_{f_{\min}}^0(L_{-10}(E_8))=L_{-4}(E_7)
$$
are nonzero and simple by \cite[Theorems~6.1--6.2]{ArakawaMoreau2017}.
The identifications with affine vertex algebras are the collapsing-level
results of Adamovi\'c, Kac, M\"oseneder Frajria, Papi and Per\v se
\cite[Theorem~3.3 and Table~4]{AdamovicEtAl2018}. The nilpotent dimensions are at most $34+18=52$ and
$58+34=92$. By conic induction (\cref{vtwo:conic-boundary}), every
nonnilpotent orbit gives rise to a nilpotent orbit of the same dimension
in the boundary of its conical saturation.
The least possible nonnilpotent orbit dimensions are respectively $54$
and $114$, so these bounds show that the two associated varieties
are contained in the nilpotent cone. Translation by integrable affine
vertex algebras, \cref{com:translation}, gives the corresponding bounds
for $L_k(E_7)$ at $k=-6,-5$ and $L_k(E_8)$ at $k=-10,-9,-8,-7$.

For $L_{-6}(E_6)$, \cite{ArakawaFutornyKrizka2024}, Theorem~3.5, gives both the $A_5$
sheet bound on nodes $1,3,4,5,6$ and the $D_4$ sheet bound on nodes
$2,3,4,5$. A nonnilpotent point of the former has centralizer
$A_5+\C$ and zero nilpotent part. If such a point also belonged to the $D_4$ sheet closure,
its centralizer would contain a subalgebra of type $D_4$. A nonzero $D_4$ representation has dimension
at least eight, whereas the defining $A_5$ representation has dimension
six. Projection to the abelian center is zero. Such an embedding is
impossible. Their intersection is therefore contained in $\N(E_6)$,
and so is $X_{L_{-6}(E_6)}$.

\subsubsection{Bounds on the nilpotent orbits and absence of nonzero semisimple points}
\label{vtwo:exceptional-caps}
The following table records unconditional bounds on $X\cap\N$, obtained
without assuming quasi-lisseness. For each orbit in the third column,
$H_f^0(L_k(\g))=0$. The closure order on nilpotent orbits gives the
bound in the last column.
\begin{center}\small
\begin{tabular}{@{}llp{.43\textwidth}r@{}}\toprule
Type&$t$&forbidden orbit(s)&bound\\\midrule
$E_6$&$7/2$&$A_4,D_4$&58\\
$E_7$&$8$&$A_3,2A_2,A_2+3A_1$&82\\
$E_7$&$11/2$&$A_4,D_4$&100\\
$E_8$&$10$&$A_3+A_2$, and $D_4$ by the level of $V^{-2}(F_4)$&176\\
$E_8$&$14$&$A_3,A_2+3A_1$&146\\
$E_8$&$15$&$A_2+2A_1$&136\\
$E_8$&$19/2$&$D_4(a_1)+A_2,D_4+A_1$&182\\\bottomrule
\end{tabular}\end{center}
The finite Weyl-group calculation, with the bounds in
Appendix~\ref{cert:scope} and the data in \hyperref[calc:ROOT-WINDOW]{Calculation~\ref*{calc:ROOT-WINDOW}} and
\hyperref[calc:WEYL-PROD]{Calculation~\ref*{calc:WEYL-PROD}}, gives constant term one for the reductions treated by character calculation. The $D_4$ entry at $E_8,t=10$ is ruled out separately below using the level of its affine vertex subalgebra. The truncated character sums $(E_6,7/2,D_4)$, $(E_7,11/2,D_4)$ and
$(E_8,19/2,D_4+A_1)$ have respectively $19,41,119$ distinct standard characters after collection. The same finite Weyl-group formula applies to all three sums. For $E_7,t=8$ and $E_8,t=14$, the same calculations determine
the coefficient of conformal degree zero directly from the Weyl
alternants and the denominator factors in
\eqref{com:standard-ds-character}. Terms of higher degree do not
contribute to this coefficient.

At $E_8,t=10$, the full closure order gives
$$
 \{\OO:A_3+A_2\notin\overline\OO\}
 \setminus\{\OO:\OO\subseteq\overline\OO_{D_4(a_1)+A_1}\}
 =\{\OO_{D_4}\}.
$$
Moreover $D_4$ is maximal in this allowed set. Were it in $X$, the
maximal-boundary lemma would apply to its reduction. Its
homomorphism $V^{\ell}(F_4)\to\mathcal W^k(E_8,D_4)$ has embedding
index one and level
$$
 \ell=(k+30)-\tfrac12(18+6)=k+18=-2,
$$
contradicting the integrability forced by bounded Cartan weights in
\cref{com:bounded-current}. Here the grade-zero
$E_6$ adjoint restricts as $52+26$, whose trace indices are $18,6$.
Here \cite[arXiv v3, Table~13, row~14]{AdamovicEtAl2026} is used only for
the centralizer and embedding data.

Suppose that $\dim(X\cap\N)\le c$. To prove $X\subseteq\N$, apply
the criterion for semisimple points and choose an orbit $Y_j$ in the
associated variety of the simple affine vertex algebra for each derived
Levi subalgebra. The corresponding induced nilpotent orbit has dimension
$\dim\g-\dim\lev_j+\dim Y_j$. A nonzero point of the corresponding
Levi center in $X$ would force this nilpotent orbit into $X\cap\N$.
If its dimension exceeds $c$, this contradicts the bound.
The following lists, in Bourbaki node order, make all inequalities strict:
\begin{center}\small
\begin{tabular}{@{}llp{.62\textwidth}@{}}\toprule
Type, $t$&bound&lower bounds for induced orbit dimensions\\\midrule
$E_6,7/2$&58&$62,60,64,68,64,62$\\
$E_7,8$&82&$84,84,94,106,100,84,86$\\
$E_7,11/2$&100&$102,108,114,120,116,110,104$\\
$E_8,10,11$&176&$178,184,196,212,208,194,188,180$\\
$E_8,14,15,16,17$&$146$&$156,184,196,212,208,194,166,148$\\
$E_8,19/2$&182&$192,184,196,212,208,194,188,184$\\\bottomrule
\end{tabular}\end{center}
Here the bound $176$ is translated from $t=10$ to $11$, and the bound $146$
from $14$ to $15,16,17$. The sharper $136$ at $15$ translates to $16,17$.
These are bounds at the stated representatives.

We justify independently the nilpotent-orbit inclusions in the associated
varieties of the simple affine vertex algebras for the Levi factors.
In $E_6,t=7/2$ the $A$ factors use admissible square-zero
orbits. The two $D_5$ factors use partition $(3^3,1)$, dimension $30$.
Its fundamental conformal-weight differences are
$(2/7,2/7,0,3/28,3/28)$, and its sole nonzero equality weight
$\omega_3$ is outside the $D_5$ root lattice. The conformal-weight inequality for the maximal ideal therefore proves nonvanishing of the reduction. In $E_7,t=8$, only the minimal
$D_6$ orbit (dimension $18$) and $2A_1$ in $E_6$ (dimension $32$)
are needed. Their conformal-weight differences are respectively
$
 (3/16,1/4,11/16,1,17/32,17/32)$ and 
 $(1/12,1/2,7/12,1,7/12,1/12).$ 
They are strictly positive. At $E_7,t=11/2$, the classical factors are
admissible: take square-zero $A$ orbits and the $D_6,D_5$ partitions
$(3,2^4,1)$ and $(3,2^2,1^3)$ of dimensions $36,24$.
For the $E_6$ factor, the orbit $A_2+2A_1$, of dimension $50$, has strict conformal-weight differences
$(5/66,2/11,1/33,4/11,1/33,5/66)$.

At $E_8,t=10,11$, nonintegrability places the minimal $D_7$ and
$E_6$ orbits, both of dimension $22$, in the corresponding associated varieties, and the cancellation of two standard characters in the constant term in \hyperref[calc:ROOT-WINDOW]{Calculation~\ref*{calc:ROOT-WINDOW}} proves the $E_7,A_2$
nonvanishing at the orbit of dimension $66$. At $t=14,15,16,17$, only the minimal
$E_7$ orbit, dimension $34$, is needed: its levels are $-4,-3,-2,-1$,
so the general minimal-orbit containment theorem
\cite[Proposition~4.7.2]{AAV} places it in the associated variety. This argument is independent of
the classification in type $E_7$. At $t=19/2$\label{exc:e8p19}, the $D_7$ square-zero orbit
$(2^4,1^6)$ has dimension $36$, the minimal $E_6$ orbit has dimension
$22$, and we use an $E_7$ orbit of dimension $70$ with adjoint decomposition
$$
 \ad(E_7)=6V_3+16V_2+20V_1+21V_0.
$$
It is obtained in $D_6+A_1$ from natural $D_6$ module
$V_2+2V_1+5V_0$ and a regular external $A_1$. The half-spin module is
$4V_2+8V_1+4V_0$. These heights are at most three, so the admissible
$q=2$ associated-variety theorem places all three orbits in the respective
associated varieties. The six rows of strict
inequalities now show that no maximal Levi center contains a nonzero
point of $X_{L_k(\g)}$ at the stated levels.

\subsubsection{The levels \texorpdfstring{$k=-12$ for $E_7$ and $k=-24$ for $E_8$}{k=-12 for E7 and k=-24 for E8}}
\label{vtwo:equivariant-exceptional}
At numerator six, the preceding dimension estimates leave some
semisimple families possible: some induced boundaries attain the bound. We show that
the relevant reductions have a one-point associated variety, which rules out
these remaining families by their Slodowy intersections.

\subsubsection{\texorpdfstring{The reduction $H^0_{D_4(a_1)}(L_{-12}(E_7))$}{The D4(a1) reduction of L-12(E7)}}
Take $f=D_4(a_1)$, diagram $0020000$, dimension $94$, and write
$U=H_f(V^{-12}(E_7))$, $I_r=H_f(J^r)$. A complete mode bound $a\le5$
gives fourteen regular root contributions through reduced degree two. After three
cancellations of standard characters their minima are one or two. Thus the
reduction is nonzero. The negative-degree subspaces for the action of $\g_0$ are
$$
 \g_{-1}=W\otimes\Lambda^2 V,\quad
 \g_{-2}=\det W\otimes\Lambda^4V,\quad
 \g_{-3}=\det W\otimes W\otimes\det V,
 \qquad \dim W=2,\ \dim V=6.
$$
Write $V=V_1\oplus V_2\oplus V_3$, with $\dim V_i=2$ and area forms
$w_i$. The pencil
$f=u\otimes w_1+v\otimes w_2+(u+v)\otimes w_3$ has an open
orbit under the degree-zero Levi subgroup: its connected stabilizer is
$A=\prod_i\mathrm{SL}(V_i)$, and its component quotient is $S_3$,
permuting the three distinct Pfaffian roots. Hence it realizes the
stated Dynkin orbit.

The $S_3$-module decomposition and the action of positive modes of
$\g_0^f$ in
\hyperref[calc:EQUIVAR]{Calculation~\ref*{calc:EQUIVAR}} give
$$
 \sum_r[(I_r)_1]=\sum_{i=1}^3[\ad_i],\qquad
 U_2^A=P_0\oplus P,\quad P_0\cong P\cong\mathbf1\oplus\mathrm{Std},
 \qquad \sum_r\dim(I_r)_2^A=6.
$$
The degree-one multiplicities imply $U_1\subset I_1$ and
$(I_r)_1=0$ for $r\ge2$. Here $U_1\cong\slg_2^{\oplus3}$, and
each associated homomorphism $V^0(\slg_2)\to U$ has level zero. Their three Sugawara vectors span $P_0$.
The joint action of positive modes of $U_1$ is injective on $P_0$
and has kernel $P$ in $U_2^A$.
Later ideals therefore lie in $P$. The conformal vector, modulo $P_0$, is the
unique invariant line in $P$. If it were absent, every residual layer
would be a sum of two-dimensional standard representations, whereas
$6-\dim P_0=3$ is odd. By \cref{v3:scalar-detection}, the conformal vector lies
in $I_1$. By \cref{com:scalar-point}, the reduction is $\C$. 

For $f'\in\OO_{(A_3+A_1)''}$, the reduced Jantzen sum for
$V^{-12}(E_7)$ has constant coefficient one. Its three signed
standard-module contributions, indexed by $(a,D,\lambda)$, are
$$
 \begin{array}{c|c|c|c}
 a&D&\lambda&\text{sign}\\\hline
 1&3&\omega_1&-\\
 1&5&2\omega_7&+\\
 3&3&\omega_1&+
 \end{array}
$$
The first and third contributions cancel. The remaining term has
constant coefficient one. Here $a$ is the positive integer indexing
the affine root and $D$ is the original degree. The complete
evaluation is recorded in
\hyperref[calc:ROOT-WINDOW]{Calculation~\ref*{calc:ROOT-WINDOW}}.
Thus \cref{com:vacuum-test} gives
$H_{f'}^0(L_{-12}(E_7))=0$, so this orbit is absent from the
associated variety. The reduction at $f\in\OO_{D_4(a_1)}$
was shown above to be $\C$. Its associated variety is a point,
which also rules out the larger orbit $D_4$ by
\cref{vtwo:point-slice}.
The closure order shows that $(A_3+A_1)''$ and $D_4$ are precisely
the minimal orbits outside $\overline\OO_{D_4(a_1)}$.
Consequently $\dim(X\cap\N)\le94$, with $D_4(a_1)$ the unique
possible orbit of that dimension. The seven Levi lower bounds are
$$
 \dim\lev_j=(67,49,39,27,33,49,79),\quad
 \dim Y_j=(28,12,0,0,0,14,42),
$$
giving induced dimensions $(94,96,94,106,100,98,96)$.
The $D_6$ orbit $(2^4,1^4)$ has positive conformal-weight differences
$(5/12,2/3,3/4,2/3,3/8,3/8)$. Nonintegrability gives the minimal
$A_6,D_5$ orbit inclusions. The two-term $E_6,A_2$ cancellation proves the required nonvanishing at dimension
$42$. The strict cases are impossible. In each equality case the center
family has dimension $95$, with $D_4(a_1)$ boundary. The intersection
of $X$ with the Slodowy slice at $f\in\OO_{D_4(a_1)}$ is therefore
positive-dimensional. This contradicts the one-point associated variety
of the reduction by \cref{vtwo:point-slice}.

\subsubsection{\texorpdfstring{The reduction $H^0_{E_8(a_7)}(L_{-24}(E_8))$}{The E8(a7) reduction of L-24(E8)}}
Now $f=E_8(a_7)$, diagram $00002000$, dimension $208$.
The centralizer of the $\slg_2$-triple is $S_5$, acting on $\g_{-5}$ as its
standard four-dimensional representation $S$. This uses
\cite{ElashviliJibladzeKac2022}, Table~1, row~18 and Example~2.12: the reducing
algebra generated by the triple and $\g_{-5}$ is all $\g$, so no
additional kernel acts trivially there. For the degree-zero Levi subalgebra with derived algebra $A_4+A_3$,
let $V,W$ be modules of dimensions five and four. Its group is
$(\mathrm{GL}(V)\times\mathrm{GL}(W))/K$, where
$K=\{(sI_5,s^{-2}I_4)\}$. The character $(\det V)^5(\det W)^3$
restricts to $s$ on $K$, so the finite group lifts. From
$$
 \g_{-1}=\Lambda^2V\otimes W,\qquad
 \g_{-2}=\Lambda^4V\otimes\Lambda^2W,\qquad
 \g_{-5}=(\det V)^2\det W\otimes W
$$
and $\g_0\cong\g_{-1}$ as $S_5$ modules, one obtains
$W=S\otimes\mathrm{sgn}$ and $V=V_{(3,2)}$ up to tensoring with a character of $K$.
The seven character classes give
$$
 \dim U_2^{S_5}=1,\qquad
 \sum_r\dim(I_r)_2^{S_5}=2-1+0-0=1.
$$
Here the complete mode bound $a\le7$ gives four nonzero standard characters
$(D,\lambda)=(22,4\omega_8),(27,2\omega_7+\omega_8),
(30,\omega_1+\omega_6+\omega_8),(31,\omega_4+\omega_8)$,
with signs $+,-,+,-$ and reduced minimum two. Thus no lower ideal term
occurs and the reduction is nonzero. The unique invariant line is the nonzero
universal conformal vector. Nestedness places it in the first reduced ideal,
so \cref{com:scalar-point} gives $H_f(L_{-24}(E_8))=\C$.

The levels of the affine vertex subalgebras, with the trace corrections
in \hyperref[calc:E8-IDEAL]{Calculation~\ref*{calc:E8-IDEAL}}, and the
closure order on nilpotent orbits leave $24$ possibilities for a maximal
nilpotent orbit in $X$, each containing $E_8(a_7)$ in its closure.
Every strictly larger candidate would contradict the zero-dimensional Slodowy slice, giving
$X\cap\N=\overline\OO_{E_8(a_7)}$. The eight Levi dimensions are
$(92,64,52,36,40,54,82,134)$. The associated varieties of the simple affine
vertex algebras for these Levi subalgebras contain nilpotent orbits of
dimensions $(52,24,12,0,0,14,42,94)$. The induced nilpotent orbits have
dimensions $(208,208,208,212,208,208,208,208)$. The required orbit
inclusions for $E_7,E_6$ were
proved just above. The $A_7$ orbit $(2^2,1^4)$ satisfies the strict conformal-weight inequalities.
For the $D_7$ orbit $(3^2,2^2,1^4)$, the sole equality is the adjoint
$\omega_2$ at Sugawara degree two. Its nonzero level is $-6$, and
\cref{vtwo:current-descendant} excludes this adjoint submodule from the
maximal ideal of $V^{-6}(D_7)$. The remaining $A_6,D_5$ minimal orbits follow from
nonintegrability. An induced orbit of dimension $212$ contradicts the
bound $\dim(X\cap\N)=208$. In each equality case, the family of
conjugacy classes obtained from the Levi center and the chosen nilpotent
orbit has a positive-dimensional Slodowy intersection, contradicting
the fact that the intersection of $X$ with the Slodowy slice at
$f\in\OO_{E_8(a_7)}$ consists of a single point, as in the $E_7$ argument.

\subsubsection{\texorpdfstring{The associated varieties of $L_{-8}(E_7)$ and $L_{-12}(E_8)$}{The associated varieties of L-8(E7) and L-12(E8)}}
\label{vtwo:exceptional-virasoro}
We prove $X_{L_{-8}(E_7)}=\overline{\OO}_{A_2}$ and
$X_{L_{-12}(E_8)}=\overline{\OO}_{A_2}$. 
Let $I_d$ be generated by the singular module of least positive degree
and set $Q_d=V/I_d$. Its reduction is a Virasoro quotient. The additional relation needed for $L_k(\g)$ is given in the following table.
\begin{center}\small
\begin{tabular}{@{}lcc@{}}\toprule
 & $E_7,t=10$ & $E_8,t=18$\\\midrule
$k$, $d$, first highest weight&$-8,3,2\omega_7$&$-12,4,\omega_2$\\
forbidden orbit(s) for $Q_d$&$(3A_1)'',4A_1$&$4A_1$\\
maximal surviving orbit $f$&$A_2$&$A_2$\\
bound on the nilpotent intersection for $Q_d$&66&114\\
affine vertex algebra at $f$&$A_5$ at level 0&$E_6$ at level 0\\
charged weight-two types&$15+15^*$&$27+27^*$\\
Virasoro central charge&$-22/5$&$-46/3$\\
proof of the relation&a nonzero mode coefficient&a composition factor\\
reduced relation degree&$4$&$8$\\
relation in the $C_2$ algebra of the reduction&$[T]^2=0$&$[T]^4=0$\\\bottomrule
\end{tabular}\end{center}
The assertion concerning the singular module of least positive degree follows from a complete ordinary
Jantzen sum with one nonzero root contribution, given in
\hyperref[calc:ROOT-WINDOW]{Calculation~\ref*{calc:ROOT-WINDOW}}. This determines the singular module in each case. At $E_8,4A_1$ the fundamental conformal-weight differences are
$$
 (1/6,0,4/9,5/6,2/3,1/2,1/3,1/6).
$$
Positive inverse-Cartan entries give equality only at $0,\omega_2$.
The latter has original degree four. At the two $E_7$ gradings the
conformal-weight differences are
$$
 (4/5,9/8,8/5,12/5,13/8,4/5,-3/40),\quad
 (3/10,1/8,3/5,9/10,5/8,3/10,-3/40).
$$
The negative difference cannot be dropped: along $m\omega_7$ the difference
is $3m(m-2)/40$. The sole negative case $\omega_7$ is outside the
inherited root lattice. Adding any other fundamental weight gives a
positive difference, including all cross terms. Equality among allowed
weights is therefore only $0,2\omega_7$, with nonzero equality degree
three. In each case, the subspace of maximal $x$-eigenvalue in the first
singular module is one-dimensional. Lemma~\ref{vtwo:scalar-vacuum} proves the asserted vanishings for $Q_d$
itself. The closure order gives the bounds in the table.

Cancellation of two standard characters independently proves nonvanishing at $A_2$
of $L_k(\g)$ at the stated level, and hence of $Q_d$. Consequently $A_2$
is a maximal nilpotent boundary of $X_{Q_d}$. The reduced Slodowy slice
has zero Poisson bracket. Lemma~\ref{com:bounded-current} applies to this
finitely strongly generated, nonnegatively integrally graded quotient
with degree-zero space $\C$. It forces the homomorphism $V^0(\mathfrak a)\to H_f(Q_d)$, with
$\mathfrak a=A_5$ or $E_6$, to factor through $L_0(\mathfrak a)=\C$.
For both centralizers the abelianization consists of one
$h$-degree-two line. Since the image of $\mathfrak a$ in $H_f(Q_d)$ is zero, \cref{v3:abelianization} gives a Virasoro surjection
$\Vir^c\twoheadrightarrow H_f(Q_d)$, with $c$ as in the table.

\subsubsection{\texorpdfstring{The Virasoro relation in $H^0_{A_2}(L_{-8}(E_7))$}{The Virasoro relation in the A2 reduction of L-8(E7)}}
The character sum gives more invariant relations than arise from the first singular
module. We must show that the extra relation belongs to the first
reduced ideal, rather than being counted only in a later ideal. The following
two elementary lemmas provide this step: Sugawara grading gives a projection
onto the subspace annihilated by $\mathfrak a[t]$, and a nonzero mode action
rules out repetition in later ideals. We give their proofs to make this
use of the character sum explicit.

\begin{lemma}
\label{v5:current-projection}
Let $W$ have a nonnegative integral grading, finite-dimensional graded pieces, and
a homomorphism $\bigotimes_iV^{\ell_i}(\mathfrak a_i)\to W$
sending $\mathfrak a=\bigoplus_i\mathfrak a_i$ to $W_1$, where
$\mathfrak a$ is semisimple. Assume that each graded subspace is a
semisimple $\mathfrak a$-module and that $\ell_i+h_i^\vee>0$. Let $S$ be
their Sugawara vector and $H=S_{(1)}$, preserving the given grading.
For $D_B=\sum_{n>0,a}J^a_{-n}W_{B-n}$, there is a canonical isomorphism
$$
 W_B^{\mathfrak a[t]}=\ker(H|W_B)
       \simeq(W_B/D_B)^{\mathfrak a},\qquad
 \pi(v)=v-(H|D_B)^{-1}Hv\quad(v\in W_B^{\mathfrak a}).
$$
Here $H|D_B$ is invertible. This projection preserves graded subspaces stable under all modes of $\mathfrak a$ and is exact on their quotients.
Suppose, in addition, that there is a graded $\mathfrak a$-stable
subspace $E\subseteq W$ such that multiplication gives an isomorphism
$$
 U(\mathfrak a\otimes t^{-1}\C[t^{-1}])\otimes E\longrightarrow W
$$
of graded $\mathfrak a$-modules and of modules over the displayed
universal enveloping algebra. Then
$\ch W^{\mathfrak a[t]}=\ch E^{\mathfrak a}$.
For the universal $W$-algebras used below, $E$ is spanned by the
normally ordered products of derivatives of the generators outside
$\mathfrak a$. The basis theorem stated in Section~\ref{vtwo:image-tools}
gives this decomposition.
\end{lemma}
\begin{proof}
Induct on degree. If a generalized $H$-eigenvalue
$z$ is not nonnegative real, every $J_n$, $n>0$, maps its generalized
eigenspace into eigenvalue $z-n$ in a smaller degree, and is therefore zero on it.
On this space Sugawara reduces to the finite Casimir, with eigenvalues
$$
 \sum_i\frac{(\lambda_i,\lambda_i+2\rho_i)}
                    {2(\ell_i+h_i^\vee)}\ge0,
$$
a contradiction. Thus the whole $H$-spectrum is nonnegative.
The same argument at eigenvalue zero shows that every $J^a_n$, $n>0$,
vanishes on that eigenspace.
The Casimir then forces the $\mathfrak a$-action to be trivial and acts as zero,
so there is no nontrivial Jordan block at zero.

The map $\bigoplus_{n>0,a}W_{B-n}\to D_B$ given by $J^a_{-n}$
intertwines $H+n$ with $H$, so $H$ has strictly positive spectrum on
$D_B$. For invariant $v$, the Sugawara formula places $Hv$ in
$D_B^{\mathfrak a}$. Subtracting its unique inverse image gives the
asserted zero-eigenvector lift. Taking $\mathfrak a$-invariants is exact on these finite-dimensional modules. On
each finite piece the zero-eigenvalue projector is a polynomial in
$H$. This proves preservation and exactness for subspaces and
quotients. Under the additional multiplication isomorphism,
$W/D\cong E$ as graded $\mathfrak a$-modules, which gives the
character assertion.
The assumptions $\ell_i+h_i^\vee>0$ ensure that the Sugawara vector
is defined and that the Casimir eigenvalues used above are nonnegative.
\end{proof}

\begin{lemma}\label{v5:mode-norepeat}
Let $W$ be a nonnegatively graded vertex algebra with
finite-dimensional graded subspaces, and let
$\mathfrak a\subset W_1$ be a finite-dimensional semisimple Lie
algebra acting on $W$ by the operators $a_{(0)}$, $a\in\mathfrak a$.
Write $W_d^{\mathfrak a}$ for the subspace annihilated by all these
operators. Let $K_1\supseteq K_2\supseteq\cdots$ be graded ideals
of $W$, and suppose
$$
 U\subseteq(K_1)_B^{\mathfrak a},\qquad
 \sum_r\dim(K_r)_B^{\mathfrak a}>\dim U.
$$
Let $\Phi$ be a finite collection of homogeneous compositions of
operators $u_{(n)}$, with $u\in W$ and $n\in\Z$, mapping degree
$B$ to degrees on which $K_2$ vanishes. If $\Phi|_U$ is injective,
then $(K_1)_B^{\mathfrak a}$ strictly contains $U$.
\end{lemma}
\begin{proof}
Every component preserves each ideal, so
$\Phi((K_2)_B)=0$ and $(K_2)_B\cap U=0$. Equality
$(K_1)_B^{\mathfrak a}=U$ would place every later invariant in that zero intersection,
contradicting the dimension inequality. The same argument applies
after restriction to the common kernel of the operators
$a(n)$, $a\in\mathfrak a$, $n>0$, and $L_n$, $n>0$, for a
specified Virasoro action, provided the induced map on each
multiplicity space is injective. The character sum is
interpreted in the completion of \cref{com:tail}.
\end{proof}

Denote $W=\mathcal W^{-8}(E_7,A_2)$, $\mathcal I=H_f(I_3)$ and $K_r=H_f(J^r)$.
The Virasoro quotient $W/\mathcal I$ is universal. Indeed, an affine
singular space $S$ generates an ideal with
$
 R_{V/\langle S\rangle}=\Sym\g/(\overline S).
$
The Poincar\'e--Birkhoff--Witt theorem for $U(\widehat\g)$
expresses that ideal as the span of products of
operators $a(-n)$, $a\in\g$, $n>0$, applied to $S$.
Modes $a_{(-n)}$, $n\ge2$, vanish in $R_V$, and $a_{(-1)}$ multiplies
by $a$. Here $S=V_{2\omega_7}$ has no $D_6$ invariants: its restriction
has nontrivial summands of dimensions $66,462,352,352,77,77,77$.
Evaluation on the center of $D_6+\C$ therefore annihilates all its
symbols. The entire Dynkin center line belongs to $X_{Q_3}$. Its
Slodowy intersection contains $f+ce$, conjugate to $\sqrt c\,h$ for
$c\ne0$. The Slodowy slice is positive-dimensional. Every nonzero proper
Virasoro quotient has point variety by \cref{v3:vir-quotient}. Hence
$W/\mathcal I=\Vir^{-22/5}$. Moreover, $\mathcal I$ is the ideal generated by $W_1$: its quotient is already Virasoro-generated and maps
identically onto this universal Virasoro quotient.

Through conformal weight four, the reduced Jantzen sum receives
contributions only from real affine roots with integer index
$a\le6$. There are twenty signed root contributions. \hyperref[calc:V5-E7-P10]{Calculation~\ref*{calc:V5-E7-P10}} evaluates the collected
ordinary and $A_5$-invariant characters, giving
$$
 \sum_r\ch(K_r)_{<4}=\ch(W/\Vir)_{<4},\qquad
 \sum_r\dim(K_r)_4^{A_5}=9,\qquad\dim \mathcal I_4^{A_5}=8.
$$
Here $W/\Vir$ denotes the vector-space quotient by its coset Virasoro
subalgebra. Exactness and positivity give $(K_2)_{<4}=0$.
It remains to show that no nonzero scalar vector in $\mathcal I_4$ lies in a later Jantzen ideal.

Let $C=V^0(A_5)\subset W$, let $S$ be its Sugawara vector,
and set $T=\omega-S$. The charged weight-two generators form the
$A_5$-module $15\oplus15^*$. Their invariant contraction
$Q=\sum:G_iG^i:$ has conformal weight four.
Appendix~\ref{app:e7-four} constructs the unique vector
$P\in \mathcal I_4^{A_5}$ with $Q$ coefficient one that is annihilated by
$a(n)$ for $a\in A_5$, $n>0$, and is Virasoro primary.
Since $(K_2)_{<4}=0$, every invariant vector in $(K_2)_4$
is annihilated by the positive modes of the $A_5$ generators and
by the positive Virasoro operators. The calculation in that appendix
shows that the entire subspace of $\mathcal I_4^{A_5}$ with these
properties is $\C P$. Thus every possible invariant in
$(K_2)_4\cap\mathcal I_4$ lies in this line.
The same calculation proves that, for a highest-weight generator
$G$ in the summand $15$,
\begin{equation}\label{v5:e7-detection}
 G_{(3)}P=-\frac{17600}{81}G\ne0.
\end{equation}
Since $K_2$ vanishes below four, \cref{v5:mode-norepeat}
excludes the line $\C P$ from $K_2$. The excess invariant belongs to
$K_1/\mathcal I$, whose degree-four singular line is therefore
$
             :TT:-\tfrac3{10}\partial^2T\ \in K_1.
$
Since $\mathcal I\subset K_1$, this vector belongs to $K_1$.
Thus the reduction has $[T]^2=0$ and point variety.

\subsubsection{\texorpdfstring{An additional relation in $H^0_{A_2}(L_{-12}(E_8))$}{An additional relation in the A2 reduction of L-12(E8)}}
The first singular module leaves a Virasoro quotient, so we need
one further relation to make its associated variety a point. We obtain
a relation in degree eight by finding a composition factor of the
vacuum module that is absent from the standard module generating the
first ideal. Exact reduction then detects this factor in the kernel
of the map to the simple affine quotient.

For the composition-multiplicity computation in
\hyperref[calc:HECKE]{Calculation~\ref*{calc:HECKE}}, use the right
antispherical Hecke module induced from the sign representation of
$W(E_8)=\langle s_1,\ldots,s_8\rangle$, with
$H_s^2=1+(v^{-1}-v)H_s$ and $N_1H_i=-vN_1$ for $1\le i\le8$.
The dominant representative has stabilizer $W_J$, $J=\{4,6,8\}$.
The three highest-weight labels use the longest elements of their cosets $wW_J$.
The affine node $0$ is adjacent to node $8$, and a word
$i_1\cdots i_r$ denotes $s_{i_1}\cdots s_{i_r}$.
By \cite[Theorem~1.4]{KT}, inversion of the simple-character matrix
and finite BGG alternation identify the canonical-to-standard
coefficients at $v=1$ with composition multiplicities, as explained in
\hyperref[com:hecke-multiplicities]{Section~\ref*{B:relations}}.
The relevant coefficients of the canonical vector in that calculation are $v$ and $0$,
respectively, giving
$$
 [V:L(16\theta)]=1,\qquad[M(\omega_2):L(16\theta)]=0.
$$
Consequently $\ker(Q_4\to L_{-12})$ contains this extra factor in a
genuine finite low-degree filtration. For $x=\theta$ and
$s=(\mu,\theta)$, Cauchy--Schwarz and
$(\rho,\mu)\ge29s/2$ give
$\delta(\mu)\ge s(s-14)/72$. Reduced degree at most eight implies
original degree at most forty. The subspace of maximal $x$-eigenvalue in the lowest-degree subspace
of $L(16\theta)$ is one-dimensional and has reduced degree eight.
The kernel of the corresponding standard module surjection has
original degree greater than forty. Its Drinfeld--Sokolov reduction is nonzero. The completed finite
filtration of \cref{vtwo:finite-filtration}, whose remaining tail cannot
contribute by the degree bound, therefore gives a nonzero kernel in
degree eight. The filtration consists of $\widehat{\g}$-submodules.
The degree bound is used only to control which composition factors can contribute.
Applying \cref{v3:vir-quotient} to the two quotient maps produces
$[T]^4=0$. The associated variety of the nonzero reduction is a single point.

For $E_7$, the seven center bounds are $(66,84,94,106,100,84,76)$:
the last uses the minimal $E_6$ orbit at level $-2$. Six are strict.
The $D_6+\C$ equality gives a dimension-$67$ center family with $A_2$
boundary, contradicting its zero-dimensional Slodowy slice.
For $E_8$ there is an even shorter conclusion requiring no center
reconstruction. If $z=s+n\in X$ is nonnilpotent, then its proper Levi
centralizer has dimension at most $134$, so
$$
 \dim Gz=248-\dim\g^s+\dim(G^s n)\ge114.
$$
Conic induction places its nilpotent boundary, of the same dimension, in
$X\cap\N$. The bound $114$ forces equality, $\g^s=E_7+\C$ and $n=0$.
But $\overline{G(\C^*s)}$ has dimension $115$ and Richardson boundary
$A_2$, so its transverse slice is positive-dimensional. This contradicts
the fact that the intersection of $X$ with the Slodowy slice at
$f\in\OO_{A_2}$ consists of a single point. Thus
$X_{L_{-12}(E_8)}=\overline\OO_{A_2}$. This direct global argument is
specific to this bound and these proper-Levi dimensions.

\subsubsection{\texorpdfstring{Quasi-lisseness of $L_{-21}(E_8)$, $L_{-18}(E_8)$ and $L_{-17}(E_8)$}{Quasi-lisseness of L-21(E8), L-18(E8) and L-17(E8)}}
\label{exc:e8-remaining}
For $L_{-21}(E_8)$, let $f$ belong to $D_4(a_1)+A_2$, of dimension $184$.
The reduced Jantzen sums at the minimal complementary orbits
$D_4,A_4,2A_3$ have constant terms $3,1,2$, respectively. Since these
sums are characters of nested ideals, the three reductions vanish.
It follows that $X\cap\N\subseteq\overline{\OO_f}$.
For the reduction at $f$, all contributions through reduced degree
two have original degree at most $31$. The two contributions have
signed dimensions $56,-28$ in degrees $19,20$, and the sum begins
$28q^2$. In particular, the reduction is nonzero.

The homomorphism
$V^6(A_2)\to\mathcal W^{-21}(E_8,D_4(a_1)+A_2)$ has Sugawara
central charge $16/3$, equal to that of the $W$-algebra. The
weight-two subspace of the universal $W$-algebra has dimension $72$.
The image of $V^6(A_2)_2$ has dimension $44$ and a nondegenerate
invariant bilinear form. Its complement annihilated by $a(n)$,
$a\in A_2$, $n>0$, is $27\oplus1$. The invariant line is spanned by
$D=\omega-\omega_{\mathrm{Sug}}$. Every proper reduced Jantzen ideal
has zero weight-one subspace. Since $28$ is not divisible by $27$,
the vector $D$ belongs to some reduced Jantzen ideal and hence to the
first. Thus $D=0$ in $H_f^0(L_{-21}(E_8))$.

The remaining strong generators of $\mathcal W^{-21}(E_8,f)$,
apart from the conformal vector and the weight-one subspace,
form the $A_2$-modules $27$ in
weight two, $10\oplus10^*$ in weight three, and $8$ in weight four.
Modulo the subspace generated by $A_2$, the only invariant strong
generator is the conformal vector. The nilpotent bound above makes
$\OO_f$ maximal, so \cref{com:maximal-boundary} implies that the
Poisson bracket on the reduced $C_2$ algebra is zero. Since
$[A_2,A_2]=A_2$, the classes of $A_2$ vanish. Their action then forces
the classes of all strong generators in nontrivial $A_2$-modules to
vanish. The relation $\omega=\omega_{\mathrm{Sug}}$ makes the remaining
class zero. By \cref{com:scalar-point}, the associated variety of the reduction is a point.
A nonnilpotent orbit would have a nilpotent boundary of the same
dimension by \cref{vtwo:conic-boundary}. The bound is $184$, and an
orbit attaining it would have boundary $\OO_f$, giving a
positive-dimensional intersection with the slice at $f$. This is incompatible
with the one-point intersection of $X$ with that Slodowy slice.
Orbit dimensions are even, so every
nonnilpotent orbit would have dimension at most $182$. Only
maximal Levi subalgebras with derived algebras $E_7,D_7,E_6+A_1$ can
satisfy this bound. The required nilpotent bounds in their derived
algebras are at most $68,26,16$, whereas the corresponding affine
vertex algebras give lower bounds $70,28,22$. Thus the
remaining semisimple points cannot occur.

For $L_{-18}(E_8)$, reduce at $2A_2$, of dimension $156$. Write $K_r=H_f^0(J^r)$ for the reduced Jantzen
ideals. The associated
homomorphism from an affine vertex algebra is
$V^0(G_2)^{\otimes2}\to\mathcal W^{-18}(E_8,2A_2)$.
The terms of reduced degree at most two in the Jantzen sum have
original degree at most $27$, and their sum is $28q+729q^2$.
Its constant term is zero. The weight-one subspace
$(14,1)\oplus(1,14)$ occurs with multiplicity one, so it belongs to
$K_1$, whereas $(K_r)_1=0$ for $r\ge2$. The weight-two subspace of the
universal $W$-algebra has dimension $484$. The subspace generated by
$G_2\oplus G_2$ contributes $434$ dimensions, and the strong generators
in $(7,7)$ contribute $49$ dimensions. Thus $(K_1)_2$ would have
dimension $483$ if it contained no invariant generator outside the
affine vertex subalgebra.

The weight-two subspaces of later ideals are sums of $(14,14)$,
$(7,7)$ and the trivial module. Indeed, the operator induced by positive
modes on $\Sym^2(14)=1+27+77$ has eigenvalues $8-C_V/2$, where
$C_V=0,28/3,20$. The resulting eigenvalues are all nonzero. In the absence of the invariant
line, the dimensions of the later ideals would be divisible by $49$.
Since $729-483=246\equiv1\pmod {49}$, \cref{v3:scalar-detection}
implies that this line belongs to $K_1$. Modulo the subspace generated
by $G_2\oplus G_2$, it is represented by the conformal vector.
The image of $G_2\oplus G_2$ in the quotient is already zero, so
$\omega=0$. By \cref{com:scalar-point}, the reduction is $\C$.
The vanishing reductions at $A_3$ and $2A_2+A_1$ exclude the minimal
complementary orbits and give $X\cap\N\subseteq\overline{\OO}_{2A_2}$. By \cref{com:bounded-current,com:maximal-boundary}, an $E_7$
center forces dimension $166$.
The pure $D_7$ sheet is excluded by the zero-dimensional Slodowy slice at $2A_2$.

\subsubsection{\texorpdfstring{A $G_2\times G_2$-invariant relation in $H^0_{2A_2}(L_{-17}(E_8))$}{A G2 x G2 invariant relation in the 2A2 reduction of L-17(E8)}}
The reduced Jantzen count first gives a relation beyond the known affine
relations. We show that the known affine invariant relations do not recur in later ideals, then project the additional
relation onto the Virasoro descendants. The Gram calculation in
Appendix~\ref{app:virasoro-four} gives a nonzero coefficient of $:TT:$,
hence a nilpotent conformal class. It uses the conformal weight for
$T$, rather than the full conformal weight.

\begin{lemma}\label{v3:no-repeat}
Let $C=\bigotimes_i V^{\ell_i}(\mathfrak a_i)$ have noncritical
levels, and let $M\subseteq C$ be graded and stable under all operators $a(n)$, $a\in\bigoplus_i\mathfrak a_i$,
$n\in\Z$. Suppose $M_e=0$ for $e\le d-2$, where $d>0$, and
$d-1\ne h_i^\vee/(\ell_i+h_i^\vee)$ for every factor.
Then $M_d^{\oplus_i\mathfrak a_i}=0$.
\end{lemma}
\begin{proof}
For $v\in M_d^{\oplus_i\mathfrak a_i}$, the map
$a\mapsto a(1)v$ has degree $d-1$ and transforms in the adjoint
representation. Every operator $b(n)$, $n>0$, annihilates its image
because $M_e=0$ for $e\le d-2$. A nonzero adjoint summand would have
Sugawara conformal weight $h_i^\vee/(\ell_i+h_i^\vee)$, contrary to
the assumption on $d-1$. Thus $a(1)v=0$. The operators $a(n)$ for
$n\ge2$ also vanish on $v$. All nonnegative modes therefore
annihilate $v$, so its Sugawara conformal weight is zero. This argument stays inside $C$,
whose grading is defined by the sum of the Sugawara vectors of
$V^{\ell_i}(\mathfrak a_i)$.
\end{proof}

\begin{lemma}
\label{v3:gram-projection}
Let $K_i$ be nested proper graded ideals in a nonnegatively graded
vertex algebra $W$, with $W_0=\C$ and normalized invariant bilinear
form. Let a reductive algebraic group $A$ act on $W$ by
grading-preserving automorphisms, and assume that every $K_i$ is
$A$-stable. In the subspace of $A$-invariants of degree $d$, suppose
$$
 W_d^A=B\oplus D\oplus\C Q,
$$
where $B\subseteq(K_1)_d^A$ lies in the radical of the invariant
bilinear form and the form on $D$ is
nondegenerate. If later layers have zero intersection with $B$ and
$\sum_i\dim(K_i)_d^A>\dim B$, then
$Q-\operatorname{proj}_D Q\in K_1$.
\end{lemma}
\begin{proof}
Otherwise, if $(K_1)_d^A=B$, nestedness and non-repetition would make
the total exactly $\dim B$. There is therefore a vector in the first ideal
outside $B$. Every proper graded ideal is contained in the radical
of the invariant bilinear form. Nondegeneracy of
$D$ forces its $Q$ coefficient to be nonzero. Normalize it to one.
Orthogonality fixes its $D$ component, and its $B$ correction already
belongs to $K_1$.
\end{proof}

\begin{corollary}
\label{v5:scalar-coefficient}
In degree four of \cref{v3:gram-projection}, write
$D=D'\oplus\C E\oplus\C F$, where $D'$ is nondegenerate and
orthogonal to $E,F$. Let $T$ be a Virasoro vector of charge $c$,
commuting with $\omega-T$, and set
$E=L^T_{-4}\mathbf1$, $F=(L^T_{-2})^2\mathbf1$.
Assume $c(5c+22)\ne0$ and
$$
 Q=\sum_{a=1}^m(G_a)_{(-1)}G^a,\qquad
 (G_a)_{(3)}G^b=N\delta_a^b\mathbf1,
$$
where the fields are primary of weight two for $\omega$ and primary
of weight $h$ for $T$. The relation, normalized to have $Q$ coefficient
one, has $F$ coefficient
\begin{equation}\label{v5:scalar-beta}
 \beta=-\frac{2mNh(5h+1)}{c(5c+22)}.
\end{equation}
If all its other terms vanish in the reduced $C_2$ algebra, then
$\beta[T]^2=0$ there. If $\beta\ne0$, then $[T]$ is nilpotent
in the unreduced $C_2$ algebra.
\end{corollary}
\begin{proof}
Existence in the first ideal follows from
\cref{v3:gram-projection}. Appendix~\ref{app:virasoro-four} computes the Gram matrix of
$E,F$ and their pairings with $Q$. Orthogonal projection onto
$\C E\oplus\C F$ gives \eqref{v5:scalar-beta}. The orthogonal
summand $D'$ does not affect this coefficient. The remaining
assertions follow by passing the relation to the reduced $C_2$ algebra.
\end{proof}

Set $W=\mathcal W^{-17}(E_8,2A_2)$, $H=H_f^0(L_{-17}(E_8))$,
$A=G_2\times G_2$ and $K_r=H_f^0(J^r)$. The complete equivariant
Jantzen sum gives
\begin{equation}\label{exc:p13-count}
 \sum_r\ch_A(K_r)_2=(77,1)+(1,77),\qquad
 \sum_r\dim(K_r)_4^A=3,
\end{equation}
with no terms in degrees zero and one. $77$ means highest weight
$2\theta$. All relevant terms have original degree at most $38$: for reduced degree at
most four, Cauchy--Schwarz requires
$(D+42)^2\le4(26D+620)$, which fails for every $D\ge39$.
The twenty signed $D_7$ characters in \hyperref[calc:P13]{Calculation~\ref*{calc:P13}} have
numerator pieces $P_2,P_3,P_4$. With
$Z=(14,1)+(1,14)+(7,1)+(1,7)+(7,7)+(1,1)$, their invariant contributions
are
$$
 \langle P_4,1\rangle=1,\qquad \langle P_3,Z\rangle=0,\qquad
 \langle P_2,Z+\Sym^2Z\rangle=2.
$$
The oscillator factor is $\prod_{n\ge1}\Sym(q^nZ)$. Hence the count is
$1+0+2=3$. Completed \emph{equivariant} exact reduction makes this a
sum of nested ideals.

The zero constant term proves that the reduction is nonzero. The
vanishing reductions at the minimal complementary orbits give $X\cap\N\subseteq\overline{\OO}_{2A_2}$, so \cref{com:bounded-current} gives an embedding of
$L_1(G_2)^{\otimes2}$ in $H$. The two singular modules generating the maximal ideals of
$V^1(G_2)$ account for the entire
degree-two sum in \eqref{exc:p13-count}. Later $K_r$ vanish through degree
two. Let $C=V^1(G_2)^{\otimes2}\subset W$. In each factor, the $G_2$-invariant subspace of degree four in the
maximal ideal is one-dimensional: its universal invariant dimension
is three, its two Sugawara descendants are nondegenerate, and the sum of
squares of a basis of its singular module of dimension $77$ has image
$(11/16)(x,x)^2\ne0$ in $R_{V^1(G_2)}=\Sym G_2$. Indeed quartic adjoint invariants are multiples
of the Casimir square. Applying the Euclidean Laplacian twice gives
$16\cdot77$ on the sum and $8\cdot14\cdot16$ on $(x,x)^2$.
Their span $B\subset(K_1)_4^A$ has dimension two.
The Sugawara argument in \cref{v3:no-repeat}, with adjoint Sugawara conformal weight $4/5\ne3$, gives $(K_r\cap C)_4^A=0$ for $r\ge2$. If $(K_1)_4^A=B$,
all later invariants would vanish, contradicting the count three.
Thus the first reduced ideal contains an additional
$G_2\times G_2$-invariant relation.

Let $S_1,S_2$ be the Sugawara vectors of the two copies of $V^1(G_2)$ and $T=\omega_W-S_1-S_2$.
Their central charges are $(14/5,14/5,-384/65)$. The $(7,7)$ generators have
respective conformal weights $(2/5,2/5,6/5)$ and pairing
$(G_a)_{(3)}G^b=-16\delta_a^b\mathbf1$, from the
independent universal Miura norm calculation in
\hyperref[calc:E8-IDEAL]{Calculation~\ref*{calc:E8-IDEAL}}. Set $Q=\sum_{a=1}^{49}:G_aG^a:$.
Counting the normally ordered products of the universal strong
generators gives $\dim W_4^A=12$: besides $B$, there are nine nondegenerate
descendants of the three commuting conformal vectors and the direction $Q$.
The three Virasoro Gram blocks have determinants
$c_i^2(5c_i+22)/2\ne0$, and the mixed norms are $c_ic_j/4\ne0$.
Thus \cref{v3:gram-projection} makes the $Q$ coefficient of the
additional relation nonzero. Normalize it to one. Only its
coefficient of $:TT:$ is needed: \cref{v5:scalar-coefficient} gives
\begin{equation}\label{exc:p13-relation}
 R=Q+\frac{1183}{4}:TT:+D\ \in K_1,\qquad
 \frac{1183}{4}=-\frac{2mNh(5h+1)}{c(5c+22)},
\end{equation}
where $(m,N,h,c)=(49,-16,6/5,-384/65)$ and
$D\in B+\C L_{-4}^T\mathbf1+D'$. Here $D'$ consists of the seven
other descendants involving $S_1,S_2$. Its form is nondegenerate and
it is orthogonal to $L_{-4}^T\mathbf1,(L_{-2}^T)^2\mathbf1$. Its
coefficients need not be evaluated. \hyperref[calc:V5-P13]{Calculation~\ref*{calc:V5-P13}} verifies the
two Ward contractions using $D_T=L^T_{-1}$ and gives this nonzero
residual coefficient. Here $D_T$ is the mode of $T$, whereas the full
translation operator is
$L_{-1}=L^{S_1}_{-1}+L^{S_2}_{-1}+L^T_{-1}$.

The character inequality above proves $R\in K_1$. In the reduced
$C_2$ algebra of $H$, the classes of $G_2\oplus G_2$ vanish by
integrability. Their Hamiltonian action then forces the classes of
the generators in nontrivial $G_2\oplus G_2$-modules to vanish.
Every term of $D$ vanishes there, as does $Q$. Hence
$(1183/4)[T]^2=0$, and the associated variety of $H$ is a point. Finite strong
generation makes $H$ $C_2$-cofinite.

The ideal generated by the two singular modules alone does not imply
this conclusion: the Miura homomorphism to
$L_1(D_7)\otimes\mathcal H_{52}$ annihilates that ideal, while
augmentation after $C_2$ sends $[T]$ to $z^2/104$. The extra Jantzen contribution is therefore essential.

The $E_7$ lower bound $52$ used here is independent of the present
$E_8$ conclusion. The known minimal reduction
$H^0_{f_{\min}}(L_{-6}(E_7))=L_{-2}(D_6)$ has associated variety
equal to the minimal nilpotent orbit closure of dimension $18$, so
$X_{L_{-6}(E_7)}$ attains dimension $34+18=52$.
Translation by the integrable affine vertex algebra $L_1(E_7)$ gives
the upper bound $52$ for $X_{L_{-5}(E_7)}$. If this bound were not attained, the only
possible maximal nilpotent boundaries would be $0,A_1$. The associated homomorphisms from $V^{-5}(E_7)$ and $V^{-1}(D_6)$
have negative levels, contradicting the integrability forced by bounded
Cartan weights in \cref{com:bounded-current}. Thus
both associated varieties contain a nilpotent orbit of dimension $52$,
which induces a nilpotent orbit of dimension $166$ in $E_8$.
The Slodowy-slice argument and these induced-orbit dimensions rule out the
remaining nonzero semisimple points as in the proof for $L_{-18}(E_8)$.
This proves quasi-lisseness of $L_{-17}(E_8)$. Simplicity of the Drinfeld--Sokolov
reduction is not required.

This proves quasi-lisseness in all $37$ cases of
\cref{exc:coverage}. The associated varieties are determined below
by the nonvanishing and vanishing of the corresponding reductions.

\section{Nilpotent orbits in the associated varieties}\label{cert:orbits}
Once $X_{L_k(\g)}\subseteq\N(\g)$ is known, its nilpotent orbits
determine it. Exact Drinfeld--Sokolov reduction detects their occurrence:
$H_f^0(L_k(\g))\ne0$ precisely when $f\in X_{L_k(\g)}$.
The proposition restates \cref{com:ds,com:completion,com:vacuum-test}
in terms of a constant coefficient of the reduced Jantzen sum, with a
proof fixing its normalization. We express this coefficient as a finite
Weyl-group sum and use the orbit-closure order to determine the variety.
The comparison of quasi-lisseness at different denominators does not
identify the associated varieties at those levels.

\begin{proposition}\label{prop:orbit-vacuum}
Let $t>0$ be rational and $x$ the dominant Dynkin grading of a nilpotent
$f$. Write $\mathcal J_{t,x}(Q)$ for the normalized specialization of the
Jantzen sum for $V^{t-h^\vee}(\g)$, in the completion of
\cref{com:tail}. Thus
$$
 \sum_{i\ge1}\ch H_f^0(J^i)
 =\mathcal J_{t,x}(Q)\,\ch\mathcal W^k(\g,f),
 \qquad k=t-h^\vee.
$$
Then
$$
 R_f(t)=[Q^0]\mathcal J_{t,x}(Q)\in\Z_{\ge0},\qquad
 f\in X_{L_{t-h^\vee}(\g)}\ \Longleftrightarrow\ R_f(t)=0.
$$
If $X_{L_{t-h^\vee}(\g)}\subseteq\N(\g)$, these equalities
determine $X_{L_{t-h^\vee}(\g)}$.
\end{proposition}
\begin{proof}
Let $J^i$ be the Jantzen filtration of $V^k(\g)$, where
$k=t-h^\vee$, and set $K_i=H_f^0(J^i)$. The specialized sum
multiplied by the character of $\mathcal W^k(\g,f)$ is
$\sum_{i\ge1}\operatorname{ch}K_i$. Since $K_i\subseteq K_1$
and each graded dimension is nonnegative, its constant coefficient
vanishes precisely when $\mathbf1\notin K_1$.
Exactness of $H_f^0$ identifies this condition with
$H_f^0(L_k(\g))\ne0$. The character of
$\mathcal W^k(\g,f)$ has constant coefficient one.
The normalized sum has no negative powers, since otherwise their product
would have a negative leading power. Its constant coefficient is therefore
$\sum_{i\ge1}\dim(K_i)_0$. Finally the associated variety is closed and
$G$-stable, and the nilpotent orbit set is finite.
\end{proof}

\subsection{Bounds on the real affine roots in the Jantzen sum}
Let $\alpha$ be a root of $\g$ and let $a\in\Z_{>0}$ be the
integer index of the real affine root $\alpha+a\delta$. Set
$$
 b=\langle\rho,\alpha^\vee\rangle+
       \frac{2at}{(\alpha,\alpha)}.
$$
Only $b\in\Z_{>0}$ occurs. The reflected finite Weyl numerator has
weight $\rho-b\alpha$ and Sugawara degree $ab$. A Weyl monomial can
contribute to the constant coefficient only if
\begin{equation}\label{eq:orbit-cutoff}
 b\bigl(a+(w\alpha,x)\bigr)=(w\rho,x)-(\rho,x)\le0.
\end{equation}
Hence $a\le\max_w(-(w\alpha,x))$. Differentiation of the Weyl
alternant at repeated grading eigenvalues does not change its exponents or
this bound. Individual signed root contributions may have
negative powers of $Q$, but their sum has no negative powers.

\subsection{A finite Weyl-group formula for \texorpdfstring{$R_f(t)$}{Rf(t)}}\label{cert:weyl-certificate}
We evaluate the specialized Weyl characters by differentiating the
alternants along the roots perpendicular to the grading element.
Let $\Phi_0^+=\{\beta\in\Phi^+:(\beta,x)=0\}$, let $W_0$ be its
Weyl group, and set
$$
 d_x=|W_0|\prod_{\beta\in\Phi_0^+}(\rho,\beta),\qquad
 c_x(\nu,D)=\frac1{d_x}
 \sum_{\substack{w\in W\\(w\nu,x)=(\rho,x)+D}}
 (-1)^{\ell(w)}\prod_{\beta\in\Phi_0^+}(w\nu,\beta).
$$
An empty product is one. The summand is constant on each left
$W_0$-coset: both its sign and the product $\prod_{\beta\in\Phi_0^+}(w\nu,\beta)$ change by the same factor. The sum over $W$ is therefore $|W_0|$ times the sum over
$W_0\backslash W$.
The resulting finite rational identity is
\begin{equation}\label{eq:finite-weyl-certificate}
 R_f(t)=\sum_{\alpha\in\Phi}\ \sum_{\substack{
 1\le a\le\lfloor\max_{w\in W}(w\alpha,x)\rfloor\\
 b=2((\rho,\alpha)+at)/(\alpha,\alpha)\in\Z_{>0}}}
 c_x(\rho-b\alpha,ab).
\end{equation}
Here $x$ is specified by the weighted Dynkin diagram $2\alpha_i(x)$
of $f$, and \eqref{eq:orbit-cutoff} gives the displayed range of $a$.

To prove the formula, apply the differential operator
$\prod_{\beta\in\Phi_0^+}\partial_\beta$ to the Weyl numerator and
denominator before restricting the torus variable to the line through
$x$. The denominator vanishes to order $|\Phi_0^+|$ along that line,
and the numerator vanishes to at least this order.
The leading denominator coefficient is $d_x$: precisely $w\in W_0$
fix the maximal exponent $(\rho,x)$, and anti-invariance makes their
contributions equal. Leibniz's rule therefore gives the specialized
character quotient, including repeated grading eigenvalues. Summing the
affine root contributions first and using the nonnegative-power conclusion in
\cref{prop:orbit-vacuum}, its constant term is the leading numerator
coefficient divided by $d_x$. The bound on $a$ is \eqref{eq:orbit-cutoff}.
No assertion of positivity is made for an individual summand.

For each pair $(p,q)$ listed in the tables and proposed maximal orbit $T$,
it is enough to verify
\begin{equation}\label{eq:table-certificate}
 R_T(p/q)=0,\qquad R_S(p/q)>0
 \quad\text{for every minimal }S\not\subseteq\overline{\OO_T}.
\end{equation}
The order in \eqref{eq:table-certificate} is the closure order on
nilpotent orbits. For classical Lie algebras it is the dominance order
on partitions, with both labels included for every very-even partition.
For exceptional types we use the weighted diagrams and closure order
in the tables of Adams \cite{AdamsOrbitTables}, based on the orbit data
of Spaltenstein and Carter. Those tables
contain finite-dimensional orbit data only, not associated varieties of affine vertex algebras.

\subsection{Classical Weyl determinants and associated varieties for large denominators}
\label{cert:classical-formula}
For classical types, the confluent Weyl alternants give determinants of
order equal to the rank. The complete entries and their integer
coefficient recurrence are in \hyperref[calc:ORBIT-CLASSICAL]{Calculation~\ref*{calc:ORBIT-CLASSICAL}}.
The Weyl group has an even number of sign changes in type $D$.
In type $C$ the pairing uses the normalization of short roots.
We use both labels of a very-even partition. The determinant
expansion gives \eqref{eq:finite-weyl-certificate}.
The irreducible components of $X_{L_k(\g)}\subseteq\N(\g)$ are
the closures of the maximal orbits among those satisfying $R_f(t)=0$.

The following bounds on the denominator imply that
$X_{L_k(\g)}=\N(\g)$ at quasi-lisse levels:
\begin{equation}\label{eq:classical-tails}
\begin{array}{c|l}
B_n&q\text{ odd}:q\ge2n+1,\quad q\text{ even}:q\ge2n+2,\\
C_n&q\text{ odd}:q\ge2n+1,\quad q\text{ even}:q\ge4n-2,\\
D_n&q\ge2n-2.
\end{array}
\end{equation}
At those denominators no integral root contribution remains, so at every quasi-lisse
level $X=\N$. For $B,D$ the additional bound $q>p$ gives
$X\cap\N=\N$, and thus $X=\N$ at every quasi-lisse level:
at $0<t<1$ the conformal-weight inequality for principal reduction is strict for every nonzero
dominant weight. In simply-laced type it is
$$
 \Delta_t(\lambda)-\lambda(\rho)
 =\frac{\|\lambda\|^2+2(1-t)(\rho,\lambda)}{2t}>0.
$$
For $B$, $2\rho-2\rho^\vee=(-1,\ldots,-1)$, and for an integral
dominant root-lattice weight the remaining term
$\sum_i\lambda_i(\lambda_i-1)$ is nonnegative. $t<1$ makes the total
strict. The type-$C$ tables use the bounds in
\eqref{eq:classical-tails}.

At admissible levels, the orbit descriptions are the published results
of \cite[Theorems~5.7.1 and~5.8.1, Tables~2--10]{AAV}, see also
\cite{AVEM} for their interpretation in terms of primitive ideals.
We recall their classical partition form for use in the tables.
For the type-$D_n$ condition below, form the strings
$d-1,d-3,\ldots,1-d$ for the parts $d$ of $\lambda$. Their positive entries in
decreasing order, padded by zeros to length $n$, form $\eta$, with
$x=\eta/2$. The conditions on the partitions labelling orbits in the associated variety are
\begin{equation}\label{eq:admissible-partitions}
\begin{array}{c|l}
A_{N-1}&\lambda_1\le q,\\
B_n&\lambda_1\le q,\\
C_n, q\text{ odd}&\lambda_1\le q,\\
C_n, q\text{ even}&\max(\lambda_1+\lambda_2-2,2\lambda_1-4)<q,\\
D_n&\eta_1+\eta_2<2q.
\end{array}
\end{equation}
The associated variety is the closure of the orbit labelled by the
maximal allowed partition. Set
$\lambda_2=0$ for a one-part partition. These admissible formulas are not
used at nonadmissible numerators.

Appendix~\ref{app:associated-varieties} lists the resulting orbits and
dimensions. The orbit comparisons include both very-even labels.

\section*{Acknowledgements}    
D.A. was partially supported by the Croatian Science Foundation under the project IP-2022-10-9006  and by the project “Implementation of cutting-edge research and its application as part of the Scientific Center of
Excellence for Quantum and Complex Systems, and Representations of Lie Algebras”, Grant No. PK.1.1.10.0004, co-financed by the European Union through the European Regional Development Fund - Competitiveness and Cohesion Programme 2021-2027. 
K.S. is supported by China NSF grant No.1250010384. 

\vspace{5mm}
\noindent
\textbf{Declaration of AI Use.} 
The manuscript grew out of authors’ long discussions with AI to determine the non-admissible quasi-lisse levels of affine vertex algebras. AI systems were used for literature searches, checking proof strategies, calculations and manuscript preparation. The authors are responsible for the final manuscript.

\appendix
\section{Character identities and vertex-algebra relations}\label{app:calculations}
The calculations listed below use the following bounds. Relations in $H_f^0(I_k)$
also require \cref{com:vacuum-test,vtwo:first-top,v5:mode-norepeat}.
Miura-image identifications require the stated strong-generator bounds.

The coefficients lie in $\Q$ or the specified algebraic extensions.
A nonzero minor modulo a prime avoiding the denominators gives a
rank lower bound over $\Q$. Kernel vectors or independent upper bounds
give equality. Vertex products include all derivative corrections:
the bracket on $\g^f$ gives only the linear part of the zeroth product.

\subsection{Bounds on root contributions to the reduced Jantzen characters}
\label{cert:scope}\label{sec:external-inputs}
If a standard module has $\delta_x(\lambda)\le B$ and
$D=\Delta_t(\lambda)$, then Cauchy--Schwarz gives
$$
 (D+(\rho,x)-B)^2\le\|x\|^2(2tD+\|\rho\|^2)
$$
whenever $D+(\rho,x)-B>0$. Beyond the larger root of this quadratic,
\cref{com:tail} shows that no term contributes through reduced degree $B$.
Each coefficient is taken in the sum of all signed root contributions.
An individual contribution need not be the character of an ideal.

A term of the Weyl numerator contributing
to reduced degree at most $B\ge0$ satisfies
$$
 b\bigl(a+(w\alpha,x)\bigr)
 \le B+(w\rho-\rho,x)\le B,
 \qquad b=\frac{2((\rho,\alpha)+at)}{(\alpha,\alpha)}\in\Z_{>0}.
$$
Since $x$ is dominant and $b\ge1$, it follows that
$$
 a\le\left\lfloor\max_{\alpha\in\Phi}(\alpha,x)+B\right\rfloor.
$$
The contributions through reduced degree $B$ are therefore supported
on this finite set of roots and integers $a$. Each summand includes the
finite Weyl alternant and the oscillator denominator. The polynomial
and determinant identities used for the classical Lie algebras hold
uniformly in the rank.

\subsection{List of calculations}
The table lists the finite identities, degree bounds and their uses
in the proofs. Detailed evaluations appear with the same numbering in
\cite{Certificates}. The cited arguments establish the resulting
vertex-algebra relations.

\newcounter{finitecalculation}
\renewcommand{\thefinitecalculation}{\thesection.\arabic{finitecalculation}}
{\small\setlength{\tabcolsep}{4pt}
\begin{longtable}{@{}p{.07\textwidth}>{\raggedright\arraybackslash}p{.89\textwidth}@{}}
\toprule No.&Identity and bound\\\midrule\endfirsthead
\toprule No.&Identity and bound\\\midrule\endhead
\midrule\multicolumn{2}{r}{Continued on the next page}\\\endfoot
\bottomrule\endlastfoot
\refstepcounter{finitecalculation}\label{calc:B-FOLDS}\thefinitecalculation&Uniform orthogonal canonical-basis identities in \cref{seed:B-data}, on the corresponding finite Bruhat intervals.\\[2pt]
\refstepcounter{finitecalculation}\label{calc:BCD}\thefinitecalculation&Classical character identities in \cref{seed:C-vacuum,seed:D-windows,seed:B-data}: the bounds $s\le3$ for $(4,1^{2n-4})$ and $s\le2$ for $(3^2,1^{2n-6})$ in type $C$, the displayed bounds in type $D$, and the kernel $K$ of \cref{com:whole-kernel} in type $B_6$ through degree ten. The uniform type-$B$ identities are proved in \cref{B:lower-integral}.\\[2pt]
\refstepcounter{finitecalculation}\label{calc:B-ODD-STRIP}\thefinitecalculation&Finite specializations of the type-$B$ character identities in \cref{B:lower-integral}, including the contributions through conformal weight two and the degree-five relation. These calculations supplement the proof for arbitrary rank given there.\\[2pt]
\refstepcounter{finitecalculation}\label{calc:B4-LIFT}\thefinitecalculation&$B_4,t=5/2$: twenty-six regular Jantzen terms through degree sixteen give the multiplicities and branching bounds in \cref{seed:B-data}.\\[2pt]
\refstepcounter{finitecalculation}\label{calc:F-CHAR}\thefinitecalculation&$F_4,t=3/2,7/2$: complete characters through degree six are $(1,0,1,1,2,2,4)$ and $(1,0,1,1,2,2,3)$, with every intervening half-degree included.\\[2pt]
\refstepcounter{finitecalculation}\label{calc:F-WINDOWS}\thefinitecalculation&The four type-$F_4$ character identities give the displayed dimensions of graded subspaces of $H_f^0(I_k)$ and constant coefficients, with the root and mode bounds of Appendix~\ref{cert:scope}.\\[2pt]
\refstepcounter{finitecalculation}\label{calc:G-PRINC}\thefinitecalculation&$G_2,t=3/2$: the thirteen terms with $D\le26$ give a one-dimensional subspace of conformal weight zero after principal reduction.\\[2pt]
\refstepcounter{finitecalculation}\label{calc:EXCEPTIONAL}\thefinitecalculation&The signed root sums, multiplied by the oscillator denominator, give the reduced Jantzen coefficients used in Section~\ref{exc:seeds}.\\[2pt]
\refstepcounter{finitecalculation}\label{calc:ROOT-WINDOW}\thefinitecalculation&Characters of $V^k(\g)$ and of the lowest nonzero graded subspaces of $I_k$ for exceptional $\g$, with modes $a\le\lfloor\max_\alpha(\alpha,x)+B\rfloor$ through reduced degree $B$.\\[2pt]
\refstepcounter{finitecalculation}\label{calc:WEYL-PROD}\thefinitecalculation&The full relative Weyl products give the additional coefficients at $E_6,7/2$, $E_7,11/2$ and $E_8,19/2$ in \cref{vtwo:exceptional-caps}.\\[2pt]
\refstepcounter{finitecalculation}\label{calc:EQUIVAR}\thefinitecalculation&The complete $S_3$- and $S_5$-equivariant characters give the dimensions of invariant subspaces at numerator six in \cref{vtwo:equivariant-exceptional}.\\[2pt]
\refstepcounter{finitecalculation}\label{calc:HECKE}\thefinitecalculation&$E_8,p=18$: fifty-six standard terms form seven eight-element blocks, each indexed by $w_i\langle s_4,s_6,s_8\rangle$, where $s_4,s_6,s_8$ are pairwise commuting simple reflections in Bourbaki numbering. Bar invariance and unitriangularity determine the coefficients of the relevant antispherical canonical-basis vector.\\[2pt]
\refstepcounter{finitecalculation}\label{calc:SCALAR}\thefinitecalculation&All $L_1,\ldots,L_4$ equations in the ten-dimensional weight-four space give the primary vector of weight four and the minor $-64r^2(5r-11)$ of \cref{vtwo:scalar-recognition}, see also Appendix~\ref{app:primary-four}.\\[2pt]
\refstepcounter{finitecalculation}\label{calc:C2-IMAGE}\thefinitecalculation&$C_2,t=3/2$: all four poles of the mixed lattice product give the stated presentation. Higher poles are impossible by grading.\\[2pt]
\refstepcounter{finitecalculation}\label{calc:V5-F-T1}\thefinitecalculation&$F_4,t=1$: three complete BRST cocycles are written as sums of 480 basis monomials in the BRST complex. Their products span the stated two-dimensional subspace of weight six and one-dimensional subspace of weight eight.\\[2pt]
\refstepcounter{finitecalculation}\label{calc:QUAD}\thefinitecalculation&Four complete fields at $F_4,t=2$ and $E_8,p=5$, from the $12\times12$ and $40\times40$ systems of \cref{com:quadratic-miura}, have the stated projections and invariant bilinear forms on the weight-one subspaces.\\[2pt]
\refstepcounter{finitecalculation}\label{calc:G-IMAGE}\thefinitecalculation&$G_2,t=1$: all singular products of the four fields in $\SF(2)$ give the subregular presentation. Pole orders are bounded by the sums of weights.\\[2pt]
\refstepcounter{finitecalculation}\label{calc:E6-LIFT}\thefinitecalculation&$E_6,p=2,3$: the $12\times12$ and $18\times20$ systems give the complete fields and invariant bilinear forms on the weight-one subspaces. At $p=3$ the formulas include nine vectors of conformal weight two and the twenty-eight-dimensional space of skew matrices.\\[2pt]
\refstepcounter{finitecalculation}\label{calc:W24}\thefinitecalculation&$E_8,p=4$: specialization of the normalized $W(2,4)$ presentation in \cite[equations~(5.10)--(5.13)]{BSW} gives the squared structure constant $(C_{44}^4)^2=-1350$ and five nonzero denominator factors.\\[2pt]
\refstepcounter{finitecalculation}\label{calc:CENTRALIZER}\thefinitecalculation&Full graded centralizers give the stated abelianizations and ranks of the bracket maps generating the relevant Lie ideals. At $E_8,p=4$, let $K_{j,q}$ denote the degree-$j$ subspace of $\g^f$ with eigenvalue $q$ for $\operatorname{ad}J$, where $J$ spans the one-dimensional reductive centralizer. The maps $K_{i,q}\otimes K_{j,-q}\to K_{i+j,0}$, $u\otimes v\mapsto[u,v]$, with $q\ne0$, together with $K_0$ determine the neutral part of the ideal generated by $K_0$ and the nonzero-charge subspaces. The $21\times21$ and $24\times24$ systems give the six quadratic generators for $E_7$ and the three quadratic generators for $E_8$.\\[2pt]
\refstepcounter{finitecalculation}\label{calc:NUM2}\thefinitecalculation&Numerator two: the quartic restrictions, exterior brackets and half-spin brackets include the seventy weights of the exterior representation and the 128 weights of the half-spin representation.\\[2pt]
\refstepcounter{finitecalculation}\label{calc:E8-TOP}\thefinitecalculation&$E_8,p=3$: nineteen root contributions with $a\le13$ give the two-dimensional lowest conformal-weight subspace, on which the three displayed matrices representing the degree-preserving modes act irreducibly.\\[2pt]
\refstepcounter{finitecalculation}\label{calc:TANGENT}\thefinitecalculation&$E_7,p=4,5$: the differentials of the $G_0$-orbit maps, $a\mapsto[a,f]$ from $\g_0$ to $\g_{-1}$, have ranks twenty-four and thirty, with diagonal $\slg_2$ kernels.\\[2pt]
\refstepcounter{finitecalculation}\label{calc:E7-FORM}\thefinitecalculation&$E_7,p=4,5$: degree-two parity and invariant forms establish PCT stability. Lowering reconstructs the entire level-one $27$-dimensional null module.\\[2pt]
\refstepcounter{finitecalculation}\label{calc:V5-E8-P7}\thefinitecalculation&$E_8,p=7$: fourteen highest vectors exhaust the $52$-dimensional centralizer. Eleven nonzero equivariant bracket maps generate the stated summands. The $2\times2$ matrix of the $e_{-\alpha_1},e_{-\alpha_2}$ coordinates of $z,f$ has determinant $-7$, proving that the two scalar directions are independent. The norm is $-240(k+23)(3k+70)$.\\[2pt]
\refstepcounter{finitecalculation}\label{calc:G-RETURN}\thefinitecalculation&$G_2,t=4/3,5/3$: the zero-weight subspaces of the singular $G_2$-modules and all Taylor terms through degrees six and ten give determinants $916,382$ modulo $1009$, with all denominators invertible.\\[2pt]
\refstepcounter{finitecalculation}\label{calc:V5-E7-P10}\thefinitecalculation&$E_7,p=10$: twenty signed terms with modes at most six give nine invariants in degree four, against eight in the known ideal. Two highest cocycles and six contractions of weight at most four give $G_{(3)}P=-(17600/81)G$. The calculation is presented in Appendix~\ref{app:e7-four}.\\[2pt]
\refstepcounter{finitecalculation}\label{calc:P13}\thefinitecalculation&$E_8,p=13$: conformal degree at most thirty-eight before Drinfeld--Sokolov reduction gives twenty signed $D_7$ characters. Two $G_2\times G_2$ projections give \cref{exc:p13-count}.\\[2pt]
\refstepcounter{finitecalculation}\label{calc:E8-IDEAL}\thefinitecalculation&The stated $E_8$ character identities and invariant subspaces show that the stated nilpotent orbits do not belong to the associated varieties at $p=6,9,12,13$. The pairing computed from the Miura homomorphism at $p=13$ is $N=-16$.\\[2pt]
\refstepcounter{finitecalculation}\label{calc:V5-P13}\thefinitecalculation&The Gram matrix on the two-dimensional Virasoro subspace and the identities involving $D_T=L^T_{-1}$, distinguished from the full translation operator $L_{-1}$, give $\beta=1183/4$ in \cref{exc:p13-relation}, with a nondegenerate complementary subspace of descendants, see Appendix~\ref{app:virasoro-four}.\\[2pt]
\refstepcounter{finitecalculation}\label{calc:ORBIT-CLASSICAL}\thefinitecalculation&Confluent Weyl determinants of order equal to the rank evaluate
\cref{cert:classical-formula}, taking account of the root-length
bounds and both very-even labels. The calculations also verify
all stated claims that nilpotent orbits are not contained in
the associated varieties.\\[2pt]
\end{longtable}
}

\section{Associated varieties in small rank}\label{app:associated-varieties}
We list associated varieties for simple Lie algebras of rank at most six.
The last subsection gives the associated varieties for $E_7$ and $E_8$.
For the noncritical quasi-lisse levels in the tables, write
$$
 k+h^\vee=p/q>0,\qquad p,q\in\Z_{>0},\qquad \gcd(p,q)=1.
$$
An entry $\lambda$ or $T$ means $X_{L_k(\g)}=\overline{\OO_\lambda}$
or $\overline{\OO_T}$, respectively. Classical orbits are labelled by
partitions in the natural representation. Exponents denote multiplicities.
Exceptional orbits carry Bala--Carter labels. In the denominator column,
a range $a\text{--}b$ denotes all integers from $a$ to $b$. The label $0$ denotes the
zero orbit. The regular orbit has closure $\N(\g)$.

At admissible levels \eqref{eq:admissible-general}, the entries follow
from \cite[Theorems~5.7.1 and~5.8.1]{AAV}. They depend on $q$, with
$\gcd(p,q)=1$, and not on the admissible numerator $p$.
The nonadmissible tables list precisely the quasi-lisse parameters with
$X_{L_k(\g)}\subsetneq\N(\g)$. At other nonadmissible quasi-lisse parameters of the
indicated type, $X_{L_k(\g)}=\N(\g)$. One must first impose
\cref{thm:classification}: omitted parameters need not be quasi-lisse.
At non-quasi-lisse levels, the tables do not determine the nonnilpotent
part of $X_{L_k(\g)}$, except when $V^k(\g)$ is simple, as described below.

The entries combine the quasi-lisse classification with
\cref{prop:orbit-vacuum} and \eqref{eq:finite-weyl-certificate}.
For each orbit, \eqref{eq:table-certificate} determines whether the orbit is contained in the associated variety and excludes
every minimal orbit outside its closure. The evaluations and closure
data are in \cite{Certificates}.

For independently normalized low-rank algebras, $B_1=C_1=A_1$ and
$D_3=A_3$. The two rank-two tables $B_2$ and $C_2$ describe the same
orbits in different natural representations. Their orbit dimensions and
levels agree. Type $D$ has two orbits for a very-even partition. None of
the maximal partitions occurring below is very even, so no ambiguity
arises in the displayed labels. Both labels are checked among
the complementary orbits.

\subsection{Associated varieties at the critical level and when \texorpdfstring{$V^k(\g)$}{Vk(g)} is simple}
For the simple graded quotient $L_{-h^\vee}(\g)$, \cref{prop:critical} gives
$X_{L_{-h^\vee}(\g)}=\N(\g)$. The dimensions are
$$
\begin{array}{c|ccccccccc}
\g&A_n&B_n&C_n&D_n&E_6&E_7&E_8&F_4&G_2\\\hline
\dim\N(\g)&n(n+1)&2n^2&2n^2&2n(n-1)&72&126&240&48&12.
\end{array}
$$
The vertex algebra $L_{-h^\vee}(\g)$ has no Sugawara conformal
vector. At levels for which $V^k(\g)$ is simple, one has
$X_{L_k(\g)}=\g^*$.

\subsection{\texorpdfstring{Types $A$, $B$, $C$ and $D$}{Types A, B, C, D}}
In type $A_{N-1}$ the partition at an admissible level is
$(q^a,b)$, where $N=aq+b$ and $0\le b<q$. A zero part is omitted.
All quasi-lisse levels in type $A$ are admissible.

{\small\setlength{\tabcolsep}{4pt}
\begin{longtable}{@{}p{.10\textwidth}p{.43\textwidth}p{.30\textwidth}p{.09\textwidth}@{}}
\caption{Associated varieties at admissible levels: $A_{1}$, $A_{2}$, $A_{3}$, $A_{4}$, $A_{5}$, $A_{6}$.}\label{tab:av-A}\\
\toprule Type&$q$&Orbit or partition&$\dim X$\\\midrule\endfirsthead
\multicolumn{4}{l}{\tablename~\thetable\ (continued)}\\\toprule
Type&$q$&Orbit or partition&$\dim X$\\\midrule\endhead
\midrule\multicolumn{4}{r}{Continued on the next page}\\\endfoot\bottomrule\endlastfoot
$A_{1}$&$1$&$(1^{2})$&$0$\\
$A_{1}$&$q\ge 2$&$(2)$&$2$\\ \midrule
$A_{2}$&$1$&$(1^{3})$&$0$\\
$A_{2}$&$2$&$(2,1)$&$4$\\
$A_{2}$&$q\ge 3$&$(3)$&$6$\\\midrule
$A_{3}$&$1$&$(1^{4})$&$0$\\
$A_{3}$&$2$&$(2^{2})$&$8$\\
$A_{3}$&$3$&$(3,1)$&$10$\\
$A_{3}$&$q\ge 4$&$(4)$&$12$\\\midrule
$A_{4}$&$1$&$(1^{5})$&$0$\\
$A_{4}$&$2$&$(2^{2},1)$&$12$\\
$A_{4}$&$3$&$(3,2)$&$16$\\
$A_{4}$&$4$&$(4,1)$&$18$\\
$A_{4}$&$q\ge 5$&$(5)$&$20$\\\midrule
$A_{5}$&$1$&$(1^{6})$&$0$\\
$A_{5}$&$2$&$(2^{3})$&$18$\\
$A_{5}$&$3$&$(3^{2})$&$24$\\
$A_{5}$&$4$&$(4,2)$&$26$\\
$A_{5}$&$5$&$(5,1)$&$28$\\
$A_{5}$&$q\ge 6$&$(6)$&$30$\\\midrule
$A_{6}$&$1$&$(1^{7})$&$0$\\
$A_{6}$&$2$&$(2^{3},1)$&$24$\\
$A_{6}$&$3$&$(3^{2},1)$&$32$\\
$A_{6}$&$4$&$(4,3)$&$36$\\
$A_{6}$&$5$&$(5,2)$&$38$\\
$A_{6}$&$6$&$(6,1)$&$40$\\
$A_{6}$&$q\ge 7$&$(7)$&$42$\\
\end{longtable}
}
For $B_n,C_n,D_n$, the admissible ranges are, respectively,
$$
\begin{array}{c|cc}
 &q\text{ odd}&q\text{ even}\\\hline
B_n&p\ge2n-1&p\ge2n\\
C_n&p\ge n+1&p\ge2n\\
D_n&p\ge2n-2&p\ge2n-2.
\end{array}
$$
The inequalities in the denominator column include all larger denominators
of the specified parity. At the nonadmissible levels it suffices to use
the bounds in \eqref{eq:classical-tails}. In types $B,D$, $q>p$ also
gives $X=\N$ at quasi-lisse levels. No such extra assertion is needed
in type $C$.

For the nonadmissible type-$B$ table, we consider reduced parameters
satisfying
$$
 \begin{cases}
  2n\le3(p-1),&p,q\text{ both odd},\\
  p>n,\quad p\ne2n-1,&\text{otherwise}.
 \end{cases}
$$
These are precisely the type-$B$ conditions in
\cref{thm:classification}. Every omitted denominator in this parameter
range gives the whole nilpotent cone.

{\small\setlength{\tabcolsep}{4pt}
\begin{longtable}{@{}p{.10\textwidth}p{.43\textwidth}p{.30\textwidth}p{.09\textwidth}@{}}
\caption{Associated varieties at admissible levels: $B_{2}$, $B_{3}$, $B_{4}$, $B_{5}$, $B_{6}$.}\label{tab:av-B-adm}\\
\toprule Type&$q$&Orbit or partition&$\dim X$\\\midrule\endfirsthead
\multicolumn{4}{l}{\tablename~\thetable\ (continued)}\\\toprule
Type&$q$&Orbit or partition&$\dim X$\\\midrule\endhead
\midrule\multicolumn{4}{r}{Continued on the next page}\\\endfoot\bottomrule\endlastfoot
$B_{2}$&$1$&$(1^{5})$&$0$\\
$B_{2}$&$2$&$(2^{2},1)$&$4$\\
$B_{2}$&$3,4$&$(3,1^{2})$&$6$\\
$B_{2}$&$\begin{gathered}q\ge 5\ (q\text{ odd}),\\q\ge 6\ (q\text{ even})\end{gathered}$&$(5)$&$8$\\ \midrule
$B_{3}$&$1$&$(1^{7})$&$0$\\
$B_{3}$&$2$&$(2^{2},1^{3})$&$8$\\
$B_{3}$&$3,4$&$(3^{2},1)$&$14$\\
$B_{3}$&$5,6$&$(5,1^{2})$&$16$\\
$B_{3}$&$\begin{gathered}q\ge 7\ (q\text{ odd}),\\q\ge 8\ (q\text{ even})\end{gathered}$&$(7)$&$18$\\ \midrule
$B_{4}$&$1$&$(1^{9})$&$0$\\
$B_{4}$&$2$&$(2^{4},1)$&$16$\\
$B_{4}$&$3$&$(3^{3})$&$24$\\
$B_{4}$&$4$&$(4^{2},1)$&$26$\\
$B_{4}$&$5,6$&$(5,3,1)$&$28$\\
$B_{4}$&$7,8$&$(7,1^{2})$&$30$\\
$B_{4}$&$\begin{gathered}q\ge 9\ (q\text{ odd}),\\q\ge 10\ (q\text{ even})\end{gathered}$&$(9)$&$32$\\ \midrule
$B_{5}$&$1$&$(1^{11})$&$0$\\
$B_{5}$&$2$&$(2^{4},1^{3})$&$24$\\
$B_{5}$&$3$&$(3^{3},1^{2})$&$36$\\
$B_{5}$&$4$&$(4^{2},3)$&$40$\\
$B_{5}$&$5,6$&$(5^{2},1)$&$44$\\
$B_{5}$&$7,8$&$(7,3,1)$&$46$\\
$B_{5}$&$9,10$&$(9,1^{2})$&$48$\\
$B_{5}$&$\begin{gathered}q\ge 11\ (q\text{ odd}),\\q\ge 12\ (q\text{ even})\end{gathered}$&$(11)$&$50$\\ \midrule
$B_{6}$&$1$&$(1^{13})$&$0$\\
$B_{6}$&$2$&$(2^{6},1)$&$36$\\
$B_{6}$&$3$&$(3^{4},1)$&$52$\\
$B_{6}$&$4$&$(4^{2},3,1^{2})$&$56$\\
$B_{6}$&$5$&$(5^{2},3)$&$62$\\
$B_{6}$&$6$&$(6^{2},1)$&$64$\\
$B_{6}$&$7,8$&$(7,5,1)$&$66$\\
$B_{6}$&$9,10$&$(9,3,1)$&$68$\\
$B_{6}$&$11,12$&$(11,1^{2})$&$70$\\
$B_{6}$&$\begin{gathered}q\ge 13\ (q\text{ odd}),\\q\ge 14\ (q\text{ even})\end{gathered}$&$(13)$&$72$\\
\end{longtable}
}
{\small\setlength{\tabcolsep}{4pt}
\begin{longtable}{@{}p{.08\textwidth}p{.06\textwidth}p{.34\textwidth}p{.26\textwidth}p{.14\textwidth}@{}}
\caption{Proper orbit closures at nonadmissible quasi-lisse levels: $B_{2}$, $B_{3}$, $B_{4}$, $B_{5}$, $B_{6}$.}\label{tab:av-B-nonadm}\\
\toprule Type&$p$&$q$&Orbit or partition&$\dim X$\\\midrule\endfirsthead
\multicolumn{5}{l}{\tablename~\thetable\ (continued)}\\\toprule
Type&$p$&$q$&Orbit or partition&$\dim X$\\\midrule\endhead
\midrule\multicolumn{5}{r}{Continued on the next page}\\\endfoot\bottomrule\endlastfoot
$B_{3}$&$3$&$1$&$(3,1^{4})$&$10$\\
$B_{3}$&$4$&$1$&$(3,1^{4})$&$10$\\
$B_{3}$&$4$&$3$&$(5,1^{2})$&$16$\\ \midrule
$B_{4}$&$5$&$1$&$(3,1^{6})$&$14$\\
$B_{4}$&$5$&$2,3$&$(5,3,1)$&$28$\\
$B_{4}$&$5$&$4$&$(7,1^{2})$&$30$\\
$B_{4}$&$6$&$1$&$(3,1^{6})$&$14$\\
$B_{4}$&$6$&$5$&$(7,1^{2})$&$30$\\ \midrule
$B_{5}$&$5$&$1$&$(3,2^{2},1^{4})$&$28$\\
$B_{5}$&$5$&$3$&$(7,3,1)$&$46$\\
$B_{5}$&$6$&$1$&$(3,2^{2},1^{4})$&$28$\\
$B_{5}$&$6$&$5$&$(9,1^{2})$&$48$\\
$B_{5}$&$7$&$1$&$(3,1^{8})$&$18$\\
$B_{5}$&$7$&$2$&$(5,3,1^{3})$&$40$\\
$B_{5}$&$7$&$3$&$(5,3^{2})$&$42$\\
$B_{5}$&$7$&$4,5$&$(7,3,1)$&$46$\\
$B_{5}$&$7$&$6$&$(9,1^{2})$&$48$\\
$B_{5}$&$8$&$1$&$(3,1^{8})$&$18$\\
$B_{5}$&$8$&$3$&$(5,3^{2})$&$42$\\
$B_{5}$&$8$&$5$&$(7,3,1)$&$46$\\
$B_{5}$&$8$&$7$&$(9,1^{2})$&$48$\\ \midrule
$B_{6}$&$5$&$1$&$(3^{3},1^{4})$&$48$\\
$B_{6}$&$5$&$3$&$(9,3,1)$&$68$\\
$B_{6}$&$7$&$1$&$(3,2^{2},1^{6})$&$36$\\
$B_{6}$&$7$&$2$&$(5,3^{2},1^{2})$&$58$\\
$B_{6}$&$7$&$3$&$(7,3^{2})$&$64$\\
$B_{6}$&$7$&$4,5$&$(9,3,1)$&$68$\\
$B_{6}$&$7$&$6$&$(11,1^{2})$&$70$\\
$B_{6}$&$8$&$1$&$(3,2^{2},1^{6})$&$36$\\
$B_{6}$&$8$&$3$&$(7,3^{2})$&$64$\\
$B_{6}$&$8$&$5$&$(9,3,1)$&$68$\\
$B_{6}$&$8$&$7$&$(11,1^{2})$&$70$\\
$B_{6}$&$9$&$1$&$(3,1^{10})$&$22$\\
$B_{6}$&$9$&$2$&$(5,3,2^{2},1)$&$56$\\
$B_{6}$&$9$&$4,5$&$(7,5,1)$&$66$\\
$B_{6}$&$9$&$7$&$(9,3,1)$&$68$\\
$B_{6}$&$9$&$8$&$(11,1^{2})$&$70$\\
$B_{6}$&$10$&$1$&$(3,1^{10})$&$22$\\
$B_{6}$&$10$&$3$&$(5,3^{2},1^{2})$&$58$\\
$B_{6}$&$10$&$7$&$(9,3,1)$&$68$\\
$B_{6}$&$10$&$9$&$(11,1^{2})$&$70$\\
\end{longtable}
}
{\small\setlength{\tabcolsep}{4pt}
\begin{longtable}{@{}p{.10\textwidth}p{.43\textwidth}p{.30\textwidth}p{.09\textwidth}@{}}
\caption{Associated varieties at admissible levels: $C_{2}$, $C_{3}$, $C_{4}$, $C_{5}$, $C_{6}$.}\label{tab:av-C-adm}\\
\toprule Type&$q$&Orbit or partition&$\dim X$\\\midrule\endfirsthead
\multicolumn{4}{l}{\tablename~\thetable\ (continued)}\\\toprule
Type&$q$&Orbit or partition&$\dim X$\\\midrule\endhead
\midrule\multicolumn{4}{r}{Continued on the next page}\\\endfoot\bottomrule\endlastfoot
$C_{2}$&$1$&$(1^{4})$&$0$\\
$C_{2}$&$2$&$(2,1^{2})$&$4$\\
$C_{2}$&$3,4$&$(2^{2})$&$6$\\
$C_{2}$&$\begin{gathered}q\ge 5\ (q\text{ odd}),\\q\ge 6\ (q\text{ even})\end{gathered}$&$(4)$&$8$\\ \midrule
$C_{3}$&$1$&$(1^{6})$&$0$\\
$C_{3}$&$2$&$(2,1^{4})$&$6$\\
$C_{3}$&$3$&$(3^{2})$&$14$\\
$C_{3}$&$4$&$(2^{3})$&$12$\\
$C_{3}$&$5,6,8$&$(4,2)$&$16$\\
$C_{3}$&$\begin{gathered}q\ge 7\ (q\text{ odd}),\\q\ge 10\ (q\text{ even})\end{gathered}$&$(6)$&$18$\\ \midrule
$C_{4}$&$1$&$(1^{8})$&$0$\\
$C_{4}$&$2$&$(2,1^{6})$&$8$\\
$C_{4}$&$3$&$(3^{2},2)$&$24$\\
$C_{4}$&$4$&$(2^{4})$&$20$\\
$C_{4}$&$5,8$&$(4^{2})$&$28$\\
$C_{4}$&$6$&$(4,2^{2})$&$26$\\
$C_{4}$&$7,10,12$&$(6,2)$&$30$\\
$C_{4}$&$\begin{gathered}q\ge 9\ (q\text{ odd}),\\q\ge 14\ (q\text{ even})\end{gathered}$&$(8)$&$32$\\ \midrule
$C_{5}$&$1$&$(1^{10})$&$0$\\
$C_{5}$&$2$&$(2,1^{8})$&$10$\\
$C_{5}$&$3$&$(3^{2},2^{2})$&$36$\\
$C_{5}$&$4$&$(2^{5})$&$30$\\
$C_{5}$&$5$&$(5^{2})$&$44$\\
$C_{5}$&$6$&$(4,3^{2})$&$40$\\
$C_{5}$&$7,10,12$&$(6,4)$&$46$\\
$C_{5}$&$8$&$(4^{2},2)$&$42$\\
$C_{5}$&$9,14,16$&$(8,2)$&$48$\\
$C_{5}$&$\begin{gathered}q\ge 11\ (q\text{ odd}),\\q\ge 18\ (q\text{ even})\end{gathered}$&$(10)$&$50$\\ \midrule
$C_{6}$&$1$&$(1^{12})$&$0$\\
$C_{6}$&$2$&$(2,1^{10})$&$12$\\
$C_{6}$&$3$&$(3^{4})$&$52$\\
$C_{6}$&$4$&$(2^{6})$&$42$\\
$C_{6}$&$5$&$(5^{2},2)$&$62$\\
$C_{6}$&$6$&$(4,3^{2},2)$&$56$\\
$C_{6}$&$7,12$&$(6^{2})$&$66$\\
$C_{6}$&$8$&$(4^{3})$&$60$\\
$C_{6}$&$9,14,16$&$(8,4)$&$68$\\
$C_{6}$&$10$&$(6,4,2)$&$64$\\
$C_{6}$&$11,18,20$&$(10,2)$&$70$\\
$C_{6}$&$\begin{gathered}q\ge 13\ (q\text{ odd}),\\q\ge 22\ (q\text{ even})\end{gathered}$&$(12)$&$72$\\
\end{longtable}
}
{\small\setlength{\tabcolsep}{4pt}
\begin{longtable}{@{}p{.08\textwidth}p{.06\textwidth}p{.34\textwidth}p{.26\textwidth}p{.14\textwidth}@{}}
\caption{Proper orbit closures at nonadmissible quasi-lisse levels: $C_{2}$, $C_{3}$, $C_{4}$, $C_{5}$, $C_{6}$.}\label{tab:av-C-nonadm}\\
\toprule Type&$p$&$q$&Orbit or partition&$\dim X$\\\midrule\endfirsthead
\multicolumn{5}{l}{\tablename~\thetable\ (continued)}\\\toprule
Type&$p$&$q$&Orbit or partition&$\dim X$\\\midrule\endhead
\midrule\multicolumn{5}{r}{Continued on the next page}\\\endfoot\bottomrule\endlastfoot
$C_{3}$&$5$&$2$&$(2^{3})$&$12$\\
$C_{3}$&$5$&$4$&$(4,2)$&$16$\\ \midrule
$C_{4}$&$3$&$1$&$(2^{4})$&$20$\\
$C_{4}$&$7$&$2$&$(2^{3},1^{2})$&$18$\\
$C_{4}$&$7$&$4$&$(4,2^{2})$&$26$\\
$C_{4}$&$7$&$6,8$&$(6,2)$&$30$\\ \midrule
$C_{5}$&$4$&$1$&$(2^{4},1^{2})$&$28$\\
$C_{5}$&$4$&$3$&$(6,4)$&$46$\\
$C_{5}$&$4$&$5$&$(8,2)$&$48$\\
$C_{5}$&$7$&$2$&$(2^{5})$&$30$\\
$C_{5}$&$7$&$4$&$(4^{2},2)$&$42$\\
$C_{5}$&$7$&$6$&$(6,4)$&$46$\\
$C_{5}$&$7$&$8$&$(8,2)$&$48$\\
$C_{5}$&$9$&$2$&$(2^{3},1^{4})$&$24$\\
$C_{5}$&$9$&$4$&$(4,2^{3})$&$38$\\
$C_{5}$&$9$&$8$&$(6,4)$&$46$\\
$C_{5}$&$9$&$10$&$(8,2)$&$48$\\ \midrule
$C_{6}$&$4$&$1$&$(2^{6})$&$42$\\
$C_{6}$&$4$&$3$&$(6^{2})$&$66$\\
$C_{6}$&$4$&$5$&$(10,2)$&$70$\\
$C_{6}$&$5$&$1$&$(2^{4},1^{4})$&$36$\\
$C_{6}$&$5$&$3$&$(6,4,2)$&$64$\\
$C_{6}$&$5$&$7$&$(10,2)$&$70$\\
$C_{6}$&$9$&$2$&$(2^{5},1^{2})$&$40$\\
$C_{6}$&$9$&$4$&$(4^{2},2^{2})$&$58$\\
$C_{6}$&$9$&$8$&$(8,4)$&$68$\\
$C_{6}$&$9$&$10$&$(10,2)$&$70$\\
$C_{6}$&$11$&$2$&$(2^{3},1^{6})$&$30$\\
$C_{6}$&$11$&$4$&$(4,2^{4})$&$52$\\
$C_{6}$&$11$&$6$&$(6,3^{2})$&$62$\\
$C_{6}$&$11$&$8$&$(6,4,2)$&$64$\\
$C_{6}$&$11$&$10,12$&$(8,4)$&$68$\\
$C_{6}$&$11$&$14,16$&$(10,2)$&$70$\\
\end{longtable}
}
{\small\setlength{\tabcolsep}{4pt}
\begin{longtable}{@{}p{.10\textwidth}p{.43\textwidth}p{.30\textwidth}p{.09\textwidth}@{}}
\caption{Associated varieties at admissible levels: $D_{4}$, $D_{5}$, $D_{6}$.}\label{tab:av-D-adm}\\
\toprule Type&$q$&Orbit or partition&$\dim X$\\\midrule\endfirsthead
\multicolumn{4}{l}{\tablename~\thetable\ (continued)}\\\toprule
Type&$q$&Orbit or partition&$\dim X$\\\midrule\endhead
\midrule\multicolumn{4}{r}{Continued on the next page}\\\endfoot\bottomrule\endlastfoot
$D_{4}$&$1$&$(1^{8})$&$0$\\
$D_{4}$&$2$&$(3,2^{2},1)$&$16$\\
$D_{4}$&$3$&$(3^{2},1^{2})$&$18$\\
$D_{4}$&$4,5$&$(5,3)$&$22$\\
$D_{4}$&$q\ge 6$&$(7,1)$&$24$\\ \midrule
$D_{5}$&$1$&$(1^{10})$&$0$\\
$D_{5}$&$2$&$(3,2^{2},1^{3})$&$24$\\
$D_{5}$&$3$&$(3^{3},1)$&$30$\\
$D_{5}$&$4$&$(5,3,1^{2})$&$34$\\
$D_{5}$&$5$&$(5^{2})$&$36$\\
$D_{5}$&$6,7$&$(7,3)$&$38$\\
$D_{5}$&$q\ge 8$&$(9,1)$&$40$\\ \midrule
$D_{6}$&$1$&$(1^{12})$&$0$\\
$D_{6}$&$2$&$(3,2^{4},1)$&$36$\\
$D_{6}$&$3$&$(3^{4})$&$44$\\
$D_{6}$&$4$&$(5,3^{2},1)$&$50$\\
$D_{6}$&$5$&$(5^{2},1^{2})$&$52$\\
$D_{6}$&$6,7$&$(7,5)$&$56$\\
$D_{6}$&$8,9$&$(9,3)$&$58$\\
$D_{6}$&$q\ge 10$&$(11,1)$&$60$\\
\end{longtable}
}
{\small\setlength{\tabcolsep}{4pt}
\begin{longtable}{@{}p{.08\textwidth}p{.06\textwidth}p{.34\textwidth}p{.26\textwidth}p{.14\textwidth}@{}}
\caption{Proper orbit closures at nonadmissible quasi-lisse levels: $D_{4}$, $D_{5}$, $D_{6}$.}\label{tab:av-D-nonadm}\\
\toprule Type&$p$&$q$&Orbit or partition&$\dim X$\\\midrule\endfirsthead
\multicolumn{5}{l}{\tablename~\thetable\ (continued)}\\\toprule
Type&$p$&$q$&Orbit or partition&$\dim X$\\\midrule\endhead
\midrule\multicolumn{5}{r}{Continued on the next page}\\\endfoot\bottomrule\endlastfoot
$D_{4}$&$4$&$1$&$(2^{2},1^{4})$&$10$\\
$D_{4}$&$4$&$3$&$(5,3)$&$22$\\
$D_{4}$&$5$&$1$&$(2^{2},1^{4})$&$10$\\
$D_{4}$&$5$&$2$&$(3^{2},1^{2})$&$18$\\
$D_{4}$&$5$&$3$&$(5,3)$&$22$\\ \midrule
$D_{5}$&$6$&$1$&$(2^{2},1^{6})$&$14$\\
$D_{5}$&$6$&$5$&$(7,3)$&$38$\\
$D_{5}$&$7$&$1$&$(2^{2},1^{6})$&$14$\\
$D_{5}$&$7$&$2$&$(3^{3},1)$&$30$\\
$D_{5}$&$7$&$3$&$(5,3,1^{2})$&$34$\\
$D_{5}$&$7$&$4,5$&$(7,3)$&$38$\\ \midrule
$D_{6}$&$6$&$1$&$(2^{4},1^{4})$&$28$\\
$D_{6}$&$6$&$5$&$(9,3)$&$58$\\
$D_{6}$&$7$&$1$&$(2^{4},1^{4})$&$28$\\
$D_{6}$&$7$&$2$&$(4^{2},3,1)$&$48$\\
$D_{6}$&$7$&$3$&$(5^{2},1^{2})$&$52$\\
$D_{6}$&$7$&$4$&$(7,5)$&$56$\\
$D_{6}$&$7$&$5$&$(9,3)$&$58$\\
$D_{6}$&$8$&$1$&$(2^{2},1^{8})$&$18$\\
$D_{6}$&$8$&$3$&$(5,3^{2},1)$&$50$\\
$D_{6}$&$8$&$5$&$(7,5)$&$56$\\
$D_{6}$&$8$&$7$&$(9,3)$&$58$\\
$D_{6}$&$9$&$1$&$(2^{2},1^{8})$&$18$\\
$D_{6}$&$9$&$2$&$(3^{3},1^{3})$&$42$\\
$D_{6}$&$9$&$4$&$(7,3,1^{2})$&$54$\\
$D_{6}$&$9$&$5$&$(7,5)$&$56$\\
$D_{6}$&$9$&$7$&$(9,3)$&$58$\\
\end{longtable}
}
\subsection{Exceptional types of rank at most six}
The admissible ranges for $E_6$, $F_4$ and $G_2$ are
$$
\begin{array}{c|l}
E_6&p\ge12,\\
F_4&p\ge9\ (q\text{ odd}),\qquad p\ge12\ (q\text{ even}),\\
G_2&p\ge4\ (3\nmid q),\qquad p\ge6\ (3\mid q).
\end{array}
$$
For the nonadmissible numerators of $E_6$, denominators $q>p$ give
$\N(E_6)$. For $F_4$ and $G_2$, the corresponding sufficient bounds
are $p/q<3/4$ and $p/q<2/3$. The remaining finite values are included
below. In particular, the nonadmissible tables may omit a denominator
smaller than these bounds when its exact calculation already gives the
regular orbit.

{\small\setlength{\tabcolsep}{4pt}
\begin{longtable}{@{}p{.10\textwidth}p{.43\textwidth}p{.30\textwidth}p{.09\textwidth}@{}}
\caption{Associated varieties at admissible levels: $E_{6}$.}\label{tab:av-E6-adm}\\
\toprule Type&$q$&Orbit or partition&$\dim X$\\\midrule\endfirsthead
\multicolumn{4}{l}{\tablename~\thetable\ (continued)}\\\toprule
Type&$q$&Orbit or partition&$\dim X$\\\midrule\endhead
\midrule\multicolumn{4}{r}{Continued on the next page}\\\endfoot\bottomrule\endlastfoot
$E_{6}$&$1$&$0$&$0$\\
$E_{6}$&$2$&$3A_{1}$&$40$\\
$E_{6}$&$3$&$2A_{2}+A_{1}$&$54$\\
$E_{6}$&$4$&$D_{4}(a_{1})$&$58$\\
$E_{6}$&$5$&$A_{4}+A_{1}$&$62$\\
$E_{6}$&$6,7$&$E_{6}(a_{3})$&$66$\\
$E_{6}$&$8$&$D_{5}$&$68$\\
$E_{6}$&$9\text{--}11$&$E_{6}(a_{1})$&$70$\\
$E_{6}$&$q\ge 12$&$E_{6}$&$72$\\
\end{longtable}
}
{\small\setlength{\tabcolsep}{4pt}
\begin{longtable}{@{}p{.08\textwidth}p{.06\textwidth}p{.34\textwidth}p{.26\textwidth}p{.14\textwidth}@{}}
\caption{Proper orbit closures at nonadmissible quasi-lisse levels: $E_{6}$.}\label{tab:av-E6-nonadm}\\
\toprule Type&$p$&$q$&Orbit or partition&$\dim X$\\\midrule\endfirsthead
\multicolumn{5}{l}{\tablename~\thetable\ (continued)}\\\toprule
Type&$p$&$q$&Orbit or partition&$\dim X$\\\midrule\endhead
\midrule\multicolumn{5}{r}{Continued on the next page}\\\endfoot\bottomrule\endlastfoot
$E_{6}$&$6$&$1$&$A_{2}$&$42$\\
$E_{6}$&$6$&$5$&$E_{6}(a_{1})$&$70$\\
$E_{6}$&$7$&$1$&$A_{2}$&$42$\\
$E_{6}$&$7$&$2$&$D_{4}(a_{1})$&$58$\\
$E_{6}$&$7$&$3$&$E_{6}(a_{3})$&$66$\\
$E_{6}$&$7$&$4$&$D_{5}$&$68$\\
$E_{6}$&$7$&$5$&$E_{6}(a_{1})$&$70$\\
$E_{6}$&$9$&$1$&$A_{1}$&$22$\\
$E_{6}$&$9$&$2$&$A_{2}+2A_{1}$&$50$\\
$E_{6}$&$9$&$4$&$D_{5}(a_{1})$&$64$\\
$E_{6}$&$9$&$5$&$E_{6}(a_{3})$&$66$\\
$E_{6}$&$9$&$7,8$&$E_{6}(a_{1})$&$70$\\
$E_{6}$&$10$&$1$&$A_{1}$&$22$\\
$E_{6}$&$10$&$3$&$D_{4}(a_{1})$&$58$\\
$E_{6}$&$10$&$7$&$E_{6}(a_{1})$&$70$\\
$E_{6}$&$11$&$1$&$A_{1}$&$22$\\
$E_{6}$&$11$&$2$&$A_{2}+2A_{1}$&$50$\\
$E_{6}$&$11$&$3$&$D_{4}(a_{1})$&$58$\\
$E_{6}$&$11$&$4$&$D_{5}(a_{1})$&$64$\\
$E_{6}$&$11$&$5$&$E_{6}(a_{3})$&$66$\\
$E_{6}$&$11$&$6$&$D_{5}$&$68$\\
$E_{6}$&$11$&$7,8$&$E_{6}(a_{1})$&$70$\\
\end{longtable}
}
{\small\setlength{\tabcolsep}{4pt}
\begin{longtable}{@{}p{.10\textwidth}p{.43\textwidth}p{.30\textwidth}p{.09\textwidth}@{}}
\caption{Associated varieties at admissible levels: $F_{4}$.}\label{tab:av-F4-adm}\\
\toprule Type&$q$&Orbit or partition&$\dim X$\\\midrule\endfirsthead
\multicolumn{4}{l}{\tablename~\thetable\ (continued)}\\\toprule
Type&$q$&Orbit or partition&$\dim X$\\\midrule\endhead
\midrule\multicolumn{4}{r}{Continued on the next page}\\\endfoot\bottomrule\endlastfoot
$F_{4}$&$1$&$0$&$0$\\
$F_{4}$&$2$&$A_{1}$&$16$\\
$F_{4}$&$3$&$\widetilde A_{2}+A_{1}$&$36$\\
$F_{4}$&$4$&$A_{2}+\widetilde A_{1}$&$34$\\
$F_{4}$&$5,6$&$F_{4}(a_{3})$&$40$\\
$F_{4}$&$7,10$&$F_{4}(a_{2})$&$44$\\
$F_{4}$&$8$&$B_{3}$&$42$\\
$F_{4}$&$9,11,12,14,16$&$F_{4}(a_{1})$&$46$\\
$F_{4}$&$\begin{gathered}q\ge13\ (q\text{ odd}),\\q\ge18\ (q\text{ even})\end{gathered}$&$F_{4}$&$48$\\
\end{longtable}
}
{\small\setlength{\tabcolsep}{4pt}
\begin{longtable}{@{}p{.08\textwidth}p{.06\textwidth}p{.34\textwidth}p{.26\textwidth}p{.14\textwidth}@{}}
\caption{Proper orbit closures at nonadmissible quasi-lisse levels: $F_{4}$.}\label{tab:av-F4-nonadm}\\
\toprule Type&$p$&$q$&Orbit or partition&$\dim X$\\\midrule\endfirsthead
\multicolumn{5}{l}{\tablename~\thetable\ (continued)}\\\toprule
Type&$p$&$q$&Orbit or partition&$\dim X$\\\midrule\endhead
\midrule\multicolumn{5}{r}{Continued on the next page}\\\endfoot\bottomrule\endlastfoot
$F_{4}$&$3$&$1$&$F_{4}(a_{3})$&$40$\\
$F_{4}$&$5$&$1$&$A_{1}+\widetilde A_{1}$&$28$\\
$F_{4}$&$5$&$2$&$F_{4}(a_{3})$&$40$\\
$F_{4}$&$5$&$3$&$F_{4}(a_{2})$&$44$\\
$F_{4}$&$5$&$4$&$F_{4}(a_{1})$&$46$\\
$F_{4}$&$6$&$1$&$\widetilde A_{1}$&$22$\\
$F_{4}$&$6$&$5,7$&$F_{4}(a_{1})$&$46$\\
$F_{4}$&$7$&$1$&$\widetilde A_{1}$&$22$\\
$F_{4}$&$7$&$2$&$B_{2}$&$36$\\
$F_{4}$&$7$&$3$&$F_{4}(a_{3})$&$40$\\
$F_{4}$&$7$&$4$&$B_{3}$&$42$\\
$F_{4}$&$7$&$5,6,8$&$F_{4}(a_{1})$&$46$\\
$F_{4}$&$8$&$1$&$\widetilde A_{1}$&$22$\\
$F_{4}$&$8$&$3$&$F_{4}(a_{3})$&$40$\\
$F_{4}$&$8$&$5,7$&$F_{4}(a_{1})$&$46$\\
$F_{4}$&$9$&$2$&$A_{2}$&$30$\\
$F_{4}$&$9$&$4$&$F_{4}(a_{3})$&$40$\\
$F_{4}$&$9$&$8,10$&$F_{4}(a_{1})$&$46$\\
$F_{4}$&$11$&$2$&$A_{2}$&$30$\\
$F_{4}$&$11$&$4$&$F_{4}(a_{3})$&$40$\\
$F_{4}$&$11$&$6$&$F_{4}(a_{2})$&$44$\\
$F_{4}$&$11$&$8,10$&$F_{4}(a_{1})$&$46$\\
\end{longtable}
}
{\small\setlength{\tabcolsep}{4pt}
\begin{longtable}{@{}p{.10\textwidth}p{.43\textwidth}p{.30\textwidth}p{.09\textwidth}@{}}
\caption{Associated varieties at admissible levels: $G_{2}$.}\label{tab:av-G2-adm}\\
\toprule Type&$q$&Orbit or partition&$\dim X$\\\midrule\endfirsthead
\multicolumn{4}{l}{\tablename~\thetable\ (continued)}\\\toprule
Type&$q$&Orbit or partition&$\dim X$\\\midrule\endhead
\midrule\multicolumn{4}{r}{Continued on the next page}\\\endfoot\bottomrule\endlastfoot
$G_{2}$&$1$&$0$&$0$\\
$G_{2}$&$2$&$\widetilde A_{1}$&$8$\\
$G_{2}$&$3$&$A_{1}$&$6$\\
$G_{2}$&$4\text{--}6,9$&$G_{2}(a_{1})$&$10$\\
$G_{2}$&$\begin{gathered}q\ge7\ (3\nmid q),\\q\ge12\ (3\mid q)\end{gathered}$&$G_{2}$&$12$\\
\end{longtable}
}
{\small\setlength{\tabcolsep}{4pt}
\begin{longtable}{@{}p{.08\textwidth}p{.06\textwidth}p{.34\textwidth}p{.26\textwidth}p{.14\textwidth}@{}}
\caption{Proper orbit closures at nonadmissible quasi-lisse levels: $G_{2}$.}\label{tab:av-G2-nonadm}\\
\toprule Type&$p$&$q$&Orbit or partition&$\dim X$\\\midrule\endfirsthead
\multicolumn{5}{l}{\tablename~\thetable\ (continued)}\\\toprule
Type&$p$&$q$&Orbit or partition&$\dim X$\\\midrule\endhead
\midrule\multicolumn{5}{r}{Continued on the next page}\\\endfoot\bottomrule\endlastfoot
$G_{2}$&$2$&$1$&$G_{2}(a_{1})$&$10$\\
$G_{2}$&$3$&$1,2$&$G_{2}(a_{1})$&$10$\\
$G_{2}$&$4$&$3$&$G_{2}(a_{1})$&$10$\\
$G_{2}$&$5$&$3$&$G_{2}(a_{1})$&$10$\\
\end{longtable}
}
\subsection{\texorpdfstring{Types $E_7$ and $E_8$}{Types E7 and E8}}
For completeness, we give the remaining exceptional types. The admissible
ranges are $p\ge18$ for $E_7$ and $p\ge30$ for $E_8$.
At the nonadmissible quasi-lisse numerators, every $q>p$ gives the full
nilpotent cone. The tables of proper orbit closures therefore require
only $q<p$, with $\gcd(p,q)=1$.

{\small\setlength{\tabcolsep}{4pt}
\begin{longtable}{@{}p{.10\textwidth}p{.43\textwidth}p{.30\textwidth}p{.09\textwidth}@{}}
\caption{Associated varieties at admissible levels: $E_{7}$.}\label{tab:av-E7-adm}\\
\toprule Type&$q$&Orbit or partition&$\dim X$\\\midrule\endfirsthead
\multicolumn{4}{l}{\tablename~\thetable\ (continued)}\\\toprule
Type&$q$&Orbit or partition&$\dim X$\\\midrule\endhead
\midrule\multicolumn{4}{r}{Continued on the next page}\\\endfoot\bottomrule\endlastfoot
$E_{7}$&$1$&$0$&$0$\\
$E_{7}$&$2$&$4A_{1}$&$70$\\
$E_{7}$&$3$&$2A_{2}+A_{1}$&$90$\\
$E_{7}$&$4$&$A_{3}+A_{2}+A_{1}$&$100$\\
$E_{7}$&$5$&$A_{4}+A_{2}$&$106$\\
$E_{7}$&$6$&$E_{7}(a_{5})$&$112$\\
$E_{7}$&$7$&$A_{6}$&$114$\\
$E_{7}$&$8$&$E_{7}(a_{4})$&$116$\\
$E_{7}$&$9$&$E_{6}(a_{1})$&$118$\\
$E_{7}$&$10,11$&$E_{7}(a_{3})$&$120$\\
$E_{7}$&$12,13$&$E_{7}(a_{2})$&$122$\\
$E_{7}$&$14\text{--}17$&$E_{7}(a_{1})$&$124$\\
$E_{7}$&$q\ge 18$&$E_{7}$&$126$\\
\end{longtable}
}
{\small\setlength{\tabcolsep}{4pt}
\begin{longtable}{@{}p{.08\textwidth}p{.06\textwidth}p{.34\textwidth}p{.26\textwidth}p{.14\textwidth}@{}}
\caption{Proper orbit closures at nonadmissible quasi-lisse levels: $E_{7}$.}\label{tab:av-E7-nonadm}\\
\toprule Type&$p$&$q$&Orbit or partition&$\dim X$\\\midrule\endfirsthead
\multicolumn{5}{l}{\tablename~\thetable\ (continued)}\\\toprule
Type&$p$&$q$&Orbit or partition&$\dim X$\\\midrule\endhead
\midrule\multicolumn{5}{r}{Continued on the next page}\\\endfoot\bottomrule\endlastfoot
$E_{7}$&$6$&$1$&$D_{4}(a_{1})$&$94$\\
$E_{7}$&$6$&$5$&$E_{7}(a_{1})$&$124$\\
$E_{7}$&$8$&$1$&$A_{2}+2A_{1}$&$82$\\
$E_{7}$&$8$&$3$&$E_{7}(a_{4})$&$116$\\
$E_{7}$&$8$&$5$&$E_{7}(a_{2})$&$122$\\
$E_{7}$&$8$&$7$&$E_{7}(a_{1})$&$124$\\
$E_{7}$&$10$&$1$&$A_{2}$&$66$\\
$E_{7}$&$10$&$3$&$E_{6}(a_{3})$&$110$\\
$E_{7}$&$10$&$7$&$E_{7}(a_{2})$&$122$\\
$E_{7}$&$10$&$9$&$E_{7}(a_{1})$&$124$\\
$E_{7}$&$11$&$1$&$A_{2}$&$66$\\
$E_{7}$&$11$&$2$&$A_{3}+A_{2}+A_{1}$&$100$\\
$E_{7}$&$11$&$3$&$E_{6}(a_{3})$&$110$\\
$E_{7}$&$11$&$4$&$E_{7}(a_{4})$&$116$\\
$E_{7}$&$11$&$5$&$E_{6}(a_{1})$&$118$\\
$E_{7}$&$11$&$6,7$&$E_{7}(a_{2})$&$122$\\
$E_{7}$&$11$&$8,9$&$E_{7}(a_{1})$&$124$\\
$E_{7}$&$12$&$1$&$2A_{1}$&$52$\\
$E_{7}$&$12$&$5$&$E_{7}(a_{4})$&$116$\\
$E_{7}$&$12$&$7$&$E_{7}(a_{3})$&$120$\\
$E_{7}$&$12$&$11$&$E_{7}(a_{1})$&$124$\\
$E_{7}$&$13$&$1$&$2A_{1}$&$52$\\
$E_{7}$&$13$&$2$&$A_{3}+2A_{1}$&$94$\\
$E_{7}$&$13$&$3$&$A_{4}+A_{2}$&$106$\\
$E_{7}$&$13$&$4$&$E_{7}(a_{5})$&$112$\\
$E_{7}$&$13$&$5$&$E_{7}(a_{4})$&$116$\\
$E_{7}$&$13$&$6,7$&$E_{7}(a_{3})$&$120$\\
$E_{7}$&$13$&$8$&$E_{7}(a_{2})$&$122$\\
$E_{7}$&$13$&$9\text{--}11$&$E_{7}(a_{1})$&$124$\\
$E_{7}$&$14$&$1$&$A_{1}$&$34$\\
$E_{7}$&$14$&$3$&$A_{3}+A_{2}+A_{1}$&$100$\\
$E_{7}$&$14$&$5$&$E_{7}(a_{5})$&$112$\\
$E_{7}$&$14$&$9$&$E_{7}(a_{2})$&$122$\\
$E_{7}$&$14$&$11,13$&$E_{7}(a_{1})$&$124$\\
$E_{7}$&$15$&$1$&$A_{1}$&$34$\\
$E_{7}$&$15$&$2$&$A_{2}+3A_{1}$&$84$\\
$E_{7}$&$15$&$4$&$D_{5}(a_{1})+A_{1}$&$108$\\
$E_{7}$&$15$&$7$&$E_{6}(a_{1})$&$118$\\
$E_{7}$&$15$&$8$&$E_{7}(a_{3})$&$120$\\
$E_{7}$&$15$&$11,13$&$E_{7}(a_{1})$&$124$\\
$E_{7}$&$16$&$1$&$A_{1}$&$34$\\
$E_{7}$&$16$&$3$&$A_{3}+A_{2}+A_{1}$&$100$\\
$E_{7}$&$16$&$5$&$E_{7}(a_{5})$&$112$\\
$E_{7}$&$16$&$7$&$E_{6}(a_{1})$&$118$\\
$E_{7}$&$16$&$9$&$E_{7}(a_{2})$&$122$\\
$E_{7}$&$16$&$11,13$&$E_{7}(a_{1})$&$124$\\
$E_{7}$&$17$&$1$&$A_{1}$&$34$\\
$E_{7}$&$17$&$2$&$A_{2}+3A_{1}$&$84$\\
$E_{7}$&$17$&$3$&$A_{3}+A_{2}+A_{1}$&$100$\\
$E_{7}$&$17$&$4$&$D_{5}(a_{1})+A_{1}$&$108$\\
$E_{7}$&$17$&$5$&$E_{7}(a_{5})$&$112$\\
$E_{7}$&$17$&$6$&$E_{7}(a_{4})$&$116$\\
$E_{7}$&$17$&$7$&$E_{6}(a_{1})$&$118$\\
$E_{7}$&$17$&$8$&$E_{7}(a_{3})$&$120$\\
$E_{7}$&$17$&$9,10$&$E_{7}(a_{2})$&$122$\\
$E_{7}$&$17$&$11\text{--}13$&$E_{7}(a_{1})$&$124$\\
\end{longtable}
}
{\small\setlength{\tabcolsep}{4pt}
\begin{longtable}{@{}p{.10\textwidth}p{.43\textwidth}p{.30\textwidth}p{.09\textwidth}@{}}
\caption{Associated varieties at admissible levels: $E_{8}$.}\label{tab:av-E8-adm}\\
\toprule Type&$q$&Orbit or partition&$\dim X$\\\midrule\endfirsthead
\multicolumn{4}{l}{\tablename~\thetable\ (continued)}\\\toprule
Type&$q$&Orbit or partition&$\dim X$\\\midrule\endhead
\midrule\multicolumn{4}{r}{Continued on the next page}\\\endfoot\bottomrule\endlastfoot
$E_{8}$&$1$&$0$&$0$\\
$E_{8}$&$2$&$4A_{1}$&$128$\\
$E_{8}$&$3$&$2A_{2}+2A_{1}$&$168$\\
$E_{8}$&$4$&$2A_{3}$&$188$\\
$E_{8}$&$5$&$A_{4}+A_{3}$&$200$\\
$E_{8}$&$6$&$E_{8}(a_{7})$&$208$\\
$E_{8}$&$7$&$A_{6}+A_{1}$&$212$\\
$E_{8}$&$8$&$A_{7}$&$218$\\
$E_{8}$&$9$&$E_{8}(b_{6})$&$220$\\
$E_{8}$&$10,11$&$E_{8}(a_{6})$&$224$\\
$E_{8}$&$12,13$&$E_{8}(a_{5})$&$228$\\
$E_{8}$&$14$&$E_{8}(b_{4})$&$230$\\
$E_{8}$&$15\text{--}17$&$E_{8}(a_{4})$&$232$\\
$E_{8}$&$18,19$&$E_{8}(a_{3})$&$234$\\
$E_{8}$&$20\text{--}23$&$E_{8}(a_{2})$&$236$\\
$E_{8}$&$24\text{--}29$&$E_{8}(a_{1})$&$238$\\
$E_{8}$&$q\ge 30$&$E_{8}$&$240$\\
\end{longtable}
}
{\small\setlength{\tabcolsep}{4pt}
\begin{longtable}{@{}p{.08\textwidth}p{.06\textwidth}p{.34\textwidth}p{.26\textwidth}p{.14\textwidth}@{}}
\caption{Proper orbit closures at nonadmissible quasi-lisse levels: $E_{8}$.}\label{tab:av-E8-nonadm}\\
\toprule Type&$p$&$q$&Orbit or partition&$\dim X$\\\midrule\endfirsthead
\multicolumn{5}{l}{\tablename~\thetable\ (continued)}\\\toprule
Type&$p$&$q$&Orbit or partition&$\dim X$\\\midrule\endhead
\midrule\multicolumn{5}{r}{Continued on the next page}\\\endfoot\bottomrule\endlastfoot
$E_{8}$&$6$&$1$&$E_{8}(a_{7})$&$208$\\
$E_{8}$&$6$&$5$&$E_{8}(a_{1})$&$238$\\
$E_{8}$&$9$&$1$&$D_{4}(a_{1})+A_{2}$&$184$\\
$E_{8}$&$9$&$2$&$D_{7}(a_{2})$&$216$\\
$E_{8}$&$9$&$4$&$E_{8}(a_{4})$&$232$\\
$E_{8}$&$9$&$5$&$E_{8}(a_{3})$&$234$\\
$E_{8}$&$9$&$7,8$&$E_{8}(a_{1})$&$238$\\
$E_{8}$&$10$&$1$&$D_{4}(a_{1})+A_{1}$&$176$\\
$E_{8}$&$10$&$3$&$E_{8}(a_{6})$&$224$\\
$E_{8}$&$10$&$7$&$E_{8}(a_{2})$&$236$\\
$E_{8}$&$10$&$9$&$E_{8}(a_{1})$&$238$\\
$E_{8}$&$11$&$1$&$D_{4}(a_{1})+A_{1}$&$176$\\
$E_{8}$&$11$&$2$&$E_{8}(a_{7})$&$208$\\
$E_{8}$&$11$&$3$&$E_{8}(a_{6})$&$224$\\
$E_{8}$&$11$&$4$&$E_{8}(a_{5})$&$228$\\
$E_{8}$&$11$&$5$&$E_{8}(a_{4})$&$232$\\
$E_{8}$&$11$&$6$&$E_{8}(a_{3})$&$234$\\
$E_{8}$&$11$&$7$&$E_{8}(a_{2})$&$236$\\
$E_{8}$&$11$&$8,9$&$E_{8}(a_{1})$&$238$\\
$E_{8}$&$12$&$1$&$2A_{2}$&$156$\\
$E_{8}$&$12$&$5$&$E_{8}(a_{5})$&$228$\\
$E_{8}$&$12$&$7$&$E_{8}(a_{3})$&$234$\\
$E_{8}$&$12$&$11$&$E_{8}(a_{1})$&$238$\\
$E_{8}$&$13$&$1$&$2A_{2}$&$156$\\
$E_{8}$&$13$&$2$&$D_{6}(a_{2})$&$204$\\
$E_{8}$&$13$&$3$&$D_{7}(a_{2})$&$216$\\
$E_{8}$&$13$&$4$&$E_{8}(a_{6})$&$224$\\
$E_{8}$&$13$&$5$&$E_{8}(a_{5})$&$228$\\
$E_{8}$&$13$&$6$&$E_{8}(a_{4})$&$232$\\
$E_{8}$&$13$&$7$&$E_{8}(a_{3})$&$234$\\
$E_{8}$&$13$&$8,9$&$E_{8}(a_{2})$&$236$\\
$E_{8}$&$13$&$10,11$&$E_{8}(a_{1})$&$238$\\
$E_{8}$&$14$&$1$&$A_{2}+2A_{1}$&$146$\\
$E_{8}$&$14$&$3$&$D_{5}+A_{2}$&$214$\\
$E_{8}$&$14$&$5$&$E_{8}(b_{5})$&$226$\\
$E_{8}$&$14$&$9$&$E_{8}(a_{2})$&$236$\\
$E_{8}$&$14$&$11,13$&$E_{8}(a_{1})$&$238$\\
$E_{8}$&$15$&$1$&$A_{2}+A_{1}$&$136$\\
$E_{8}$&$15$&$2$&$A_{4}+2A_{1}$&$192$\\
$E_{8}$&$15$&$4$&$E_{8}(b_{6})$&$220$\\
$E_{8}$&$15$&$7,8$&$E_{8}(a_{4})$&$232$\\
$E_{8}$&$15$&$11$&$E_{8}(a_{2})$&$236$\\
$E_{8}$&$15$&$13,14$&$E_{8}(a_{1})$&$238$\\
$E_{8}$&$16$&$1$&$A_{2}+A_{1}$&$136$\\
$E_{8}$&$16$&$3$&$E_{8}(a_{7})$&$208$\\
$E_{8}$&$16$&$5$&$E_{8}(a_{6})$&$224$\\
$E_{8}$&$16$&$7$&$E_{8}(a_{4})$&$232$\\
$E_{8}$&$16$&$9$&$E_{8}(a_{3})$&$234$\\
$E_{8}$&$16$&$11$&$E_{8}(a_{2})$&$236$\\
$E_{8}$&$16$&$13$&$E_{8}(a_{1})$&$238$\\
$E_{8}$&$17$&$1$&$A_{2}+A_{1}$&$136$\\
$E_{8}$&$17$&$2$&$A_{4}+2A_{1}$&$192$\\
$E_{8}$&$17$&$3$&$E_{8}(a_{7})$&$208$\\
$E_{8}$&$17$&$4$&$E_{8}(b_{6})$&$220$\\
$E_{8}$&$17$&$5$&$E_{8}(a_{6})$&$224$\\
$E_{8}$&$17$&$6$&$E_{8}(a_{5})$&$228$\\
$E_{8}$&$17$&$7,8$&$E_{8}(a_{4})$&$232$\\
$E_{8}$&$17$&$9$&$E_{8}(a_{3})$&$234$\\
$E_{8}$&$17$&$10,11$&$E_{8}(a_{2})$&$236$\\
$E_{8}$&$17$&$12\text{--}14$&$E_{8}(a_{1})$&$238$\\
$E_{8}$&$18$&$1$&$A_{2}$&$114$\\
$E_{8}$&$18$&$5$&$E_{8}(b_{6})$&$220$\\
$E_{8}$&$18$&$7$&$E_{8}(a_{5})$&$228$\\
$E_{8}$&$18$&$11$&$E_{8}(a_{3})$&$234$\\
$E_{8}$&$18$&$13$&$E_{8}(a_{2})$&$236$\\
$E_{8}$&$18$&$17$&$E_{8}(a_{1})$&$238$\\
$E_{8}$&$19$&$1$&$A_{2}$&$114$\\
$E_{8}$&$19$&$2$&$A_{3}+A_{2}+A_{1}$&$182$\\
$E_{8}$&$19$&$3$&$E_{6}(a_{3})+A_{1}$&$204$\\
$E_{8}$&$19$&$4$&$D_{5}+A_{2}$&$214$\\
$E_{8}$&$19$&$5$&$E_{8}(b_{6})$&$220$\\
$E_{8}$&$19$&$6$&$E_{8}(b_{5})$&$226$\\
$E_{8}$&$19$&$7$&$E_{8}(a_{5})$&$228$\\
$E_{8}$&$19$&$8$&$E_{8}(b_{4})$&$230$\\
$E_{8}$&$19$&$9$&$E_{8}(a_{4})$&$232$\\
$E_{8}$&$19$&$10,11$&$E_{8}(a_{3})$&$234$\\
$E_{8}$&$19$&$12,13$&$E_{8}(a_{2})$&$236$\\
$E_{8}$&$19$&$14\text{--}17$&$E_{8}(a_{1})$&$238$\\
$E_{8}$&$20$&$1$&$2A_{1}$&$92$\\
$E_{8}$&$20$&$3$&$A_{4}+A_{2}+A_{1}$&$196$\\
$E_{8}$&$20$&$7$&$E_{8}(a_{6})$&$224$\\
$E_{8}$&$20$&$9$&$E_{8}(b_{4})$&$230$\\
$E_{8}$&$20$&$11$&$E_{8}(a_{4})$&$232$\\
$E_{8}$&$20$&$13$&$E_{8}(a_{2})$&$236$\\
$E_{8}$&$20$&$17,19$&$E_{8}(a_{1})$&$238$\\
$E_{8}$&$21$&$1$&$2A_{1}$&$92$\\
$E_{8}$&$21$&$2$&$A_{3}+2A_{1}$&$172$\\
$E_{8}$&$21$&$4$&$E_{8}(a_{7})$&$208$\\
$E_{8}$&$21$&$5$&$D_{7}(a_{2})$&$216$\\
$E_{8}$&$21$&$8$&$E_{8}(a_{5})$&$228$\\
$E_{8}$&$21$&$10,11$&$E_{8}(a_{4})$&$232$\\
$E_{8}$&$21$&$13$&$E_{8}(a_{2})$&$236$\\
$E_{8}$&$21$&$16,17,19$&$E_{8}(a_{1})$&$238$\\
$E_{8}$&$22$&$1$&$2A_{1}$&$92$\\
$E_{8}$&$22$&$3$&$A_{4}+A_{2}+A_{1}$&$196$\\
$E_{8}$&$22$&$5$&$D_{7}(a_{2})$&$216$\\
$E_{8}$&$22$&$7$&$E_{8}(a_{6})$&$224$\\
$E_{8}$&$22$&$9$&$E_{8}(b_{4})$&$230$\\
$E_{8}$&$22$&$13,15$&$E_{8}(a_{2})$&$236$\\
$E_{8}$&$22$&$17,19$&$E_{8}(a_{1})$&$238$\\
$E_{8}$&$23$&$1$&$2A_{1}$&$92$\\
$E_{8}$&$23$&$2$&$A_{3}+2A_{1}$&$172$\\
$E_{8}$&$23$&$3$&$A_{4}+A_{2}+A_{1}$&$196$\\
$E_{8}$&$23$&$4$&$E_{8}(a_{7})$&$208$\\
$E_{8}$&$23$&$5$&$D_{7}(a_{2})$&$216$\\
$E_{8}$&$23$&$6$&$D_{7}(a_{1})$&$222$\\
$E_{8}$&$23$&$7$&$E_{8}(a_{6})$&$224$\\
$E_{8}$&$23$&$8$&$E_{8}(a_{5})$&$228$\\
$E_{8}$&$23$&$9$&$E_{8}(b_{4})$&$230$\\
$E_{8}$&$23$&$10,11$&$E_{8}(a_{4})$&$232$\\
$E_{8}$&$23$&$12$&$E_{8}(a_{3})$&$234$\\
$E_{8}$&$23$&$13\text{--}15$&$E_{8}(a_{2})$&$236$\\
$E_{8}$&$23$&$16\text{--}19$&$E_{8}(a_{1})$&$238$\\
$E_{8}$&$24$&$1$&$A_{1}$&$58$\\
$E_{8}$&$24$&$5$&$E_{8}(a_{7})$&$208$\\
$E_{8}$&$24$&$7$&$E_{8}(b_{6})$&$220$\\
$E_{8}$&$24$&$11$&$E_{8}(b_{4})$&$230$\\
$E_{8}$&$24$&$13$&$E_{8}(a_{4})$&$232$\\
$E_{8}$&$24$&$17$&$E_{8}(a_{2})$&$236$\\
$E_{8}$&$24$&$19,23$&$E_{8}(a_{1})$&$238$\\
$E_{8}$&$25$&$1$&$A_{1}$&$58$\\
$E_{8}$&$25$&$2$&$A_{2}+3A_{1}$&$154$\\
$E_{8}$&$25$&$3$&$D_{4}(a_{1})+A_{2}$&$184$\\
$E_{8}$&$25$&$4$&$D_{5}(a_{1})+A_{2}$&$202$\\
$E_{8}$&$25$&$6$&$D_{7}(a_{2})$&$216$\\
$E_{8}$&$25$&$7$&$E_{8}(b_{6})$&$220$\\
$E_{8}$&$25$&$8$&$E_{8}(a_{6})$&$224$\\
$E_{8}$&$25$&$9$&$E_{8}(b_{5})$&$226$\\
$E_{8}$&$25$&$11$&$E_{8}(b_{4})$&$230$\\
$E_{8}$&$25$&$12,13$&$E_{8}(a_{4})$&$232$\\
$E_{8}$&$25$&$14$&$E_{8}(a_{3})$&$234$\\
$E_{8}$&$25$&$16\text{--}18$&$E_{8}(a_{2})$&$236$\\
$E_{8}$&$25$&$19,21\text{--}23$&$E_{8}(a_{1})$&$238$\\
$E_{8}$&$26$&$1$&$A_{1}$&$58$\\
$E_{8}$&$26$&$3$&$D_{4}(a_{1})+A_{2}$&$184$\\
$E_{8}$&$26$&$5$&$E_{8}(a_{7})$&$208$\\
$E_{8}$&$26$&$7$&$E_{8}(b_{6})$&$220$\\
$E_{8}$&$26$&$9$&$E_{8}(b_{5})$&$226$\\
$E_{8}$&$26$&$11$&$E_{8}(b_{4})$&$230$\\
$E_{8}$&$26$&$15$&$E_{8}(a_{3})$&$234$\\
$E_{8}$&$26$&$17$&$E_{8}(a_{2})$&$236$\\
$E_{8}$&$26$&$19,21,23$&$E_{8}(a_{1})$&$238$\\
$E_{8}$&$27$&$1$&$A_{1}$&$58$\\
$E_{8}$&$27$&$2$&$A_{2}+3A_{1}$&$154$\\
$E_{8}$&$27$&$4$&$D_{5}(a_{1})+A_{2}$&$202$\\
$E_{8}$&$27$&$5$&$E_{8}(a_{7})$&$208$\\
$E_{8}$&$27$&$7$&$E_{8}(b_{6})$&$220$\\
$E_{8}$&$27$&$8$&$E_{8}(a_{6})$&$224$\\
$E_{8}$&$27$&$10$&$E_{8}(a_{5})$&$228$\\
$E_{8}$&$27$&$11$&$E_{8}(b_{4})$&$230$\\
$E_{8}$&$27$&$13$&$E_{8}(a_{4})$&$232$\\
$E_{8}$&$27$&$14$&$E_{8}(a_{3})$&$234$\\
$E_{8}$&$27$&$16,17$&$E_{8}(a_{2})$&$236$\\
$E_{8}$&$27$&$19,20,22,23$&$E_{8}(a_{1})$&$238$\\
$E_{8}$&$28$&$1$&$A_{1}$&$58$\\
$E_{8}$&$28$&$3$&$D_{4}(a_{1})+A_{2}$&$184$\\
$E_{8}$&$28$&$5$&$E_{8}(a_{7})$&$208$\\
$E_{8}$&$28$&$9$&$E_{8}(b_{5})$&$226$\\
$E_{8}$&$28$&$11$&$E_{8}(b_{4})$&$230$\\
$E_{8}$&$28$&$13$&$E_{8}(a_{4})$&$232$\\
$E_{8}$&$28$&$15$&$E_{8}(a_{3})$&$234$\\
$E_{8}$&$28$&$17$&$E_{8}(a_{2})$&$236$\\
$E_{8}$&$28$&$19,23$&$E_{8}(a_{1})$&$238$\\
$E_{8}$&$29$&$1$&$A_{1}$&$58$\\
$E_{8}$&$29$&$2$&$A_{2}+3A_{1}$&$154$\\
$E_{8}$&$29$&$3$&$D_{4}(a_{1})+A_{2}$&$184$\\
$E_{8}$&$29$&$4$&$D_{5}(a_{1})+A_{2}$&$202$\\
$E_{8}$&$29$&$5$&$E_{8}(a_{7})$&$208$\\
$E_{8}$&$29$&$6$&$D_{7}(a_{2})$&$216$\\
$E_{8}$&$29$&$7$&$E_{8}(b_{6})$&$220$\\
$E_{8}$&$29$&$8$&$E_{8}(a_{6})$&$224$\\
$E_{8}$&$29$&$9$&$E_{8}(b_{5})$&$226$\\
$E_{8}$&$29$&$10$&$E_{8}(a_{5})$&$228$\\
$E_{8}$&$29$&$11$&$E_{8}(b_{4})$&$230$\\
$E_{8}$&$29$&$12,13$&$E_{8}(a_{4})$&$232$\\
$E_{8}$&$29$&$14,15$&$E_{8}(a_{3})$&$234$\\
$E_{8}$&$29$&$16\text{--}18$&$E_{8}(a_{2})$&$236$\\
$E_{8}$&$29$&$19\text{--}23$&$E_{8}(a_{1})$&$238$\\
\end{longtable}
}

\bigskip

\section{Calculations in conformal weight four}\label{app:weight-four}
This appendix gives the three coefficient calculations used in
\cref{com:sf-recognition,v5:e7-detection,v5:scalar-coefficient},
including the derivative terms in the complete vectors.

\subsection{Primary vectors in the symplectic-fermion realization}
\label{app:primary-four}
We use the common lattice realization and notation of
\cref{com:sf-recognition}. In particular,
$\mathcal A_r=\SF(r)^{\Sp_{2r}}$, $T=T_a+T_r$,
$p_j=a^j+\sum_i h_i^j$, and $[a_\lambda a]=\lambda\mathbf1$.
We show that the $C_2$ projection is injective on the weight-four
primary subspace and that a primary vector with symbol
$\alpha p_4+\beta p_2^2$, $\alpha\ne0$, is a multiple of the
weight-four primary generator of $\mathcal A_{r+1}$.

Replacing $a,h_i$ by $-a,-h_i$, respectively, gives
$T_a=(:a^2:+\partial a)/2$ and
$T_r=\sum_i(:h_i^2:+\partial h_i)/2$.
This substitution fixes every polynomial of even total degree. Write $L'_n$ for the modes of
$T_r$, $D=5r-11$, $S=\sum h_i^2$, and $R=\sum h_i^4$.
Through weight four the tensor factor $\mathcal A_r$ has only its Virasoro
states and, when $r\ge2$, one primary vector $W_r$ independent of Virasoro descendants.
The weight-four space of $\mathcal H_a\otimes\Vir_{-2r}$ has
dimension ten. The $L_1,L_2$ equations have an eight-column minor
$-64r^2(5r-11)\ne0$. Its two kernel vectors $P,Q$ also satisfy
$L_3P=L_4P=L_3Q=L_4Q=0$ and have leading symbols
$$
 P:\ a^4+\frac9{r(5r-11)}S^2,\qquad
 Q:\ 3a^4+\frac6r a^2S-\frac{21}{2r(5r-11)}S^2.
$$
The coefficients of the complete vectors and the minor are recorded in
\hyperref[calc:SCALAR]{Calculation~\ref*{calc:SCALAR}}. These explicit equations establish
both the kernel dimension and the displayed symbols.
They are independent. When $r\ge2$, the additional primary $W_r$ has
symbol $R+b_rS^2$ and adds an independent $R$ term. This normalization exists because the invariant quadratic fermion of weight four bosonizes with
nonzero $R$ coefficient, and primary correction adds only $S^2$.
Consequently the map to $R_{\mathcal H_a\otimes\mathcal H^{\otimes r}}$
is injective on the space of primary vectors of conformal weight four
in $\mathcal H_a\otimes \mathcal A_r$.

For $r=1$, comparison of the three quartic coefficients shows that the
unique primary with symbol $\alpha(a^4+S^2)+\beta(a^2+S)^2$ has
$\beta=-6\alpha$. It is the primary of the concrete $\mathcal A_2$.
Inductively, if $b_r=6/(5r-11)$, write the primary as
$uP+vQ+\alpha W_r$. Comparing the $a^4,a^2S,S^2$ coefficients gives
$$
 v=\frac{r\beta}{3},\qquad u=\alpha+(1-r)\beta,
 \qquad (5r-6)(2r+3)\beta=6(2r+3)\alpha.
$$
Thus $\beta=6\alpha/(5r-6)$, and the primary is unique. The concrete
inclusion $\mathcal A_{r+1}\subseteq\mathcal H_a\otimes \mathcal A_r$ realizes it:
take invariants under the subgroup $\Sp_{2r}$ fixing the first fermion
pair and its torus. This proves both uniqueness and the induction.

\subsection{The product \texorpdfstring{$G_{(3)}P$}{G(3)P} in \texorpdfstring{$\mathcal W^{-8}(E_7,A_2)$}{W-8(E7,A2)}}
\label{app:e7-four}
We use the notation of the $E_7,t=10$ argument in
Section~\ref{exc:seeds}: $W=\mathcal W^{-8}(E_7,A_2)$ and
$\mathcal I=H_f^0(I_3)$, where $I_3$ is the ideal of $V^{-8}(E_7)$
generated by its singular subspace of degree three.
The quotient $W/\mathcal I$ is $\Vir^{-22/5}$. We construct the invariant
weight-four vector $P\in \mathcal I$ used in that argument and compute
its product with a charged generator $G$.

Let $C=V^0(A_5)\subset W$, $S$ its Sugawara vector,
$T=\omega-S$, and $H=L_0^S$. The nontrivial $A_5$-modules of strong generators are $15+15^*$,
both of conformal weight two. Their invariant contraction is $Q=\sum:G_iG^i:$,
using the dual of $([e,u_i],[e,u_j])$ for their finite lowest vectors.
Choose an $A_5$-invariant generator $w\in \mathcal I_3$ by subtracting its
Virasoro quotient image, and set $w_0=(1-H/3)w$. It is annihilated by $a(n)$ for $a\in A_5$, $n>0$,
and is Virasoro primary: the subspace of $A_5$-invariants in
$\mathcal I_2$ has no nonzero vector annihilated by all these operators.
The basis of normally ordered products of the universal strong
generators gives
$$
 \mathcal I_4^{A_5}=C_4^{A_5}\oplus\C:ST:\oplus\C\partial w_0\oplus\C Q,
 \qquad\dim C_4^{A_5}=5.
$$
The vectors $a(n)Q$, $a\in A_5$, $n>0$, lie in the adjoint
submodules of $C$ and in the span of $:aT:$.
Terms involving $15$ or $15^*$ have nontrivial central character for
the simply connected group of type $A_5$. Thus
$HQ\in C_4^{A_5}+\C:ST:$, where $H$ has eigenvalues four and two.
By \cref{v5:current-projection},
$$
 Q_0=(1-H/2)(1-H/4)Q,\qquad
 P=Q_0-\tfrac16\partial(L_1Q_0).
$$
The vector $P$ spans the unique $A_5$-invariant line in $\mathcal I_4$ annihilated by all
operators $a(n)$, $a\in A_5$, $n>0$, and consisting of Virasoro-primary
vectors. The $Q$ coefficient is one.
$L_1\partial w_0=6w_0$ proves uniqueness.

The charged fields are reconstructed from two complete highest-weight
BRST cocycles of $A_5$-weights $\omega_2$ and $\omega_4$. Their
lowering zero modes give the two fifteen-dimensional modules. These
zero modes preserve BRST closure and primarity. Set
$\Phi(H)=(1-H/2)(1-H/4)$ and choose the derivative coordinate $\rho$
on the highest charged field $G$ of weight $\omega_2$ in the
equivariant decomposition
$\Sym^2\g_0\oplus\partial\g_0$ of the weight-two subspace of
the affine vertex algebra in the Miura homomorphism.
Since $W_1$ contains no charged module, translation gives
$$
 \rho(G_{(3)}P)=
 \rho\bigl(G_{(3)}\Phi(H)Q-\tfrac12G_{(2)}L_1\Phi(H)Q\bigr).
$$
The Weyl subgroup $S_2\times S_4$ fixing the highest weight fixes $G$
and its nonzero derivative projection, which spans the
one-dimensional weight-$\omega_2$ space in $\partial\g_0$. Thus
$\rho$ is invariant under this subgroup. The ordered dual contractions split into six orbits,
of sizes $1,8,6$ in each orientation, from the decomposition
$\Lambda^2(\C^2\oplus\C^4)$. The two orientations give the same
three scalar values. In the normalization $\rho(G)=3/2$, their
weighted sum is
$$
 \rho(G_{(3)}P)=2\left(-\frac{2849}{216}
           -8\frac{9889}{864}-6\frac{1397}{144}\right)
       =-\frac{8800}{27}.
$$
All vectors entering this identity have weight at most four. The complete highest fields,
finite pairings and six scalar evaluations are reproduced in
\hyperref[calc:V5-E7-P10]{Calculation~\ref*{calc:V5-E7-P10}}.
All calculations take place before passing from $V^0(A_5)$
to its simple quotient.
The $A_5$-module containing $G$ occurs with multiplicity one in $W_2$,
so this nonzero coordinate of the Miura image determines the scalar in $W$. Consequently
$$
                  G_{(3)}P=-\frac{17600}{81}G\ne0,
$$
which is the identity \eqref{v5:e7-detection} used in the main proof.

\subsection{The coefficient of \texorpdfstring{$:TT:$}{:TT:} in a product of primary vectors}
\label{app:virasoro-four}
We use the hypotheses and notation of
\cref{v5:scalar-coefficient}. Thus $T$ is a Virasoro vector of
central charge $c$, commuting with $\omega-T$,
$E=L^T_{-4}\mathbf1$, $F=(L^T_{-2})^2\mathbf1$ and
$Q=\sum_{a=1}^m(G_a)_{(-1)}G^a$. The vectors $G_a,G^a$ have
conformal weight two for $\omega$ and weight $h$ for $T$, with
$(G_a)_{(3)}G^b=N\delta_a^b\mathbf1$.
The relation is normalized to have $Q$ coefficient one.

The Gram matrix on the span of $E,F$ and its determinant are
$$
 M=\begin{pmatrix}5c&3c\\3c&c(c+8)/2\end{pmatrix},\qquad
 \det M=c^2(5c+22)/2.
$$
The column vector of pairings with $E,F$, obtained from the Virasoro commutation
relations, is $(3mNh,mNh(h+2))^{\mathsf T}$. The coefficients
are its negative image under $M^{-1}$. To verify this column even
when $T$ is only part of the conformal vector, set $D_T=L^T_{-1}$.
The full translation operator is
$L_{-1}=L^T_{-1}+L^{\omega-T}_{-1}$. The full invariant adjoint identity gives, for a dual pair $G,G'$,
$$
 (D_TG)_{(4)}G'=-2hN,\qquad
 (D_T^2G)_{(5)}G'=2h(2h+1)N.
$$
Indeed $L_1D_TG=2hG$ and
$L_1^2D_T^2G=4h(2h+1)G$. The second adjoint expansion has its
factor $1/2$. Using
$[L^T_n,G_{(r)}]=(D_TG)_{(n+r+1)}+(n+1)hG_{(n+r)}$
now gives $\langle E,Q\rangle=3mNh$ and
$$
 \langle F,Q\rangle/m
 =2h(2h+1)N-(6h+3)2hN+(6h+9h^2)N=Nh(h+2).
$$
Since $D'$ is orthogonal to $E,F$, its coefficients do not enter
the calculation of $\beta$.

Inverting the displayed matrix gives
$$
               \beta=-\frac{2mNh(5h+1)}{c(5c+22)}.
$$
For the application to $\mathcal W^{-17}(E_8,2A_2)$,
$(m,N,h,c)=(49,-16,6/5,-384/65)$, and therefore
$\beta=1183/4$, as in \eqref{exc:p13-relation}.

\section{Twisting functors and the equivalence of positive blocks}
\label{app:twisting-compatibility}
We prove Lemma~\ref{com:twisting-compatibility}. The notation and the categories
are those of Section~\ref{com:transport}. The semiregular twisting
functors and their adjunction are those of \cite[Sections~4--5]{ABGG}.
For the conjugation construction for twisting functors associated with
roots simple only in the integral root system, we use
\cite[Corollary~5.6 and Lemma~5.8]{ABGG}. The BGG reciprocity statements
used for projectives are \cite[Theorem~2.2]{Fiebig} and
\cite[Theorem~3.4]{Liu2022}. We include the details to establish
compatibility on homomorphisms in the specified large positive blocks,
including singular highest weights.

\begin{proof}[Proof of Lemma~\ref{com:twisting-compatibility}]
For a simple root $\alpha_i$ of $\widehat\g$ that belongs to the integral
root system, let $T_i$ be localization
twisting and $Z_i$ the maximal rank-one-integrable quotient.
On a small projective $P$ we first establish the natural identification
\begin{equation}\label{com:projective-twist}
                      T_iP=\ker(P\longrightarrow Z_iP).
\end{equation}
For an integral $h_i$-weight
vector $m$ of weight $\mu$, the localization vector $f^{-j}\otimes m$
is sent, in the usual Weyl-lift convention, to
\begin{equation}\label{com:localization-map}
 \sum_{\substack{l\ge0\\\mu+j+l\ge0}}
 \frac{(-1)^l f^{\mu+j+l}e^lm}
 {(j-1)!\,l!\,(\mu+j+l)!},\qquad j\ge1.
\end{equation}
This is finite because $e$ is locally nilpotent. The relations
$[e,f^a]=a f^{a-1}(h-a+1)$ verify localization compatibility and
rank-one covariance. Expanding the same finite binomial identity on
$E\otimes M$, for a finite-dimensional rank-one module $E$, proves
tensor compatibility. Apply it to each finite adjoint rank-one
summand of the affine algebra to obtain full affine covariance.
Let $\lev^{(i)}=\widehat{\h}+\mathfrak{sl}_2^{(i)}$ and let
$\mathfrak p_i=\lev^{(i)}+\mathfrak u_i^+$ be the corresponding
minimal standard parabolic subalgebra. Write $\mathfrak u_i^-$ for
its opposite nilradical. Choose $H_i$ with
$\alpha_i(H_i)=0$ and $\alpha_j(H_i)=1$ for $j\ne i$.
Each fixed $H_i$-degree of $U(\mathfrak u_i^\pm)$ is a finite
sum of finite-dimensional integrable $\mathfrak{sl}_2^{(i)}$-modules.
Indeed finite real-root strings bound the coefficient of $\alpha_i$
once the other simple-root coefficients are bounded, and the Poincar\'e--Birkhoff--Witt theorem then
gives the assertion for the enveloping algebra.

For a weight generator $v$ of a small projective $P$, take a finitely generated
projective $Q$ in the rank-one category $\mathcal O$ surjecting onto $U(\lev^{(i)})v$. The common
highest-weight bound makes all sufficiently large positive
$H_i$-degrees act trivially on its image. For such $N$, the
$\mathfrak p_i$-module
$$
 Q_N=\bigl(U(\mathfrak u_i^+)/U(\mathfrak u_i^+)_{>N}\bigr)
                 \otimes Q
$$
maps to $P$. Induction to the affine algebra has a finite filtration by ordinary
Verma modules, obtained by filtering $Q_N$ by $H_i$-degree and its
rank-one Verma filtrations. Its restriction to $\lev^{(i)}$ is
$U(\mathfrak u_i^-)\otimes Q_N$, a direct sum of projective
rank-one modules: tensoring with a finite-dimensional integrable
module preserves projectives because tensoring with its dual is
an exact adjoint. The projection to the block of $P$ lies in
$\mathcal O^{\leq\lambda^+}_\Lambda$ and is a direct summand of this
module. Taking the sum over a finite generating set of $P$ gives a
surjection onto $P$. This surjection splits by projectivity.
Thus the rank-one restriction of $P$ is projective, and the
complementary Cartan grading is preserved.

The natural map \eqref{com:localization-map} is injective on each
rank-one projective. On $M(m)$, $m\ge0$, it has image
$M(-m-2)\subset M(m)$ and integrable cokernel $L(m)$.
On the singular projective $M(-1)$ it is an isomorphism.
The other projective in the regular block is the corresponding
block summand of $L(m+1)\otimes M(-1)$. Tensor compatibility
makes the map an isomorphism on this summand as well. These are
all rank-one projectives, so the map on $P$ is injective and its
cokernel is $i$-integrable. Every map from $P$ to an $i$-integrable
module kills its image, by naturality and $T_iN=0$ for such $N$.
The cokernel is therefore exactly $Z_iP$. This proves
\eqref{com:projective-twist} as a natural identification of
$\widehat\g$-module functors on small projectives, including
singular positive blocks.

The integrable quotient is intrinsic to its simple labels. In fact a
finitely generated object all of whose constituents are $i$-integrable
has reflection-invariant rank-one character. Each weight string is
bounded above, hence also below, so $f_i$ is locally nilpotent.
Pass to unions for the large block. The labelled equivalence $\mathcal E$
identifies the maximal quotients determined by the corresponding simple
labels and their natural quotient maps. If a corresponding root is
simple only in the integral root system, use the conjugated quotient
described below. Exactness then identifies the kernels in
\eqref{com:projective-twist}. Right exactness, algebraic direct sums
and projective presentations extend the identification to all objects
and their homomorphisms.

If $\alpha$ is simple only in the integral root system, conjugate
through equivalences $Q$ across nonintegral walls to a simple root of $\widehat\g$.
The cited twisting lemmas identify
$$
 T_{s_\alpha}=Q^{-1}T_iQ,
 \qquad \widetilde Z_\alpha=Q^{-1}Z_iQ.
$$
Its quotient labels are those with
$\langle\mu+\widehat\rho,\alpha^\vee\rangle\in\Z_{>0}$.
They are preserved by \eqref{com:root-dilation}. This is a
\emph{conjugated} integrable quotient. Literal integrability for
the $\mathfrak{sl}_2$ associated with the original root $\alpha$ would be incorrect.
Taking compositions along reduced expressions in the integral
Weyl group gives the required compatibility between $\mathcal E$
and the twisting functors, with the full singular stabilizer
taken into account. 
\end{proof}

\bigskip

\end{document}